\documentclass{amsart}
\usepackage{graphicx} % Required for inserting images
\numberwithin{equation}{section}

\title{Translation functors, branching laws and the restriction  of coherent cohomology}
\date{September 21, 2025}
\usepackage{graphicx}
\usepackage{amssymb, amscd}
\usepackage{epstopdf}
\usepackage{color}
\usepackage{framed, hyperref}
\usepackage[titletoc]{appendix}
\usepackage{bbm}
\usepackage{lipsum}
\usepackage{adjustbox}
\usepackage{amsmath, amsthm,latexsym, verbatim}
\usepackage{here}
\usepackage[all]{xy}
\definecolor{darkgreen}{rgb}{0,0.35,0}

\numberwithin{table}{section}

\newtheorem{theorem}{Theorem}[section]
\newtheorem{prop}[theorem]{Proposition}
\newtheorem{proposition}[theorem]{Proposition}

\newtheorem{lemma}[theorem]{Lemma}
\newtheorem{cor}[theorem]{Corollary}
\newtheorem{definition}[theorem]{Definition}
\newtheorem{example}[theorem]{Example}
\newtheorem{conj}[theorem]{Conjecture}
\newtheorem{vconj}[theorem]{Vague Conjecture}
\newtheorem{conv}[theorem]{Convention}

\newtheorem{question}[theorem]{Question}
\newtheorem{setting}[theorem]{Setting}

\newcommand{\btheorem}{\begin{theorem}}
\newcommand{\etheorem}{\end{theorem}}
\newcommand{\bprop}{\begin{prop}}
\newcommand{\eprop}{\end{prop}}
\newcommand{\blemma}{\begin{lemma}}
\newcommand{\elemma}{\end{lemma}}
\newcommand{\bcor}{\begin{cor}}
\newcommand{\ecor}{\end{cor}}

\theoremstyle{remark}
\newtheorem{remark}[theorem]{Remark}

\newcommand{\db}{{\overline{\partial}}}
\newcommand{\ttl}{{\tilde{\lambda}}}

\newcommand{\af}{{\mathbf{A}_f}}

\newcommand{\bA}{{\mathbf A}}

\newcommand{\CK}{{\mathcal K}}
\newcommand{\CKp}{{\mathcal{K}^+}}
\newcommand{\CA}{{\mathcal A}}

\newcommand{\Qbar}{{\bar{\mathbb Q}}}
\newcommand{\C}{{\mathbb C}}
\newcommand{\CC}{{\mathbb C}}
\newcommand{\RR}{{\mathbb R}}

\newcommand{\QQ}{{\mathbb Q}}
\newcommand{\CV}{{\mathcal{V}}}

\newcommand{\Z}{{\mathbb Z}}
\newcommand{\IS}{{\mathbb S}}

\newcommand{\G}{{\mathbb G}}

\newcommand{\IP}{{\mathbb P}^-}

\newcommand{\ra}{\rightarrow}

\newcommand{\fg}{\mathfrak g}

\newcommand{\fk}{\mathfrak k}

\newcommand{\fp}{\mathfrak p}
\newcommand{\fP}{\mathfrak P}
\newcommand{\ft}{\mathfrak t}
\newcommand{\op}{\oplus}
\newcommand{\oi}{\ominus}

\newcommand{\isoarrow}{{~\overset\sim\longrightarrow~}}

\definecolor{lightgreen}{rgb}{0.6,1.0,0.6}
\definecolor{lemon}{rgb}{1.0,1.0,0.5}
\definecolor{lightblue}{rgb}{0.6,1.0,1.0}
\newcommand{\tlambda}{\lambda}
\newcommand{\poslmd}{\lambda^+}
\newcommand{\poslmdp}{\lambda^{\prime, +}}
\newcommand{\mulmd}{\mu_{\lambda}}
\newcommand{\mulmdp}{\mu^\prime_{\lambda'}}
\newcommand{\Wlmdp}{W^\prime_{\lambda'}}
\newcommand{\pilmdp}{\pi'_{\lambda'}}

\usepackage[parfill]{parskip}

\begin{document}

\author{Michael Harris}
\author{Toshiyuki Kobayashi}
\author{Birgit Speh}

\address[Michael HARRIS] {%
Department of Mathematics,
Columbia University,
New York, NY  10027
USA}

\address[Toshiyuki KOBAYASHI]{%
Graduate School of Mathematical Sciences,
The University of Tokyo;
%3-8-1 Komaba, %Meguro, Tokyo, 153-8914, Japan;
French-Japanese Laboratory in Mathematics and its Interactions,
FJ-LMI CNRS IRL2025, %Tokyo,
Japan.}

\address[Birgit SPEH]{%
Department of Mathematics,
Cornell University,
Ithaca, NY 14850,
USA.
}

\title[Translation functors and restriction of coherent cohomology]{Translation functors, branching problems, and applications to the restriction of coherent cohomology of Shimura varieties}

\thanks{M.H. was partially supported by NSF Grant DMS-2302208. 
T.K. was partially supported by the JSPS under Grant-in Aid for Scientific Research (A) (JP23H00084). The research of B.S. was partially supported by a grant from the Simons foundation}

\begin{abstract}
We study properties of the restriction of discrete series representations of $G=U(p,q)$ to $G'= U(p-1,q)$ and the corresponding symmetry breaking operators in $\operatorname{Hom}_{G'}(\pi|_{G'}, \pi')$. This leads to the introduction of elementary and coherent pairs of discrete series representations and their classification.

Translations of symmetry breaking operators are defined via tensor products with finite-dimensional representations, which leads to the study of the coherent cohomology of discrete series representations under restriction and translations.  This is applied to the study of cup products of coherent cohomology of associated Shimura varieties, and to the arithmetic of central values of certain Rankin--Selberg $L$-functions of $GL(n+1)\times GL(n)$.
\end{abstract}

\maketitle


\setcounter{tocdepth}{2}
\tableofcontents

\section{Introduction}
\label{sec:intro}
Let $V$ be an $(n+1)$-dimensional complex vector space equipped with a non-degenerate hermitian inner product, and let $V' \subset V$ be a (non-degenerate) 
hermitian subspace of codimension one.
Denote by $G = U(V)$ and $G' = U(V')$ the corresponding unitary groups.
The \emph{branching law} for the restriction of discrete series representations of $G$ to the subgroup $G'$ describes the pairs $(\pi,\pi')$ of discrete series representations of $G$ and $G'$, respectively, such that
\begin{equation}\label{GGP1} \text{Hom}_{G'}(\pi|_{G'},\pi') \neq 0.
\end{equation}
Here, $\pi$ and $\pi'$ are understood as smooth Fr\'echet representations of moderate growth, as in \cite{GGP}.
 In \cite{KO} a complete classification is given for pairs of reductive groups $(G,G')$ for which the dimension of \eqref{GGP1} is uniformly bounded for all irreducible representations $(\pi, \pi')$. Moreover, for unitary groups as above, it is known that this space in \eqref{GGP1} is one-dimensional \cite{BP20, He, SZ}. The  branching law  was  conjectured by W.~Gan, B.~Gross and D.~Prasad for tempered representations, and is thus commonly referred to as \emph{Gan--Gross--Prasad (GGP) conjecture}. It was proved independently by Hongyu He and Raphaël Beuzart-Plessis.  In this article, we use He's description in \cite{He}, which provides a simple combinatorial criterion for \eqref{GGP1} in terms of Harish-Chandra parameters.  To $\pi$ (resp. $\pi'$)  H.He
assigns a sequence of $n+1$ (resp. $n$) symbols $+, -$ (resp. $\oplus,\ominus$) that we refer to as {\it decorations}, representing the positive chambers for the given Harish-Chandra parameters, and \eqref{GGP1} holds if and only if these symbols can be interleaved according to certain combinatorial rules.   We denote these decorations $\delta(\pi)$ and $\delta'(\pi')$ respectively, and write $ [\delta(\pi) \delta'(\pi' )] $ for the interleaved decorations and refer to $ [\delta(\pi) \delta'(\pi' )] $ as the signature of the pair $\pi$, $\pi'$. For a different perspective on the branching law for real forms of orthogonal groups, see \cite[Chapter 3]{KS}.

Now suppose that $H$, $H'$ are unitary groups of hermitian spaces over an imaginary quadratic field \footnote{The discussion applies 
to more general CM fields.} with $H(\RR) = G$,
$H'(\RR) = G'$, and let $Sh(H)$ and $Sh(H')$ be the corresponding Shimura varieties \cite{H21}.\footnote{We work with these Shimura varieties of {\it abelian type}, rather than the more familiar PEL type Shimura varieties attached to unitary similitude groups.}   Each of these Shimura varieties has a canonical family of {\it automorphic vector bundles}, attached to representations of the maximal compact subgroups $K \subset G$ and $K' \subset G'$, respectively.    If $[W]$ is such an automorphic vector bundle, then we denote by $H^*_!(Sh(H),[W])$ the interior coherent cohomology of $Sh(H)$ with coefficients in $[W]$, defined as in \cite{H14}. See Section~\ref{avbH} for more details.  There are also automorphic vector bundles $[F], [F']$ attached to finite-dimensional (algebraic) representations $F, F'$ of $H$ and $H'$; these are endowed with flat connections, and we will need them in the sequel.

To the representations $\pi$ of $G$ and $\pi'$ of $G'$ we can canonically associate automorphic vector bundles $[W_\pi]$ and $[W_{\pi'}]$ on $Sh(H)$ and $Sh(H')$, respectively, in such a way that, if $\Pi$ (resp. $\Pi'$) is 
a cuspidal automorphic representation of $H$ (resp. $H'$) such that
$$ \Pi \isoarrow \pi \otimes \Pi_f; ~~ \Pi' \isoarrow \pi' \otimes \Pi'_f$$
with $\Pi_f$ (resp. $\Pi'_f$) a representation of the finite ad\`eles $H(\af)$ of $H$ (resp. $H'(\af)$), then
\begin{equation}\label{degree}
\operatorname{Hom}_{H(\af)}(\Pi_f^{\vee},H^j_!(Sh(H),[W_\pi]))  \neq 0 \iff j = q(\pi),
\end{equation}
with a similar relation for $H'$.  Here $q(\pi)$ is an integer between $0$ and $\dim Sh(H)$, which can be computed explicitly in terms of the decoration $\delta(\pi)$ attached to $\pi$ (Definition~\ref{def:q}).

There is a natural inclusion of Shimura varieties 
$Sh(H') \hookrightarrow Sh(H)$\footnote{This is not quite correct -- one has to replace $Sh(H')$ by its product with a $0$-dimensional Shimura variety, as in \eqref{inclusionmap}, but this is close enough for the purposes of the introduction.} and we can consider the pullback of $[W_\pi]$ to $Sh(H')$. This is the first of a series of papers that aim to determine when the restriction maps of \eqref{GGP1} correspond to non-trivial restrictions of coherent cohomology of Shimura varieties.  It is easy to see that a necessary condition is that $q(\pi) = q(\pi')$, but it is not sufficient as we discuss in Section~\ref{sec:coherent}. 

\medskip

Our first main result (Theorem~\ref{thm:elementary_pair_equiv}) is the classification of those pairs $(\sigma,\sigma')$ for which there is a natural map $[W_\sigma] \ra [W_{\sigma'}]$, and such that the map in \eqref{GGP1} is induced from a non-zero $K'$-intertwining operator from the minimal $K$-type of $\sigma$ to the minimal $K'$-type of $\sigma'$.  We refer to such pairs $(\sigma,\sigma')$ as \emph{elementary pair} (Definition~\ref{def:elementary}), 
a notion that is an outgrowth of the work~\cite{H14}.
For a discussion of subtle differences among the various related properties, see 
Proposition~\ref{prop:23022216} and
Remark~\ref{rem:equiv_coherent_sign}.

For an elementary pair the inclusion of $Sh(H')$ in $Sh(H)$ defines a canonical homomorphism
\begin{equation}\label{resH}
H^q_!(Sh(H),[W_\sigma]))  \ra H^q_!(Sh(H'),[W_{\sigma'}])),
\end{equation}
and this can be shown to be non-trivial, in most cases, when $q = q(\sigma) = q(\sigma')$.

Now we consider pairs $(\pi,\pi')$ such that $\delta(\pi) = \delta(\sigma)$, $\delta'(\pi') = \delta'(\sigma')$ for some elementary pair $(\sigma,\sigma')$.  We call such a pair {\it coherent} (Definition~\ref{def:coh_deco}).
For any coherent pair $(\pi, \pi')$, its associated signature $[\delta(\pi)\delta'(\pi')]$
satisfies the GGP interleaving relations (Definition~\ref{def:interlace}), since the signature $[\delta(\sigma)$ $\delta'(\sigma')]$ associated with the elementary pair $(\sigma, \sigma')$ does.  Our ultimate objective is to show that the $1$-dimensional space $\operatorname{Hom}_{G'}(\pi|_{G'},\pi')$, in the case of coherent pairs, also gives rise to a non-trivial homomorphism of coherent cohomology of the corresponding Shimura varieties.   However, for a general coherent pair $(\pi,\pi')$ there is no intertwining map $[W_\pi] \ra [W_{\pi'}]$ of automorphic vector bundles.  Following an idea of Loeffler, Pilloni, Skinner, and Zerbes for the  analogous inclusion of a product of modular curves in the Siegel modular threefold, we replace $[W_\pi]$ and $[W_{\pi'}]$ by automorphic vector bundles of the form $[F]\otimes [W_\pi^\vee]$ and $[F']\otimes [W_{\pi'}^\vee]$, where $[F]$ and $[F']$ are the flat automorphic vector bundles mentioned above, in such a way that there are natural maps
\begin{equation}\label{auv}  [F] \otimes [W_\sigma^\vee] \ra [F']\otimes [W_{\sigma'}^\vee]
\end{equation}
that can be used to globalize the homomorphism from $\pi$ to $\pi'$.   We conjecture that such pairs $(F,F')$ always exist.  

Our second main result (Theorem \ref{thm:23081510}) is that any coherent pair $(\sigma,\sigma')$ can be linked to an elementary pair with the same decorations by a chain of diagrams of the form
\begin{equation}
\label{eqn:Fdiag20}
\xymatrix{
   \sigma \otimes F
     \ar[r]
     &
     \sigma' \otimes F'
     \ar[d]
\\
     \pi=\psi_{\lambda}^{\lambda+\nu}(\sigma)
     \ar@{^{(}-_>}[u]
     &\psi_{\lambda'}^{\lambda'+\nu'}(\sigma')=\pi'
}
\end{equation}
where $F$ and $F'$ are irreducible finite-dimensional representations of $H$ and $H'$, respectively, such that $\dim \operatorname{Hom}_{H'}(F,F') = 1$ and the left-hand (resp. right-hand) vertical map is an inclusion (resp. a projection) and the top horizontal map is the tensor product of the natural map
$\sigma \ra \sigma'$ of an elementary pair with an $H'$-equivariant homomorphism from $F$ to $F'$.  

When $F$ is the standard representation of $H$ or its dual, we prove that the corresponding map \eqref{auv} gives rise to a non-trivial map of interior coherent cohomology in degree $q(\pi) = q(\pi')$ that is compatible with the  branching rule at the archimedean components in the form proved by He; this is our third main result.   More precisely  we state in Conjecture~\ref{conj:250629}
%\ref{mainconj} 
a precise version of the following
\begin{vconj}\label{vague}  The chain of diagrams constructed in Theorem \ref{thm:23081510} that links a coherent pair $(\pi,\pi')$ to an elementary pair $(\sigma,\sigma')$ similarly gives rise to a non-trivial map of coherent cohomology.
\end{vconj}
Then we prove in Theorem~\ref{thm:non-van} that this conjecture holds when the chain is of length $1$,
by introducing a new framework for translation functors that applies to \emph{symmetry breaking} arising from the restriction $G \downarrow G'$,
(Theorem~\ref{thm:23081404}).
We plan to address this conjecture in general in subsequent work.

Finally, we apply these considerations to central critical values of $L$-functions.  The Ichino--Ikeda identity, in the version worked out for unitary groups by N. Harris \cite{II,NH} and proved in a series of articles, expresses the central value of the global automorphic $L$-function $L(s,\Pi\otimes (\Pi')^\vee)$ (or more precisely the Rankin--Selberg $L$-function for the base change of $\Pi\otimes (\Pi')^\vee$ to $GL(n+1)\times GL(n)$) in terms of period integrals over $H(\QQ)\backslash H'(\bA)$ of the products of automorphic forms in $\Pi\otimes (\Pi')^\vee$.  When the period integral can be expressed as a cup product pairing in coherent cohomology, we obtain an expression for the central value of $L(s,\Pi\otimes (\Pi')^\vee)$ in terms of automorphic periods, and this expression can be compared to Deligne's conjecture on critical values of motivic $L$-functions (see \cite{GHL}, for example).  Conjecture \ref{vague} implies that this is always the case when the archimedean components of $\Pi$ and $\Pi'$ form a coherent pair.  Our final result works out the consequences of our verification of this conjecture for chains of length one.

Most pairs $(\pi,\pi')$ that satisfy \eqref{GGP1} are not coherent.  For such pairs we do not know how to interpret the restriction maps in terms of algebraic geometry.  Nevertheless, the  Ichino--Ikeda identity holds whenever there is a global restriction map from $\Pi$ to $\Pi'$, and thus the period integrals attached to such $(\pi,\pi')$ must have an arithmetic meaning.

\bigskip
{\em An outline of the article:} The article is divided into 4 parts . In  part one and two we discuss some problems in representation theory of  unitary groups and thus they do not require any knowledge of Shimura varieties. The first part (Sections~\ref{sec:convent}--\ref{sec:coherent}) concerns  discrete series representations of the unitary groups $U(p,q)$ and their restriction to subgroups $U(p-1,q)$. The second part (Sections~\ref{sec:SBO_translation}--\ref{sec:Pk} introduces a new framework for translation functor that applies to \emph{symmetry breaking} arising from the restriction of  representations $U(p,q)\downarrow U(p-1,q)$
and discusses its applications. The third part (Sections~\ref{Shimura7}--\ref{specialvalues}) is devoted to applications of the results of the previous parts to global questions: Shimura varieties, periods and special values of L-functions. Part~\ref{part:4} consist of examples which illustrate the considerations in Parts~\ref{part:1} and~\ref{part:translation}. We close with a list of notation to help the reader.

After the introduction in Section~\ref{sec:intro}  we introduce in Section~\ref{sec:convent}  the Harish-Chandra parameter $\lambda$ of discrete series representations $\pi$, a reference parameter $\lambda ^+$, a decoration $\delta$, and in Definition~\ref{def:W_lambda} a {\it coherent parameter} $\Lambda$. 

\medskip
In Section~\ref{sec:resds}, we review He's combinatorial description of the GGP \emph{interleaving relations} \cite{He}, which characterize the non-vanishing of the spaces of symmetry breaking operators arising from the restriction of discrete series representations in the setting $(G, G')=(U(p,q), U(p-1,q))$.  In Definition~\ref{def:string_dd}, we define the signature $[\delta\delta']$ associated with a pair $(\pi, \pi')$ of discrete series representations of $G$ and $G'$, respectively,
where $\delta$ and $\delta'$ are the decorations of $\pi$ and $\pi'$, respectively.

Among the GGP interleaving patterns,
{\em coherent signatures}  are introduced in Definition~\ref{def:coh_deco}, and their equivalent characterizations and properties are discussed
in Proposition~\ref{prop:23022216}.
We call a pair $(\pi, \pi')$ of discrete series representations a {\em coherent pair}  (Definition~\ref{def:coherent_pair}) if \begin{equation} \text{Hom}_{G'}(\pi|_{G'},\pi') \neq 0
\end{equation}
and its signature $[\delta \delta']$ is coherent. 

We call a coherent pair $(\pi, \ \pi')$ an {\em elementary pair} (Definition~\ref{def:elementary}) if the minimal $K$-type of $\pi'$ is a direct summand of the restriction of the minimal $K$-type of $\pi$. 
In Theorem~\ref{thm:elementary_pair_equiv}, we give several equivalent conditions for a coherent pair to be elementary. In particular, the pair $(\pi,\pi')$ is elementary if the following conditions hold: $\operatorname{Hom}_{G'}(\pi|_{G'},\pi') \neq 0$, the minimal $K$-type of $\pi'$ is a direct summand of the restriction of the minimal $K$-type of $\pi$, and the finite-dimensional representation of $G'$ with highest weight $\Lambda'$ of $\pi'$ is a direct summand of the restriction of the finite-dimensional representation of $G$ with highest weight $\Lambda$.

Conversely, when viewed from the standpoint of an elementary pair,
Proposition~\ref{prop:cohpair} ensures that
a pair $(\pi,\pi')$ of discrete series representations of $G$ and $G'$ forms a coherent pair if and only if the following two conditions hold:
\begin{enumerate}
\item  $\operatorname{Hom}_{G'}(\pi|_{G'},\pi') \neq 0$; i.e. the decorations $\delta$ and $\delta'$ satisfy the GGP interleaving relation.
\item There exists an elementary pair $(\sigma,\sigma')$ with $\delta(\sigma) = \delta$ and $\delta(\sigma') = \delta'$.  
\end{enumerate}

The main theorems in Section~\ref{sec:coherent} are the construction of elementary pairs $(\pi, \pi')$. 
Theorem~\ref{thm:230816} addresses the problem of finding an elementary pair when a discrete series representation $\pi'$ of the subgroup $G'$ is given, by finding a corresponding discrete series representation $\pi$ of $G$. In contrast, 
Theorem~\ref{thm:23022420b} provides an elementary pair when a discrete series representation $\pi$ of the group $G$ is given, by finding a corresponding discrete series representation $\pi'$ of the subgroup $G'$. The latter theorem relies on a specific assumption about the parameters of $\pi$. In both cases, we also prove the uniqueness of the constructed pairs.

Furthermore we show in (iv) of Theorem~\ref{thm:elementary_pair_equiv} that for an elementary pair $\pi=\pi_\lambda$ and $\pi' =\pi'_{\lambda'}$ 
 $$
    H^q(\mathfrak{P}, K; \pi \otimes W_{\pi}^\vee) \longrightarrow
    H^q(\mathfrak{P}', K'; \pi^{\prime} \otimes W_{{\pi}'}^\vee)
    $$
    (see Section~\ref{sec:Pk} for further details of {\em coherent cohomology})
    is non-trivial.

    \medskip
    In Section \ref{sec:SBO_translation} we introduce a framework for translation functor that applies to symmetry breaking. 
    Although the translation functor serves as a powerful tool for studying families of representations of a \emph{single} group, no general theory has been developed to analyze \emph{symmetry breaking}---that is, morphisms arising from restriction of representations---between representations of a \emph{pair of groups} $G \supset G'$ under translation. This is rather delicate: for example, the vanishing condition for symmetry breaking is not preserved under translations even within the strict dominant chambers of the product group $G \times G'$.

\bigskip

    Given a coherent pair $(\pi, \pi')$, their coherent cohomology groups are concentrated in the same degree $q=q(\pi)=q(\pi')$.
In Section \ref{sec:Pk}, we investigate whether the symmetry breaking operator $S\colon \pi \to \pi'$ induces a natural homomorphism between the corresponding coherent $(\fP,K)$-cohomology groups.
Such a morphism does exist and is, in fact, an isomorphism, when $(\pi, \pi')$ forms an elementary pair (Theorem~\ref{thm:elementary_pair_equiv}.
However, for a non-elementary pair $(\pi, \pi')$, 
the existence of an analogous morphism cannot generally be expected, as already mentioned in the case of a general coherent pair, for which there exists no intertwining map $[W_\pi] \ra [W_{\pi'}]$ of automorphic vector bundles. 

We propose a conjectural framework involving the \emph{translation of symmetry breaking} in Conjecture~\ref{conj:250629}, and provide supporting evidence in Theorem~\ref{thm:non-van}, based on
 the new method of symmetry breaking under translation developed in Section~\ref{sec:SBO_translation}.

Section~\ref{Shimura7} briefly recalls the theory of coherent cohomology of Shimura varieties, in the special case of varieties attached to unitary groups.  The main results of the section (Theorem~\ref{cohc} and Proposition~\ref{redu}) describe the global calculation of the (interior) coherent cohomology of sufficiently regular automorphic vector bundles in terms of $(\fP,K)$-cohomology.  This calculation is then combined in Section~\ref{sec:rest_cohomology} with the results of Part~\ref{part:translation} to compute the restriction of coherent cohomology for an inclusion of Shimura varieties.

To end Part \ref{global}, Section~\ref{specialvalues} reviews the definition of automorphic periods in terms of coherent cohomology. By combining the relation of restriction maps for coherent pairs to cup products in coherent cohomology, on the one hand, with the Ichino--Ikeda--N. Harris formula for central values of Rankin--Selberg $L$-functions, on the other hand, we obtain an explicit expression of these central values in terms of automorphic periods.  This is stated in Theorem \ref{IIperiods} and more explicitly in Corollary~\ref{rewrite}.  The section also includes a brief section recalling Deligne's conjecture on critical values of motivic $L$-functions and references to papers that explain how the period invariant defined by Deligne can also be expressed in terms of automorphic periods, thus justifying the claim that our results on the central values of Rankin--Selberg $L$-functions are consistent with Deligne's conjecture.  

Readers are warned that Part~\ref{global} makes use of conventions different from those in the rest of the paper, to conform to those standard in the theory of Shimura varieties. 

\medskip
In  Part~\ref{part:4}, Section~\ref{sec:10}, we enumerate GGP interleaving patterns and then identify the conditions under which each of the various associated properties is satisfied in the low-dimensional case.  In particular, in Section~\ref{subsec:Part4_Branching},
we illustrate with examples which parts of GGP interlacing conditions give rise to \emph{coherent pairs} $(\pi, \pi')$, or more strongly, \emph{elementary pairs} $(\pi, \pi')$. We close by examining the symmetry breaking of non-holomorphic discrete series representations, which first appears in the case $(G,G') =(U(2,1),U(1,1))$.

 Section~\ref{sec:10}, the last section, is a guide to the notation used in the article.
We summarize first some notation for discrete series representations introduced in Section~\ref{sec:convent}
 of the main paper. 
 
We then establish a correspondence between the conventions adopted in Part~\ref{part:translation}, based on $(\mathfrak{p}^+, K)$-cohomology, and those in Part~\ref{global}, which are based on $(\mathfrak{p}^-, K)$-cohomology. More generally, several conventions for $(\mathfrak{p}, K)$-cohomology are found in the literature.
 These conventions vary depending on the choice of complex structure on the hermitian symmetric space $G/K$ (in our case, the holomorphic tangent bundle is identified with $G\times_K {\fp^+}$),  the definition of a maximal parabolic subalgebra $\fP$ (whether one takes $\fP^+=\mathfrak{k}+\fp^+$ or $\fP^-=\mathfrak{k}+\fp^-$), and whether the coefficients are given by a representation $\pi$ or its contragredient $\pi^\vee$ as coefficients.
Thus, in theory, there are 
$2\times 2 \times 2=8$ conceivable conventions.
Needless to say, the results proved in Part~\ref{part:translation} hold under all eight conceivable conventions.
For the reader’s convenience, we provide a correspondence among the different conventions.
To clarify our own conventions, we consider in Section~\ref{sec:10} the simplest case: a holomorphic discrete series representation.

\bigskip
{\bf{Acknowledgement.}}
\newline
The authors would like to thank the College of William and Mary for its hospitality.
T.~Kobayashi and B.~Speh would also like to thank  
 Mathematisches Forschungsinstitut, Oberwolfach (MFO), the Institut des Hautes {\'E}tudes Scientifiques (Bures-sur-Yvette)
  and Institut Henry Poincar{\'e} (Paris)  for the opportunity to work together under pleasant conditions.

\bigskip

\newpage
\part{Local theory}
\label{part:1}
\medskip

\medskip
\section{Conventions and parameters}
\label{sec:convent}

In this section, we present some conventions used throughout this article for the parametrization of discrete series representations (three equivalent parametrizations are reviewed in Proposition~\ref{prop:discpara}), as well as the associated combinatorial data such as decorations and coherent parameters.

\subsection{Unitary Groups and Subgroups}\label{subsec:ugps}

Let $V$ be an $(n+1)$-dimensional vector space over $\CC$ endowed with a non-degenerate hermitian form $h_V$ of signature $(p,q)$, $p + q = n+1$, and let $V' \subset V$ be a codimension $1$ subspace, where $h_V$ restricts to a form of signature $(p-1,q)$.  Let $G = U(V)$, $G^\flat = U(V')\times U((V')^\perp)$, viewed as the subgroup of $G$
that stabilizes $V'$ and its orthogonal complement in $V$, and $G' = U(V')$ the image of the action of $G^\flat$ on $V'$.
Then we may identify $G$ and $G'$ with the indefinite unitary group $U(p,q)$ and $U(p-1,q)$, respectively.

We write $V = V^+ \oplus V^-$, where $V^+$ and $V^-$ are maximal subspaces on which $h_V$ is positive-definite and negative-definite, respectively. For simplicity, we may assume that $V' = V^{\prime,+} \oplus V^-$, where $V^{\prime,+} = V^+ \cap V'$.
Let $K \subset G$ be the stabilizer of $V^+$ and $V^-$, $K' = K \cap G'$; thus
we may identify $K \simeq U(p) \times U(q)$ and $K' \simeq U(p-1)\times U(q)$, in the standard notation for compact unitary groups. 

We choose a basis $f_1, \dots, f_{p}$ of $V^+$
$V^+$ and a basis $f_{p+1}, \dots, f_{p+q}$ of $V^-$, such that  $f_1, \dots, f_{p-1}$ form a basis of $V^{\prime,+}$.
Let $T \subset K$ be the maximal torus that acts diagonally with respect to the basis $\{f_i\}$, $T' = T \cap K'$.   
We let $W_G = W(G,T)$, $W_{G'} = W(G',T')$ denote the (absolute) Weyl groups.

Complexified Lie algebras are denoted by gothic letters $\fg$, $\fk$, $\ft$, etc.  
We let $e_i \in \ft^*$ be the character of $\ft$ acting on the vector $f_i$, $i = 1, \dots, p+q$.  
The basis of $\mathfrak{t}'^*$ is given by
$\{e_1, \dots, \widehat{e_p}, \dots, e_{p+q}\}$.

We identify $W_G$ with ${\mathfrak{S}}_{p+q}$
 and ${\mathfrak{t}}^{\ast}$ 
 with ${\mathbb{C}}^{p+q}$
 by letting $w \in {\mathfrak{S}}_{p+q}$ act on $\nu \in {\mathbb{C}}^{p+q}$
 by 
\[
  w\nu=(\nu_{w^{-1}(1)}, \dots, \nu_{w^{-1}(n+1)}).  
\]
This means that the action of the basis $e_j \in \mathfrak{t}^*$ is
 given by $w e_j=e_{w(j)}$
 $(1 \le j \le p+q)$.

\subsubsection{Root Systems and the Harish-Chandra Decomposition}
\label{subsec:HCdeco}
We fix our choice of positive systems $\Phi^{+}$ and $\Phi^{c,+}$  for the pairs $({\mathfrak{g}}, {\mathfrak{t}})$ and $({\mathfrak{k}}, {\mathfrak{t}})$, respectively, as follows:
\begin{alignat*}{2}
\Phi^{c,+} \equiv\, &\Delta^+({\mathfrak{k}})
:=&&\{e_i-e_j:\text{$1 \le i<j \le p$ or $p+1 \le i < j \le p+q$}\}, 
\\
\Phi^{nc,+} \equiv\, &\Delta({\mathfrak{p}}^+)
:=&&\{e_i-e_j:\text{$1 \le i \le p$ and $p+1 \le j \le p+q$}\},
\\
\Phi^+ \equiv\, &\Delta^+({\mathfrak{g}})
:=&&\Phi^{c,+} \cup \Phi^{nc,+}.  
\end{alignat*}

We further define 
\[
\Phi^-:=-\Phi^{+}, \quad 
\Phi^{c,-}:=-\Phi^{c,+}, \quad 
\Phi^{nc,-}\equiv \Delta({\mathfrak{p}}^-):=
-\Phi^{nc,+}.
\]
In the context of the pair $(\fg',\ft')$, we define
 $\Phi'^+ \subset \Phi^{+}$ to be the subset of $e_i-e_j$ where neither $i$ nor $j$ equals $p$,
 and denote the corresponding subsets of $\Phi^{c,+}$ and  $\Phi^{nc,+}$ by $\Phi^{\prime,c,+}$ and  $\Phi^{\prime,nc,+}$, respectively.  
Accordingly, We obtain compatible Harish-Chandra decompositions:
\begin{equation}\label{eqn:HCD} \fg = \fk \oplus \fp^+ \oplus \fp^-~;  \fg' = \fk' \oplus \fp^{\prime,+} \oplus \fp^{\prime,-},
\end{equation}
with
$\fp^{\prime,+} \subset \fp^+; ~ \fp^{\prime,-} \subset \fp^-$.

Let $X = G/K$ denote the Hermitian symmetric space equipped with the complex structure on $G/K$
 via the Borel embedding
\begin{equation}
\label{eqn:Borel_embed}
G/K \hookrightarrow G_{\mathbb{C}}/{P_{\mathbb{C}}^-},
\end{equation}
 where ${P_{\mathbb{C}}^-} \subset G_{\mathbb{C}}$ is 
 a parabolic subgroup with complex Lie algebra
\begin{equation}
\label{eqn:GK_embed}
{\mathfrak{P}}^-:={\mathfrak{k}}+{\mathfrak{p}}^-.  
\end{equation}
Then the holomorphic tangent space
$T_{X,o}$ (respectively, the anti-holomorphic tangent space $T_{X,o}^-$ at the base point $o=eK$)
is identified with $\mathfrak{p}^+$
(respectively $\mathfrak{p}^-$.

The same setup applies to  $X' = G'/K'$.
In particular, the natural embedding $X' \hookrightarrow X$ is holomorphic.

The Harish-Chandra decomposition~\eqref{eqn:HCD} can be also interpreted as follows: if $h \in X$ is a homomorphism
$$h:  \C^\times \ra K_\C$$ 
with values in the center of $K_\C$, then the adjoint action of $h(z)$ induces a Hodge decomposition of $\fg$.  By Deligne's conventions, 
$\fp^- = \fg^{1,-1}$ is the eigenspace for 
$$z \mapsto z^{-1}\cdot\bar{z},$$
while $\fp^+ = \fg^{-1,1}$ is the eigenspace for
$$z \mapsto z\cdot\bar{z}^{-1}.$$

\subsubsection{Root Systems and the Shimura Data}\label{subsec:Shimura_data}

In applications to Shimura varieties we will want to distinguish the real Lie group $G$ from the algebraic group 
$\G = U(V)$ over $\QQ$ such that $G = \G(\RR)$.  Then $(\G,X)$ is given
the structure of Shimura datum as in \cite{H21}.   We recall that there the datum is written $(U(V),Y_V)$, where $Y_V$ has the base point 
\begin{equation}\label{yV} y_{V}(z) = \begin{pmatrix}  (z/\bar{z})I_{p} & 0 \\ 0 & I_{q} \end{pmatrix} \end{equation}
in the standard basis $\{e_i: 1 \leq i \leq n+1=p+q\}$ of $V$.  Note that the eigenvalue of $y_V(z)$ on $e_i-e_j$ is $z/\bar{z}$
if $1 \leq i \leq p$ and $p+1 \leq j \leq n+1$; thus $\fp^+$ as defined above is the sum of the corresponding root spaces.  Similarly, $\fp^-$ is the sum of the root spaces for $e_i-e_j$ with $1 \leq j \leq p$ and $p+1 \leq i \leq n+1$.  The Shimura datum for $U(q,p)$ is obtained from that for $U(p,q)$ by complex conjugation.

\subsection{Weyl Groups, Decoration}
\label{subsec:Weylgp}

Associated with the positive system $\Phi^+=\Delta^+({\mathfrak{g}})$
 and 
$\Phi^{c,+}=\Delta^+({\mathfrak{k}})$, 
 the regular dominant chambers are given in the standard coordinates by 
\begin{align}
    \label{eqn:C+}
   C_+:=&\{x \in {\mathbb{R}}^{p+q}:x_1>\cdots>x_{p+q}\}, 
   \\
\label{eqn:C+K}
%\notag
   C_+^K:=&\{x \in {\mathbb{R}}^{p+q}:
   \text{$x_1>\cdots>x_{p}$ and $x_{p+1}>\cdots>x_{p;q}$}\}.  
\end{align}

We write $\rho_{\mathfrak g}$
for $\Phi^+$
 and $\rho_c$ for $\Phi^{c,+}$.  
One has 
\begin{alignat}{3}
\label{eqn:rho_G}
   \rho_{\mathfrak g}
   :=& \tfrac 1 2 \sum_{\alpha \in \Phi^+}\alpha
   &=&\tfrac 1 2 (p+q-1, \dots, 1-p-q)) &\in& C_+.  
\\
\label{eqn:rho_c}
  \rho_c
  :=& \tfrac 1 2 \sum_{\alpha \in \Phi^{c,+}}\alpha
   &=&\tfrac 1 2(p-1, \dots, 1-p;q-1, \dots, 1-q)
   &\in& C_+^K.  
\end{alignat}

\subsubsection{Harish-Chandra Isomorphism and ${\mathfrak{Z}}({\mathfrak{g}})$-Infinitesimal Characters}
Let ${\mathfrak{Z}}({\mathfrak{g}})$ denote
 the center 
  of the universal enveloping algebra
  $U({\mathfrak{g}})$.  
When ${\mathfrak{g}}={\mathfrak{gl}}_N$, 
 we normalize the Harish-Chandra isomorphism
\[
   \operatorname{Hom}_{{\mathbb{C}}\operatorname{-alg}}({\mathfrak{Z}}({\mathfrak{g l}}_N), {\mathbb{C}}) \simeq {\mathbb{C}}^N/{\mathfrak{S}}_N
\]
 such that the trivial one-dimensional representation 
 ${\bf{1}}$
 of ${\mathfrak{g l}}_{N}$ has the infinitesimal character $\rho_{{\mathfrak{gl}}_N}=\frac 1 2 (N-1, N-3, \dots, 1-N)$.

\subsubsection{Irreducible Finite-Dimensional Representations}
Let $N$ be a positive integer.
We define the subsets
\[
{\mathbb{Z}}_{\operatorname{int}, >}^N
 \subset
 {\mathbb{Z}}_{\operatorname{int}, \mathrm{reg}}^N
  \subset
  {\mathbb{Z}}_{\operatorname{int}}^N
\]
as follows:
\begin{align}
\label{eqn:Zint}
{\mathbb{Z}}_{\operatorname{int}}^N
:=&({\mathbb{Z}}+\tfrac{N+1}2)^N, 
\\
\label{eqn:Z_N_reg}
   {\mathbb{Z}}^N_{\operatorname{int}, \mathrm{reg}}
   :=&
   \{x \in {\mathbb{Z}}_{\operatorname{int}}^N :
x_i \ne x_j \text{ if } i \ne j\},  
\\
%\notag
{\mathbb{Z}}^N_{\operatorname{int}, >}
:=&\{x \in ({\mathbb{Z}}_{\operatorname{int}})^N:
   x_1> \cdots >x_N\}, 
\end{align}

Then $\rho_{{\mathfrak{gl}}_N}
\in \mathbb{Z}^N_{\operatorname{int,>}}$.  
For each $\mu=(\mu_1, \dots, \mu_N) \in {\mathbb{Z}}^N$
 with $\mu_1 \ge \mu_2 \ge \cdots \ge \mu_N$, 
  we denote by $F^G(\mu)$ 
 the irreducible finite-dimensional representation of $G$  with highest weight $\mu $, 
where $G$ is a real form of $GL(N,{\mathbb{C}})$.  
The ${\mathfrak{Z}}({\mathfrak{g}})$-infinitesimal character of $F^G(\mu)$ is given by
$\lambda=\mu+ \rho_{{\mathfrak{gl}}_N}
\bmod{\mathfrak{S}}_N$.  
We note that $\lambda$ belongs to ${\mathbb{Z}}^N_{\operatorname{int}, >}$.

\subsubsection{Kostant Representative $W^{{\mathfrak{k}}}$
 for $W_K \backslash W_G$}
Let $\ell \colon W_G \to {\mathbb{N}}$ denote the standard length function on the Weyl group $W_G$, 
   defined with respect to the generating set of simple reflections.  
Equivalently, 
\begin{equation*}
\ell(w)\equiv \ell^G(w) = \# \Delta^+(w)
\quad\text{for $w \in W_G$, }
\end{equation*}
where we define 
\begin{equation*}
\Delta^+(w):=\Phi^+ \cap w\Phi^-.  
\end{equation*}

Similar notation is applied to $K$:
for $\sigma \in W_K \simeq {\mathfrak{S}}_p \times {\mathfrak{S}}_q$, 
$\ell^K(\sigma)=\sharp \Delta^+(\sigma)$.

The set of complete representatives for the left coset space $W_K \backslash W_G$,
known as the {\it{Kostant representatives}},
is defined by
\begin{align}
\label{eqn:Kostant_W}
W^{\mathfrak{k}}:={}&\{w \in W_G:
                     \text{$w \nu$ is $\Phi^{c,+}$-dominant
whenever $\nu$ is $\Phi^+$-dominant}\}
\\
={}&\{w \in W_G :\Delta^+(w) \subset \Phi^{nc, +}\}.  
\end{align}

Consider the decomposition 
 $\Phi^+=\Phi^{c, +} \cup \Phi^{nc, +}$
 into compact and noncompact roots.
Let $\sigma \in W_K$ and $w \in W^{\mathfrak{k}}$.
We claim that the set $\Delta^+(\sigma w)$ decomposes as a disjoint union:
\begin{equation}
\label{eqn:sigma_w}
   \Delta^+(\sigma w) = \Delta_K(\sigma) \amalg 
   \sigma \Delta^+(w).  
\end{equation}
Let us verify this.

First, suppose $\alpha \in \Phi^{c, +}$. 
Then $\alpha \in \Delta^+(\sigma w)$ if and only if
$\sigma^{-1} \alpha \in w \Phi^{-}$,
which is equivalent to $\langle \sigma^{-1} \alpha, w \nu \rangle <0$ for any $\Phi^{+}$-dominant $\nu$.
Since $w \nu$ is $\Phi^{c, +}$-dominant when $w \in W^{\mathfrak{k}}$, this inequality holds if and only if the compact root $\sigma^{-1} \alpha$ lies in $\Phi^{c, -}$; that is, $\alpha\in \Delta_K(\sigma)$.
Next, suppose $\alpha \in \Phi^{nc,+}$. Then
$\alpha \in \Delta^+(\sigma w)$ if and only if $\sigma^{-1} \alpha \in w \Phi^-$. On the other hand, since $\sigma$ preserves $\Phi^{nc,+}$,we have $\sigma^{-1}\alpha \in \Phi^+$. Thus, $\alpha \in \Delta^+(\sigma w)$ if and only if $\sigma^{-1} \alpha \in \Phi^+ \cap w \Phi^{-}$, which means that $\alpha \in \sigma \Delta^+(w)$.
Thus, \eqref{eqn:sigma_w} is verified.

The decomposition \eqref{eqn:sigma_w} gives rise to a natural bijection
\[
  W_K \times W^{\mathfrak{k}} \overset \sim \longrightarrow 
  W_G, 
\qquad
 (\sigma, w) \mapsto \sigma w,
\]
such that 
\[
   \ell^G(\sigma w) = \ell^K(\sigma) + \ell^G(w).  
\]

In particular, we obtain a natural
 bijection:
\[
W^{\mathfrak{k}} \overset \sim \longrightarrow
 W_K \backslash W_G 
 \simeq 
 ({\mathfrak{S}}_p \times {\mathfrak{S}}_q)\backslash {\mathfrak{S}}_{p+q}.  
\]

We now introduce the finite set
 of \emph{$\Phi^{c,+}$-compatible Weyl chambers}
 by
\begin{align*}
  {\mathcal{C}}_K:=&\{\text{Weyl chambers for $({\mathfrak{g}}, {\mathfrak{t}})$
 for which $\rho_c$ is dominant}\}
 \\
 =&\{\text{Weyl chambers contained in $C_+^K$}\}.  
\end{align*}
For instance, 
the original dominant chamber $C_+$, 
defined in \eqref{eqn:C+}, 
belongs to ${\mathcal{C}}_K.$

Given $w \in W_G\simeq {\mathfrak{S}}_{p+q}$, 
 we define the corresponding Weyl chamber by
\begin{equation}
\label{eqn:w_to_CK}
  C_w :=
  \{x \in {\mathbb{R}}^n:
   x_{w(1)}>x_{w(2)}>\cdots >x_{w(p+q)}
\}
=w C_+.  
\end{equation}

By definition, an element $w \in W_G $ lies in $W^{\mathfrak{k}}$ if and only if 
\[
w^{-1}(i) < w^{-1}(j) \quad \text{for all } 1 \le i < j \le p \text{ or } p+1 \le i < j \le p+q.
\] 
In particular,  $w \in W^{\mathfrak{k}}$
if and only if
$\langle \rho_K, w\nu\rangle >0$
  for any $\nu \in C_+$.

 This yields a natural bijection
\begin{equation}
\label{eqn:WkCK}
W^{\mathfrak{k}} \overset \sim \to
{\mathcal{C}}_K, \quad w \mapsto C_w.  
\end{equation}

\subsubsection{Decoration $\operatorname{Deco}(p,q)$}
We consider another piece of combinatorial data:
\begin{definition}[Decoration]
\label{def:deco}
A {\bf{decoration}} of the signature $(p,q)$ is defined as a string consisting of $p$
  plus signs $(+)$ and $q$ minus signs $(-)$,
  arranged in arbitrary order.  
  The set of all such decorations with signature $(p,q)$ is denoted by $\operatorname{Deco}(p,q)$.  
\end{definition}

\begin{example}
   $\operatorname{Deco}(2,1) 
   =\{(++-), (+-+), (-++)\}$.   
\end{example}

The cardinality of $\operatorname{Deco}(p,q)$ is  $(p+q)!/p! \, q!$,
 which coincides with that of both $W^{\mathfrak{k}}$ and ${\mathcal{C}}_K$.
We now describe explicit bijections between these sets.

Let $\delta \in \operatorname{Deco}(p,q)$ be a decoration.
We associate to $\delta$ a total order of the variables $x_1$, $\dots$, $x_{p+q}$
 as follows.  
Each of the $p$ plus signs $(+)$ in $\delta$
 is replaced, from left to right, 
 with $x_{1}, x_{2}, \dots, x_{p}$, 
 and each of the $q$ minus signs $(-)$ is replaced, 
 also from left to right, 
 with $x_{p+1}, \dots, x_{p+q}$.  
Let $(x_{i_1}, \dots, x_{i_{p+q}})$
 denote the resulting sequence.  
Then the chamber $C(\delta)$ is defined by 
\begin{equation}
\label{eqn:Cdelta}
  C(\delta):=\{x \in {\mathbb{R}}^{p+q}: x_{i_1} > x_{i_{2}} > \cdots > x_{i_{p+q}}\}.
\end{equation}
By construction, 
$C(\delta)$ lies in ${\mathcal{C}}_K$.  
Reversing the above procedure yields
 the inverse map, 
 establishing the following bijection:
\[
\operatorname{Deco}(p,q)
  \overset \sim \longrightarrow
 {\mathcal{C}}_K, \ 
 \delta \mapsto C(\delta).
\]

Composing this map with the bijection $\mathcal{C}_K \overset{\sim}{\longleftarrow} W^{\mathfrak{k}}$ given in \eqref{eqn:WkCK}, we obtain the following result.
 
\begin{lemma}
\label{lem:decobij}
There exists a one-to-one correspondence
 between the following three sets:
\begin{equation}
\label{eqn:decobij}
  \operatorname{Deco}(p,q)
  \,\,\overset \sim \longrightarrow\,\,
  \hphantom{i} {\mathcal{C}}_K
  \,\,\overset \sim \longleftarrow\,\, 
  W^{\mathfrak{k}}, 
  \quad
   \delta \mapsto C(\delta)=C_w\,\,
\raisebox{.1em}{\rotatebox[origin=c]{180}{$\mapsto$}}\,\, w
\end{equation}
where the maps 
$\delta \mapsto C(\delta)$ and
$C_w\,\, 
\raisebox{.1em}{\rotatebox[origin=r]{180}{$\mapsto$}}
\,\, w$
are given by \eqref{eqn:Cdelta} and \eqref{eqn:WkCK}, respectively.  
\end{lemma}

For later argument, 
 we also describe the direct correspondence
\[
W^{\mathfrak{k}} \overset \sim \longrightarrow \operatorname{Deco}(p,q),
\quad 
w \mapsto \delta:=w^{-1}(+\cdots+-\cdots-),
\]
as follows.  
\begin{equation}
\label{eqn:w_delta}
\begin{split}
\text{Given } w \in W^{\mathfrak{k}}, \text{ replace } x_{w(i)} \text{ by } + \text{ if } w(i) \le p, \\
\text{and by } - \text{ if } w(i) > p \text{ in the sequence } (x_{w(1)}, \dots, x_{w(p+q)}).
\end{split}
\end{equation}

\begin{example}
\label{ex:deco_w_lmd}
We illustrate the bijection \eqref{eqn:decobij} between 
$\operatorname{Deco}(p,q)$, $\mathcal{C}_K$, and $W^{\mathfrak{k}}$.

\smallskip

\noindent
\textrm{(1)}\quad
\[
\delta = (++ \cdots + - \cdots -)
\quad \longleftrightarrow \quad
C_+ 
\quad \longleftrightarrow \quad
e.
\]

\smallskip

\noindent
\textrm{(2)}\quad
\[
\delta = (-- \cdots - + \cdots +)
\quad \longleftrightarrow \quad
C(\delta)
\quad \longleftrightarrow \quad
w,
\]
where
\[
C(\delta) = \{ x \in \mathbb{R}^{p+q} : x_{p+1} > \cdots > x_{p+q} > x_1 > \cdots > x_p \},
\]
and
\[
w =
\begin{pmatrix}
1 & 2 & \cdots & q & q+1 & \cdots & p+q \\
p+1 & p+2 & \cdots & p+q & 1 & \cdots & p
\end{pmatrix}
\in W^{\mathfrak{k}}.
\]

\smallskip

\noindent
\textrm{(3)}\quad
Let $(p,q) = (3,2)$, and take $\delta = (-+-++) \in \operatorname{Deco}(3,2)$.  
Then
\[
w = 
\begin{pmatrix}
1 & 2 & 3 & 4 & 5 \\
4 & 1 & 5 & 2 & 3
\end{pmatrix}
\in W^{\mathfrak{k}}, \quad
\delta = w^{-1}(+++--) = (-+-++),
\]
\[
C_w = w C_+ = \{ x \in \mathbb{R}^5 : x_4 > x_1 > x_5 > x_2 > x_3 \}.
\]
\end{example}

For each $\delta \in \operatorname{Deco}(p,q)$, 
 we define
\begin{equation}
\label{eqn:rhodelta}
\rho(\delta):=w \rho_{\mathfrak g} \in C(\delta)\subset{\mathfrak{t}}^{\ast},
\end{equation}
where $w$ is determined by the correspondence \eqref{eqn:w_delta} in Lemma~\ref{lem:decobij}.  
We keep the notation from Lemma \ref{lem:decobij}, 
 and define a Borel subalgebra
\begin{equation}
\label{eqn:Boreldelta}
   {\mathfrak{b}} \equiv {\mathfrak{b}}(\delta)
                  = {\mathfrak{t}}+{\mathfrak{u}}
\end{equation}
 of ${\mathfrak{g}}={\mathfrak{g l}}(p+q,{\mathbb{C}})$
 such that the nilpotent radical 
 ${\mathfrak{u}}$ is given by its weights
\[
  \Delta({\mathfrak{u}})
  :=
  w \Phi^+.  
\]
We note that $\rho(\delta)$ in \eqref{eqn:rhodelta}
 coincides with
$
   \frac12 \sum_{\alpha \in \Delta({\mathfrak{u}})} \alpha.  
$
The Borel subalgebra $\mathfrak{b}(\delta)$ will be used in Section~\ref{subsec:Zuckerman} to describe the underlying $(\mathfrak{g},K)$-module of the discrete series representation $\pi_{\lambda}^\delta$ via cohomological parabolic induction.

\subsection{Parameters of Discrete Series Representations}
\label{subsec:disc_para}

Let $G$ be a real reductive group. We denote by $\widehat{G}$ the unitary dual, i.e., the set of equivalence classes of irreducible unitary representations of $G$. The subset 
$\operatorname{Disc}(G)$ $\subset G$
consists of the discrete series representations---those equivalence classes of irreducible unitary representations that occur as subrepresentations of the regular representation on $L^2(G)$.  
When $G$ is compact, 
 $L^2(K)=L^2(G)$.

Now we describe $\operatorname{Disc}(G)$ for $G=U(p,q)$.

We recall from \eqref{eqn:Z_N_reg}
 that the set 
$
 {\mathbb{Z}}_{\operatorname{int}, \mathrm{reg}}^{p+q}
  \subset
  {\mathbb{Z}}_{\operatorname{int}}^{p+q}
   =
({\mathbb{Z}}+\frac{p+q+1}{2})^{p+q}
$ 
is the set of {\emph{regular integral}} elements. 
We recall from \eqref{eqn:C+K} that $C_+^K$ denotes the strict dominant Weyl chamber corresponding to the positive system $\Phi^{c,+}$ for $(\mathfrak{k}, \mathfrak{t})$.
By Lemma~\ref{lem:decobij}, we have the following decomposition into a disjoint union:
\begin{align}
\label{eqn:Disj}
{\mathbb{Z}}^{p+q}_{\operatorname{int},\operatorname{reg}}
  \cap C_+^K
  =& \{x \in {\mathbb{Z}}^{p+q}_{\operatorname{int},\operatorname{reg}}:
    \text{$x_1 > \cdots >x_p\,$ and $\,x_{p+1}>\cdots>x_{p+q}$}\}
\\
\notag
  =&
   \coprod_{\delta \in \operatorname{Deco}(p,q)} {\mathbb{Z}}_{\operatorname{int}}^{p+q}(\delta),
\end{align}
where, for each decoration $\delta \in \operatorname{Deco}(p,q)$, 
 we define
\begin{equation}
\label{eqn:Z_int_delta}
  {\mathbb{Z}}_{\operatorname{int}}^{p+q}(\delta)
:=
  {\mathbb{Z}}_{\operatorname{int}}^{p+q}
  \cap 
  C(\delta)
  \subset {\mathbb{Z}}_{\mathrm{int, reg}}^{p+q}.  
\end{equation}
Therefore, we obtain the following one-to-one correspondence as a consequence of Lemma~\ref{lem:decobij}.

\begin{proposition}
[Parameter Space of Discrete Series Representation of $U(p,q)$]
\label{prop:discpara}
There exists a one-to-one correspondence
 between the following three sets:
\newline
{\rm{(i)}}\enspace
${\mathbb{Z}}^{p+q}_{\operatorname{int},\operatorname{reg}} \cap C_+^K
 \ni \tlambda$, 
\newline
{\rm{(ii)}}\enspace
${\mathbb{Z}}^{p+q}_{\operatorname{int},>} \times W^{\mathfrak{k}}
\ni (\poslmd, w)$, 
\newline
{\rm{(iii)}}\enspace
${\mathbb{Z}}_{\operatorname{int},>}^{p+q} \times \operatorname{Deco}(p,q)
\ni (\poslmd, \delta)$.

Here, 
 the correspondence
 $\lambda \leftrightarrow (\lambda^+, w)
 \leftrightarrow (\lambda^+, \delta)$
 in (i)---(iii) is given as follows.  
\begin{enumerate}
\item[$\bullet$]
$(\lambda^+, w) \mapsto \lambda$:
$\lambda=w \lambda^+$;

\item[$\bullet$]
$\lambda \mapsto (\lambda^+, w)$:
$w \in W^{\mathfrak{k}}$ is the unique
element of $W_G$
 such that $\lambda \in C_w$, 
 and $\lambda^+ = w^{-1}\lambda$;

\item[$\bullet$]
$(\lambda^+, w) \leftrightarrow (\lambda, \delta)$:
see \eqref{eqn:w_delta}, 
 and Lemma \ref{lem:decobij};

\item[$\bullet$]
$(\poslmd, \delta) \leftrightarrow \tlambda$:
Given $\poslmd \in {\mathbb{Z}}^{p+q}_{\operatorname{int},>}$
 and a signature $\delta \in \operatorname{Deco}(p,q)$, 
 we define
$\tlambda \in {\mathbb{Z}}^{p+q}_{\operatorname{int}, \mathrm{reg}}$
 such that $\tlambda \in C(\delta)$, 
 namely, 
 $\tlambda \in {\mathbb{Z}}^{p+q}_{\operatorname{int}}(\delta)$.  
\end{enumerate}
\end{proposition}

\begin{definition}
\label{def:ref_deco}
We say that $\poslmd$ is the {\bf reference parameter} for $\tlambda$, and $\delta\equiv \delta(\lambda)$ is the {\bf decoration} of $\lambda$.
Thus, the coordinates of the reference parameter $\lambda^+$ satisfy
\[
\lambda_1^+ > \lambda_2^+ > \cdots > \lambda_{p+q}^+.
\]
\end{definition}

By Harish-Chandra's classification of discrete series representations, each discrete series representation is determined by its Harish-Chandra parameter. Using the bijection described in Proposition~\ref{prop:discpara}, we fix a parametrization of $\operatorname{Disc}(U(p, q))$ as follows. 
\begin{definition}[Parametrization of Discrete Series]
\label{def:HC_para}
Let $G=U(p,q)$. Let  $\lambda
\in {\mathbb{Z}}^{p+q}_{\operatorname{int},\operatorname{reg}} \cap C_+^K$, and let
 $\lambda \leftrightarrow (\lambda^+, w)
 \leftrightarrow (\lambda^+, \delta)$
 be the correspondence as given in Proposition~\ref{prop:discpara}.
 \newline
  {\rm{(1). (Parametrization of $\operatorname{Disc}(U(p,q))$)}}.
 We denote by $\pi_\lambda$ the underlying smooth representation of the unique irreducible unitary representation of $G$ with the following properties:
\begin{enumerate}
\item[$\bullet$] it is a discrete series representation of $G$;
\item[$\bullet$] its $\mathfrak{Z}(\mathfrak{g})$-infinitesimal character is $\lambda \bmod \mathfrak{S}_{p+q} = \lambda^+ \bmod \mathfrak{S}_{p+q}$;
\item[$\bullet$] its minimal $K$-type is given by $F^K(\mulmd)$;
\end{enumerate}
where we define
\begin{equation}
\label{eqn:min_K_type_delta}
\mulmd := \lambda + \rho(\delta) - 2\rho_c
= w(\lambda^+ + \rho_{\mathfrak g}) - 2\rho_c.
\end{equation}
The parameter $\lambda$ is called
 the {\bf {Harish-Chandra parameter}}, 
  and $\mu_{\lambda}$ is the {\bf {Blattner parameter}}.  
Following Definition~\ref{def:ref_deco}, we say that $\poslmd$ is the {\bf reference parameter} for the discrete series representation $\pi_\tlambda$, and 
the string of signatures $\delta$ is the {\bf decoration} by $\delta(\lambda)$.

In addition to these, the {\bf{coherent parameter}} $W_\lambda$ associated with $\pi_\lambda$, which plays a role in $(\mathfrak{P}, K)$-cohomologies, will be given in Definition~\ref{def:W_lambda}.
\newline
\rm{(2). (Family associated with decoration $\delta$)}.
In accordance with the decomposition
\eqref{eqn:Disj}, the set of discrete series representations of $U(p,q)$ decomposes as a disjoint union indexed by decorations $\delta \in \operatorname{Deco}(p,q)$:
\begin{equation}
\label{eqn:Disc_delta}
\operatorname{Disc}(U(p,q))
=   \coprod_{\delta \in \operatorname{Deco}(p,q)} 
\operatorname{Disc}(U(p,q))(\delta)
\end{equation}
as given in Proposition~\ref{prop:discpara}.
We define $\operatorname{Disc}(U(p,q))(\delta)$ to be the subset corresponding to the decoration $\delta \in \operatorname{Deco}(p,q)$.

To emphasize the family to which it belongs, we write $\pi_{\lambda}^\delta$ for $\pi_\lambda$, using the decoration $\delta$ as a superscript.

\end{definition}

\begin{example}
\label{ex:23021209}
We illustrate the correspondence
$\lambda \leftrightarrow (\lambda^+, w) \leftrightarrow (\lambda^+, \delta)$
as described in Proposition \ref{prop:discpara} (2)--(4)
  with the following example, see also Example~\ref{ex:deco_w_lmd}:  
\begin{alignat*}{2}
&\tlambda&=&(15,5,0;20,10)
\\
\leftrightarrow
&(\poslmd, w)&=&((20,15,10,5,0), \begin{pmatrix} 1 & 2 & 3 & 4 & 5 \\ 4 & 1 &5 & 2 & 3 \end{pmatrix})
\\
\leftrightarrow
&(\poslmd, \delta)&=&
((20,15,10,5,0), (- + - + +))
\end{alignat*}
The computation of the associated parameters proceeds as follows.
Using our permutation rule
\begin{equation*}
w(\nu_1,\nu_2,\nu_3,\nu_4,\nu_5)
=(\nu_{w^{-1}(1)},\nu_{w^{-1}(2)},\nu_{w^{-1}(3)},\nu_{w^{-1}(4)},\nu_{w^{-1}(5)})
=(\nu_2,\nu_4,\nu_5,\nu_1,\nu_3), 
\end{equation*}
we obtain
\begin{equation*}
\tlambda=
w \poslmd
=(\poslmd_2, \poslmd_4, \poslmd_5, \poslmd_1, \poslmd_3)
=(15, 5, 0, 20, 10).
\end{equation*}
Accordingly, the chamber containing $\lambda$ is given by
\begin{equation*}
 C_w:=w C_+= \{x \in {\mathbb{R}}^5:x_4>x_1>x_5>x_2>x_3\}=C(\delta).  
\end{equation*}
\end{example}

 \begin{remark}
     He's sign conventions
(see Section~\ref{sec:resds})
     are not identical to the conventions  used by Adams--Barbasch--Vogan~\cite{ABV} and in the Atlas software.
 \end{remark}

\subsection{Equivalent Definitions of \texorpdfstring{$\Phi^{nc,+}(\delta)$}{Phi-nc,+(delta)} and a Degree \texorpdfstring{$q$}{q}}
\label{subsec:q}

As recalled in Lemma~\ref{lem:decobij}, 
there is a one-to-one correspondence
\[
\operatorname{Deco}(p,q)\simeq{\mathcal{C}}_K \simeq W^{\mathfrak{k}}.
\]
We now assign to each of these objects a \emph{degree} $q$, following the notation used for the $\bar\partial$-cohomoogy in Narasimhan--Okamoto~\cite{NO70}.    

First, 
we define
\begin{equation}
\label{eqn:q_delta}
q \colon \operatorname{Deco}(p,q) \to {\mathbb{N}}, 
\qquad
\delta \mapsto q(\delta)  
\end{equation}
as the number of pairs $(+,-)$ in $\delta$
 such that $+$ appears to the left of $-$.  
Equivalently, this is the cardinality of the set
\begin{equation}
\label{eq:Phi_nc_delta}
\Phi^{nc,+}(\delta):=
\{e_i-e_j: \text{$+_i$ appears to the left of $-_j$ in $\delta$}\},
\end{equation}
where $+_i$ denotes the $i$-th $+$ and $-_j$ denotes the $j$-th $-$ in $\delta$.

Second, 
for each $w\in W^{\mathfrak{k}}$,
we define
\[
\Phi^{nc,+}(w):=\Phi^{nc,+} \cap w \Phi^+,
\]
\[
q \colon W^{\mathfrak{k}} \to {\mathbb{N}},\quad
  w \mapsto q(w)
  :=\# (\Phi^{nc,+}(w)).  
\]

Third, for any Weyl chamber $C$ in ${\mathcal{C}}_K$, 
 we choose $\lambda \in C$ and
 define
\[\Phi^{n c, +}(\tlambda):=\{\alpha \in \Phi^{nc,+}:
\langle \alpha, \tlambda \rangle >0\}, 
\]
\begin{equation}
\label{eqn:qnc}
  q(\tlambda):=\# \Phi^{n c, +}(\tlambda).  
\end{equation}

We note that $\Phi^{n c, +}(\tlambda)$ is determined
 solely
 by the Weyl chamber $C \in {\mathcal{C}}_K$, 
  which $\lambda$ belongs to, 
  since
\begin{equation}
\label{eqn:ncdelta}
   \Phi^{n c, +}(\tlambda)=
   \{\alpha \in \Phi^{nc,+}:
     \langle \alpha, x \rangle >0\quad
\text{for any $x \in C$}\}.
\end{equation}

\begin{lemma}\label{lem:q}
Let \[
\tlambda \leftrightarrow 
(\poslmd, w) \leftrightarrow
(\poslmd, \delta)
\]
 be the correspondence given as in Proposition \ref{prop:discpara}.
\newline
 (1) 
 $\Phi^{nc,+}(\lambda) =\Phi^{nc,+}(\delta)=\Phi^{nc,+}(w)$. 
\newline
 (2) 
$
  q(\lambda) = q(\delta)=q(w).
$
Moreover, 
 $q(w)=pq-\ell(w)$.  
\end{lemma}
\begin{definition}[Degree $q$]
\label{def:q}
Let $\pi$ be a discrete series representation of $G$, and let $\tlambda \leftrightarrow 
(\poslmd, w) \leftrightarrow (\poslmd, \delta)$ be the associated parameters as given in Definition~\ref{def:HC_para}.
We define the degree $q(\pi)$ of $\pi$ by
\[
q(\pi):= q(\delta) = q(\lambda)=q(w).
\]
\end{definition}
\begin{proof}[Proof of Lemma~\ref{lem:q}]
Suppose that
$w \in W^{\mathfrak{k}}$ is defined by the decoration $\delta \in \operatorname{Deco}(p,q)$
  as in Lemma \ref{lem:decobij}.
  Then the $i$-th $+$ appears in the $w^{-1}(i)$-th position from the left, and the $j$-th $-$ appears in the $w^{-1}(j)$-th position from the left.
Hence, we have
\[
\Phi^{nc,+}(\delta)=
\{e_i-e_j: i\le p <j \text{ and } w^{-1}(i)<w^{-1}(j) \}.
\]

On the other hand, 
 since $\Delta^+(w)  \subset \Phi^{nc,+}$
 for $w \in W^{\mathfrak{k}}$, 
we have the following disjoint union:
\[
\label{eqn:Phi_nc_decom}
    \Phi^{nc,+}
=(\Phi^{nc,+} \cap w \Phi^{-}) 
  \amalg 
  (\Phi^{nc,+} \cap w \Phi^{+})
\\
=\Delta^+(w) \amalg \Phi^{nc,+}(w). 
\]
In view of
$
  w \Phi^-
  =\{e_{w(i)}-e_{w(j)}:1 \le j<i \le p+q\}, 
$
 we have
\begin{align*}
  \Delta^+(w)
=&\{e_{w(i)}-e_{w(j)}:1 \le j <i \le p+q, \quad w(i) \le p < w(j)\}, 
\\
\Phi^{nc,+}(w)
=&\{e_{w(i)}-e_{w(j)}
   :
1 \le i< j \le p+q, \,\,
w(i) \le p < w(j)\}.  
\end{align*}

Hence, we have shown that $\Phi^{nc,+}(\delta)=\Phi^{nc,+}(w)$.

If $\tlambda=w \poslmd$, 
 then it follows from \eqref{eqn:ncdelta} that
\begin{equation*}
\Phi^{n c, +}(\tlambda)
=\{\alpha \in \Phi^{nc, +}:
\langle w^{-1} \alpha, \poslmd \rangle >0\}
=\Phi^{nc,+} \cap
 w  \Phi^+ =\Phi^{nc,+}(w). 
 \end{equation*}
 Thus the first statement is proven.

 As the cardinality of the set $\Phi^{nc,+}(\lambda) =\Phi^{nc,+}(\delta)=\Phi^{nc,+}(w)$,
 we obtain
 $q(\lambda)=q(\delta)=q(w)$.
Finally, it follows from \eqref{eqn:Phi_nc_decom} that
\[
  q(w) 
  =
   \# \Phi^
   {nc, +}-\sharp (\Delta^+(w))
   =pq-\ell(w).  
\] 
\end{proof}

\begin{example}
[See also Example \ref{ex:23021209}]
Let $\delta=(-+-++)\in \operatorname{Deco}(3,2)
$.
The corresponding element
$w \in W^{\mathfrak{k}}$,
 defined in Lemma \ref{lem:decobij},
is given by  $w=\begin{pmatrix} 1 & 2 & 3 & 4 & 5 \\ 4 & 1 & 5 & 2 & 3 \end{pmatrix}$.
Then
\[
\Phi^{nc,+}(\delta) = \Phi^{nc,+}(w) = \{ e_1 - e_5 \}, \quad 
\ell(w) = 5, \quad \text{and} \quad q(w) = q(\delta) = 1.
\]
\end{example}

\subsection{Coherent Parameter}
\label{sebsec:Wlmbd}
Associated to the data 
\[
  \tlambda \leftrightarrow (\poslmd, w)
  \leftrightarrow (\poslmd, \delta)
\]
of the discrete series representation $\pi_\lambda^\delta$,
as described in Definition~\ref{def:HC_para} and in
Proposition \ref{prop:discpara}, 
  there are two specific  $K$-types.  
One is the minimal $K$-type
 with highest weight $\mu_{\lambda}$,
 as defined in \eqref{eqn:min_K_type_delta} .

The other $K$-type is $W_{\tlambda}$, which plays a crucial role in the $(\mathfrak{P}, K)$-cohomologies:
\begin{definition}[Coherent parameter $\Lambda$ and $W_\lambda$]
\label{def:W_lambda}
Let $W_\lambda$ denote the finite-dimensional irreducible representation of $K=U(p) \times U(q)$ with  highest weight 
\[
\Lambda=\tlambda-\rho_{\mathfrak g}=w \poslmd -\rho_{\mathfrak g}
\]
\begin{equation}
\label{eqn:cohlmd}
W_{\tlambda} := F^K(\tlambda-\rho_{\mathfrak g}).  
\end{equation}
For a discrete series representation $\pi=\pi_\lambda$ with Harish-Chandra parameter $\lambda$, we  write $W_\pi$ for $W_\lambda$, and refer to $\Lambda$  as the {\textbf {coherent parameter}} associated with $\pi_{\lambda}$. 
\end{definition}

See Theorem~\ref{thm:elementary_pair_equiv} and \ref{prop:dbc} for the role of $W_\pi$ in the computation of the $(\mathfrak{P}, K)$-cohomology of $\pi$.

\begin{remark}
Note that $\Lambda$ is a dominant weight with respect to the positive system $\Phi^{c,+}$ for $\Delta(\mathfrak{k}, \mathfrak{t})$.  
  We emphasize that, in our convention,
  $\rho_{\mathfrak g}$ is defined as
 $\rho_{\mathfrak g}=\rho_c+\rho({\mathfrak{p}}^+)$,
 as in \eqref{eqn:rho_G},
 rather than $\rho_c +\rho({\mathfrak{p}}^-)$.  
\end{remark}

\begin{remark}
The degree $q(\pi)$ and the coherent parameter $W_\pi$ are sensitive to the choice of the set of noncompact roots $\Phi^{nc,+}=\Delta(\mathfrak{p}^+)$, or, equivalently, to the complex structure of the hermitian symmetric space $G/K$. 
The dependence is summarized in Table~\ref{tab:250503} in Part~\ref{part:4}.
\end{remark}

The following lemma provides a relation 
between the Blattner parameter $\mu_\lambda$ of the minimal $K$-type and the coherent parameter $\Lambda = \lambda -\rho_\mathfrak{g}$ of $W_{\lambda}$.  
\begin{lemma}
\label{lem:24042109}
Let $\mu_\lambda$ denote the Blattner parameter, and let $\Lambda = \lambda - \rho_{\mathfrak g}$ be the coherent parameter of a discrete series representation $\pi$ with Harish-Chandra parameter $\lambda \in \mathbb{Z}_{\mathrm{int, reg}}^{p+q} \cap \mathcal C_K$,
(see Definitions~\ref{def:HC_para} and~\ref{def:W_lambda}).
Let 
\[
\lambda \leftrightarrow (\lambda^+, w) \leftrightarrow (\lambda^+, \delta)
\]
denote the correspondence as in Proposition~\ref{prop:discpara}.
Then we have
\begin{align}
\label{eqn:rho+wrho}
 \rho_{\mathfrak g} +  \rho(\delta)
   =& \sum_{\alpha \in  \Phi^{nc, +}(\delta)} \alpha  + 2 \rho_c, 
\\
\label{eqn:mu_Lambda}
  \mu_{\lambda}
   =& \Lambda + \sum_{\alpha \in  \Phi^{nc, +}(\delta)} \alpha.  
\end{align}
\end{lemma}
\begin{proof}
By definition, we have
\[
\rho_{\mathfrak g}+\rho(\delta)=\rho_{\mathfrak g}+w \rho_{\mathfrak g}
=
\frac{1}{2}\sum_{\alpha \in \Phi^+}  \alpha
+
\frac{1}{2}\sum_{\alpha \in \Phi^+} w \alpha
=\sum_{\alpha \in w \Phi^+ \cap \Phi^+}  \alpha.
\]
The identity \eqref{eqn:rho+wrho} then follows from
\[
w \Phi^+ \cap \Phi^+=(w \Phi^+ \cap \Phi^{nc,+}) \cup (w \Phi^+ \cap \Phi^{c,+})=\Delta^+(w) \cup \Phi^{c,+}.
\]
The identity \eqref{eqn:mu_Lambda} is straightforward from the definitions of $\mu_\lambda$ and $\Lambda$.
\end{proof}

  \begin{example}
[See also Example~\ref{ex:23021209}]
Consider the discrete series representation of \( G = U(3,2) \) with Harish-Chandra parameter \( \tlambda = (15, 5, 0;\, 20,10) \).

As shown in Example~\ref{ex:23021209}, the correspondence
\[
\lambda \leftrightarrow (\lambda^+, w) \leftrightarrow (\lambda^+, \delta)
\]
in Proposition~\ref{prop:discpara} (2)--(4) gives
\[
\lambda^+ = (20,15,10,5,0), \quad
w = \begin{pmatrix} 1 & 2 & 3 & 4 & 5 \\ 4 & 1 & 5 & 2 & 3 \end{pmatrix}, \quad
\delta = (- + - + +).
\]

The discrete series representation $\pi_\lambda
= \pi_{(15,5,0;\, 20,10)}$, has infinitesimal character
\[
(15, 5, 0;\, 20,10)
\bmod \mathcal{S}_5
=
(20,15,10,5,0)
\bmod \mathcal{S}_5,
\]
while the Blattner parameter 
$\mu_\lambda=\tlambda + \rho(\tlambda) - 2\rho_c = (14,4,0;\, 21,11)$ 
and the coherent parameter $\Lambda=\lambda-\rho_{\mathfrak{g}}$ yield 
the following two $K$-types:
\begin{alignat*}{2}
&\text{(minimal $K$-type)} &\quad \mu_\lambda &:= F^{GL_3}(14,4,0) \boxtimes F^{GL_2}(21,11), \\
&\text{(coherent parameter)} &\quad W_{\tlambda} &:= F^{GL_3}(13,4,0) \boxtimes F^{GL_2}(21,12).
\end{alignat*}

The relevant quantities in this computation are:
\begin{align*}
\rho_{\mathfrak g} &= (2,1,0,-1,-2), \\
\rho(\tlambda) &= (1,-1,-2,2,0) = w \rho_{\mathfrak g}, \\
2\rho_c &= (2,0,-2;\, 1,-1).
\end{align*}
\end{example}

\subsection{The Extended \texorpdfstring{$L$}{L}-packet}

The extended $L$-packet---also referred to as the Vogan $L$-packet $\Phi(\poslmd)$ with infinitesimal parameter $\poslmd$---contains $2^{n+1}$ discrete series representations of the $n+1$ groups $U(p,q)$ with $p+q = n+1$. 
This count is reflected in the identity
\[
\sum_{\substack{p+q = n+1}} \# 
\operatorname{Deco}(p, q) 
=
\sum_{\substack{p+q = n+1}} 
\frac{(p+q)!}{p! q!}
= 2^{n+1}
\]

The members of the $L$-packet $\Phi(\poslmd)$ are  indexed by permutations of the coefficients of $\poslmd$, as follows.
Fix the signature $(p,q)$ and let $\Phi_{p,q}(\poslmd)$ denote the permutations 
\begin{equation}\label{mu}
\begin{aligned}
\tlambda & = (\tlambda_1 > \tlambda_2 > \dots > \tlambda_p; \tlambda_1 > \dots > \tlambda_q) \\ & = (\poslmd_{w^{-1}(1)} > \dots > \poslmd_{w^{-1}(p)}; \poslmd_{w^{-1}(p+1)} > \dots > \poslmd_{w^{-1}(p+q)})
\end{aligned}
\end{equation}

The permutation $\poslmd$ of $\tlambda$ with coefficients in decreasing order is called the {\it reference parameter} for $\tlambda$ (Definition~\ref{def:HC_para}),
where $w \in W_G$ has the property that
$w^{-1}(a)<w^{-1}(a+1)$
if $a \neq p$; that is, $w\in W^\mathfrak{k}$ for $(\mathfrak{g}, \mathfrak{k})=(\mathfrak{gl}(p+q,\mathbb C), \mathfrak{gl}(p,\mathbb C)\oplus \mathfrak{gl}(q,\mathbb C))$

On the other hand, $W_G$ acts on $(\pm 1)^{n+1}$ by permutation, and thus on $\Phi(\poslmd)$, preserving each $\Phi(\poslmd)_{p,q}$.

\subsection{Cohomological Parabolic Induction}
\label{subsec:Zuckerman}
This section gives an expression of discrete series representations
in terms of Zuckerman's derived functor modules.

Let ${\mathfrak {q}}={\mathfrak {l}} + {\mathfrak{u}}$
 be a $\theta$-stable parabolic subalgebra.  
 (We shall treat the case 
 that ${\mathfrak {q}}$ is a $\theta$-stable Borel subalgebra.)  
As an algebraic analogue of the Dolbeault cohomology
 with coefficients in a $G$-equivariant holomorphic vector bundle
 over the complex manifold $G/L$, 
Zuckerman introduced a derived functor 
$
   {\mathcal{R}}_{\mathfrak {q}}^j \equiv 
  ({\mathcal{R}}_{\mathfrak {q}}^{\mathfrak{g}})^j
$
 ($j \in {\mathbb{N}}$) as a cohomological parabolic induction.  
We use the normalization so that 
 ${\mathcal {R}}_{{\mathfrak {q}}}^j$ preserves 
 the infinitesimal characters.  
More precisely, 
 we regard ${\mathcal{R}}_{\mathfrak {q}}^j$
 as a covariant functor
 from the category of metaplectic $({\mathfrak {l}}, 
(L \cap K)\!\widetilde{\hphantom{m}})$-modules
 to the category of $({\mathfrak {g}}, K)$-modules
 so that 
 if $\nu \in {\mathfrak {t}}^*/W_L$
 is the ${\mathfrak{Z}}({\mathfrak {l}})$-infinitesimal character
 of an $({\mathfrak{l}}, (L \cap K)\!\widetilde{\hphantom{m}})$-module $V$
 then the ${\mathfrak{Z}}({\mathfrak {g}})$-infinitesimal character of ${\mathcal{R}}_{{\mathfrak {q}}}^{S}(V)$
 equals 
$
   \nu \in {\mathfrak {t}}^*/W_G
$.

\begin{lemma}
\label{lem:23082617}
For each decoration $\delta \in \operatorname{Deco}(p,q)$, we take a Borel subalgebra 
\[
{\mathfrak{b}}\equiv{\mathfrak{b}}(\delta)={\mathfrak{t}} + {\mathfrak{u}}
\]
 of ${\mathfrak{g}}={\mathfrak{g l}}(p+q, {\mathbb{C}})$
 as in \eqref{eqn:Boreldelta}.  
Suppose $\tlambda \in {\mathbb{Z}}_{\operatorname{int}}^{p+q} \cap C(\delta)$.  Then the underlying $({\mathfrak{g}}, K)$-module
 of the discrete series representation $\pi_{\tlambda}^{\delta}$
 is isomorphic to ${\mathcal{R}}_{\mathfrak{b}}^S({\mathbb{C}}_{\tlambda})$, 
 where $S=\frac 1 2 (p(p-1)+q(q-1))$.  
\end{lemma}

We recall the standard argument of spectral sequence, 
 which is useful
 in computing the translation of discrete series representations
 in terms of Zuckerman's derived functor modules as follows:
 
\begin{lemma}
\label{lem:Rspecseq}
Let 
$
   {\mathfrak{b}}\equiv{\mathfrak{b}}(\delta)
  ={\mathfrak{t}} + {\mathfrak{u}}
$
 be a Borel subalgebra associated with the signature $\delta$
 as in Lemma \ref{lem:23082617}.  
Let $F$ be a finite-dimensional representation of $G$, 
 and 
\[
  0 = F_0 \subset F_1 \subset \cdots \subset F_m = F
\]
 be a ${\mathfrak{b}}$-stable filtration
 such that ${\mathfrak{u}}$
 acts trivially on $F_i/F_{i-1}$.  
Then there is a natural spectral sequence
\begin{equation}
\label{eqn:Rbspec}
   {\mathcal{R}}_{\mathfrak{b}}^j({\mathbb{C}}_{\lambda} \otimes F_i/F_{i-1})
   \Rightarrow
   {\mathcal{R}}_{\mathfrak{b}}^j({\mathbb{C}}_{\lambda}) \otimes F
\end{equation}
of $({\mathfrak{g}}, K)$-modules.  
\end{lemma}

\subsection{Convention}

We recall from Section~\ref{subsec:ugps} that the complex structure on $G/K$ is defined via the Borel embedding
$
G/K \hookrightarrow G_{\mathbb{C}} / P_{\mathbb{C}}^-,
$
see \eqref{eqn:Borel_embed}.
Accordingly, 
 the holomorphic tangent (respectively, cotangent) bundle of $G/K$ is identified with the associated bundle $G \times_K \mathfrak{p}^+$ (respectively, $G \times_K \mathfrak{p}^-$). 
 
An irreducible unitary representation realized in the space of holomorphic $L^2$-sections of a Hermitian holomorphic vector bundle over this complex manifold $G/K$ is called a \emph{holomorphic discrete series representation}.
The distinction between holomorphic and antiholomorphic structures is encoded into our conventions as the following example.
Furthermore, we provide the value of $q(\delta)=q(\lambda) \in \{0, 1, \dots, pq\}$ for $\pi_{\tlambda} = \pi_{\poslmd}^\delta \in \operatorname{Disc}(U(p,q))$, see Lemma~\ref{lem:q} for equivalent definitions.

 \begin{conv}\label{holo}  \begin{enumerate}
 \item The decoration $$\delta=+~ \dotsb ~+ ~ - \dotsb ~ -$$ denotes an antiholomorphic representation with $q(\delta) = pq$.
  \item The decoration $$\delta=-~ \dotsb ~- ~ + \dotsb ~ +$$ denotes a holomorphic representation with $q(\delta) = 0$.
  \item The decoration $$\delta= + \dotsb +$$ is a decoration of an irreducible finite-dimensional representation of $U(p,0)$ with $q(\delta)=0$.
  \item  When $p = q$, the decorations $$\delta=+~ -~ +~ - \dotsb ~+ ~ - ~\dotsb ~+ ~ - \text{ and }
  -~ +~ -~ +~  \dotsb ~+ ~ - ~\dotsb ~- ~ +$$ denote generic discrete series,
   with $q(\delta) = \frac{p(p+1)}{2}$
  and $q(\mu) = \frac{p(p-1)}{2}$, respectively.
\end{enumerate}
\end{conv}

\section{Preliminaries: Restriction of discrete series}\label{sec:resds}

In this section, we review He's combinatorial description of the GGP \emph{interleaving relations} \cite{He}, which characterize the non-vanishing of the spaces of symmetry breaking operators arising from the restriction of discrete series representations 
in the setting $(G, G')=(U(p,q), U(p-1,q))$.
\subsection{H.He's Description}
 \label{subsec:He_description}
  
Every discrete series representation $\pi=\pi_{\lambda}$ of $G = U(p,q)$ 
is indexed by its Harish-Chandra parameter $\lambda\in {\Z}^{p+q}_{\operatorname{int},\operatorname{reg}} \cap C_+^K$,
as stated in Definition~\ref{def:HC_para}.
As described in Proposition~\ref{prop:discpara}, it may also be parametrized by a pair consisting of a reference parameter $\poslmd$ 
and a decoration $\delta$.

Similarly, a discrete series representation $\pi'$
of $G'=U(p-1,q)$ is determined by its Harish-Chandra parameter $\lambda'\in {\Z}^{p+q-1}_{\operatorname{int},\operatorname{reg}} \cap C_+^{K'}$, or equivalently by a pair $(\lambda'^+,\delta')$.

Recall that the decorations $\delta$ and $\delta'$ are $(p+q)$-tuples and $(p+q-1)$-tuples of signs, respectively, with exactly $p$ (respectively, $p-1$) positive entries and $q$ negative entries.
Following the convention of He~\cite{He}, we use the symbols $+$ and $-$ (respectively, $\op$ and $\oi$)  
to denote the entries of the decoration $\delta$ for $G$ (respectively, the decoration $\delta'$ for $G'$).

\medskip

Since $p+q$ and $(p-1)+q$ have opposite parities, no component $\lambda_i$ can be equal to any $\lambda'_j$.
Let the components $\lambda_1, \dots, \lambda_{p+q}$ of $\lambda$ and $\lambda'_1, \dots, \lambda'_{p-1},\widehat{\lambda'_p}, \lambda'_{p+1}, \lambda'_{p+q}$ of $\lambda'$ be rearranged in descending order.

Assign symbols to the entries as follows:
\begin{itemize}
  \item the entries corresponding to $\lambda_1, \dots, \lambda_p$ are replaced by $+$,
  \item those corresponding to $\lambda_{p+1}, \dots, \lambda_{p+q}$ are replaced by $-$,
  \item those corresponding to $\lambda'_1, \dots, \lambda'_{p-1}$ are replaced by $\op$,
  \item those corresponding to $\lambda'_{p+1}, \dots, \lambda'_{p+q}$ are replaced by $\oi$.
\end{itemize}
 The resulting sequence of symbols is denoted by $[\delta\delta'] \in \{+, -,\op, \oi\}^{2p+2q-1}$.
  
\begin{definition}[Signature ${[\delta \delta']}$ associated with $(\pi, \pi')$]
\label{def:string_dd}
Let $(\pi, \pi')$ be a pair of discrete series representations of $G=U(p,q)$ and $G'=U(p-1,q)$,
and let $(\lambda, \lambda')\in (
 {\Z}^{p+q}_{\operatorname{int},\operatorname{reg}} \cap C_+^K, {\Z}^{p+q-1}_{\operatorname{int},\operatorname{reg}} \cap C_+^{K'}$)
 be their Harish-Chandra parameters.
We define
\[
[\delta\delta'] \in \{+, -, \op,\oi\}^{2p+2q-1}
\] to be the sequence of $(2n+1)$ symbols, obtained as described above,
and refer to it as the {\bf{associated signature}} of the pair $(\pi, \pi')$ or $(\lambda, \lambda')$.

\end{definition}

We note that $\delta$ and $\delta'$ are the decorations associated with $\pi$ and $\pi'$, respectively;  
however, the signature $[\delta \delta']$ is not uniquely determined by the pair $(\delta, \delta')$.

\begin{example}
Let $(G, G') = (U(2,1), U(1,1))$.  
In the following two cases, the decorations are the same: $\delta = +-+$ and $\delta' = \op \oi$,  
yet the associated signature $[\delta \delta']$ differs.

\begin{itemize}
  \item[(1)] Let $\lambda = (10, 1; 2)$ and $\lambda' = \left(\frac{9}{2}, \frac{7}{2}\right)$.  
  The associated decorations are $\delta = +-+$ and $\delta' = \op \oi$.  
  The signature $[\delta \delta']$ is $+\, \op\, \oi\, -\, +$, since  
  $\lambda_1 > \lambda_1' > \lambda_3' > \lambda_3 > \lambda_2$.
  
  \item[(2)] Let $\lambda = (10, 8; 9)$ and $\lambda' = \left(\frac{3}{2}, \frac{1}{2}\right)$.  
  Again, the decorations are $\delta = +-+$ and $\delta' = \op \oi$.  
  However, the signature $[\delta \delta']$ is $+\, -\, +\, \op\, \oi$, since  
  $\lambda_1 > \lambda_3 > \lambda_2 > \lambda'_1 > \lambda'_2$.
\end{itemize}
\end{example}

In what follows, we use the same letter $\pi$ to denote irreducible Fr\'echet 
smooth representations of moderate growth,  obtained by considering smooth vectors in the unitary representations $\pi$ of $G$. 
Similarly, we denote by $\pi'$ the corresponding smooth representations of the subgroup $G'$.

Let
$\operatorname{Hom}_{G'}(\pi|_{G'},\pi')$ denote
the space of continuous $G'$-equivariant homomorphisms with respect to the Fr\'echet topology, as in \cite{GGP}.
Equivalently, we consider continuous $G'$-invariant bilinear forms
$$\pi \otimes \overline{\pi'} \ra \CC,$$
where $\overline{\pi'}$ denotes the complex conjugate representation of $\pi'$.

 It has been shown in \cite{SZ, JSZ} that $\dim \operatorname{Hom}_{G'}(\pi|_{G'},\pi') \leq 1$.   

\begin{definition}[GGP Interleaving Relation]\label{def:interlace}
Let $\pi$ and $\pi'$ be discrete series representations of $G$ and $G'$, respectively,
and let $[\delta \delta']$ denote the associated signature (Definition~\ref{def:string_dd}).  
We say that the pair $(\pi,\pi')$ satisfies the {\bf{GGP interleaving relation}}
if every adjacent pair of symbols in the associated signature
$[\delta\delta']$ belongs to the following list:
\begin{equation}
\label{eqn:GGP_pattern}
(\op +),\quad (+ \op),\quad (- \oi),\quad (\oi -),\quad 
(+ -),\quad (- +),\quad (\op \oi),\quad (\oi \op).
\end{equation}
We say that the {\it signature } $[\delta \delta']$ is a {\bf {GGP interleaving pattern}} if all adjacent pairs of symbols in $[\delta\delta']$ belong to this list.
\end{definition}
Here, instead of using the standard term \emph{interlacing pattern}, we deliberately adopt the term \emph{interleaving pattern} to accommodate a broader framework. See also \cite{KS}.
\begin{example}
Consider the following examples of signatures $[\delta \delta']$:
\begin{itemize}
  \item[(1)] The signature $[\delta \delta'] = +\, \op\, \oi\, -\, +$ is a {\bf GGP interleaving pattern},  
  since each adjacent pair belongs to the list \eqref{eqn:GGP_pattern}.
 
  \item[(2)] The signature $[\delta \delta'] = +\, -\, \op\, +\, \oi$ is {\bf not} a GGP interleaving pattern, since
   the pair $( -\, \op )$ is not included in the allowed list.
\end{itemize}
\end{example}

He's interpretation of the Gan--Gross--Prasad conjecture for discrete series $L$-packets is the following.

\begin{theorem}[He]\label{thm:GGP_He}   The space $\operatorname{Hom}_{G'}(\pi|_{G'},\pi')$ is of dimension $1$ if and only if 
$(\pi,\pi')$ satisfies the GGP interleaving relation.
\end{theorem}

\begin{remark}
In Corollary~\ref{cor:25061806}, 
We provide a simple and novel proof that the multiplicity
\[
[\pi|_{G'} : \pi'] = 1,
\]
assuming that the associated signature $[\delta \delta']$ is coherent (Definition~\ref{def:coh_deco}), 
which is the case we focus on in Part~\ref{part:translation} and Part~\ref{global}.
\end{remark}

It is convenient to introduce the following notation.
\begin{definition}[Support of $\pi|_{G'}$ at Decoration $\delta'$]
\label{def:Supp_delta}
Let $\pi$ be a discrete series representation of $G=U(p,q)$,
and let $G'=U(p-1,q)$.
We define a subset of $\widehat{G'}$ by
\[
\operatorname{Supp}(\pi|_{G'})
:=\left\{\pi'\in \operatorname{Disc}(G'):\operatorname{Hom}_{G'}(\pi|_{G'},\pi')\ne \{0\}\right\},
\]
where morphisms are taken in the category of smooth Fr\'echet representations of moderate growth.
For each decoration $\delta'\in \operatorname{Deco}(p-1,q)$, we define
\[
\operatorname{Supp}(\pi|_{G'})(\delta')
\]
to be the subset corresponding to the decomposition \eqref{eqn:Disc_delta}.
\end{definition}

See Section~\ref{subsec:Part4_Branching} in Part~\ref{part:4} for some examples.

R. Beuzart-Plessis proved the multiplicity one part of the GGP conjecture for all tempered $L$-packets of unitary groups by a variant of the proof for
non-archimedean local fields.  H. He's proof, which is specific to discrete series representations, confirms the prediction of \cite{GGP} that identifies the member of the packet that carries a $G'$-invariant linear form; it relies on the combinatorics of the theta correspondence.

\subsection{Examples}

\subsubsection{Compact unitary groups}\label{cpt}
When $q = 0$, all signs are $+$, and the only signature of the form $[\delta\delta']$ that satisfies the condition in Theorem~\ref{thm:GGP_He} is
\begin{equation}\label{classical}
+ ~\op ~ + ~\op~ \dots ~ \op ~ +.
\end{equation}

The interleaving relation associated with \eqref{classical} corresponds to the interlacing
\[
\poslmd_1 > \poslmdp_1 > \poslmd_2 > \dots > \poslmdp_{p-1} > \poslmd_p,
\]
which precisely matches Weyl's branching rule for the restriction of representations from the compact group $U(p)$ to $U(p-1)$.

\subsubsection{Holomorphic Discrete Series Representations}\label{holorep}
For $p, q\ge 1$, a discrete series representation $\pi$ of $G=U(p,q)$  decomposes \emph{discretely} upon restriction to the subgroup $G'=U(p-1,q)$ if and only if $\pi$ is either  holomorphic or antiholomorphic (\cite{K98}).

As noted in Convention~\ref{holo}, $\pi_{\tlambda}$ is a holomorphic discrete series representation if and only if its decoration begins with
$q$ minus signs followed by $p$ plus signs.  The representation is described by the reference parameter
$$\poslmd = (\lambda_{p+1}, \lambda_{p+2}, \dots,   \lambda_{p+q}, \lambda_1, \lambda_2, \dots \lambda_p),$$
with the associated decoration $[\delta]=-  \dotsb - + \dots +$.
It is easy to see that the only GGP interleaving pattern $[\delta\delta']$ compatible with this decoration is
\begin{equation}\label{holodeco}
\oi ~ - ~ \oi ~\dotsb ~ - ~ \oi ~ - ~ + ~ \op ~ + ~ \op ~ +~ \dotsb ~\op ~ +
\end{equation}
In particular, the corresponding decoration $[\delta']=\oi \dotsb \oi \op \dotsb \op$ indicates that $\operatorname{Supp}(\pi|_\lambda)$ consists
entirely of holomorphic discrete series representations.   

In this case, the explicit (discretely decomposable) branching laws can be found in \cite[Theorem~8.10]{K08}, and the decomposition is described by the \emph{interleaving condition} as follows:
\begin{equation}
\label{eqn:holo_branch}
\lambda_{p+1}'>\lambda_{p+1}>\dotsb>\lambda_{p+q}'>
\lambda_{p+q}>\lambda_1>\lambda_1'>\dotsb>\lambda_{p}'>\lambda_{p},
\end{equation}
which is consistent with
the inequalities implicit in \eqref{holodeco}.

The restriction of an antiholomorphic discrete series representation $\pi$ to $G'$ is also a direct sum of antiholomorphic representations.
This is natural since the complex conjugation that acts on the hermitian symmetric spaces associated with $G$ and $G'$ exchanges
holomorphic and antiholomorphic representations.  The relevant interleaving pattern $[\delta\delta']$ is obtained by reflecting \eqref{holodeco} horizontally:
$$+ ~ \op ~ + ~\dotsb ~ \op ~ + ~ - ~ \oi ~\dotsb ~ - ~\oi.$$

\section{Coherent Pairs and Elementary Pairs}\label{sec:coherent}

We continue our study of pairs of discrete series representations $\pi$ of $G = U(p,q)$ and $\pi'$ of $G' = U(p-1,q)$. Among the pairs $(\pi, \pi')$ for which $\operatorname{Hom}_{G'}(\pi|_{G'}, \pi') \neq 0$, we introduce two distinguished classes,  \emph{elementary pairs} (Definition~\ref{def:elementary}) and
\emph{coherent pairs} (Definition~\ref{def:coherent_pair}). 

The interleaving patterns associated with coherent and elementary pairs impose stronger conditions than those given by the GGP interleaving relation (Definition~\ref{def:interlace}); we refer to this refined condition as \emph{coherent signature} (Definition~\ref{def:coh_deco}). We provide six combinatorially equivalent formulations of the coherent signature in Proposition~\ref{prop:23022216}. Building on this framework, we present four equivalent characterizations of \emph{elementary pairs} $(\pi, \pi')$ in Theorem~\ref{thm:elementary_pair_equiv}.
These results set the stage for Theorems~\ref{thm:23022420b} and~\ref{thm:230816}, which describe how every coherent pair $(\pi, \pi')$ is related to a corresponding elementary pair that exhibits more favorable properties, one of which includes the non-vanishing of the natural morphisms between $(\fP,K)$-cohomologies
(Proposition~\ref{prop:23022216}).

In Section~\ref{sec:rest_cohomology} of Part~\ref{global}, we will discuss morphisms between their $(\mathfrak{P}, K)$-cohomologies, focusing on non-elementary coherent pairs.

\subsection{Various Conditions for Signatures}
\label{subsec:various_dd}

In this section, we describe several combinatorial conditions on a pair of decorations $\delta$ and $\delta'$, as well as their interleaving patterns $[\delta \delta']$. These conditions are formulated combinatorially in Definition~\ref{def:23022211}.

The relationships among these conditions are clarified in Proposition~\ref{prop:23022216} in the next subsection, which highlight the structural core of this section. This proposition establishes six mutually equivalent conditions on the signature $[\delta \delta']$, 
referred to as \emph{coherent signature} (Definition~\ref{def:coh_deco}).
These conditions provide a characterization of the interleaving pattern for a pair of discrete series representations of $G = U(p, q)$ and $G' = U(p-1, q)$, when the pair $(\pi, \pi')$ forms an elementary pair (Definition~\ref{def:elementary}), or more generally, a coherent pair (Definition~\ref{def:coherent_pair}).

\medskip

Suppose that $(\delta, \delta') \in (\operatorname{Deco}(p,q), \operatorname{Deco}(p-1,q))$
 satisfies the GGP interleaving relations.

We recall the following from earlier sections:
from Lemma~\ref{lem:q}, the equivalent definitions of degree $q(\delta)$ and the set of non-compact roots $\Phi^{nc,+}(\delta)$ associated with a decoration $\delta \in \operatorname{Deco}(p,q)$;
from \eqref{eqn:min_K_type_delta}, the minimal $K$-type parameter $\mu_\lambda$;
and from Definition~\ref{def:W_lambda}, the coherent parameter $W_\lambda$.

.
\begin{definition}
\label{def:23022211}
Suppose that a signature $[\delta\delta']$ satisfies the GGP interleaving pattern (Definition~\ref{def:interlace}), 
 where $\delta \in \operatorname{Deco}(p, q)$, 
$\delta' \in \operatorname{Deco}(p-1,q)$.  
We consider the following set of conditions.
Note that the first three conditions---$(q)$, 
 $(\delta+)$, and $(nc)$---pertain to
 the pair of decorations $\delta$ and $\delta'$, 
 whereas the last three conditions---$(\mu)$, $(W)$, and $(-\oi -)$---concern the interleaving pattern $[\delta\delta']$.

\begin{enumerate}
\item[$(q)$]
$q(\delta)=q(\delta')$.  

\item[$(\delta+)$]
$\delta$ is obtained by converting 
 every $\op$ to $+$
 and every $\oi$ to $-$, 
 and by adding $+$
 at the rightmost.  
 
 \item[$(nc)$]
$\Phi^{nc,+}(\delta)=\Phi^{nc,+}(\delta')$,  
(cf.\ \cite[(4.9)]{H14}).  

\item[$(\mu)$]
The interleaving pattern $[\delta\delta']$ arises from some $(\tlambda, \tlambda')
\in \left(\mathbb{Z}^{p+q}_{\mathrm{int,reg}}\cap C_+^K\right) \times  \left(\mathbb{Z}^{p+q-1}_{\mathrm{int,reg}}\cap C_+^{K'}\right)$
(see Proposition~\ref{prop:discpara}) that
 satisfies the GGP interleaving condition
 and also fulfills $[\mulmd:\mulmdp]\ne 0$;
 that is, $\operatorname{Hom}_{K'}({\mu_\lambda|}_{K'},\mu_{\lambda'}) \neq 0$.  
 
\item[$(W)$]
The interleaving pattern $[\delta\delta']$ arises from some $(\tlambda, \tlambda')
\in \left(\mathbb{Z}^{p+q}_{\mathrm{int,reg}}\cap C_+^K\right) \times  \left(\mathbb{Z}^{p+q-1}_{\mathrm{int,reg}}\cap C_+^{K'}\right)$
(see Proposition~\ref{prop:discpara}) that
 satisfies the GGP interleaving condition
 and also fulfills
$[W_{\tlambda}:W_{\tlambda'}]\ne 0$;
  that is, $\operatorname{Hom}_{K'}({W_\lambda|}_{K'},\Wlmdp) \neq 0$.  
 \item[{\rm{($-\oi -$)}}]
 The interleaving pattern $[\delta\delta']$
satisfies the following two conditions:
\begin{enumerate}
\item  For $1 \leq i \leq p-1$, the  symbol $\op_i$ appears to the {\bf right} of $+_i$ and  to the {\bf left} of $+_{i+1}$. 
This is expressed as 
$$+_i > \op_i > +_{i+1}.$$
\item  For $1 \leq j \leq q$, the  symbol $\oi_j$ appears to the {\bf left} of $-_j$ and to the {\bf right} of $-_{j-1}$.  
This is expressed as
$$-_{j-1} > \oi_j > -_j.$$
In the case $j=1$, we disregard $-_{j-1}$, and consider only the condition ${\oi}_{j} > -_j$.
\end{enumerate}
   
\end{enumerate}
\end{definition}

The following equivalence and implication are immediate from Lemma~\ref{lem:q} 
on equivalent definitions of the set $\Phi^{nc,+}(\delta)$ and the degree $q(\delta)$:
\begin{equation}
\label{eqn:nc_q}
    (\delta +) \Longleftrightarrow (nc) \Longrightarrow (q).
\end{equation}

We also consider reduction procedures for interleaving patterns of the form $[\delta \delta']$,
which are used in conditions (iv) and (v) of Proposition~\ref{prop:23022216} below.

\begin{definition}[Reduction]
\label{def:reduction}
We introduce the following reduction procedures:  
\begin{equation}
\label{eqn:reduction}
\text{
Delete pairs $+\op$, $\oi -$ from the left, 
 or $\op +$ from the right.  
}
\end{equation}
We also consider a weaker version of this reduction:  
\begin{equation}
\label{eqn:wreduction}
\text{
Delete pairs $+\op$, $\oi -$ from the left.  
}
\end{equation}
\end{definition}

\subsection{Six Equivalent Conditions for Coherent Signatures}

We examine the relationships among the definitions involving signature $[\delta \delta']\in \{+, -, \op, \oi\}^{2p+2q-1}$, where $\delta \in \operatorname{Deco}(p,q)$ and $\delta' \in \operatorname{Deco}(p-1,q)$, as introduced in Definition~\ref{def:23022211}. We then establish the equivalence of six associated properties, stated as Proposition~\ref{prop:23022216}, which constitutes the foundation of the main results of this section.

\begin{proposition}[Coherent Signature, see Definition~\ref{def:coh_deco} below]
\label{prop:23022216}
Let $[\delta \delta']$ be a signature satisfying the GGP interleaving pattern
(Definition~\ref{def:interlace}),
 where $\delta \in \operatorname{Deco}(p,q)$, 
 $\delta' \in \operatorname{Deco}(p-1,q)$. 

\noindent
\textrm{(1)}
The following six conditions on $[\delta\delta']$ are then equivalent. 
\begin{enumerate}
\item[{\rm{(i)}}]
Condition $(W)$ is satisfied by the signature $[\delta \delta']$.  

\item[{\rm{(ii)}}]
The signature $[\delta \delta']$ reduces to a single \lq\lq{$+$}\rq\rq\
 under the reduction process described in
 \eqref{eqn:reduction}.  

\item[{\rm{(iii)}}]
The signature $[\delta \delta']$ reduces to a single \lq\lq{$+$}\rq\rq\
 under the weaker reduction process described in
 \eqref{eqn:wreduction}.  

\item[{\rm{(iv)}}]
Condition ${\rm (-\oi -)}$ is satisfied.

\item[{\rm{(v)}}]
Condition ${\rm{(-\oi -)}}$ is satisfied, 
and $+$ appears as the rightmost symbol in $[\delta\delta']$.  

\item[{\rm{(vi)}}]
Both conditions (nc) and ($\mu$) are satisfied.  
\end{enumerate}

\textrm{(2)}
If any (and thus all) of these equivalent conditions hold on the signature $[\delta\delta']$, then the interleaving pattern $[\delta \delta']$ satisfies condition $(\mu)$, and the pair $(\delta, \delta')$ satisfies conditions $(q)$ and $(\delta+)$, as given in Definition~\ref{def:23022211}.

\end{proposition}

\medskip
We provide a proof of Proposition~\ref{prop:23022216} in Section~\ref{subsec:pf_prop:23022216}.

Before presenting the proof, we clarify the mutual relationships among the six conditions formulated in the proposition.
For subtle examples below, see Table~
\ref{tab:U21} in Part~\ref{part:4}.
\begin{remark}
\label{rem:equiv_coherent_sign}

 \par\noindent
{\rm{(1)}}\enspace 
In contrast to the obvious implication $(\delta +) \Rightarrow (q)$ stated in \eqref{eqn:nc_q},
 the converse implication $(q) \Rightarrow (\delta+)$ 
is not true.  
A counterexample is 
\[  [\delta \delta']=+-+\op\oi\]
 for $U(2,1) \downarrow U(1,1)$, 
 for which $q(\delta)=q(\delta')=1$, 
 whereas the pair $(\delta, \delta')=(+-+,\op\oi)$ does not satisfy the condition $(\delta+)$.  
 In particular, $(q) \Rightarrow (W)$ does not hold in general, since $(W) \rightarrow (\delta +)$ as stated in Proposition~\ref{prop:23022216}.
 See Case III in Table~\ref{tab:U21}.

{\rm{(2)}}\enspace
In contrast to the implication $(W) \Rightarrow (nc)$ stated in (i) and~(vi) of Proposition~\ref{prop:23022216}, the reverse implication $(nc) \Rightarrow (W)$ does not hold in general.  
For instance, 
 both $[+-+\op \oi]$ and $[+ \op \oi -+]$ satisfy
 $(nc)$ for $(G, G')=(U(2,1), U(1,1))$, 
but only $[+\op\oi-+]$ satisfies $(W)$.

\par\noindent
{\rm{(3)}}\enspace
The implication $(\delta +) \Rightarrow (W)$ is not valid.  
This claim is immediate from (1) and 
the equivalence \eqref{eqn:nc_q} between $(\delta +)$ and $(nc)$.

\par\noindent
{\rm{(4)}}\enspace
 In contrast to the implication $(W) \Rightarrow (\mu)$ stated in (i) and (vi) of Proposition~\ref{prop:23022216}, the reverse implication $(\mu) \Rightarrow (W)$ does not generally hold.

As a counterexample, consider the interleaving pattern $[\delta \delta'] := + - \oi \op +$, which corresponds to the restriction $U(2,1) \downarrow U(1,1)$, (see Case IV in Table~\ref{tab:U21}). 
This interleaving pattern satisfies condition $(\mu)$ but not $(q)$.

Indeed, suppose $(\lambda_1, \lambda_2, \lambda_3) \in \mathbb{Z}^3$ satisfies $\lambda_1 > \lambda_3 > \lambda_2 + 1$, and choose $(\lambda_1', \lambda_3') \in (\mathbb{Z}+\tfrac{1}{2})^2$ such that $\lambda_3 - \tfrac{1}{2} > \lambda_1' > \lambda_2$ and define $\lambda_3' := \lambda_3 - \tfrac{1}{2}$. Then we have $[\mulmd : \mulmdp] = 1$, while $q(+ - +) = 1$ and $q(\oi \op) = 0$.

Thus, by the implication $(W) \Rightarrow (nc) \Rightarrow (q)$ stated in Proposition~\ref{prop:23022216} and equation~\eqref{eqn:nc_q}, it follows that $(\mu)$ does not imply $(W)$ in this case.

 \par\noindent
{\rm{(5)}}\enspace
Condition (vi) in Proposition~\ref{prop:23022216} clarifies that neither $(nc)$ nor $(\mu)$ alone is sufficient to ensure the condition on coherent parameters $(W)$, as demonstrated their counterexamples in (2) and~(4). However, taken together, they are complementary and jointly equivalent to condition $(W)$.

\par\noindent
{\rm{(6)}}\enspace
The equivalence (vi) $\Leftrightarrow$ (i); that is, $(nc)+(\mu) \Leftrightarrow (W)$ does not imply that
the same pair $(\lambda, \lambda')$ used to verify $(W)$ also satisfies $(\mu)$, even under Condition $(nc)$.
This observation motivates the precise definition 
 of {\it{elementary pairs}} in Definition~\ref{def:elementary}
 as well as four equivalent characterizations of \emph{elementary pairs} $(\pi, \pi')$
 presented in Theorem~\ref{thm:elementary_pair_equiv} below.

\par\noindent
{\rm{(7)}}\enspace
By (ii) of Proposition~\ref{prop:23022216}, 
the only adjacent symbol pairs in $[\delta\delta']$ that satisfy any (hence, all) of conditions in Proposition~\ref{prop:23022216} are 
$$(+ ~ \op), (\op ~ \oi), (\oi ~- ), (- ~ +)$$ 
along with one possible $(\op ~ +)$ occuring at the very end on the right.    In particular, the adjacent pairs 
$$ (-~ \oi),  (+ ~ -),  (\oi ~ \op)$$ 
never occur, even though they are compatible with a GGP interleaving pattern (Definition~\ref{def:interlace}).

\end{remark}

\subsection{Coherent Pairs and Elementary Pairs}

In this section, we introduce the notions of {\emph{coherent pair}} and {\emph{elementary pair}} for a pair of discrete series representations $\pi$ and $\pi'$ of $G$ and $G'$, respectively, where $(G,G')=(U(p,q),U(p-1,q))$. Proposition~\ref{prop:23022216} from the previous section, which provides equivalent conditions on the signature $[\delta \delta']$, 
clarifies the structural background of these notions.

\begin{definition}[Coherent 
Signature]\label{def:coh_deco}
A signature $[\delta\delta']$ is said to be \textbf{coherent} 
if it satisfies any (and hence all) of the equivalent conditions listed in Proposition~\ref{prop:23022216}.
\end{definition}
As a direct consequence of Proposition~\ref{prop:23022216},
 it follows that $q(\delta)=q(\delta')$ whenever $[\delta\delta']$ forms a coherent pair.

\begin{example}
Here are a few examples of coherent signatures $[\delta \delta']$. The first case is when $(p,q)=(5,2)$, and the second two cases is when $(p,q)=(4,3)$.
We also provide the values of $q(\delta)$ and $q(\delta')$:
$$+ ~ \op ~ + ~ \op ~ + ~ \op ~ \oi ~ - ~ \oi ~ - ~ + ~ \op ~ +; ~~  q(\delta) = q(\delta') = 6.$$
$$+ ~ \op ~ \oi ~ - ~ + ~ \op ~ + ~ \op ~ \oi ~ - ~ \oi ~ -  ~ +; ~ ~ q(\delta) = q(\delta') = 7.$$
$$ \oi ~ - ~ + ~ \op ~ + ~ \op ~ \oi ~ - ~ \oi ~ -  ~ + ~ \op ~ +; ~~  q(\delta) = q(\delta') = 4.$$
\end{example}

\begin{definition}[Coherent Pair]\label{def:coherent_pair}
Let $\pi = \pi_{\tlambda}$ and $\pi' = \pi'_{\tlambda'}$ be discrete series representations of $G$ and $G'$, respectively, with Harish-Chandra parameters $\lambda$ and $\lambda'$.
Assume that
\[
\operatorname{Hom}_{G'}(\pi|_{G'}, \pi') \neq \{0\},
\]
and that the signature $[\delta\delta']$ associated with the pair of Harish-Chandra parameters $(\lambda, \lambda')$ as in Section~\ref{subsec:He_description} is coherent.
Then the pair $(\pi, \pi')$ is called a \textbf{coherent pair}.
\end{definition}

By definition, a coherent pair $(\pi, \pi')$ satisfies the GGP interleaving condition.  However, not every pair $(\pi,\pi')$ such that $\operatorname{Hom}_{G'}(\pi|_{G'},\pi') \neq \{0\}$ is coherent.

For examples illustrating explicit conditions on parameters, see Section~\ref{subsec:Part4_Branching} in Part~\ref{part:4}.

\begin{definition}[Elementary Pair]
\label{def:elementary}
Let $(\pi, \pi') = (\pi_\lambda, \pi'_{\lambda'})$
be a coherent pair as defined in Definition~\ref{def:coherent_pair}.
We say that $(\pi, \pi')$ is an \textbf{elementary pair} if their minimal $K$-types satisfy
\[
[\mu_\lambda : \mulmdp] \neq 0.
\]
\end{definition}
\begin{remark}
In Theorem~\ref{thm:elementary_pair_equiv},
we present four equivalent characterizations of \emph{elementary pairs} $(\pi, \pi')$,
one of which includes the non-vanishing of the natural morphisms between $(\fP,K)$-cohomologies.
\end{remark}

\subsection{Proof of Proposition~\ref{prop:23022216}
}
\label{subsec:pf_prop:23022216}
~~~
\newline
This section provides a proof of Proposition~\ref{prop:23022216}.

Obviously, Conditions (iii), (iv), and (v) in Proposition~\ref{prop:23022216} are preserved under the reduction procedures \eqref{eqn:reduction} defined in Definition~\ref{def:reduction}. Once the equivalence among conditions (i)‐‐(vi) is established, it follows that the remaining conditions in Proposition~\ref{prop:23022216} are also preserved under the reduction procedures. However, this fact is used as part of the proof of Proposition~\ref{prop:23022216} itself.
In particular, we will prove in  Lemma~\ref{lem:23022211} that Condition (i) in Proposition~\ref{prop:23022216} 
is preserved under the reduction procedures \eqref{eqn:reduction}.

\begin{proof}
[Proof of Proposition \ref{prop:23022216}]

{\bf Step 1. (iii) $\Rightarrow$ (ii) $\Rightarrow$ (v) $\Rightarrow$ (iv).}
The implications (iii) $\Rightarrow$ (ii) $\Rightarrow$ (v) $\Rightarrow$ (iv) are straightforward.

{\bf Step 2.  (iv) $\Rightarrow$ (iii).}
Since properties (iii) and (iv) are 
preserved under the reduction process of deleting $+\op$ or $\oi -$ from the left (see \eqref{eqn:wreduction}),  the implication (iv) $\Rightarrow$ (iii) follows.

{\bf Step 3. (i) $\Rightarrow$ (iii).}
We now turn to the most substantial implication, namely, (i) $\Rightarrow$ (iii), which requires a more detailed argument.

Let $\Lambda=(\Lambda_1, \dots, \Lambda_{p+q})=\tlambda-\rho_{\mathfrak g}$
and $\Lambda'=(\Lambda'_1, \dots, \widehat{\Lambda_p}, \dots,{\Lambda}'_{p+q})=\lambda'-\rho_{\mathfrak{g}'}$ 
denote the coherent parameters.
Since $K \simeq U(p) \times U(q)$
and $K' \simeq U(p-1) \times U(q)$ as Lie groups,  
the condition $[W_{\tlambda} : \Wlmdp] \ne 0$ implies that the highest weights $\Lambda$ and $\Lambda'$ 
must satisfy the inequalities
\[
\Lambda_1 \ge \Lambda_1' \quad \text{ and } \quad \Lambda_{p+1}=\Lambda'_p.
\]

A direct computation of the coherent parameters $\Lambda=\tlambda-\rho_{\mathfrak g}$
and $\Lambda'=\lambda'-\rho_{\mathfrak{g}'}$ yields:
\begin{alignat*}{2}
\Lambda_{1} &= \lambda_{1} - \frac{p + q - 1}{2}, \quad
&\Lambda_{p+1} &= \lambda_{p+1} + \frac{p - q + 1}{2}, \\
\Lambda'_{1} &= \lambda_1' - \frac{p + q - 2}{2}, \quad
&\Lambda'_{p+1} &= \lambda_{p+1}' + \frac{p - q}{2}.
\end{alignat*}

Hence, under condition (W) holds, we obtain
\begin{align}
\label{eqn:Lambda_1}
\tlambda_{1} &\ge \lambda'_1+\frac{1}{2},
\\
\label{eqn:Lambda_p}
\tlambda_{p+1} &= \lambda'_{p+1}-\frac{1}{2}
\end{align}

Assume first that the string $[\delta\delta']$ begins with $+-$ on the left.  
By definition, this implies that among the \( 2(p+q)-1 \) numbers  
\[
\{ \lambda_1, \dots, \lambda_{p+q},\ \lambda'_1, \dots, \widehat{\lambda'_{p}}, \dots \lambda'_{p+q} \},
\]  
\(\lambda_1\) is the largest and \(\lambda_{p+1}\) is the second largest.
In particular, we have the strict inequalities:
 $\lambda_1 > \lambda_{p+1} > \lambda_{p+1}'$.
However,  the inequality $\tlambda_{p+1} > \lambda_{p+1}'$ contradicts \eqref{eqn:Lambda_p}.
Therefore, under condition (W), the string $[\delta\delta']$ cannot begin with $+-$.

Similarly, if $[\delta\delta']$ begins with either $- \oi$ or $- +$, then again $\lambda_{p+1}>\lambda_p'$ contradicts the same equation \eqref{eqn:Lambda_p}.
Therefore, these cases are excluded under condition (W).

Now consider $[\delta\delta']$ starting with $\op \oi$, $\op +$ or $\oi \op$. In these cases, one finds inequality $\lambda_1'>\lambda_1$, which contradicts \eqref{eqn:Lambda_1}.

From Definition~\ref{def:interlace} concerning the GGP interleaving relation, we conclude that
if condition $(W)$ holds, 
 then the string $[\delta\delta']$ must begin with either $+\op$ or $\oi -$. 

In turn, owing to Lemma~\ref{lem:23022211}, which we will prove in Section~\ref{subsec:reduction_K_types} and which concerns a \emph{reduction process},
 we can iteratively apply this argument,
 deleting $+\op$ or $\oi -$ from the left each step,
 until the string reduces to the singleton $\{+\}$.  
This establishes the implication (i) $\Rightarrow$ (iii).  

{\bf Step 4. (iii) $\Rightarrow$ (i).}
The implication (iii) $\Rightarrow$ (i) is straightforward, 
 as Lemma~\ref{lem:23022211}~(1) and (2), to be stated
 in Section~\ref{subsec:reduction_K_types},
will include this implication.

{\bf Step 5. (iv) $\Rightarrow$ (iii).}
To prove the implication (iv) $\Rightarrow$ (iii), 
 we claim the leftmost pair
 in $[\delta\delta']$ is either $+\op$ or $\oi-$.  
Indeed, 
 the first symbol in $[\delta\delta']$ must be
either $+$ or $\oi$, since $+_1 > \op_1$
 and $\oi_1 > -_1$.  
If it starts with $+$, 
 then, by the GGP interleaving relation, the next symbol must be either $\op$ or $-$; 
however $-$ is ruled out because $\oi_1 > -_1$.  
Similarly, 
 if it starts with $\oi$, 
 then the next one must be either $\op$ or $-$, 
 but $\op$ is excluded 
 because $+_1 > \op_1$.  
This confirms our claim.

{\bf Step 6. (iii) $\Rightarrow$ (vi).}
The implication (iii) $\Rightarrow$ (vi) is verified at the end of the proof.  

{\bf Step 7. (vi) $\Rightarrow$ (iii).}
Suppose that $[\delta \delta']$ satisfies the conditions $(nc)$ and $(\mu)$.

For $\nu=(\nu_1, \dots, \nu_p)$, we use the notation $\nu(\hat{p})$ to denote the truncated vector $(\nu_1, \dots, \nu_{p-1})$ by omitting the $p$-th entry.
Then we have 
\[
\rho_{\mathfrak{g}'}=\rho_{\mathfrak g}(\hat{p})+ (-\frac{1}{2} {\bf{1}}_{p-1} \oplus \frac{1}{2} {\bf{1}}_q).
\]

Since $\mu_{\lambda}
   = (\lambda-\rho_{\mathfrak g}) + \sum_{\alpha \in  \Phi^{nc, +}(\delta)} \alpha$ by Lemma~\ref{lem:24042109},
by condition $(nc)$, we have 
\[
\mulmdp = \mu_\lambda(\hat{p})+\lambda' -\lambda(\hat{p})+
(\frac{1}{2} {\bf{1}}_{p-1} \oplus (-\frac{1}{2}) {\bf{1}}_q).
\]

We write $\lambda=\lambda^{(p)}\oplus \lambda^{(q)}
\in \mathbb{C}^p \oplus \mathbb{C}^q$
and $\mu_{\lambda}=\mu_{\lambda}^{(p)} \oplus \mu_{\lambda}^{(q)}$.
Similarly, we write $\lambda'=\lambda'^{(p-1)} \oplus \lambda'^{(q)} = (\lambda'_1, \dots, \widehat{\lambda'_p}, \dots, \lambda'_{p+q})$.

Then $[\mu_{\lambda}:\mulmdp] \ne 0$ holds
if and only if
 then  the parameters
$\lambda$ and $\lambda'$ satisfy
\begin{align}
\label{eqn:mup}
\lambda_i \ge \lambda_i' + \frac 1 2, 
&\quad
\lambda_i' -\lambda_{i+1} \ge a_{i+1} - a_i  + \frac{1}
{2}, \quad (1 \le i \le p-1), 
\\
\label{eqn:muq}
\lambda^{(q)} &= {\lambda'}^{(q)} - \frac 1 2 {\bf{1}}_q,  
\end{align}
where we set
\[
 (a_1, \dots a_{p+q}):=\sum_{\alpha \in \Phi^{nc,+}(\delta)} \alpha.
\]
By the GGP interleaving relation, 
 we deduce from \eqref{eqn:muq}
 that $\oi$ and $-$ always occurs as the adjacent pair
 $\oi -$ in $[\delta \delta']$.  

We now consider a maximal connected substring in $[\delta \delta']$
 consisting only of $\op$ and $+$.  
By condition $(nc)$, 
\begin{equation*}
  \#\{\text{$+$s in this string}\}
  =
  \#\{\text{$\op$s in this string}\} 
 \ (+1)
\end{equation*}
where we add $+1$ 
 if there is no adjacent pair 
  $\oi -$ to the right. 
In such a string, 
 one has $\rho(\delta)_{i}-\rho(\delta)_{i+1}=1$.  
Hence, if $[\mu_{\lambda}:\mulmdp] \ne 0$, then  by \eqref{eqn:mup} we obtain
\[
\lambda_i \ge \lambda_i'+\frac 1 2
\ \text{ and } \ \lambda_i' \ge \lambda_{i+1}+\frac 1 2
\]
  
This determines the structure of the string as an iteration
 of adjacent pairs $+\op$ 
 (and ends at $+$ 
 if there is no $\oi -$ pair to the rightmost position).  
Therefore, condition (vi)--that is,
 (nc) and $(\mu)$--implies (iii).

For the last statement, 
 one sees readily
\begin{alignat*}{2}
&\text{(iii) $\Rightarrow$ (q)}
&&\text{by Lemma \ref{lem:23020112}, }
\\
&\text{(iii) $\Rightarrow$ ($\mu$)}
&&\text{by Lemma \ref{lem:23022211}, }
\\
&\text{(iii) $\Rightarrow$ ($\delta+$)}
\,\,
&&\text{by definition. }
\end{alignat*}
Since $\Phi^{nc, +}(\lambda)$ determines $\delta$, and vice versa, by \eqref{eqn:ncdelta}, 
 the condition $(\delta+)$ is equivalent to $(nc)$
 provided that the rightmost entry of $\delta$ is $+$.  

Hence, the proof of 
 Proposition~\ref{prop:23022216} completes,
 postponing the proof for related reduction procedures 
 in Section~\ref{subsec:reduction_K_types}.
\end{proof}

\subsection{Comparison of Two \texorpdfstring{$K$}{K}-types \texorpdfstring{$\mu_\lambda$}{mulambda} and \texorpdfstring{$W_\lambda$}{Wlambda}}

We prove that for any elementary pair $(\pi_\lambda, \pilmdp)$ (Definition~\ref{def:elementary}) the multiplicity $[W_\lambda:\Wlmdp]$ is non-zero.

For later purpose, we explicitly compute 
the conditions on the parameters $\lambda$ and $\lambda'$ under which the multiplicities $[W_\lambda:\Wlmdp] \neq 0$
and
$[\mu_\lambda: \mulmdp] \neq 0$
hold 
for a coherent pair $(\pi_\lambda, \pilmdp)$.

\begin{lemma}
\label{lem:25060510}
Let $\pi$ be a discrete series representation of $G$ with decoration $\delta$, and let $\pi'$ be a discrete series representation of $G'$ with decoration $\delta'$ satisfying 
\[
\operatorname{Hom}_{G'}(\pi|_{G'}, 
\pi') \neq 0 \quad \text{ and } \quad
\Phi^{nc, +}(\delta)=\Phi^{nc,+}(\delta').
\]
Let $\lambda$ and $\lambda'$ be the Harish-Chandra parameters of $\pi$ and $\pi'$, respectively.
Then the following hold:
\newline
(1).  $[W_\lambda:\Wlmdp] \neq 0$ if and only if 
\begin{equation}
\label{eqn:W_cond}
    \lambda'_j=\lambda_j+\frac{1}{2} \text{ for all } p+1 \le j \le p+q.
\end{equation}
\newline
(2).
Define
\[
 (a_1, \dots a_{p+q}):=\sum_{\alpha \in \Phi^{nc,+}(\delta)} \alpha.
\]
Then 
$[\mu_\lambda: \mu'_{\lambda'}] \neq 0$
if and only if condition \eqref{eqn:W_cond} holds and 
\begin{equation}
\label{eqn:mu_cond}
    \lambda'_i\ge \lambda_{i+1}+ a_{i+1} - a_{i} \text{ for all } 1 \le i \le p-1.
\end{equation}
(3)
If $[\mu_\lambda: \mulmdp] \neq 0$,
then $[W_\lambda:\Wlmdp] \neq 0$.

\end{lemma}
\begin{proof}
It follows from Remark~\ref{rem:equiv_coherent_sign}~(7) that, under condition $(nc)$, 
the Harish-Chandra parameters~$\lambda$ and~$\lambda'$
satisfy the following inequalities:
\begin{equation}
\label{eqn:nc_lmd_ineq}
  \lambda_1>\lambda'_1> \dotsb > \lambda'_{p-1} >\lambda_p; \quad
  \lambda_{p+1}' > \lambda_{p+1} > \dotsb >\lambda'_{p+q}>\lambda_{p+q}.
\end{equation}
\newline
(1) Let $\Lambda=(\Lambda_1, \dots, \Lambda_{p+q})=\tlambda-\rho_{\mathfrak g}$
and $\Lambda'=(\Lambda'_1, \dots, \widehat{\Lambda_p}, \dots,{\Lambda}'_{p+q})=\lambda'-\rho_{\mathfrak{g}'}$ 
denote the coherent parameters.
A straightforward computation shows
\begin{alignat*}{2}
\Lambda_{i} &= \lambda_{i} - \frac{p + q + 1 - 2i}{2}
&\quad &(1 \le i \le p+q), \quad
 \\
\Lambda'_{j} &= \lambda_j' - \frac{p + q - 2j}{2}
&\quad &(1 \le j \le p-1),
\\
\Lambda'_{j} &= \lambda_{j}' + \frac{p +q+2-2j}{2}
&\quad &(p+1 \le j \le p+q).
\end{alignat*}
Since $K \simeq U(p) \times U(q)$
and $K' \simeq U(p-1) \times U(q)$ are isomorphisms as Lie groups,  
the condition $[W_{\tlambda} : \Wlmdp] \ne 0$ holds if and only if the highest weights $\Lambda$ and $\Lambda'$ of $W_\tlambda$ and $\Wlmdp$, respectively satisfy the following equalities and inequalities:
\begin{alignat*}{1}
&\Lambda_1 \ge \Lambda'_1 \ge \cdots \ge \Lambda'_{p-1}
\ge \Lambda_p,
\\
&\Lambda_j=\Lambda'_j \quad (p+1 \le j \le p+q).
\end{alignat*}
These conditions are equivalent to
\begin{alignat*}{2}
&\lambda_i \ge \lambda'_i +\frac{1}{2},  \quad
\lambda_i'\ge \lambda_{i+1}+\frac{1}{2} &\quad& (1 \le i \le p-1),
\\
&\lambda'_{j}=\lambda_j+\frac{1}{2} &\quad & (p+1 \le j \le p+q).
 \end{alignat*}
The first two sets of inequalities follow from \eqref{eqn:nc_lmd_ineq} since the difference between any components of
$\lambda$ and any component of $\lambda'$ is a half-integer.
This completes the proof of the first statement.
\newline
(2).
This was established in Step 7 of the proof of Proposition~\ref{prop:23022216}.
\newline
(3).
The claim follows from the explicit conditions given in (1) and (2).
\end{proof}

\subsection{Equivalent Conditions for Elementary Pairs}
\label{subsec:equiv_elementary}
In this section, we summarize several equivalent definitions of elementary pairs. In particular, condition (ii) below corresponds to the original definition introduced in Definition~\ref{def:elementary}.

In (iv) below, the relative $(\fP, K)$-cohomology
 (see Section~\ref{sec:Pk})
 is defined for the maximal parabolic subalgebra
\begin{equation}
\label{eqn:P_+}
\mathfrak{P}\equiv\mathfrak{P}^+:={\mathfrak{p}}_+ + {\mathfrak{k}}_{\mathbb{C}}.
\end{equation}

\begin{theorem}[Equivalent Conditions for Elementary Pairs $(\pi, \pi')$]
\label{thm:elementary_pair_equiv}
Let $\pi = \pi_\lambda$ and $\pi' = \pilmdp$
be discrete series representations of $G$ and $G'$, respectively, with Harish-Chandra parameters $\lambda$
and $\lambda'$.
Let $\delta$ and $\delta'$ denote the decorations of $\pi$ and $\pi'$, respectively.
Let $W_\lambda$ and $W'_{\lambda'}$ denote the coherent parameters of $\pi$ and $\pi'$, respectively
(Definition~\ref{def:W_lambda}).
Assume that the signature $[\delta \delta']$ associated with $\lambda$ and $\lambda'$ is coherent
(see Proposition~\ref{prop:23022216} for six equivalent definitions).
Then, the following four conditions are equivalent:
\begin{itemize}
    \item[{\rm{(i)}}]
    $[\mu_\lambda : \mulmdp] \neq 0$.
    
    \item[{\rm{(ii)}}]
    In addition to (i), $\operatorname{Hom}_{G'}(\pi|_{G'}, \pi') \neq \{0\}$.
    
    \item[{\rm{(iii)}}]
    In addition to (ii), $[W_\lambda : \Wlmdp] \neq 0$.
    
    \item[{\rm{(iv)}}]
    In addition to (iii), the restriction map
    $T\otimes \operatorname{pr}_{W} $
    induces an isomorphism
    $$ 
 H^\ell(\mathfrak{P}, K; \pi \otimes W_{\lambda}^\vee) \overset{\sim}\longrightarrow
    H^\ell(\mathfrak{P}', K'; \pi' \otimes {W'_{{\lambda}'}}^\vee),
    $$
where $T \in \operatorname{Hom}_{G'}(\pi|_{G'}, \pi')$ and $\operatorname{pr}_{W} \in \operatorname{Hom}(W_\lambda^\vee|_{K'}, {W'_{\lambda'}}^\vee)$ are non-zero homomorphisms as guaranteed by (ii) and (iii), respectively, each unique up to scalar multiplication.
The cohomology groups are one-dimensional when
$\ell=q(\pi) = q(\pi')$ (Definition~\ref{def:q}), and vanish in all other degrees.
\end{itemize}
\end{theorem}

\begin{remark}
(1)
The implication (ii) $\Rightarrow$ (iii) in Theorem~\ref{thm:elementary_pair_equiv} does not necessarily hold without the assumption that the associated signature $[\delta\delta']$ is a coherent signature. See Remark~\ref{rem:equiv_coherent_sign}.

(2)
Since both pairs $(G, G')$ and $(K, K')$ are strongly Gelfand pairs, the multiplicities
\[
\dim \operatorname{Hom}_{G'}(\pi|_{G'}, \pi'), \quad
\dim \operatorname{Hom}_{K'}(\mu_\lambda|_{K'},\mulmdp), 
\quad
\dim \operatorname{Hom}_{K'}(W_\lambda|_{K'}, \Wlmdp)
\]
can only take values $0$ or $1$. Therefore,
\[ 
[\pi : \pi']
= [\mu_\lambda : \mulmdp] 
= [W_\lambda : \Wlmdp] = 1,
\]
provided that $(\pi, \pi')$ is an elementary pair---that is, any (hence, all) of the conditions (i)--(iv) in Theorem~\ref{thm:elementary_pair_equiv} are satisfied.
\end{remark}

\begin{proof}[Proof of Theorem~\ref{thm:elementary_pair_equiv}]
{\bf{(i)$ \Rightarrow$ (ii).}}
The argument is the same as in \cite{HL}, as recalled in the proof of \cite[Prop.~2.5]{H14}:  By
Flensted-Jensen's theorem on the positivity of the matrix coefficients of the minimal types, it follows that $\pi'$ is weakly contained in $\pi$, provided that the restriction to $G'$ of the matrix coefficients of $\pi$ lies in $L^{2+\varepsilon}(G')$ for all $\varepsilon$.  
But this was verified in \cite[Prop.~2.3.4]{HL}.
\newline
{\bf{(ii)$ \Rightarrow$ (iii).}}
By Definition~\ref{def:coherent_pair}, the pair $(\pi, \pi')$ is coherent.
Then implication (ii)$ \Rightarrow$ (iii) was established in Lemma~\ref{lem:25060510}
\newline
{\bf{(iii)$ \Rightarrow$ (iv).}}
The $K$-module $\wedge^j({\mathfrak{p}}^+)$ is not irreducible in general. 
We define
\[
 t(\delta)
   := \sum_{\alpha \in  \Phi^{nc, +}(\delta)} \alpha 
   = \sum_{\alpha \in \Delta({\mathfrak{p}}^+) \cap w \Delta^+({\mathfrak{g}})} \alpha.
\]
Then the $K$-module  $\wedge^j({\mathfrak{p}}^+)$
decomposes into a diret sum of irreducible representations as follows.
\begin{align}
\label{eqn:25050326}
 \wedge^j({\mathfrak{p}}^+)
   \simeq& 
   \bigoplus_{\substack{\delta \in \operatorname{Deco}(p,q)\\ q(\delta)=j}} 
   F^K(t(\delta)),
\\
\notag
  \wedge({\mathfrak{p}}^+)
   \simeq& \bigoplus_{\delta \in \operatorname{Deco}(p,q)} F^K(t(\delta)).
\end{align}
A small combinatorial  computation using the Blattner formula included after this proof  shows that
\begin{equation}
\label{eqn:25060908}
\begin{aligned}
\operatorname{Hom}_K(\wedge^\ell \mathfrak{p}^+, \pi \otimes W_\lambda^\vee)
&\simeq\operatorname{Hom}_K(
W_\lambda \otimes \wedge^\ell \mathfrak{p}^+, \pi)
\\
&\simeq
\begin{cases}
0 & \ell \neq q(\lambda),
\\
    \operatorname{Hom}_K(W_\lambda \otimes F^K(t(\delta)), \mu_\lambda) \simeq \mathbb{C} &\ell = q(\lambda),
\end{cases}
\end{aligned}
\end{equation}
which induces the following equalities:
\begin{equation}
\label{eqn:pK_vanish}
H^\ell(\fP, K; \pi \otimes W_\lambda^\vee)
\simeq
\begin{cases}
0 & \ \ell \neq q(\lambda),
\\
 \mathbb{C} &\ \ell = q(\lambda).
\end{cases}
\end{equation}

The $(\mathfrak{P}, K)$-cohomologies of the representation $\pi$ can be computed by
the complex 
$\{ \operatorname{Hom}_K(\wedge^\ell \mathfrak{p}^+, \pi \otimes W_\lambda^\vee)\}_{\ell \in \mathbb{N}}$.
The identities \eqref{eqn:25060908} implies that
there is a canonical bijection
\[
\operatorname{Hom}_K(
W_\lambda \otimes F^K(t(\delta)), \mu_\lambda)
\overset\sim\rightarrow
H^{q(\lambda)}(\mathfrak{P}, K; \pi \otimes W_{\lambda}^\vee). 
\]
Similarly, there is a canonical bijection for $G'$:
\[
\operatorname{Hom}_{K'}(
\Wlmdp\otimes 
F^{K'}(t^\prime(\delta')), \mulmdp)
\overset\sim\rightarrow
H^{q(\lambda')}(\mathfrak{P}', K'; \pi' \otimes {\Wlmdp}^{\vee}). 
\]
Thus the assertion follows from the commutative diagram below.
\begin{equation}
\begin{CD}\label{eqn:rest_PK}
\operatorname{Hom}_K(
W_\lambda \otimes F^K(t(\delta)), \mu_\lambda)
 @>>> 
H^{q(\lambda)}(\mathfrak{P}, K; \pi \otimes W_{\lambda}^\vee)
  \\
@VVV     @VVV \\
\operatorname{Hom}_{K'}(
{\Wlmdp} \otimes 
F^{K'}(t^\prime(\delta')), \mulmdp)
 @>>>  H^{q(\lambda')}(\mathfrak{P}', K'; \pi' \otimes {\Wlmdp}^{\vee})
\end{CD}
\end{equation}
\end{proof}

For completeness, we provide a short proof for \eqref{eqn:25060908}. 
The same argument remains valid when $\pi_\lambda$ is a limit of discrete series representations, provided that the minimal $K$-type parameter $\mu_\lambda$ is dominant with respect to $\Phi^{c,+}$.

\begin{proof}[Proof of \eqref{eqn:25060908}]
    
Let $U$ be any irreducible $K$-module occurring in the tensor product representation $W_\lambda \otimes \wedge^* \mathfrak{p}^+$. Then its highest weight
must be of the form
\[
\lambda - \rho_{\mathfrak g} +\sum_{\beta \in J} \beta,
\]
for some subset $J \subset \Phi^{nc,+}$.
Let $w \in W^\mathfrak{k}$ be the corresponding element to the decoration $\delta$. Then the set $\Phi^{nc,+}$ decomposed into a disjoint union:
\[
\Phi^{nc,+}=(\Phi^{nc,+} \cap w\Phi^+) \sqcup (\Phi^{nc,+} \cap w\Phi^-)
=\Phi^{nc,+}(w) \sqcup \Delta^+(w). 
\]
Accordingly, we write a decomposition of $J$ as $J={J_+} \sqcup {J_-}$.

On the other hand, by the Blattner formula,
the highest weight of any $K$-type in $\pi_\lambda$ must be of the form
\[
\mu_\lambda + \sum_{\alpha \in 
\Delta(\mathfrak{p}) \cap w \Phi^+}
n_\alpha \alpha,
\]
for some non-negative integers $n_\alpha$ for $\alpha \in \Delta(\mathfrak{p}) \cap w \Phi^+$.
The set $
\Delta(\mathfrak{p}) \cap w \Phi^+$ is expressed as the following disjoint union:
\[
\Delta(\mathfrak{p}) \cap w \Phi^+
=\Phi^{nc,+}(\delta)\sqcup (\Phi^{nc,-}\cap w \Phi^+).
\]
Therefore, if $\operatorname{Hom}_{K}(W_\lambda \otimes \wedge^*\mathfrak{p}^+, \pi_\lambda|_K) \neq \{0\}$, then the following identity holds:
\[
\lambda - \rho_{\mathfrak g} +\sum_{\beta \in J_+} \beta + \sum_{\beta \in J_-} \beta
=
\mu_\lambda + 
\sum_{\alpha \in \Phi^{nc,+}(\delta)}
n_\alpha \alpha
-\sum_{\alpha \in \Phi^{nc,+}\cap w \Phi^{+}} n_\alpha \alpha.
\]
Since
$\mu_{\lambda}
   = \lambda -\rho_{\mathfrak g} + \sum_{\alpha \in  \Phi^{nc, +}(\delta)} \alpha. 
$ by \eqref{eqn:mu_Lambda},
we obtain:
\[
\sum_{\beta \in \Phi^{nc,+}(\delta)\setminus J_+} \beta
+\sum_{\alpha\in\Phi^{nc,+(\delta)}} n_\alpha \alpha=
\sum_{\beta \in J_-} \beta + \sum_{\alpha \in w \Phi^{-} \cap \Phi^{nc,+}} n_\alpha \alpha.
\]
Taking the inner product against $\rho(\delta)$, we see that $J_+=\Phi^{nc,+}(\delta)$, $J_-=\emptyset$, and $n_\alpha=0$
for all $\alpha \in \Delta(\mathfrak{p}) \cap w \Phi^+$.

Hence, \eqref{eqn:25060908} is proved.
\end{proof}
 
\begin{remark}
\label{rem:elementary_SBO_surjects_minK}
(1) 
As the above proof shows, the implication (i) $\Rightarrow$ (ii) does not require the condition that the signature $[\delta\delta']$ is coherent; however, this condition is necessary in order to apply the proposition to restrictions of coherent cohomology. This point will be discussed in Part~\ref{global}.
\newline
(2)
It follows from the proof of  Theorem~\ref{thm:elementary_pair_equiv} that
if $(\pi, \pi')=(\pi_, \pi'_{\lambda'})$ is elementary, then any element $T$ of $\operatorname{Hom}_{G'}(\pi|_{G'}, \pi')$ induces
a non-trivial homomorphism of minimal $K$-types
$\mu_\lambda \to \mu'_{\lambda'}$.
Indeed, the proof for (i) $\Rightarrow$ (ii)
constructs such $T$ that induces a non-trivial map
$\mu_\lambda \to \mu'_{\lambda'}$. Then the claim follows since $\dim \operatorname{Hom}_{G'}(\pi|_{G'}, \pi') =1$.
\end{remark}

\subsection{Unique Existence of Elementary Pairs}

In this section, we present two theorems related to the construction of elementary pairs $(\pi, \pi')$. 
Theorem~\ref{thm:230816} addresses the problem of finding an elementary pair when a discrete series representation $\pi'$ of the subgroup $G'$ is given, by finding a corresponding discrete series representation $\pi$ of $G$. In contrast, 
Theorem~\ref{thm:23022420b} provides an elementary pair when a discrete series representation $\pi$ of the group $G$ is given, by finding a corresponding discrete series representation $\pi'$ of the subgroup $G'$. The latter theorem relies on a specific assumption about the parameters of $\pi$.

In both cases, we also prove the uniqueness of the constructed pairs.

Let us begin with some combinatorial preparations concerning two decorations $\delta, \delta'$,
which follows readily from
Proposition~\ref{prop:23022216}.

\begin{lemma}
\label{lem:23022420}
Suppose that $\delta \in \operatorname{Deco}(p,q)$
 ends with $+$ in the rightmost entry.  
Then there exists a unique coherent signature $[\delta\delta']$ (see Definition~\ref{def:coh_deco}),
 where $\delta' \in \operatorname{Deco}(p-1,q)$.  
\end{lemma}

\begin{proof}
[Proof of Lemma~\ref{lem:23022420}]
Construct the signature $[\delta \delta']$ by replacing each $+$ in $\delta$ with $+\ \op$,  
and each $-$ with $\oi\ -$.  
This defines an element $\delta' \in \operatorname{Deco}(p-1,q)$.
By (iii) of Proposition~\ref{prop:23022216} that
the resulting signature $[\delta \delta']$ is coherent,
and this coherent signature is unique.
\end{proof}

Given a coherent signature $[\delta \delta']$, the following lemma provides a precise condition under which the support $\operatorname{Supp}(\pi_{\lambda}^{\delta}|_{G'})(\delta')$, as defined in Definition~\ref{def:Supp_delta}, is non-empty.

\begin{lemma}
\label{lem:23022509}
Let $\tlambda \leftrightarrow (\poslmd, w) \leftrightarrow (\poslmd, \delta)$ denote the correspondence described in Proposition~\ref{prop:discpara}.
Let $[\delta \delta']$ be the coherent signature,
as given in Lemma~\ref{lem:23022420},
where $\delta' \in \operatorname{Deco}(p-1,q)$. 
\newline
(1).
The following three conditions are equivalent:
\begin{itemize}
  \item[{\rm (i)}] $\operatorname{Supp}(\pi_{\lambda}^{\delta}|_{G'})(\delta') \neq \emptyset$.
  \item[{\rm (ii)}] $\tlambda_i - \tlambda_j \ge 2$ holds for all indices $i \le p < j$ such that $w^{-1}(i) + 1 = w^{-1}(j)$.
  \item[{\rm (iii)}] $(\poslmd)_a - (\poslmd)_{a+1} \ge 2$ holds for all $a$ such that $w(a) \le p < w(a+1)$.
\end{itemize}

(2).
Assume that any (and hence all) of the above equivalent conditions hold.
We define $\lambda'\in\mathbb{Z}_{\mathrm{int}}^{p+q-1}$ by
\begin{equation}
\begin{cases}
\lambda_i':=\tlambda_i-\frac  1 2\qquad &(1 \le i \le p-1), 
\\
\lambda_j':=\tlambda_{j+1}+\frac  1 2\qquad &(p \le j \le p+q-1).
\end{cases}
\end{equation}
We express this definition as
\begin{equation}
\label{eqn:lmdtolmd1}
  \lambda' := (\tlambda^{(p)} - \tfrac{1}{2} \mathbf{1}_p)(\widehat{p}) \oplus (\tlambda^{(q)} + \tfrac{1}{2} \mathbf{1}_q).
\end{equation}
Then $\lambda'$ is regular; that is, it
belongs to $\mathbb{Z}_{\mathrm{int, reg}}^{p+q-1}$,
defined by
\[
   \{x \in {\mathbb{Z}}_{\operatorname{int}}^{p+q-1}
:
x_i \ne x_j\quad\text{if $i \ne j$}\},
\]
see \eqref{eqn:Z_N_reg} for the definition.
\end{lemma}

\begin{remark}
\label{rem:strong_regularity}
Condition (iii) of Lemma~\ref{lem:23022509} is weaker
 than the assumption
$(\poslmd)_a-(\poslmd)_{a+1} \ge 2$
 for all $1 \le a \le p+q-1$,
as stated in \cite[Hypothesis 4.6]{H14}.  
\end{remark}

\begin{proof}[Proof of Lemma~\ref{lem:23022509}]
{\bf{(ii) $\Leftrightarrow$ (iii).}}
The equivalence between (ii) and (iii) follows directly from the definitions $\tlambda=w \poslmd$,
where $w\in W^\mathfrak{k}$.
\paragraph{\bf (i) $\Rightarrow$ (ii).}  
Assume that condition (i) holds.
We take $\lambda'$ in the support  
$\operatorname{Supp}(\pi_\lambda^\delta|_{G'})(\delta')$.
Let $i, j$ be a pair of indices such that
$i \le p < j$ and $w^{-1}(i)+1=w^{-1}(j)$.
Then
the decoration $\delta$ contains
the  $+-$ pattern at the positions $(w^{-1}(i), w^{-1}(j))$ by
the correspondence $w \leftrightarrow \delta$ in
Proposition~\ref{prop:discpara}.

It follows from condition (ii) in Proposition \ref{prop:23022216} that this pattern appears in the interleaving string $[\delta\delta']$ as a subsequence of the form 
\[
+ \op \oi +.
\]
By the GGP interleaving condition, the corresponding parameter $\lambda'$ in the support must satisfy
\[
\lambda_i > \lambda_i' > \lambda_{j-1}' > \lambda_j.
\]
Taking into account the parity shifts between $\lambda
\in ({\mathbb{Z}}_{\operatorname{int}})^{p+q}
=({\mathbb{Z}}+\tfrac{p+q+1}2)^{p+q}$
and $\lambda'
\in ({\mathbb{Z}}_{\operatorname{int}})^{p+q-1}
=({\mathbb{Z}}+\tfrac{p+q}2)^{p+q-1}
$, this implies 
\[
\lambda_i - \lambda_j \ge \tfrac{1}{2} + 1 + \tfrac{1}{2} = 2.
\]
\paragraph{\bf (ii) $\Rightarrow$ (i).}  
Conversely, if condition (ii) holds, then the parameter $\lambda'$ defined in \eqref{eqn:lmdtolmd1} is regular and satisfies the GGP interleaving condition.  
Therefore, the discrete series representation $\pi'$ of $G'$
with Harish-Chandra parameter $\lambda'$ and
 decoration $\delta'$ appears in the restriction of $\pi_\tlambda^\delta$ to $G'$, that is,
\[ 
\pi' \in \operatorname{Supp}(\pi_{\tlambda}^{\delta}|_{G'})(\delta').
\]
\end{proof}

The following is also a straightforward consequence
 of the simpler part of Proposition \ref{prop:23022216}, 
 but we give a precise formulation
 for future reference.  The choice of $\lambda'$ below is essentially the same as the construction in \cite[\S 4.2]{H14}
 under a slightly stronger regularity condition.

\begin{theorem}[Construction of Elementary Pair: from $G$ to $G'$]
\label{thm:23022420b}
Let $\pi=\pi_\lambda$ be a discrete series representation of $G=U(p,q)$ with Harish-Chandra parameter $\lambda \in {\mathbb{Z}}_{\operatorname{int}}^{p+q}$.
Let $\delta \in \operatorname{Deco}(p,q)$ denote its decoration.
Assume that $\delta$ ends with $+$ at the rightmost and that
 $\lambda$ satisfies any (hence, all) of the equivalent
 \lq\lq{regularity condition}\rq\rq\
 as in Lemma \ref{lem:23022509}~(1). 
 
Then there exists a discrete series representation $\pi'=\pi'_{\lambda'}$ of $G'$ such that
$(\pi, \pi')$ is an elementary pair (Definition~\ref{def:elementary}).

More precisely, let $[\delta\delta']$ be the unique coherent signature, as given in Lemma~\ref{lem:23022420}, where $\delta' \in \operatorname{Deco}(p-1,q)$.
Then there exists
 $\lambda' \in ({\mathbb{Z}}_{\operatorname{int}}^{p+q-1})(\delta')$ (see \eqref{eqn:Z_int_delta} for definition)
 such that $[\mulmd:\mu'_{\lambda'}]\ne 0$
 and $[W_{\lambda}:W'_{\lambda'}] \ne 0$.  
 
\end{theorem}

\begin{remark}
The same conclusion holds
 if a stronger assumption of the regularity 
$(\poslmd)_a-(\poslmd)_{a+1} \ge 2$
 $(1 \le i \le p+q-1)$
 is satisfied.  
See Remark~\ref{rem:strong_regularity}.
\end{remark}

\begin{proof}
[Proof of Theorem~\ref{thm:23022420b}]
We take
\[
   \lambda':=
   (\tlambda^{(p)}-\tfrac 1 2 {\bf{1}}_p)(\widehat p) 
   \oplus
   (\tlambda^{(q)}+\tfrac 1 2 {\bf{1}}_q), 
\]
as in \eqref{eqn:lmdtolmd1}.
By Lemma \ref{lem:23022509}, 
 $\lambda'$ is regular.

By the choice of $\delta'$, 
 one has
\[
  \rho(\delta')=
  (\rho(\delta)^{(p)}-\tfrac 1 2 {\bf{1}}_p)(\widehat p)
  \oplus
  (\rho(\delta)^{(q)}-\tfrac 1 2 {\bf{1}}_q).  
\]
Combined with 
\[
  2\rho_c'=
  (2\rho_c^{(p)}+{\bf{1}}_p)(\widehat p)
  \oplus
  2\rho_c^{(q)}, 
\]
one has
\[
  \text{$\mulmdp=\mu_{\lambda}(\widehat{p})$, 
namely, 
 ${\mulmdp}^{(p-1)}={\mulmd}^{(p)} (\widehat p)$
 and ${\mulmdp}^{(q)}= \mulmd^{(q)}$.  
}
\]
Hence $[\mulmd:\mu'_{\lambda'}] \ne 0$.  

It then follows from Theorem~\ref{thm:elementary_pair_equiv} that
 $[W_{\tlambda}:W'_{\lambda'}]\ne 0$.  
\end{proof}

\medskip
As a counterpart to Theorem~\ref{thm:23022420b} with the roles of \( G \) and \( G' \) reversed, the following theorem guarantees that for any discrete series representation \( \pi' \) of \( G' \), there exists a discrete series representation \( \pi \) of \( G \) such that  $(\pi, \pi')$ forms an elementary pair.

\begin{theorem}[Construction of an Elementary Pair: from $G'$ to $G$]
\label{thm:230816}
Let $(G, G')=(U(p,q), U(p-1,q))$.
Let $\delta' \in \operatorname{Deco}(p-1, q)$.

\begin{itemize}
  \item[(1)] There exists a unique coherent signature $[\delta\delta']$ with $\delta \in \operatorname{Deco}(p, q)$.

  \item[(2)] Let $\delta \in \operatorname{Deco}(p, q)$ be the decoration determined in (1).
  Then for any discrete series representation $\pi'=\pi'_{\lambda'}$ of $G'$ with Harish-Chandra parameter $\lambda'$ belonging to the decoration $\delta'$, there exists
  a discrete series representation $\pi=\pi_\lambda$ of $G$ with Harish-Chandra parameter $\lambda$ belonging to the decoration $\delta$ such that $(\pi, \pi')$ forms an elementary pair.
  In other words, for any
  $\lambda' \in {\mathbb{Z}}_{\operatorname{int}}^{p+q-1}(\delta')$, there exists $\lambda \in {\mathbb{Z}}_{\operatorname{int}}^{p+q}(\delta)$ such that the following three conditions are satisfied:
   \[
  \pi'_{\lambda'}\in \operatorname{Supp} \left( \left. \pi_{\lambda} \right|_{G'} \right), 
  \quad
  [\mu_{\lambda} : \mu'_{\lambda'}] \ne 0, 
  \quad
  \text{and }
  \quad
  [W_{\lambda} : W'_{\lambda'}] \ne 0.
  \]
\end{itemize}
\end{theorem}

\begin{proof}[Proof of Theorem~\ref{thm:230816}]
(1)\enspace
Construct the string $[\delta \delta']$
 by replacing 
 $\op$ by $+\op$, 
 $\oi$ by $\oi-$, 
 and by adding $+$ at the rightmost.  
Then this is the unique coherent signature
 by Proposition \ref{prop:23022216}.  
\par\noindent
(2)\enspace
We write $\lambda'=({\lambda'}^{(p-1)}, {\lambda'}^{(q)})$.  
Take any $\lambda_p \in {\mathbb{Z}}+\frac 1 2 (p+q+1)$
 such that 
\[
  \lambda_p< \operatorname{min}\{\lambda_1'+\tfrac 1 2 , \dots, \lambda_{p-1}'+\tfrac 1 2, \lambda_p'-\tfrac 1 2, \dots, \lambda_{p+q-1}'-\tfrac 12\}, 
\]
and set
\begin{equation}
   \lambda:=((\lambda')^{(p-1)}+\tfrac 1 2 {\bf{1}}_p, \lambda_p, (\lambda')^{(q)}-\tfrac 1 2 {\bf{1}}_q) \in ({\mathbb{Z}}^{p+q})_{\operatorname{int}}.  
\end{equation}
Then, by the choice of the coherent signature $[\delta \delta']$ in (1), we have 
\[
\lambda \in {\mathbb{Z}}_{\operatorname{int}}^{p+q} \cap C(\delta) = ({\mathbb{Z}}_{\operatorname{int}}^{p+q})(\delta).
\]
Hence, $\pi_{\lambda}$ is a discrete series representation of $G$ with Harish-Chandra parameter $\lambda$ and decoration $\delta$.  
It then follows from Theorem~\ref{thm:23022420b}
 that $[\mu_{\lambda}:\mu'_{\lambda'}] \ne 0$ and $[W_{\lambda}:W'_{\lambda'}]\ne 0$.  
\end{proof}

We may take $\sigma' = \pi'$ in the second statement of Theorem~\ref{thm:230816},
as seen from the above proof
together with the definition of coherent pairs (Definition~\ref{def:coherent_pair}).
This statement is formulated precisely as follows:

\begin{cor}\label{cor:coh_to_elem}
 Given any coherent pair $(\pi, \pi')$ with associated signature $[\delta\delta']$,
 there exists a discrete series representation $\sigma$ of $G$ such that $(\sigma, \pi')$ constitutes an elementary pair having the same signature $[\delta\delta']$.
\end{cor}

The following proposition provides an equivalent definition of coherent pair $(\pi, \pi')$:
\begin{prop}
\label{prop:cohpair} The pair $(\pi,\pi')$ is a coherent pair (Definition~\ref{def:coherent_pair}) if 
\begin{enumerate}
\item  $\operatorname{Hom}_{G'}(\pi|_{G'},\pi') \neq 0$; i.e. the decorations $\delta$ and $\delta'$ satisfy the GGP interleaving relation;
\item there is an elementary pair $(\sigma,\sigma')$ with $\delta(\sigma) = \delta$ and $\delta(\sigma') = \delta'$.  
\end{enumerate}
\end{prop}

\begin{proof}
The implication \textup{(1)} $\Rightarrow$ \textup{(2)} is immediate from Corollary~\ref{cor:coh_to_elem}.

The reverse implication \textup{(2)} $\Rightarrow$ \textup{(1)} follows directly from the definition of coherent pairs.
\end{proof}

\vskip 1pc

The method of \emph{symmetry breaking under translation},
developed in Section~\ref{sec:SBO_translation},
together with Theorem~\ref{thm:230816},
yields a simple and alternative proof of He's theorem in the special case where the signature is coherent.

\begin{cor}
\label{cor:25061806}
Let $[\delta \delta']$ be a coherent signature (see Definition~\ref{def:coh_deco}).
Then, for any $\lambda \in \mathbb{Z}_{\mathrm{int}}^{p+q}(\delta)$ and any $\lambda' \in \mathbb{Z}_{\mathrm{int}}^{p+q-1}(\delta')$ (see \eqref{eqn:Z_int_delta} for notation), we have
\[
\operatorname{Hom}_{G'}(\pi_\lambda|_{G'}, \pi'_{\lambda'}) \neq 0,
\]
where $\pi_\lambda$ and $\pi'_{\lambda'}$ are the representations associated to $\lambda$ and $\lambda'$, respectively.
\end{cor}

\begin{proof}
We set $\pi := \pi_\lambda$ and $\pi' := \pi'_{\lambda'}$, and also write $\sigma' := \pi'$.
It follows from Theorem~\ref{thm:230816} (see also Corollary~\ref{cor:coh_to_elem}) that there exists $\xi \in \mathbb{Z}_{\mathrm{int}}^{p+q}(\delta)$ such that the pair $(\sigma, \sigma')$, with $\sigma := \pi_\xi$, forms an elementary pair.

Since $(\lambda, \lambda')$ and $(\xi, \lambda')$ define the same interleaving pattern $[\delta \delta']$, a result from Section~\ref{sec:SBO_translation}---specifically, the \emph{stability theorem for multiplicities} in branching laws
within the same \emph{fence}, as proved in \cite[Theorem 3.3]{KS} and recalled here in Remark~\ref{rem:stable_multiplicity}---implies
\[
  \operatorname{Hom}(\pi|_{G'}, \pi') \simeq \operatorname{Hom}(\sigma|_{G'}, \sigma').
\]

Because $\operatorname{Hom}(\sigma|_{G'}, \sigma') \neq 0$, we conclude that $\operatorname{Hom}(\pi|_{G'}, \pi') \neq 0$.
\end{proof}
\begin{remark}
The result of Corollary~\ref{cor:25061806} can be extended to the case of \emph{non-coherent signatures}, as long as the interleaving pattern $[\delta \delta']$ satisfies condition~$(\mu)$ (Definition~\ref{def:23022211}).

For example, consider the symmetry breaking of non-holomorphic discrete series representations for the pair $(G, G') = (U(2,1), U(1,1))$ (see Subsection~\ref{subsubsec:SB_u21}). The interleaving pattern $+ - \oi \op +$ satisfies condition~$(\mu)$, but it does not form a coherent signature, since $q(\delta)\neq q(\delta')$ (Case~IV in Table~\ref{tab:U21}).
Nevertheless, the proof of Corollary~\ref{cor:25061806}, which relies on
using the \emph{stability theorem for multiplicities} in branching laws,
still applies in this case and yields the explicit branching law for this interleaving pattern.
\end{remark}

\subsection{Reduction and Properties of \texorpdfstring{$K$}{K}-types}
\label{subsec:reduction_K_types}

We recall the following from earlier sections:
from Lemma~\ref{lem:q}, the equivalent definitions of the degree $q(\delta)$ and the set $\Phi^{nc,+}(\delta)$ associated with decoration $\delta \in \operatorname{Deco}(p,q)$;
from equation \eqref{eqn:min_K_type_delta}, the minimal $K$-type parameter $\mu_\lambda$;
and from Definition~\ref{def:W_lambda}, the coherent parameter $W_\lambda$.
In this section, we examine the extent to which conditions $(q)$, $(\mu)$, 
 and $(W)$ from Definition \ref{def:23022211}
are preserved under reduction procedures---or equivalently, under expansion procedures.
The main result of this section, Lemma~\ref{lem:23022211} below, completes Step 3 
in Section~\ref{subsec:pf_prop:23022216} 
in the proof of implication (i) $\Rightarrow$ (iii) in Proposition~\ref{prop:23022216}.

We now consider the operation of appending $+$ to the decoration $\delta$ and $\op$ to the decoration $ \delta'$. 
We use the following lemma, which is a direct consequence of the definition of degree $q(\delta)$ given by \eqref{eqn:q_delta}.
\begin{lemma}
\label{lem:23020112}
\begin{enumerate}
\item[{\rm{(1)}}]
Let $\delta \in \operatorname{Deco}(p,q)$.  

We consider $\widetilde\delta =(-\delta)$ or 
 $(\delta-) \in \operatorname{Deco}(p,q+1)$.  
Then 
\[
  q(\widetilde \delta) - q(\delta)
  =
  \begin{cases}
  0  \qquad &\text{if }\quad \widetilde \delta=-\delta, 
\\
    p  \qquad &\text{if }\quad \widetilde \delta=\delta-.  
\end{cases}
\]

\item[{\rm{(2)}}]
We consider $\widetilde\delta =(+\delta)$ or 
 $(\delta+) \in \operatorname{Deco}(p+1,q)$.  
Then 
\[
  q(\widetilde \delta) - q(\delta)
  =
  \begin{cases}
  q  \qquad &\text{if }\quad \widetilde \delta=+\delta, 
\\
    0  \qquad &\text{if }\quad \widetilde \delta=\delta+.  
\end{cases}
\]

\item[{\rm{(3)}}]
Let $\delta \in \operatorname{Deco}(p,q)$
 and $\delta' \in \operatorname{Deco}(p-1,q)$.     
Then 
\[
  q(\delta)-q(\delta')=q(+\delta)-q(\op\delta')=q(-\delta)-q(\oi \delta')
  =q(\delta+)-q(\delta'\op).  
\]

\end{enumerate}
\end{lemma}

There are multiple ways to extend a given signature
 $[\delta\delta']$ via such operations.
 The following lemma describes how the conditions~$(\mu)$ and~(W) are affected under these extensions, in the specific context we are concerned with.

\begin{lemma}[Expansion of Parameters]
\label{lem:23022104}
Let $[\delta\delta']$ be a signature satisfying the GGP interleaving pattern, 
 where $\delta \in \operatorname{Deco}(p, q)$, 
$\delta' \in \operatorname{Deco}(p-1,q)$.  

Let $(\tlambda, \lambda')\in (\Z_\mathrm{int}^{p+q}(\delta), \Z_\mathrm{int}^{p+q-1}(\delta'))$ be a pair of discrete series parameters of $G$ and $G'$ as in Proposition~\ref{prop:discpara} such that the associated signature is $[delta \delta']$.
We assume that condition (W) is satisfied in 
Definition~\ref{def:23022211}; that is,
$[W_{\tlambda}:\Wlmdp]\ne 0$. 

\begin{enumerate}
\item[{\rm(1)}]
($G=U(p,q) \rightsquigarrow \widetilde{G}=U(p+1,q)$).
We define 
\[
\widetilde \delta := +\delta \in \operatorname{Deco}(p+1, q), \quad
\widetilde \delta' := \op \delta' \in \operatorname{Deco}(p, q).
\]
Let
\[
\widetilde \lambda := \left( \lambda_0, \tlambda - \tfrac{1}{2} \mathbf{1}_{p+q} \right), \quad
\widetilde \lambda' := \left( \lambda_0', \lambda' - \tfrac{1}{2} \mathbf{1}_{p+q-1} \right),
\]
subject to the conditions given in Proposition~\ref{prop:discpara}, namely,
\begin{align*}
\lambda_0 &\in \mathbb{Z} + \tfrac{1}{2}(p+q), 
& \lambda_0 &> \max_{1 \le i \le p+q} \tlambda_i, \\
\lambda_0' &\in \mathbb{Z} + \tfrac{1}{2}(p+q-1),
& \lambda_0' &> \max_{1 \le j \le p+q-1} \lambda_j'.
\end{align*}

Then, the following equivalences hold:
\begin{align*}
[\mu_{\widetilde \lambda} : \mu_{\widetilde \lambda'}] \ne 0
&\iff 
[\mu_{\tlambda} : \mulmdp] \ne 0 
\quad \text{and} \quad 
\lambda_0 > \lambda_0' > \lambda_1 + \rho(\tlambda)_1 - \tfrac{p+q}{2}, \\
[W_{\widetilde \lambda} : W_{\widetilde \lambda'}] \ne 0
&\iff 
[W_{\tlambda} : \Wlmdp] \ne 0 
\quad \text{and} \quad 
\lambda_0 > \lambda_0' > \tlambda_1.
\end{align*}

\item[{\rm(2)}]
($G=U(p,q) \rightsquigarrow \widetilde{G}=U(p,q+1)$).
We define
\[
\widetilde \delta := -\delta \in \operatorname{Deco}(p, q+1), \quad
\widetilde \delta' := \oi \delta' \in \operatorname{Deco}(p-1, q+1).
\]
Let
\begin{align*}
\widetilde \lambda 
&:= \left( \tlambda^{(p)} - \tfrac{1}{2} \mathbf{1}_p;\ \lambda_0,\ \tlambda^{(q)} - \tfrac{1}{2} \mathbf{1}_q \right), \\
\widetilde \lambda' 
&:= \left( \lambda'^{(p-1)} - \tfrac{1}{2} \mathbf{1}_{p-1};\ \lambda_0',\ \lambda'^{(q)} - \tfrac{1}{2} \mathbf{1}_q \right),
\end{align*}
subject to the condition given in Proposition~\ref{prop:discpara}, namely,
\begin{align*}
\tlambda_0 &\in \mathbb{Z} + \tfrac{1}{2}(p+q), 
& \tlambda_0 &> \max_{1 \le i \le p+q} \tlambda_i, \\
\lambda_0' &\in \mathbb{Z} + \tfrac{1}{2}(p+q-1),
& \lambda_0' &> \max_{1 \le j \le p+q-1} \lambda_j'.
\end{align*}

Then the following equivalences hold:
\begin{align*}
[\mu_{\widetilde \lambda} : \mu_{\widetilde \lambda'}] \ne 0 
&\iff 
[\mu_{\tlambda} : \mulmdp] \ne 0
\quad \text{and} \quad 
\lambda_0 = \lambda_0' - \tfrac{1}{2}, \\
[W_{\widetilde \lambda} : W_{\widetilde \lambda'}] \ne 0 
&\iff 
[W_{\tlambda} : \Wlmdp] \ne 0 
\quad \text{and} \quad 
\tlambda_0 = \lambda_0' + \tfrac{1}{2}.
\end{align*}

\item[{\rm(3)}]
($G=U(p,q) \rightsquigarrow \widetilde{G}=U(p+1,q)$).
We define
\[
\widetilde \delta := \delta+ \in \operatorname{Deco}(p+1, q), \quad
\widetilde \delta' := \delta' \op \in \operatorname{Deco}(p, q).
\]
Let
\[
\widetilde \lambda := \left( \lambda_0,\ \tlambda + \tfrac{1}{2} \mathbf{1}_{p+q} \right), \quad
\widetilde \lambda' := \left( \lambda_0',\ \lambda' + \tfrac{1}{2} \mathbf{1}_{p+q-1} \right),
\]
subject to the condition given in Proposition~\ref{prop:discpara}, namely,
\begin{align*}
\tlambda_0 &\in \mathbb{Z} + \tfrac{1}{2}(p+q), 
& \tlambda_0 &< \min_{1 \le i \le p+q} \tlambda_i, \\
\lambda_0' &\in \mathbb{Z} + \tfrac{1}{2}(p+q-1),
& \lambda_0' &< \min_{1 \le j \le p+q-1} \lambda_j'.
\end{align*}

Then the following equivalences hold:
\begin{align*}
[\mu_{\widetilde \lambda} : \mu_{\widetilde \lambda'}] \ne 0 
&\iff 
[\mu_{\tlambda} : \mulmdp] \ne 0
\quad \text{and} \quad 
\lambda_p + \rho(\lambda)_p + \tfrac{1}{2}(p+q-1) > \lambda_0' > \lambda_0, \\
[W_{\widetilde \lambda} : W_{\widetilde \lambda'}] \ne 0 
&\iff 
[W_{\tlambda} : \Wlmdp] \ne 0
\quad \text{and} \quad 
\tlambda_p > \lambda_0' > \lambda_0.
\end{align*}
\end{enumerate}
\end{lemma}

Prior to the proof of Lemma~\ref{lem:23022104},
we prepare some necessary computation.

We write
\begin{alignat*}{2}
\rho_c=&(\rho_c^{(p)}; \rho_c^{(q)})
=&&(\tfrac{p-1}{2}, \ldots, \tfrac{1-p}{2};\tfrac{q-1}{2}, \ldots, \tfrac{1-q}{2}),
\\
\rho_{\mathfrak g}=&(\rho_{\mathfrak g}^{(p)}; \rho_{\mathfrak g}^{(q)})
=&&(\tfrac{p+q-1}{2}, \ldots, \tfrac{1+q-p}{2};\tfrac{q-p-1}{2}, \ldots, \tfrac{1-p-q}{2}).
\end{alignat*}

Let $\delta \in \operatorname{Deco}(p,q)$
 denote the decoration of $\tlambda \in {\mathbb{Z}}_{\operatorname{int},>}^{p+q}$, as given in Proposition~\ref{prop:discpara}.  
We extend the parameter space by one dimension.
In the first example, the new component is placed in the first coordinate; in the second, it is placed in the $(p+1)$st coordinate; and in the third, it is placed in the last coordinate.

As a preliminary step, we derive the following identities directly from the definitions.

(1)\enspace
($G=U(p,q) \rightsquigarrow \widetilde{G}=U(p+1,q)$).
Let $\widetilde \delta := +\delta \in \operatorname{Deco}(p+1, q)$, and let 
$\lambda_0 \in \mathbb{Z} + \tfrac{1}{2}(p+q)$ be such that
\[
\widetilde{\lambda} := \left( \lambda_0,\ \tlambda - \tfrac{1}{2} \mathbf{1}_{p+q} \right) 
\in \mathbb{Z}_{\operatorname{int}, >}^{p+q+1}.
\]
Then we have
\begin{equation}
\begin{aligned}
\mu_{\widetilde{\lambda}} 
&= \left( \lambda_0 + \tfrac{q - p}{2}, \ \mu_{\tlambda}^{(p)};\ \mu_{\tlambda}^{(q)} \right), \\
W_{\widetilde{\lambda}} 
&= \left( \lambda_0 - \tfrac{p+q}{2}, \ W_{\tlambda}^{(p)};\ W_{\tlambda}^{(q)} \right),
\label{eqn:+Wlmd}
\end{aligned}
\end{equation}
since the following identities hold:
\begin{align*}
\rho(\widetilde{\lambda}) 
&= \left( \tfrac{p+q}{2},\ \rho(\tlambda) - \tfrac{1}{2} \mathbf{1}_{p+q} \right), \\
\widetilde{\rho}_c 
&= \left( \tfrac{p}{2},\ \rho_c^{(p)} - \tfrac{1}{2} \mathbf{1}_p \right) \oplus \rho_c^{(q)}, \\
\widetilde{\rho}_{\mathfrak{gl}_{p+q+1}}
&= \left( \tfrac{p+q}{2}, \ \rho_{\mathfrak g} - \tfrac{1}{2} \mathbf{1}_{p+q} \right).
\end{align*}

(2)\enspace
($G=U(p,q) \rightsquigarrow \widetilde{G}=U(p,q+1)$).
Let $\widetilde \delta = -\delta \in \operatorname{Deco}(p, q+1)$, and 
\[
\widetilde{\lambda} = \bigl(\tlambda^{(p)} - \tfrac{1}{2} \mathbf{1}_p;\ \lambda_0;\ \tlambda^{(q)} - \tfrac{1}{2} \mathbf{1}_q \bigr).
\]

Then we have
\begin{equation}
\begin{aligned}
\mu_{\widetilde{\lambda}} 
&= \bigl(\mu_{\tlambda}^{(p)} - \mathbf{1}_p;\ \lambda_0 + \tfrac{p - q}{2};\ \mu_{\tlambda}^{(q)} \bigr), \\
W_{\widetilde{\lambda}} 
&= \bigl(W_{\tlambda}^{(p)};\ \lambda_0 + \tfrac{p - q}{2};\ W_{\tlambda}^{(q)} \bigr),
\label{eqn:-Wlmd}
\end{aligned}
\end{equation}
since the following identities hold:
\begin{align*}
\rho(\widetilde{\lambda})
&= \bigl(\rho^{(p)}(\tlambda) - \tfrac{1}{2} \mathbf{1}_p;\ \tfrac{p+q}{2};\ \rho^{(q)}(\tlambda) - \tfrac{1}{2} \mathbf{1}_q \bigr), \\
\widetilde{\rho}_c
&= \bigl(\rho_c^{(p)};\ \tfrac{q}{2};\ \rho_c^{(q)} - \tfrac{1}{2} \mathbf{1}_q \bigr), \\
\widetilde{\rho}_{\mathfrak{gl}_{p+q+1}}
&= \bigl(\rho_{\mathfrak g}^{(p)} - \tfrac{1}{2} \mathbf{1}_p;\ \tfrac{p+q}{2};\ \rho_{\mathfrak g}^{(q)} - \tfrac{1}{2} \mathbf{1}_q \bigr).
\end{align*}

(3)\enspace
($G=U(p,q) \rightsquigarrow \widetilde{G}=U(p+1,q)$).
Let \(\widetilde \delta = \delta+ \in \operatorname{Deco}(p+1, q)\), and
\[
\widetilde{\lambda} = \bigl(\tlambda - \tfrac{1}{2} \mathbf{1}_p;\ \lambda_0;\ \tlambda^{(q)} - \tfrac{1}{2} \mathbf{1}_q \bigr).
\]

Then we have
\begin{equation}
\begin{aligned}
\mu_{\widetilde{\lambda}} 
&= \bigl(\mu_{\tlambda}^{(p)} - \mathbf{1}_p;\ \lambda_0 + \tfrac{p - q}{2};\ \mu_{\tlambda}^{(q)} \bigr),
\\
W_{\widetilde{\lambda}} 
&= \bigl( W_{\tlambda}^{(p)} - \mathbf{1}_p;\ \lambda_0 + \tfrac{p - q}{2};\ W_{\tlambda}^{(q)} \bigr).
\end{aligned}
\label{eqn:Wlmd+}
\end{equation}
since the following identity holds:
\begin{align*}
\widetilde{\rho}_{\mathfrak{gl}_{p+q+1}}
&= \bigl( \rho_{\mathfrak g}^{(p)} + \tfrac{1}{2} \mathbf{1}_p;\ \tfrac{-p + q}{2};\ \rho_{\mathfrak g}^{(q)} - \tfrac{1}{2} \mathbf{1}_q \bigr).
\end{align*}

\begin{proof}[Proof of Lemma~\ref{lem:23022104}]
\begin{enumerate}
\item[{\rm{(1)}}]
Under the assumption that $[W_{\tlambda}:\Wlmdp]\ne 0$, 
it follows from the formula  \eqref{eqn:+Wlmd} that \[
[W_{\widetilde\lambda}:W_{\widetilde\lambda'}]\ne 0
\]
holds if and only if
\[\lambda_0-\tfrac{p+q}{2} \ge \lambda_0'-\tfrac{p+q-1}{2} \ge \tlambda_1-\tfrac{p+q-1}{2},
\]
 equivalently,
 \[\lambda_0-\tfrac{1}{2} \ge \lambda_0'\ge \lambda_1.
 \]
\item[{\rm{(2)}}]
We use \eqref{eqn:-Wlmd}.  
Look at the $U(q+1)$-modules.  
\[
   \lambda_0-\tfrac{p+q}{2} = \lambda_0'-\tfrac{p+q-1}{2}
\Leftrightarrow
   \lambda_0=\lambda_0'+\tfrac{1}{2}.  
\]
\item[{\rm{(3)}}]
We use \eqref{eqn:Wlmd+}.  
$\lambda_p-\frac{-p+q+1}{2}-1 \ge \lambda_0'+\frac{p-q-1}{2} \ge \lambda_0+\frac{p-q}{2}$
$\Leftrightarrow$
 $\tlambda_p \ge \lambda_0'+1 \ge \lambda_0+\frac 3 2$.
Taking the parity condition into account, 
 this is equivalent to $\tlambda_p > \lambda_0' > \lambda_0$.  
\end{enumerate}
\end{proof}

The following lemma is derived from Lemma~\ref{lem:23022104}.

\begin{lemma}[Expansion of Interleaving Pattern]
\label{lem:23022211}
Suppose that a signature $[\delta\delta']$ is a
 GGP interleaving pattern (Definition~\ref{def:interlace}), where
 $\delta \in \operatorname{Deco}(p,q)$ and
 $\delta' \in \operatorname{Deco}(p-1,q)$.   
\begin{enumerate}
\item[{\rm{(1)}}] \rm{(Case $+\op [\delta \delta']$
for $G=U(p,q) \rightsquigarrow \widetilde{G}=U(p+1,q)$
)}.
We prepend $+\op$ to the signature $[\delta \delta']$ to form a new signature 
\[
[\widetilde \delta \widetilde \delta'] := +\op [\delta \delta'],
\]
where $\widetilde \delta := +\delta \in \operatorname{Deco}(p+1, q)$  
and $\widetilde \delta' := \op \delta' \in \operatorname{Deco}(p, q)$.  
Then, the conditions $(\mu)$ and $(W)$ are preserved:
\begin{align*}
\text{$(\mu)$ holds for $[\delta \delta']$}
&\iff \text{$(\mu)$ holds for $+\op [\delta \delta']$,}
\\
\text{$(W)$ holds for $[\delta \delta']$}
&\iff \text{$(W)$ holds for $+\op [\delta \delta']$.}
\end{align*}

\item[{\rm{(2)}}]
 \rm{(Case $\oi - [\delta \delta']$
 for $G=U(p,q) \rightsquigarrow \widetilde{G}=U(p,q+1)$)}.
We prepend $\oi -$ to the signature $[\delta \delta']$ to form a new signature 
\[
[\widetilde \delta \widetilde \delta'] := \oi - [\delta \delta'],
\]
where
$\widetilde \delta := -\delta \in \operatorname{Deco}(p, q+1)$  
and $\widetilde \delta' := \oi \delta' \in \operatorname{Deco}(p-1, q+1)$.  
Then, the conditions $(\mu)$ and $(W)$ are preserved:
\begin{align*}
\text{$(\mu)$ holds for $[\delta \delta']$}
&\iff\text{$(\mu)$ holds for $\oi- [\delta\delta']$,}
\\
\text{$(W)$ holds for $[\delta \delta']$}
&\iff\text{$(W)$ holds for $\oi-[\delta\delta']$.}
\end{align*}

\item[{\rm{(3)}}]
\rm{(Case $[\delta \delta'] \op +$ for
$G=U(p,q) \rightsquigarrow \widetilde{G}=U(p+1,q)$)}.
We append $\op +$ to the signature $[\delta \delta']$ to form a new signature
\[
[\widetilde \delta \widetilde \delta']:=[\delta \delta']\op+,
\]
where
$\widetilde \delta := \delta + \in \operatorname{Deco}(p+1, q)$  
and $\widetilde \delta' :=  \delta' \op \in \operatorname{Deco}(p, q)$.  
Then,  the conditions $(\mu)$ and $(W)$ are preserved:
\begin{align*}
\text{$(\mu)$ holds for $[\delta \delta']$}
&\iff\text{$(\mu)$ holds for $[\delta\delta']\op+$,}
\\
\text{$(W)$ holds for $[\delta \delta']$}
&\iff\text{$(W)$ holds for $[\delta\delta']\op+$.}
\end{align*}

Thus, Step~3 in Section~\ref{subsec:pf_prop:23022216} is now complete.  
With this, the proof of Proposition~\ref{prop:23022216} is fully established.

\end{enumerate}
\end{lemma}

\vskip 1pc

\medskip

  \medskip

\medskip

\newpage
                         
\part{A Translation Functor and Restrictions of Coherent Cohomology}
\label{part:translation}

\section{Symmetry Breaking under Translation}
\label{sec:SBO_translation}
~~~
\newline
In representation theory of reductive groups
the \emph{translation functor} has been a powerful tool for studying families of representations of a \emph{single} group.
However, a systematic framework for applying it to \emph{symmetry breaking} phenomena---arising from the restriction of representations between a \emph{pair of groups} $G \supset G'$---has so far been lacking. 
For example, the vanishing condition for symmetry breaking is rather delicate, as it is not preserved under translation functors for the product group $G \times G'$, even when the parameters lie within the regular dominant chambers.
 
In this section, we introduce a new framework for translation functors that applies not only to representations of individual groups
 but also to \emph{symmetry breaking} arising from the restriction  $G \downarrow G'$, where $(G,G')=(U(p,q), U(p-1,q))$.  
 
The tool developed here is formulated 
in terms of the pair of complexified Lie algebras 
 $({\mathfrak{g}}, {\mathfrak{g}'})$.
For convenience, we adopt the following convention
 throughout this section:
 \begin{equation}
 \label{eqn:n=pq1}
     n=p+q-1.
 \end{equation}
With this convention, $({\mathfrak{g}}, {\mathfrak{g}'})=({\mathfrak{gl}}(n+1,\C), {\mathfrak{gl}}(n,\C))$.

The main result of this section is Theorem~\ref{thm:23081404}.

\subsection{Translation Functor and Symmetry Breaking---Motivation from Coherent Cohomology}
\label{subsec:MFO23}

~~~

Motivated by the main theme of this paper, we begin by considering the following question.

\begin{question}
\label{q:mfo1}
Let $(G, G')=(U(p,q), U(p-1,q))$.
Let $(\pi,\pi')$ be a coherent pair
with associated signature $[\delta \delta']$ as in Definition~\ref{def:coherent_pair}.
Assume that their Harish-Chandra parameters 
are sufficiently far from the walls. 
Can we find an elementary pair 
$(\sigma, \sigma')$ having the same associated signature $[\delta \delta']$,
and a pair of 
finite-dimensional irreducible representations $(F,F')$ of $(G,G')$, 
such that the following three conditions are satisfied?
\begin{enumerate}
\item $q(\pi) = q(\pi') = q(\sigma) = q(\sigma')$.
\item $\operatorname{Hom}_{G'}(F|_{G'},F') \neq 0$.
\item The composition 
 of the following three maps in \eqref{eqn:Fdiag2} is non-zero:
\begin{equation}
\label{eqn:Fdiag2}
\xymatrix{
   \sigma \otimes F
     \ar[r]^{T \otimes \operatorname{pr}_F}
     &
     \sigma' \otimes F'
     \ar[d]
\\
     \pi
     \ar@{^{(}-_>}[u]
     &\pi' 
     }
\end{equation}

Here, the top arrow is the tensor product of non-vanishing elements of $T\in \operatorname{Hom}_{G'}(\sigma|_{G'},\sigma')$ and $\operatorname{pr}_{F} \in \operatorname{Hom}_{G'}(F|_{G'},F')$.
\end{enumerate}

\end{question}

Condition (1) of Question~\ref{q:mfo1} is motivated by the study of morphisms between $(\fP,K)$-cohomologies arising from symmetry breaking (see Section~\ref{sec:Pk}).
In the context we consider here, this coincidence of degrees indeed holds, as shown in Proposition~\ref{prop:23022216} and Definition~\ref{def:coh_deco}.

Question~\ref{q:mfo1} is also concerned with finding an elementary pair $(\sigma, \sigma')$.
This part was addressed in the previous sections, specifically in Theorems~\ref{thm:23022420b} and~\ref{thm:230816}.

Condition~(3) of Question~\ref{q:mfo1}, however, presents significantly a serious problem.
To analyze this, we consider a more general setting by dropping the assumption that $\pi, \sigma, \pi', \sigma'$ are discrete series representations.
Furthermore, we set aside the specific pair $(G,G')=(U(p,q), U(p-1,q)$,
and instead work purely with the complexified Lie algebras $(\mathfrak{g}, \mathfrak{g'})=(\mathfrak{gl}_{n+1}, \mathfrak{gl}_n)$, where $n=p+q-1$.
This conceptual shift allows us to formulate the problem in a broader and more algebraic framework,
which, to the best of our knowledge, has not been previously explored in the literature.  
This broader approach is developed from Section~\ref{subsec:SBO_translate} to Section~\ref{subsec:pf_thm:23081343}.
Theorem~\ref{thm:23081404} provides a positive result for Condition~(3) in the special setting where $F$ is the standard representation or its dual, and $F'$ is the trivial representation.

In Section~\ref{subsec:string}, we discuss iterated applications of suitable chosen pairs of \emph{small} representations $(F, F')$
to connect two pairs $(\pi, \pi')$ and $(\sigma, \sigma')$,
see Question~\ref{q:mfo2}.

\medskip
\subsection{Non-Vanishing Translation of Symmetry Breaking}
\label{subsec:SBO_translate}
~~~
\newline

We keep convention in \eqref{eqn:n=pq1}, which means that the complexified Lie algebras corresponding to the pair $(G,G')=(U(p,q), U(p-1,q))$ are given by
\[
(\mathfrak{g}, \mathfrak{g}')=({\mathfrak{g l}}_{n+1}({\mathbb{C}}), {\mathfrak{g l}}_{n}({\mathbb{C}})).
\]
Accordingly, in this section, we denote characters of the Cartan subalgebras of $\mathfrak{g} \supset \mathfrak{g}'$ by
\[
\lambda=(\lambda_0, \lambda_1, \dots, \lambda_n),
\quad
\lambda'=(\lambda'_1, \dots, \lambda'_n),
\]
with respect to the standard bases
$\{e_0, e_1, \dots, e_n\}$ and
$\{e'_1, \dots, e'_n\}$,
rather than
\[
\lambda=(\lambda_1, \dots, \lambda_{p+q}),
\quad
\lambda'=(\lambda'_1, \dots, \widehat{\lambda'_p}, \dots, \lambda'_{p+q})
\]
as in the previous sections.
This change in convention should not cause confusion, since $\lambda$ and $\lambda'$  are used here basically as parameters in the Harish-Chandra isomorphisms, where they are considered only up to Weyl group action, as we now recall.

Let $N=n+1$ or $N=n$.
We normalize the Harish-Chandra isomorphism
\[
   \operatorname{Hom}_{{\mathbb{C}}\operatorname{-alg}}({\mathfrak{Z}}({\mathfrak{g l}}_N), {\mathbb{C}}) \simeq {\mathbb{C}}^N/{\mathfrak{S}}_N
\]
 such that the trivial one-dimensional representation 
 of $\mathfrak{g}={\mathfrak{g l}}_{N}$ has 
 the infinitesimal character 
 $\rho_{\mathfrak g}=(\frac{N-1}{2}, \frac{N-3}{2}, \dots, \frac{1-N}{2})$.

 Suppose that $\tau\in \operatorname{Hom}_{{\mathbb{C}}\operatorname{-alg}}({\mathfrak{Z}}({\mathfrak{g l}}_N), {\mathbb{C}}) \simeq {\mathbb{C}}^N/{\mathfrak{S}}_N$.
 For a $\mathfrak{g}$-module $V$, let 
 $P_\tau(V)$ denote the $\tau$-\emph{primary component} of V; that is, 
 \begin{equation}
 \label {eqn:pr_gen_inf}
  P_{\tau}(V):=\bigcup_{k=1}^{\infty}\, 
               \bigcap_{z \in {\mathfrak{Z}}(G)} \operatorname{Ker}(z - \tau(z))^k.  
 \end{equation}

\begin{setting}
\label{set:Fpi}
Suppose that
 $(\mathfrak{g}, \mathfrak{g}')=
({\mathfrak{g l}}_{n+1}({\mathbb{C}}), {\mathfrak{g l}}_{n}({\mathbb{C}}))$.  
Let $\pi \in \operatorname{Irr}(G)$
 and $\pi' \in \operatorname{Irr}(G')$.  
Let $\lambda=(\lambda_0, \dots, \lambda_n)$ and $\lambda'=(\lambda_1', \dots, \lambda_n')$ denote
 the infinitesimal characters of $\pi$ and $\pi'$, respectively, via the Harish-Chandra isomorphism.  
Let $T \colon \pi \to \pi'$ be a non-zero $G'$-homomorphism.  
Consider the standard representation of $G$
 on $F={\mathbb{C}}^{n+1}$, or its contragredient representation $({\mathbb{C}}^{n+1})^{\vee}$.
We decompose $F=F' \oplus F''$,
 where $F'$ is the unique trivial one-dimensional $G'$-submodule and $F''$ is the $n$-dimensional submodule.  
Let $\{f_0, \dots, f_n\}$ be the standard basis of $F$ such that $F'={\mathbb{C}}f_0$
 and $F''=\operatorname{span}_{\mathbb{C}}\{f_1, \dots, f_n\}$.  
\end{setting}

We denote by $\operatorname{pr}_{F\to F'} \colon F \to F/F'' \simeq F'$ the natural projection, 
 and identify $\pi' \otimes F'$ with $\pi'$ as $G'$-modules via the isomorphism $F' \simeq {\mathbb{C}}$.  
(Hereafter, we simply write $\operatorname{pr}_F$ for $\operatorname{pr}_{F \to F'}$.)

The following theorem gives a sufficient condition on the generalized ${\mathfrak{Z}}(\mathfrak{g})$-eigenspaces in $\pi \otimes F$
 for the restriction of the linear map $T \otimes \operatorname{p r}_{F \to F'} \colon \pi \otimes F \to \pi'\otimes F'$ to be nonzero.

\begin{theorem}
\label{thm:23081404}
Suppose that any generalized eigenspaces of ${\mathfrak{Z}}(\mathfrak{g})$
 in $\pi \otimes F$ are eigenspaces.  
\par\noindent
{\rm{(1)}}\enspace
Let $F={\mathbb{C}}^{n+1}$.  
For any $i$ such that 
$\lambda_i \not \in \{\lambda_1'-\frac 1 2, \lambda_2'-\frac 1 2, \dots, \lambda_n'-\frac 1 2\}$, 
 $T \otimes \operatorname{pr}_{F}$ does not vanish on 
 the primary component $P_{\lambda+e_i}(\pi \otimes {\mathbb{C}}^{n+1})$.  
\par\noindent
{\rm{(2)}}\enspace
Let $F=({\mathbb{C}}^{n+1})^{\vee}$.  
For any $i$ such that 
 $\lambda_i \not \in \{\lambda_1'+\frac 1 2, \lambda_2'+\frac 1 2, \dots, \lambda_n'+\frac 1 2\}$, 
 $T \otimes \operatorname{pr}_{F}$ does not vanish on  the primary component
 $P_{\lambda-e_i}(\pi \otimes ({\mathbb{C}}^{n+1})^{\vee})$.  
\end{theorem}

\begin{remark}
See Proposition \ref{prop:prireg}
 for a setting where any generalized eigenspaces of ${\mathfrak{Z}}(\mathfrak{g})$
 in $\pi \otimes F$ are eigenspaces.  
\end{remark}

In this paper, we apply Theorem~\ref{thm:23081404} to the pair $(G,G')=(U(p,q), U(p-1,q))$.
It is worth noting that the theorem also applies to other real forms, such as $(G L_{n+1}({\mathbb{R}}), G L_{n}({\mathbb{R}}))$; see \cite{KS} for further applications.

The rest of this section is devoted to the proof of Theorem \ref{thm:23081404}.  
\subsection{Separation of \texorpdfstring{$\mathfrak{Z}(\mathfrak{g})$}{Z(g)}-Infinitesimal Eigenspaces}
\label{subsec:powerCG}
~~~
\newline
The primary objective of this section is to present Theorem~\ref{thm:23081343}, which subsequently leads to Theorem~\ref{thm:23081404}.

Our key idea is to harness only the power of the Casimir element in order to perform the translation functor. At first glance, this may appear unrealistic, since the center $\mathfrak{Z}(\mathfrak{g})$ is a polynomial ring with $n+1$ algebraically independent elements, while the Casimir element accounts for only one of them. However, when $F$ is a small representation, it turns out that the $\mathfrak{Z}(\mathfrak{g})$-character of the tensor product with $F$ can be distinguished solely through the Casimir element. Moreover, this method enables us to trace the action of symmetry breaking operators, which leads to the formulation of
Theorem~\ref{thm:23081343}.

We define the Casimir element $c_G$
 for ${\mathfrak{g}}={\mathfrak{g l}}_{n+1}({\mathbb{C}})$
 with respect to the non-degenerate symmetric bilinear form
\[
  {\mathfrak{g l}}_{n+1} \times {\mathfrak{g l}}_{n+1} \to {\mathbb{C}}, 
\quad
  (X,Y) \mapsto \operatorname{Trace}(XY).  
\]
 by 
\begin{equation}
\label{eqn:Casimir}
   c_G:= \sum_{j=0}^n \sum_{i=0}^n E_{i j} E_{j i} \in U(\mathfrak{g}).  
\end{equation}

Likewise, 
 we define the Casimir element for $\mathfrak{g}'={\mathfrak{g l}}_n({\mathbb{C}})$ by
\[
     c_{G'}:= \sum_{j=1}^n \sum_{i=1}^n E_{i j} E_{j i} \in U(\mathfrak{g}').  
\]

Then $c_G$ acts on $\pi_{\lambda}$ and on $F={\mathbb{C}}^{n+1}$ or its contragredient $({\mathbb{C}}^{n+1})^{\vee}$ 
 as scalar multiplication by
\begin{align}
\label{eqn:A_lmd}
 A\equiv&A(\lambda):= \sum_{i=0}^n \lambda_i^2-|\rho_{\mathfrak g}|^2,
\quad
\text{and \quad $n+1$,}
\intertext{respectively. On the other hand, $c_{G'}$ acts on $\pilmdp$ by}
\label{eqn:B_lmd}
 B \equiv & B(\lambda'):= \sum_{i=1}^n {\lambda'}_i^2-|\rho_{\mathfrak{g}'}|^2.  
\end{align}

For $0 \le i \le n$ and $\varepsilon \in \{+, -\}$, 
 we define linear operators 
 \[
 \operatorname{pr}_{i, \varepsilon} \colon \pi \otimes F \to \pi \otimes F
 \]
 by 
\begin{align}
\label{eqn:pri+}
\operatorname{pr}_{i, +}
:=&\prod_{j \ne i} (c_G - |\lambda+e_j|^2+|\rho_{\mathfrak g}|^2)
\quad
\text{for $F={\mathbb{C}}^{n+1}$}, 
\\
\label{eqn:pri-}
\operatorname{pr}_{i, -}
:=&\prod_{j \ne i} (c_G - |\lambda-e_j|^2+|\rho_{\mathfrak g}|^2)
\quad
\text{for $F=({\mathbb{C}}^{n+1})^{\vee}$}, 
\end{align}
where the product is taken over $j \in \{0, 1, \dots, n\} \setminus \{i\}$.

Here we describe the meaning of $\operatorname{pr}_{i, \pm}$
 in the case 
 where $(G,G')=(U(p,q), U(p-1,q))$
 and $(\pi, \pi')$ are discrete series representations
 of $G$ and $G'$, 
 respectively.   

We begin by recalling some general terminology.
Let $\mu$ be an extremal weight of an irreducible finite-dimensional representation $F$ of $\mathfrak{g}$.
We define the {\emph{translation functors}} 
in the category of $(\mathfrak{g}, K)$-modules
or 
 in the  category of smooth admissible representations of finite length by
 \begin{equation}
\label{eqn:tr_functor}
    \psi_{\tau}^{\tau+\mu}(\cdot)
    :=
    P_{\tau+\mu}(P_{\tau}(\cdot) \otimes F).
\end{equation}

\begin{proposition}
\label{prop:prireg}
Let $\pi=\pi_{\lambda} \in \operatorname{Irr}(G)$ be the smooth representation associated with a discrete series representation of $G=U(p,q)$
 with Harish-Chandra parameter $\lambda$,
 and let $\delta$ be its decoration, as in Definition~
\ref{def:HC_para}.
\par\noindent
{\rm{(1)}}\enspace
Let $F={\mathbb{C}}^{p+q}$.  
Then there is a direct sum decomposition 
 of $G$-modules:
\begin{equation}
\label{eqn:Fdirctsum}
  \pi_{\lambda} \otimes F
  \simeq
 \bigoplus_{i=0}^n \operatorname{p r}_{i,+}(\pi_{\lambda} \otimes F).  
\end{equation}
Moreover, 
 for any $i$ such that $\lambda+e_i$ is regular, 
 there is a natural $G$-isomorphism:
\[
  \pi_{\lambda+e_i} \simeq \psi_{\lambda}^{\lambda+e_i}(\pi_{\lambda})=\operatorname{p r}_{i,+}(\pi_{\lambda} \otimes F),
\]
where $\pi_{\lambda+e_i}$ denotes the smooth representation associated to the discrete series representation with Harish-Chandra parameter $\lambda+e_i$, and with the same decoration $\delta$.
\par\noindent
{\rm{(2)}}\enspace
Let $F=({\mathbb{C}}^{p+q})^{\vee}$.  
Then there is a direct sum decomposition 
 of $G$-modules:
\[
  \pi_{\lambda} \otimes F
  \simeq
 \bigoplus_{i=0}^n \operatorname{p r}_{i,-}(\pi_{\lambda} \otimes F).  
\]
Moreover, 
 for any $i$ such that $\lambda-e_i$ is regular,  there is a natural $({\mathfrak{g}}, K)$-isomorphism:
\[
  \pi_{\lambda-e_i} \simeq \psi_{\lambda}^{\lambda-e_i}(\pi_{\lambda})
 =\operatorname{p r}_{i,-}(\pi_{\lambda} \otimes F),  
\]
where $\pi_{\lambda-e_i}$ denotes the smooth representation associated to the discrete series representation with Harish-Chandra parameter $\lambda-e_i$, and with the same decoration $\delta$.
\end{proposition}

\begin{proof}
Although the argument is standard ({\it{cf}}. \cite{KV, K92})), we include a brief proof here for the reader's convenience.

We realize $\pi_{\lambda}$
 as a cohomological parabolic induction
 ${\mathcal{R}}_{\mathfrak{b}}^S({\mathbb{C}}_{\lambda})$
 as in Lemma \ref{lem:23082617}, 
 where $S=\frac 1 2 ((p-1)p+(q-1)q)$ and ${\mathfrak{b}}={\mathfrak{b}}(\delta)$ is the $\theta$-stable Borel subalgebra associated with signature $\delta$.

\par\noindent
(1)\enspace
Since $\lambda$ belongs to the regular Weyl chamber $C(\delta)$ (see Proposition~\ref{prop:discpara} for notation), we have $\lambda + e_j \in \overline{C(\delta)}$ for any $1 \le j \le p+q$.  
Moreover, since $\lambda$ is regular, we have $\lambda + e_j \not\equiv \lambda + e_k$ in $\mathbb{C}^{p+q}/\mathfrak{S}_{p+q}$ for any $j \ne k$.  
Hence, the spectral sequence in Lemma~\ref{lem:Rspecseq} collapses and yields the following direct sum decomposition as a $({\mathfrak{g}}, K)$-module:
\begin{equation}
\label{eqn:Rfjsum}
 (\pi_{\lambda})_K \otimes  F\simeq
  \bigoplus_{j=1}^{p+q} {\mathcal{R}}_{\mathfrak{b}}^S({\mathbb{C}}_{\lambda+e_j}).  
\end{equation}
Since ${\mathfrak{Z}}(\mathfrak{g})$ acts on each $(\mathfrak{g},K)$-module
 ${\mathcal{R}}_{\mathfrak{b}}^S({\mathbb{C}}_{\lambda+e_j})$ by scalar, 
 each ${\mathfrak{Z}}(\mathfrak{g})$-generalized eigenspace in $\pi_{\lambda} \otimes F$ is actually an eigenspace.

The Casimir element $c_G$ acts on ${\mathcal{R}}_{\mathfrak{b}}^S({\mathbb{C}}_{\lambda+e_j})$ by the scalar $|\lambda+e_j|^2-|\rho_{\mathfrak g}|^2$.
Since $|\lambda+e_j|^2\ne |\lambda+e_k|^2$ 
 for any $j \ne k$,
 the direct summands in \eqref{eqn:Rfjsum} are separated by the action of $c_G$ alone.

Therefore, we have
\begin{equation}
\label{eqn:pr_j_pm}    
\operatorname{p r}_{j, +}(\pi_{\lambda} \otimes F)
 = P_{\lambda+e_j}(\pi_{\lambda} \otimes F)
 (= \psi_{\lambda}^{\lambda+e_j}(\pi_{\lambda} \otimes F)).
 \end{equation}
\par\noindent
(2)\enspace
The proof is similar and is omitted.  
\end{proof}

We now consider the behavior of the operator $T \otimes \operatorname{pr}_F$ on the subspace \eqref{eqn:pr_j_pm}.
\begin{theorem}
\label{thm:23081343}
Suppose we are in Setting \ref{set:Fpi}.  
Then there exist constants $c_\pm \in \mathbb{C}^\times$, independent of $\lambda$ and $\lambda'$ such that 
\[
  (T \otimes \operatorname{pr}_{F})
  \operatorname{pr}_{i,\pm}(u \otimes f_0)
  =
  c_{\pm} \prod_{j=1}^n (\lambda_i - \lambda_j' \pm \frac 1 2 ) T u
\]
for any $u \in \pi$.    
\end{theorem}

\begin{example}
\label{ex:23081337}
For $n=2$, $c_{\pm}=4$.  
Indeed, 
 by using the notation in the proof of Theorem \ref{thm:23081343}, 
 one has 
\[
 \widetilde b^{(2)}
 -2 (\lambda_1+\lambda_2+1) \widetilde b^{(1)}
 +(2\lambda_1+1)(2\lambda_2+1)
=
 4(\lambda_0-\lambda_1'+\frac 12)(\lambda_0-\lambda_2'+\frac 12).  
\]
\end{example}

\begin{remark}
\label{rem:stable_multiplicity}
Building on Theorem~\ref{thm:23081343}, one can prove that 
\[
\dim \operatorname{Hom}_{G}(\pi|_{G'}, \pi')
\]
remains constant for all $\pi$ in the same coherent family, provided that the pair of parameters $(\pi, \pi')$ lie within the same \emph{fence}.  
See [KS, Theorem 3.3] for a precise formulation of the stability theorem for multiplicities in branching laws.
\end{remark}
We shall give a proof of Theorem \ref{thm:23081343}
 in the next section.  
We first show
 that the proof of Theorem \ref{thm:23081404} can be reduced to Theorem \ref{thm:23081343}.  

\begin{proof}
[Proof of Theorem~\ref{thm:23081404}]
(1) By Proposition~\ref{prop:prireg}, 
 one has 
\[
\psi_{\lambda}^{\lambda+e_i}(\pi \otimes F)=
\prod_{j \ne i}(c_G-|\lambda+e_j|^2+|\rho_{\mathfrak g}|^2)(\pi \otimes F)
=\operatorname{p r}_{i, +}(\pi \otimes F), 
\]
where $F={\mathbb C}^{p+q}$.
Hence, it suffices to prove
\[
  (T \otimes \operatorname{p r}_{F \to F'}) \circ \operatorname{p r}_{i,+}
(\pi \otimes F) \ne \{0\}.  
\]
By Theorem \ref{thm:23081343}, 
the non-vanishing follows provided that
 $\lambda_i \not \in \{\lambda_1'-\frac 1 2, \lambda_2'-\frac 1 2, \dots, \lambda_n'-\frac 1 2\}$
 because $T \ne 0$.  

(2)  The proof in the case of $\psi_{\lambda}^{\lambda - e_i}$, 
where $F = (\mathbb{C}^{n+1})^\vee$, 
is analogous and therefore omitted.  
\end{proof}

\subsection{Proof of Theorem \ref{thm:23081343}}
\label{subsec:pf_thm:23081343}
~~~
\newline
This section studies the diagonal action
 of the power of the Casimir element
 on $\pi \otimes F$, 
 and gives a proof of Theorem \ref{thm:23081343}.

We set 
\begin{equation}
\label{eqn:a}
   a \equiv a(\lambda, \lambda'):= \sum_{i=0}^n \lambda_i - \sum_{i=1}^n {\lambda'}_i.  
\end{equation}

\begin{lemma}
\label{lem:23081140}
{\rm{(1)}}\enspace
For any $u \in \pi$, 
 one has
\[  c_G (u \otimes f_0)
 \equiv (A+n+1+2 a)(u \otimes f_0)
 + 2 \sum_{j=1}^n E_{0 j} u \otimes f_j
\]
$\bmod \operatorname{Ker} T \otimes F$.  
\par\noindent
{\rm{(2)}}\enspace
For any $w=\underset{j=1}{\overset n{\sum}} u_j \otimes f_j$
 with $u_1, \dots, u_n \in \pi$, 
 one has
\[
 c_G w
 =
 (A+n+1)w 
 +2(\sum_{j=1}^n E_{j0} u_j) \otimes f_0
 +2 \sum_{j=1}^n (\sum_{i=1}^n E_{i j} u_i) \otimes f_j.  
\]
\end{lemma}

\begin{proof}
By \eqref{eqn:Casimir}, 
 one may express the Casimir element $c_G$ as 
\begin{align}
\notag
c_G=&E_{00}^2 + c_{G'} +\sum_{j=1}^n(E_{j0}E_{0j}+E_{0j}E_{j0})
\\
\label{eqn:23081017}
=& E_{00}^2 + (n+1)E_{00}-\sum_{i=0}^n E_{ii} + c_{G'} 
+ 2 \sum_{j=1}^n E_{j0}E_{0j}.  
\end{align}

For any $u \in \pi$ and $f \in F$, the diagonal action of the Casimir element is given by
\begin{align*}
   c_G(u \otimes f)
  =&
  (c_G u) \otimes f
  +2 \sum_{i=0}^n \sum_{j=0}^n (E_{i j}u) \otimes (E_{j i}f)
  + u \otimes c_G f
\\
  =&(A+n+1) u \otimes f
  + 2 \sum_{i=0}^n \sum_{j=0}^n (E_{i j}u) \otimes (E_{j i}f)
\end{align*}
by the Leibniz rule.  

Observe that $E_{00}u-a u \in \operatorname{Ker} T$, since 
\begin{align*}
T(E_{00}u)
=& T((\sum_{i=0}^nE_{ii}-\sum_{j=1}^nE_{jj})u)
\\
=& (\sum_{i=0}^n \lambda_{i})T u-(\sum_{j=1}^nE_{jj})Tu
\\
=& (\sum_{i=0}^n \lambda_{i}-\sum_{j=1}^n \lambda_j')Tu
=au.  
\end{align*}
Since $E_{i j} f_k = \delta_{j k} f_i$
 $(0 \le i, j , k \le n)$, 
 the lemma follows.  
\end{proof}

In order to analyze the $k$-th power $c_G^k$
 of the diagonal action of $c_G$ on $\pi \otimes F$
 followed by the symmetry breaking operator 
 $T \otimes \operatorname{p r}_{F \to F'}$, 
 we focus on the difference between the powers
 of $c_G$ 
 and those of $c_G'$.  
The above proof of Lemma \ref{lem:23081140}, 
 which is the first step
 in the inductive argument
 for Proposition \ref{prop:23081142} below, 
 has already used
 such a formula \eqref{eqn:23081017}.  
For the general case, 
 we begin with a set-up on the level of the enveloping algebra
 $U({\mathfrak{g l}}_{n+1})$.

For $\ell \in {\mathbb{N}}$
 and $1 \le j \le n$, 
 we define $E_j^{(\ell)}, E^{(\ell)} \in U({\mathfrak{g l}}_{n+1})$
 inductively by 
\begin{align*}
   E_j^{(0)}:=&E_{0j}, 
\\
   E_j^{(\ell)}:=&\sum_{i=1}^n E_{i j} E_i^{(\ell-1)}, 
\qquad
\text{for $\ell \ge 1$}
\\
E^{(\ell)} :=& \sum_{j=1}^n E_{j 0} E_j^{(\ell)}.  
\end{align*}
In particular, we obtain 

\begin{align*}
  E_j^{(\ell)} =&\sum_{1 \le i_1, \dots, i_{\ell} \le n}
  E_{i_{\ell} j} E_{i_{\ell-1} i_{\ell}} \cdots E_{i_1 i_2} E_{0 i_1}, 
\\
  E^{(\ell)} =&
  \sum_{1 \le i_1, \dots, i_{\ell}, i_{\ell+1} \le n}
  E_{i_{\ell+1} 0} E_{i_{\ell} i_{\ell+1}} E_{i_{\ell-1} i_{\ell}}\cdots E_{i_1 i_2} E_{0 i_1}.  
\intertext{For example, }
E^{(0)}=& \sum_{j=1}^n E_{j0} E_{0 j}, 
\\
   E_j^{(1)}=& \sum_{i=1}^n E_{i j} E_{0 i}
\quad\text{and 
$E^{(1)}=\sum_{i=1}^n \sum_{j=1}^n E_{j0} E_{i j} E_{0 i}$}.  
\end{align*}

We recall
 that $(G L_{n+1} \times G L_{n})/\operatorname{diag}(G L_{n})$
 is a spherical variety, 
 and the algebra $U({\mathfrak{g l}}_{n+1})^{{\mathfrak{g l}}_{n}}$
 of ${\mathfrak{g l}}_{n}$-invariant elements
 is commutative.

\begin{lemma}
\label{lem:Hirshen0812}
For any $\ell \in {\mathbb{N}}$, 
 one has 
$E^{(\ell)} \in U({\mathfrak{g l}}_{n+1})^{{\mathfrak{g l}}_{n}}$.  
\end{lemma}

\begin{proof}
[Proof of Lemma \ref{lem:Hirshen0812}]
Let $1 \le a, b \le n$.  
We verify the vanishing of the following.

\begin{equation*}
\operatorname{a d}(E_{a b}) E^{(\ell)}
=\operatorname{a d}(E_{a b}) 
 \sum_{1 \le i_1, \dots, i_{\ell}, i_{\ell+1} \le n} 
 E_{i_{\ell+1}0} E_{i_{\ell} i_{\ell+1}} E_{i_{\ell-1}i_{\ell}}
 \cdots E_{i_1 i_2} E_{0i_1}.  
\end{equation*}
To illustrate the computation, we consider the case $\ell=0$.  
In this case the right-hand side is equal to 
\begin{align*}
&\operatorname{a d}(E_{a b}) \sum_{j=1}^n \sum_{i=1}^n E_{j0}E_{i j} E_{0i}
\\
=& E_{a 0} \sum_{i=1}^n E_{i b} E_{0 i}
   + (\sum_{j=1}^n E_{j0} E_{a j} E_{0 b} - E_{a 0} \sum_{i =1}^n E_{i b} E_{0 i}) - \sum_{j=1}^n E_{j0} E_{a j} E_{0 b}
\\
=&0.  
\end{align*}
The computation for general $\ell$ goes similarly. 
\end{proof}

We now consider an analogue of the Harish-Chandra homomorphism for the algebra $U({\mathfrak{g l}}_{n+1})^{{\mathfrak{g l}}_{n}}$.
Since the natural homomorphism
${\mathfrak{Z}}({\mathfrak{g l}}_{n+1}) \otimes {\mathfrak{Z}}({\mathfrak{g l}}_{n})
\to U({\mathfrak{g l}}_{n+1})^{{\mathfrak{g l}}_{n}}$ is known to be bijective, we obtain the followng two lemmas.
\begin{lemma}
For $\lambda \in {\mathbb{C}}^{n+1}/{\mathfrak{S}}_{n+1}
 \simeq \operatorname{Hom}_{\mathbb{C}\operatorname{-alg}}({\mathfrak{Z}}({\mathfrak{g l}}_{n+1}), {\mathbb{C}})$ 
 and $\lambda' \in {\mathbb{C}}^{n}/{\mathfrak{S}}_{n}
 \simeq \operatorname{Hom}_{\mathbb{C}\operatorname{-alg}}({\mathfrak{Z}}({\mathfrak{g l}}_{n}), {\mathbb{C}})$, 
 one has an algebra homomorphism
 $\chi_{\lambda,\lambda'}$ 
 with the commutative diagram below.  
\begin{equation*}
\xymatrix{
   {\mathfrak{Z}}({\mathfrak{g l}}_{n+1}) \otimes {\mathfrak{Z}}({\mathfrak{g l}}_{n})
     \ar[d]\ar[rd]
     &
\\
     U({\mathfrak{g l}}_{n+1})^{{\mathfrak{g l}}_{n}}
     \ar[r]_{\chi_{\lambda,\lambda'}}
     &{\mathbb{C}}
}
\end{equation*}
\end{lemma}

We note that $\chi_{\lambda, \lambda'}(E^{(\ell)})$ is invariant 
by ${\mathfrak{S}}_{n+1} \times {\mathfrak{S}}_n$.  
\begin{lemma}
\label{lem:23081228}
$\chi_{\lambda, \lambda'}(E^{(\ell)})$ is a polynomial 
 of $\lambda_0$, $\dots$, $\lambda_n$, 
$\lambda_1'$, $\dots$, $\lambda_n'$
 of degree at most $\ell+2$.  
\end{lemma}

\begin{example}
By \eqref{eqn:23081017}, 
\[
E^{(0)} = \sum_{j=1}^n E_{j 0} E_{0 j}
= \frac 1 2(c_G - E_{00}^2 - (n+1)E_{00}
  + \sum_{i=0}^n E_{i i} - c_{G'}).  
\]
In particular, 
 $\chi_{\lambda, \lambda'}(E^{(0)})$ is a symmetric polynomial in the variables $\lambda$ and $\lambda'$, 
 given explicitly by 
\begin{align*}
 2\chi_{\lambda, \lambda'}(E^{(0)})
 =&A-a^2-(n+1)a+\sum_{i=0}^n \lambda_i -B
\\
=& |\lambda|^2-|\rho_{{\mathfrak{gl}}_{n+1}}|^2
  -(\sum_{i=0}^n \lambda_i - \sum_{j=1}^n \lambda_j')^2
  -(n+1)(\sum_{i=0}^n \lambda_i - \sum_{j=1}^n \lambda_j')
\\
  &+ \sum_{i=0}^n \lambda_i 
  - |\lambda'|^2
  +|\rho_{{\mathfrak{gl}}_n}|^2.  
\end{align*}
\end{example}

We keep the notation $A\equiv A(\lambda), B\equiv B(\lambda)$, and $a \equiv a(\lambda, \lambda')$
from \eqref{eqn:A_lmd}, \eqref{eqn:B_lmd}, and \eqref{eqn:a}.
We introduce polynomials $b^{(k)} \equiv b^{(k)}(\lambda, \lambda')$
 and $B_{\ell}^{(k)}\equiv B_{\ell}^{(k)}(\lambda, \lambda')$
 $(0 \le \ell \le k-1)$
as polynomials in the variables $\lambda$ and $\lambda'$
 inductively for each $k \in {\mathbb{N}}$
 by starting with $b^{(0)}:=1$ and
\begin{align}
b^{(k+1)}:=& (A+2a + n + 1) b^{(k)} 
  + 2 \sum_{\ell=0}^{k-1} B_{\ell}^{(k)}\chi_{\lambda, \lambda'}(E^{(\ell)}), 
\label{eqn:bk}
\\
B_{\ell}^{(k+1)}:=&
\begin{cases}
(A+n+1) B_{0}^{(k)} + 2 b^{(k)}
\qquad
&\text{if $\ell=0$, }
\\
(A+n+1) B_{\ell}^{(k)} + 2 B_{\ell-1}^{(k)}
\qquad
&\text{if $1 \le \ell \le k-1$, }
\label{eqn:Bk}
\\
2 B_{k-1}^{(k)}
&\text{if $\ell=k$}.  
\end{cases}
\end{align}

Here, we regard $B_{0}^{(0)} \equiv B_{-1}^{(0)} \equiv 0$.  
\begin{example}
\label{ex:23081236}
{\rm{(1)}}\enspace
For $k=1$, 
 one has 
\[
  b^{(1)} = A + 2a +n +1, 
\quad
 B_{0}^{(1)} =2.  
\]
\par\noindent
{\rm({2)}}\enspace
For $k=2$, 
 one has 
\begin{align*}
  b^{(2)} =& (A + 2a +n +1)^2+2\chi_{\lambda,\lambda'}(E^{(1)}), 
\\
 B_{0}^{(2)} =&3(A+2a+n+1), 
\quad
B_1^{(2)}=2.  
\end{align*}
\end{example}

We now describe the diagonal action of powers of the Casimir element, modulo the kernel of the operator $T \otimes \operatorname{pr}_F$:
\begin{proposition}
\label{prop:23081223}
Suppose we are in Setting~\ref{set:Fpi}.
For any $k \in {\mathbb{N}}$ and $u\in \pi
$, 
 we have
\[
  (T \otimes \operatorname{pr}) c_G^k (u \otimes f_0)
 =b^{(k)}(\lambda,\lambda')\, T u.  
\]
\end{proposition}
This proposition follows directly from the next result:
\begin{proposition}
\label{prop:23081142}
Retain Setting~\ref{set:Fpi}.  
For any $k \in {\mathbb{N}}$, 
 we have
\begin{equation*}
c_G^k(u \otimes f_0)
\equiv  b^{(k)}(\lambda, \lambda') \, u \otimes f_0 
+ \sum_{\ell=0}^{k-1} B_{\ell}^{(k)} \left(\sum_{j=1}^n E_j^{(\ell)} u \otimes f_j\right)
\end{equation*}
$\bmod \operatorname{Ker}T \otimes F$ $(\subset \pi \otimes F)$.  
\end{proposition}

\begin{proof}
[Proof of Proposition~\ref{prop:23081142}]
We prove the formula by induction on $k$.  
The case $k=0$ is clear, since $b^{(0)}=1$.  

Assume the result holds for some $k \in \mathbb{N}$. Then, modulo $\operatorname{Ker} T \otimes F, we have$:
\begin{equation*}
   c_G^{k+1}(u \otimes f_0)
   \equiv
   b^{(k)} c_G (u \otimes f_0)
   + \sum_{\ell=0}^{k-1} B_{\ell}^{(k)} c_G(E_j^{(\ell)} u \otimes f_j).  
\end{equation*}
By Lemma \ref{lem:23081140} (1), 
 the first term on the right-hand side equals
\begin{equation}
\label{eqn:cG1}
(A+n+1+2a) b^{(k)} u \otimes f_0 + 2 b^{(k)} \sum_j E_{0j} u \otimes f_j.  
\end{equation}
On the other hand, 
 by Lemma \ref{lem:23081140} (2), each term
 $c_G(E_j^{(\ell)} u \otimes f_j)$ is equal to 
\begin{align*}
&2 (\sum_{j=1}^n E_{j 0} E_j^{(\ell)} u) \otimes f_0
+ (A+n+1) (E_j^{(\ell)} u \otimes f_j)
+ 2 (\sum_{1 \le j, k \le \ell} E_{j k} E_j^{(\ell)}u) \otimes f_k
\\
=\,&2 E^{(\ell)} u \otimes f_0
+ (A+n+1) (E_j^{(\ell)} u \otimes f_j)
+ 2 \sum_{j=1}^n E_j^{(\ell+1)} u \otimes f_j.  
\end{align*}

By Lemma \ref{lem:Hirshen0812}, we have
\[
   E^{(\ell)}u -\chi_{\lambda,\lambda'}(E^{(\ell)})u \in \operatorname{Ker}T.
\]
Hence, modulo $\operatorname{Ker} T \otimes F$,
\begin{equation}
\label{eqn:cG2}
   c_G(\sum_{j=1}^n E_j^{(\ell)}u \otimes f_j)
\equiv
 2 \chi_{\lambda,\lambda'}(E^{(\ell)})u \otimes f_0
 +(A+n+1) \sum_{j=1}^n (E_j^{(\ell)}u \otimes f_j)
 + 2 \sum_{j=1}^n E_j^{(\ell+1)}u \otimes f_j.  
\end{equation}

Combining~\eqref{eqn:cG1} and~\eqref{eqn:cG2}, we obtain
\begin{align*}
c_G^{k+1}(u \otimes f_0)
 \equiv\, &
\left( (A + n + 1 + 2a)\, b^{(k)} + 2 \sum_{\ell=0}^{k-1} B_{\ell}^{(k)} \chi_{\lambda, \lambda'}(E^{(\ell)}) \right) u \otimes f_0 \\
&+ \sum_{j=1}^n \left( 2 b^{(k)} E_{0j}
+ \sum_{\ell=0}^{k-1} (A + n + 1) B_{\ell}^{(k)} E_j^{(\ell)}
+ 2 \sum_{\ell=0}^{k-1} B_{\ell}^{(k)} E_j^{(\ell+1)} \right) u \otimes f_j.
\end{align*}

By the inductive definitions~\eqref{eqn:bk} and~\eqref{eqn:Bk} of $b^{(k+1)}$ and $B_{\ell}^{(k+1)}$, respectively, this simplifies to
\[
c_G^{k+1}(u \otimes f_0)
\equiv
b^{(k+1)} u \otimes f_0
+ \sum_{\ell=0}^{k} B_{\ell}^{(k+1)} \left( \sum_{j=1}^{n} E_j^{(\ell)} u \otimes f_j \right),
\]
as desired.
\end{proof}

Similarly to $b^{(k)}$ and $B_{\ell}^{(k)}$, 
 we define $\widetilde b^{(k)} \equiv \widetilde b^{(k)}(\lambda,\lambda')$
 for $k \ge 0$
 and $\widetilde B_{\ell}^{(k)}\equiv \widetilde B_{\ell}^{(k)}(\lambda,\lambda')$
 for $0 \le \ell \le k-1$, 
 inductively on $k$,
 by setting $\widetilde b^{(0)}:=1$ and $\widetilde{B}^{(1)}_0:=2$, and defining
\begin{align*}
\widetilde b^{(k+1)}:=& (2a + n + 1) \widetilde b^{(k)} 
  + 2 \sum_{\ell=0}^{k-1} \widetilde B_{\ell}^{(k)}\chi_{\lambda, \lambda'}(E^{(\ell)}), 
\\
\widetilde B_{\ell}^{(k+1)}:=&
\begin{cases}
(n+1) \widetilde B_{0}^{(k)} + 2 \widetilde b^{(k)}
\qquad
&\text{if $\ell=0$, }
\\
(n+1) \widetilde B_{\ell}^{(k)} + 2 \widetilde B_{\ell-1}^{(k)}
\qquad
&\text{if $1 \le \ell \le k-1$, }
\\
2 \widetilde B_{k-1}^{(k)}
&\text{if $\ell=k$}.  
\end{cases}
\end{align*}

An advantage of this renormalization is the following:
\begin{lemma}
\label{lem:23081341}
$\widetilde b^{(k)}(\lambda, \lambda')$ is 
 a polynomial in $\lambda$ and $\lambda'$
 of degree at most $k$, 
 and $\widetilde B_{\ell}^{(k)}(\lambda, \lambda')$ is 
 that of $k-\ell-1$.  
Moreover, 
both are invariant under the action of ${\mathfrak{S}}_{n+1} \times {\mathfrak{S}}_{n}$.  
\end{lemma}

\begin{proof}
Recall that $a=\sum_{i=0}^n \lambda_i - \sum_{j=1}^n \lambda_j'$.  
As seen in Example \ref{ex:23081236}, 
 $\widetilde b^{(1)} = 2a+n+1$
 and $B_0^{(1)}=2$.  
Hence the lemma holds for $k=1$.  

By Lemma \ref{lem:23081228},
 $\chi_{\lambda, \lambda'}(E^{(\ell)})$ is a symmetric polynomial in $\lambda$ and $\lambda'$
 of degree at most $\ell+2$.
Therefore, by induction on $k$ using the recursive definitions of $\widetilde b^{(k)}$ and $\widetilde B_{\ell}^{(k)}$, the claim follows.
\end{proof}

Instead of considering the diagonal action of $c_G$
 on $\pi \otimes F$, 
 we consider 
 \begin{equation}
 \widetilde c_G:=c_G - A \operatorname{id},
 \end{equation}
 where we recall $\lambda=(\lambda_0, \dots, \lambda_n)$
 is the infinitesimal character 
 of $\pi$, and $A=|\lambda|^2-|\rho_{\mathfrak g}|^2$, as before.

The proof of the following proposition is parallel
 to Proposition \ref{prop:23081223}, and will therefore be omitted.

\begin{proposition}
\label{prop:23081236}
For any $k \in {\mathbb{N}}$ and any $u \in \pi$, 
we have
\[
  (T \otimes \operatorname{pr}) \, \widetilde c_G^k
 (u \otimes f_0)
  =\widetilde b^{(k)} (\lambda,\lambda') T u.  
\]
\end{proposition}

We are ready to complete the proof of Theorem~\ref{thm:23081343}.
\begin{proof}
[Proof of Theorem~\ref{thm:23081343}]
~~~
Observe that since $c_G-|\lambda \pm e_j|^2+|\rho_{\mathfrak g}|^2=\widetilde c_G \mp(1+2\lambda_j)$, 
 it follows that
\begin{align*}
\operatorname{p r}_{i, \pm}
=&\prod_{j \ne i} (\widetilde c_G \mp (1+2\lambda_j))
\\
=&\sum_{k=0}^n (\mp 1)^k \,\,\widetilde c_G^k e_k(1+2 \lambda_0, \dots, \widehat{1+2\lambda_i}, \dots, 1+2\lambda_n)
\end{align*}
where $e_k(x_1, \dots, x_n)$ denotes the $k$-th symmetric polynomial defined via the identity
\[
  \prod_{j=1}^n (t-x_j)
  =
 \sum_{k=0}^n (-1)^k t^{n-k}e_k(x_1, \dots, x_n).  
\]
Applying Proposition~\ref{prop:23081236}, 
 we obtain 
\[
  (T \otimes \operatorname{p r}_{F \to F'})\operatorname{p r}_{i, \pm}
  (u \otimes f_0)
  =g_{i, \pm}(\lambda,\lambda')T u
\]
where the function $g_{i, \pm} (\lambda,\lambda')$
is defined by
\[
  g_{i, \pm} (\lambda,\lambda')
  :=\sum_{k=0}^n (\mp 1)^k
   \tilde{b}^{(n-k)}(\lambda,\lambda')
   e_k(1+2\lambda_0,\dots, \widehat{1+2\lambda_i}, \dots, 1+2\lambda_n).  
\]

By Lemma~\ref{lem:23081341},   
 $g_{i, \pm}(\lambda,\lambda')$ is a polynomial in $\lambda$, $\lambda'$
 of degree 
 at most $n$, 
 and is symmetric with respect to $\lambda_1'$, $\dots$, $\lambda_n'$.

Let us find the formula of $g_{i, \pm}(\lambda,\lambda')$.

By the definition, 
\begin{equation}
\label{eqn:gsigma}
 g_{\sigma(i)}(\sigma^{-1} \lambda, \lambda')
 =g_{i}(\lambda, \lambda')
\end{equation}
 for any $i$ and any $\sigma \in {\mathfrak{S}}_{n+1}$.

Since these formul{\ae} are applicable also any representation 
in Setting \ref{set:Fpi}, 
 we evaluate them to pairs of irreducible finite-dimensional representations.  
For this, 
 let $\pi$ and $\pi'$ be the irreducible finite-dimensional representations
 of $G L_{n+1}$ and $G L_n$
 having highest weights, 
 respectively,  
 $\mu=(\mu_0, \mu_1, \dots, \mu_n) \in {\mathbb{Z}}^{n+1}$
 and $\mu'=(\mu_1', \dots, \mu_n') \in {\mathbb{Z}}^{n}$
 satisfying the interlacing relation
 $\mu_0 \ge \mu_1' \ge \mu_1 \ge \cdots \ge \mu_n' \ge \mu_n$.  
The infinitesimal characters of $\pi$ and $\pi'$
 are given respectively 
 by $\lambda = \mu + \rho_{\mathfrak g}$ and $\lambda' = \mu' + \rho_{\mathfrak{g}'}$.

We note that $\sigma_{i,+}:=\operatorname{p r}_{i,+}(\pi \otimes F)$
 is an irreducible representation of $G L_{n+1}$
 with highest weight $\mu+e_i$
 as far as $i=0$ or $\mu_{i-1} > \mu_i$.  
Likewise, 
 $\sigma_{i,-}:=\operatorname{p r}_{i,-}(\pi \otimes F)$ is the one 
with highest weight $\mu-e_i$
 as far as $i=n$ or $\mu_i > \mu_{i+1}$.  

We now assume $\mu_0 = \mu_1' >\mu_1$.  
Then $[\sigma_{0, -}:\pi']=0$, 
and thus $(T \otimes \operatorname{p r}_{F \to F'})\circ \operatorname{p r}_{0, -}(\pi \otimes F)=0$.  
This means 
 that the polynomial $g_{0, -}(\lambda,\lambda')$ vanishes
 if $\mu_0=\mu_1'$, 
 that is, 
 if $\lambda_0-\frac n 2=\lambda_1'-\frac{n-1}2$.  
Therefore, 
 $g_{0,-}(\lambda, \lambda')$ has $\lambda_0-\lambda_1'-\frac 1 2$
 as an irreducible factor.  
Since $g_{0,-}(\lambda, \lambda')$ is symmetric 
 with respect to $\lambda_1'$, $\dots$, $\lambda_n'$, 
 $g_{0,-}(\lambda,\lambda')$ must be
 of the form 
 $c_- \prod_{j=1}^n (\lambda_0-\lambda_j'-\frac 1 2)$
 for some constant $c_-$.  
By \eqref{eqn:gsigma}, one has
\[
   g_{i,-}(\lambda,\lambda')=c_-\prod_{j=1}^n(\lambda_i-\lambda_j'-\frac 1 2)
\]
for all $0 \le i \le n$.

On the other hand, 
 we consider the case $i=n$. 
Take a Borel subalgebra
 of ${\mathfrak{g l}}_{n+1}({\mathbb{C}})$
 such that $f_0$ is the lowest vector 
in $F=({\mathbb{C}}^{n+1})^{\vee}$.   
Then $(T \circ \operatorname{p r}_{F \to F'}) \operatorname{p r}_{n,-}(u \otimes f_0) \ne 0$
 if $u$ is the lowest weight vector in $\pi$.  
Hence $c_- \ne 0$.

Similarly, 
 we obtain the formula for $(T \circ \operatorname{p r}_{F \to F'})(\operatorname{p r}_{i,+}(u \otimes f_0))$
 by starting with $i=n$ and $\mu_n=\mu_n'$.

Thus Theorem \ref{thm:23081343} has been proved.  
\end{proof}

Hence, the proof of Theorem~\ref{thm:23081510} and Theorem \ref{thm:23081404}
 is completed.

\subsection{String of Translations of Morphisms}
\label{subsec:string}

We now return to Question~\ref{q:mfo1}, which was introduced at the beginning of this section. As we have seen, Theorem~\ref{thm:23081404} provides a positive result for Condition~(3) of Question~\ref{q:mfo1} in the special setting where $F$ is the standard representation or its dual. and $F'$ is the trivial representation.

To bridge a given coherent pair $(\pi, \pi')$ with an elementary pair $(\sigma, \sigma')$,
one may consider iterated applications of appropriately chosen pairs of \emph{small} representations $(F, F')$.

This leads us to ask the following:
\begin{question}
\label{q:mfo2}
To what extent can such sequences of \lq\lq{translating morphisms}\rq\rq\ be continued?
\end{question}

\begin{definition}
\label{def:230808}
Let $\pi, \sigma \in \operatorname{Irr}(G)$
 and $\pi', \sigma' \in \operatorname{Irr}(G')$
 be such that
\[
  \dim \operatorname{Hom}_{G'}(\pi|_{G'},\pi') = \dim\operatorname{Hom}_{G'}(\sigma|_{G'},\sigma') =1.  
\]
We say that $(\sigma,\sigma')\to(\pi,\pi')$
 if the composition of the three maps in \eqref{eqn:Fdiag2} is non-zero, obtained by taking
$(F, F')=({\mathbb{C}}^{n+1}, {\mathbb{C}})$
 or $(({\mathbb{C}}^{n+1})^{\vee}, {\mathbb{C}})$.

Let $\sim$ denote the equivalence relation generated by $\longrightarrow$.  
\end{definition}

For example, $(\pi, \pi') \sim (\sigma, \sigma')$ if there exist pairs $(\pi_i, \pi_i')$ for $i=1,2$ such that
\[
  (\pi, \pi') \ \rightarrow \ (\pi_1, \pi_1') \ \leftarrow \ (\pi_2, \pi_2') \ \leftarrow \ (\sigma, \sigma').
\]

\begin{theorem}
\label{thm:23081510}
For any coherent pair $(\pi, \pi')$
with associated signature $[\delta\delta']$
 (see Definition \ref{def:coherent_pair}), 
 there exists an elementary pair $(\sigma, \sigma')$
 with the same coherent signature $[\delta \delta']$
 such that $(\pi, \pi') \sim (\sigma, \sigma')$.  
 
\end{theorem}

\begin{proof}
[Proof of Theorem \ref{thm:23081510}]

Let $[\delta \delta']$ be the associated signature
of the coherent pair $(\pi, \pi')$ 
 (see Definition~\ref{def:coh_deco}).
Let $\lambda$ and $\lambda'$ be the Harish-Chandra parameters of $\pi$ and $\pi'$, respectively,
where $\lambda \in {\mathbb{Z}}_{\operatorname{int}}^{p+q}(\delta) ={\mathbb{Z}}_{\operatorname{int}}^{p+q}\cap C(\delta)$ (see \eqref{eqn:Z_int_delta}) and $\lambda' \in {\mathbb{Z}}_{\operatorname{int}}^{p+q-1}(\delta')$.
  
We define
$\widetilde \lambda=(\widetilde\lambda_1, \dots, \widetilde \lambda_{p+q}) \in {\mathbb{Z}}_{\operatorname{int}}^{p+q}$ 
 by the following formula:
\begin{equation}
\label{eqn:coh_to_elem}
\widetilde\lambda_i
:=
\begin{cases}
\lambda_i
\quad
&(1 \le i \le p), 
\\
\lambda_{i-1}'-\frac 1 2
\quad
&(p+1 \le i \le p+q).  
\end{cases}
\end{equation}
It follows from Corollary~\ref{cor:coh_to_elem} that
$\widetilde \lambda \in C(\delta)$ and $(\sigma, \sigma'):=(\pi_{\widetilde \lambda}^{\delta}, \pi')$ form an elementary pair with the same associated signature $[\delta\delta']$.  

We now define
\[
   \ell_j:=\widetilde \lambda_{p+j}-\lambda_{p+j} \ge 0
  \quad\text{for $1 \le j \le q$.}  
\]
Then the string of weights
\begin{align*}
  &\widetilde \lambda \rightarrow \widetilde \lambda -e_{p+1}
  \rightarrow \widetilde \lambda -2 e_{p+1}
  \rightarrow \cdots
  \rightarrow \widetilde \lambda -\ell_1 e_{p+1}
  \rightarrow \widetilde \lambda -\ell_1 e_{p+1}-e_{p+2}
\\
  \rightarrow 
  &\cdots
  \rightarrow \widetilde \lambda -\ell_1 e_{p+1}-\ell_2 e_{p+2}
  \rightarrow \cdots
  \rightarrow \lambda=\widetilde \lambda-\sum_{j=1}^q \ell_j e_{p+j} 
\end{align*}
starts from the elementary pair $(\sigma, \sigma')$
 and ends with the given coherent pair $(\pi, \pi')$.  
All the intermediate weights lie in the strict dominant Weyl chamber $C(\delta)$ associated with the decoration $\delta$, 
 and fulfill the assumption of Theorem~\ref{thm:23081404}~(1). 
In view of Proposition~\ref{prop:prireg}, 
 one can proceed by applying $\operatorname{p r}_{p+j, +}$ $(1 \le j \le q)$ along the string.  
\end{proof}

\begin{remark}
Conversely, 
 the iteration of the translation functors along the string 
\begin{align*}
  &\lambda \rightarrow \lambda + e_{p+1}
  \rightarrow \lambda +2 e_{p+1}
  \rightarrow \cdots
  \rightarrow \lambda + \ell_1 e_{p+1}
  \rightarrow \lambda +\ell_1 e_{p+1}+e_{p+2}
\\
  \rightarrow 
  &\cdots
  \rightarrow \lambda +\ell_1 e_{p+1}+\ell_2 e_{p+2}
  \rightarrow \cdots
  \rightarrow \widetilde \lambda=\lambda+\sum_{j=1}^q \ell_j e_{p+j} 
\end{align*}
 allows us to start from the given coherent pair $(\pi, \pi')$
 and ends the elementary pair $(\sigma, \sigma')$.  
\end{remark}

\begin{example}
[$U(2,2) \downarrow U(1,2)$]
\label{ex:23081410}
Consider the interleaving pattern \linebreak $[\delta \delta']=\oi - + \op \oi - +$.  
Then the pair of discrete series representations $(\pi_\lambda, \pi'_{\lambda'})$ of $G\supset G'$
has the coherent signature $[\delta \delta']$
if and only if the Harish-Chandra parameters $\lambda=(\lambda_1, \lambda_2, \lambda_3, \lambda_4)$ and $\lambda'=(\lambda_1', \lambda_2', \lambda_3')$
satisfy the interleaving condition
\[\lambda_2' > \lambda_3 > \lambda_1 > \lambda_1' > \lambda_3' > \lambda_4 > \lambda_2.
\]
Moreover, $(\pi_\lambda, \pi'_{\lambda'})$ forms an elementary pair if and only if $\lambda_2' = \lambda_3+\frac 1 2$
 and $\lambda_3'=\lambda_4+\frac 1 2$.

For instance, 
 $(\pi, {\pi'})
 = (\pi_{\frac{11}2, \frac{1}2, \frac{13}2, \frac{3}2}, 
 \pi_{5, 8 , 4})$
 forms a non-elementary coherent pair
 and 
 $(\sigma, \sigma')= (\pi_{\frac{11}2, \frac{1}2, \frac{15}2, \frac{7}2}, {\pi}_{5, 8, 4})$ forms an elementary pair with the same associated signature $[\delta \delta']=\oi - + \op \oi - +$.  .

The condition in Theorem~\ref{thm:23081404} (2) 
 guarantees that $\lambda$ can be translated 
 to $\lambda-e_i$ provided $\lambda_i \not \in \{\frac{11}2, \frac 7 2, \frac 92\}$.  
An iteratative procedure shows that the string
\[
  \frac 1 2 (11, 1, 15, 7) \overset 3\longrightarrow
  \frac 1 2 (11, 1, 13, 7) \overset 4\longrightarrow
  \frac 1 2 (11, 1, 13, 5) \overset 4\longrightarrow
  \frac 1 2 (11, 1, 13, 3)
\]
demonstrates that $(\pi, \pi') \sim (\sigma, \sigma')$.  
\end{example}

\section{\texorpdfstring{$({\mathfrak{P}}, K)$}{(P, K)}-cohomologies and Coherent Pairs}
\label{sec:Pk}

Given a coherent pair $(\pi, \pi')$, their coherent $(\fP,K)$-cohomology groups are concentrated in the same degree $q=q(\pi)=q(\pi')$, as recalled from Section~\ref{sec:coherent}.

In this section, we investigate whether a symmetry breaking operator $S\colon \pi \to \pi'$ induces a natural homomorphism between the corresponding coherent $(\fP,K)$-cohomology groups.
Such a morphism does exist and is, in fact, an isomorphism, when $(\pi, \pi')$ forms an elementary pair (Theorem~\ref{thm:elementary_pair_equiv}).
However, for a non-elementary pair, the existence of an analogous morphism cannot generally be expected.

We propose a conjectural framework involving \emph{translation of symmetry breaking} in Conjecture~\ref{conj:250629}, and establish a supporting evidence in Theorem~\ref{thm:non-van}, based on
 the method of symmetry breaking under translation developed in Section~\ref{sec:SBO_translation}.

\subsection{Preliminaries on \texorpdfstring{$(\fP,K)$}{(P,K)}-cohomology}
\label{subsec:basic_Pk_cohomology}
This subsection is devoted to summarizing the $(\fP,K)$-cohomology for discrete series representations, as preparation for the subsequent discussion.
We also recall from Section~\ref{sec:coherent} certain natural isomorphisms of $(\fP,K)$-cohomology spaces, induced by symmetry breaking in the restriction of discrete series representations, within the framework of \emph{elementary pairs}.

As recalled from \eqref{eqn:P_+},
$\fP \equiv \mathfrak{P}^+ = \fk \oplus \fp^+$
is a maximal parabolic subalgebra of $\mathfrak{g}=\mathfrak{gl}_{p+q}$.

This is a $\theta$-stable maximal parabolic subalgebra of $\fg$ with
unipotent radical $\fp^+$ and Levi quotient $\fk$.

Any smooth admissible representation $\pi$ gives rise to a $(\fg,K)$-module $\pi_K$ by taking $K$-finite vectors, which in turn defines a $(\fP, K)$-module by restriction.
The action of $\mathfrak{p}^+$ on any irreducible $(\fP,K)$-module is trivial, so such a module can be identified with an irreducible $(\fk,K)$-module, and hence with an irreducible representation of the compact connected group $K$. 

Let $W$ be a finite-dimensional $K$-module, regarded as
a $(\fP,K)$-module by letting $\fp^+$ act trivially.  
\begin{definition}   The $\db$-cohomology (or $(\mathfrak{P}, K)$-cohomology) of $\pi$ with coefficients in $W^\vee$ is defined to be the relative Lie algebra cohomology $H^*(\fP, K; \pi_K \otimes W^\vee)$, which we simply denote by
$$H^*(\fP,K;\pi\otimes W^\vee).$$
\end{definition}

The $\db$-cohomology of discrete series representations $\pi$ with coefficients in $K$-modules $W$ originates in the work of Narasimhan and Okamoto \cite{NO70}, who studied the geometric realization of $\pi$ for generic parameters. 
The $(\mathfrak{P},K)$-cohomology can be regarded as a \emph{cohomological Frobenius reciprocity} of the geometric realization of discrete series representations via Dolbeault cohomology.
This cohomology was fully calculated by W.~Schmid, and was later extended to
 non-degenerate limits of discrete series by F.~Williams. 
 The following statement is reproduced from \cite[Theorem 3.2.1]{BHR}. 

% 1000
% See \eqref{eqn:pK_vanish}.

\begin{prop}\label{prop:dbc}  
Let $\pi$ be a discrete series representation of $G$ with Harish-Chandra parameter $\tlambda$, and let $q(\pi)$ denote the degree of $\pi$, as defined in Definition~\ref{def:q}.  
Let $W_\pi$ be the irreducible representation of $K$ with highest weight $\Lambda = \tlambda - \rho_\mathfrak{g}$, referred to as the \emph{coherent parameter}, as defined in Definition~\ref{def:W_lambda}.  
Here, $\rho_\mathfrak{g}$ denotes half the sum of the positive roots in the root system containing those in $\fp^+$.

\begin{enumerate}
    \item For any irreducible representation $W$ of $K$,
    \[
    \dim_{\C} H^i(\fP, K ; \pi \otimes W^\vee) = 
    \begin{cases}
    1, & \text{if } i = q(\pi) \text{ and } W \simeq W_\pi, \\
    0, & \text{otherwise}.
    \end{cases}
    \]

    \item Suppose that $\sigma$ is another discrete series representation of $G$, and that
    \[
    H^i(\fP, K ; \sigma \otimes W_\pi^\vee) \neq 0.
    \]
    Then $\sigma \simeq \pi_{\tlambda}$ and $i = q(\pi)$.
\end{enumerate}
\end{prop}

%1000

As shown in Theorem~\ref{thm:elementary_pair_equiv}, 
there are four equivalent characterizations of elementary pairs $(\pi, \pi')$, one of which involves natural isomorphisms between $(\fP,K)$- and $(\mathfrak{P}',K')$-cohomologies.
We now recall this characterization.

Let $(\sigma, \sigma')$ be an elementary pair (Definition~\ref{def:elementary}),  
with coherent parameter $$(W, W') := (W_\sigma, W'_{\sigma'}).$$  
In particular, parts (ii) and (iii) of Theorem~\ref{thm:elementary_pair_equiv} assert the existence of non-zero homomorphisms
\[
T \in \operatorname{Hom}_{G'}(\sigma|_{G'}, \sigma')
\quad \text{and} \quad
\operatorname{pr}_W \in \operatorname{Hom}_{K'}(W^\vee|_{K'}, {W'}^\vee),
\]
each unique up to scalar multiplication.  
Let
\[
\fP' := \fk' \oplus \fp^{\prime,+},
\]
which is a subalgebra of $\fP$, and set $q := q(\sigma) = q(\sigma')$  
(see Definition~\ref{def:q} for the convention, and Proposition~\ref{prop:23022216} for the equality).

As recalled in part (iv) of Theorem~\ref{thm:elementary_pair_equiv}, 
the morphism $T \otimes \operatorname{pr}_{W}$ induces a bijection on $({\mathfrak{P}}, K)$-cohomology:
\begin{equation}
\label{eqn:TprW_iso}
\xymatrix@C=30pt{
H^{\ast}({\mathfrak{P}}, K;\, \sigma \otimes W^\vee)
 \ar[r]^-{\sim}
& 
H^{\ast}({\mathfrak{P}}', K';\, {\sigma}'\otimes {W'}^\vee).
}
\end{equation}

\subsection{A Conjecture for Coherent Pairs $(\texorpdfstring{\pi}{pi}, \texorpdfstring{\pi'}{pi'})$}
\label{conj_PK}
~
\newline
We now turn to the more general setting, by considering morphisms between
$H^*(\fP,K;\pi\otimes W_\pi^\vee)$
and
$H^*(\fP',K';\pi^{\prime}\otimes W_{\pi'}^\vee)$
in the setting where $(\pi, \pi')$ is a
\emph{non-elementary coherent pair}.

In this section, we present a conjecture (see Conjecture~\ref{conj:250629}), which
will be addressed in a future article.  In the subsequent sections, we shall state and prove a special case of conjecture.

To this end, we begin by fixing the relevant notation and framework.

Let $(\pi, \pi')$ be a coherent pair, 
 as defined in Definition~ \ref{def:coherent_pair},
 and let 
 \[
 S \in \operatorname{Hom}_{G'}(\pi|_{G'}, \pi')
 \]
 be a non-zero symmetry breaking operator.  
Let $(W_{\pi}, W_{\pi'}')$ denote the associated coherent parameter  (see Definition \ref{def:W_lambda}).  
It is important to note that, in general, $\operatorname{Hom}_{K'}(W_{\pi}^\vee|_{K'}, {W'_{\pi'}}^\vee) =0$ 
 when the pair $(\pi, \pi')$ is non-elementary.
Consequently, we cannot expect to define an induced morphism such as \eqref{eqn:TprW_iso} in this case. 

On the other hand, as recalled in Theorem~\ref{thm:230816},
for any coherent pair $(\pi, \pi')$
 with associated signature $[\delta \delta']$, 
 there exists an elementary pair $(\sigma, \sigma')$
 with the same signature $[\delta\delta']$.

\smallskip

In order to relate a coherent pair $(\pi, \pi')$ to an elementary pair $(\sigma, \sigma')$, we consider the following setup.
Choose a pair of irreducible {\it{finite-dimensional}} representations $F$ and $F'$
 of $G$ and $G'$, 
 respectively, 
 such that 
 \begin{equation}
\operatorname{Hom}_{G'}(F|_{G'},F') \ne 0,
 \quad 
 \operatorname{Hom}_{G}(\pi \otimes F, \sigma)\neq 0, \quad
  \operatorname{Hom}_{G'}(\sigma', \pi' \otimes F')\neq 0.
 \end{equation}
Let $\operatorname{pr}_F \colon F \to F'$ denote a non-zero $G'$-homomorphism, 
 which is unique up to scalar multiplication.  

If there exist morphisms
\begin{equation}
\label{eqn:pi_F_sigma}
  A \in \operatorname{Hom}_{G}(\sigma, \pi \otimes F) 
\quad\text{and}\quad
  A' \in \operatorname{Hom}_{G'}(\sigma', \pi' \otimes F')
\end{equation}
satisfying 
\begin{equation}
\label{eqn:ST_comm}
(S \otimes \operatorname{pr}_F)\circ A= A' \circ T, 
\end{equation}
then the morphisms $A$ and $A'$ induce the following commutative diagram 
 of $({\mathfrak{P}}, K)$-cohomologies:
\begin{equation}
\label{eqn:fP_TS_diagram}
\xymatrix@C=75pt{
H^{\ast}({\mathfrak{P}}, K;\sigma \otimes W^\vee)
 \ar[r]^{T \otimes \operatorname{pr}_{W}}
\ar[d]_{A_{\ast}}
& 
H^{\ast}({\mathfrak{P}}', K';{\sigma}'\otimes {W'}^\vee)
 \ar[d]^{A_{\ast}'}
\\
H^{\ast}({\mathfrak{P}}, K;\pi \otimes F \otimes W^\vee)
 \ar[r]^{S \otimes \operatorname{pr}_{F} \otimes \operatorname{pr}_{W}} 
&
H^{\ast}({\mathfrak{P}}', K';\pi' \otimes F' \otimes {W'}^\vee)
}
\end{equation}

Suppose further that we are given morphisms
\begin{equation}
\label{eqn:Wpi_F}
  B \in \operatorname{Hom}_{K}(W_{\pi}^{\vee},  F \otimes W^{\vee})
\quad\text{and}\quad
  B' \in \operatorname{Hom}_{K'}(W'^{\vee} \otimes F', W_{\pi'}'^{\vee}).  
\end{equation}
Passing to $({\mathfrak{P}}, K)$-cohomologies, 
 the morphisms $B$ and $B'$ give rise to the following diagram:  
\begin{equation}
\label{eqn:B_fP}
\xymatrix@C=75pt{
H^{\ast}({\mathfrak{P}}, K;\pi \otimes F \otimes W^\vee)
 \ar[r]^{S \otimes \operatorname{pr}_{F} \otimes \operatorname{pr}_{W}} 
&
 H^{\ast}({\mathfrak{P}}', K';\pi' \otimes F' \otimes {W'}^\vee)
 \ar[d]^{B_{\ast}'}
\\
H^{\ast}({\mathfrak{P}}, K;\pi \otimes W_{\pi}^\vee)
 \ar[u]^{B_{\ast}}
& 
H^{\ast}({\mathfrak{P}}', K';\pi' \otimes {W'_{\pi'}}^\vee)
}
\end{equation}

Since an induced morphisms such as \eqref{eqn:TprW_iso} is not available in the non-elementary case,
we instead consider the composite map
\[
  B_{\ast}' \circ (S \otimes \operatorname{pr}_F \otimes \operatorname{pr}_{W}) \circ B_{\ast}.
\]
The following conjecture asserts, in particular, that this composion map is an isomorphism 
 between 
$
 H^{\ast}({\mathfrak{P}}, K;\pi \otimes W_{\pi}^\vee)
$
 and 
$
   H^{\ast}({\mathfrak{P}}', K';\pi' \otimes {W'_{\pi'}}^\vee)
$
for a suitable choice of data, 
 formulated as below.  

\begin{conj}
\label{conj:250629}
Let $(\pi,\pi')$ be a coherent pair whose
 Harish-Chandra parameters $\lambda$ and $\lambda'$
 are sufficiently far from the walls.   
 Then there exist:
\begin{itemize}
  \item an elementary pair $(\sigma, \sigma')$,
  \item a pair of irreducible finite-dimensional representations $F$ and $F'$ of $G$ and $G'$, respectively,
  \item and morphisms $A$, $A'$, $B$, $B'$, $T$, and $S$ as above,
\end{itemize}
such that the following conditions hold:
\begin{enumerate}
\item[$\bullet$]
 (Symmetry breaking under translation) \eqref{eqn:ST_comm} commutes;

\item[$\bullet$]
 All the induced morphisms in cohomology (see \eqref{eqn:3_diagram}) are isomorphisms.  
\end{enumerate}
In particular, all six cohomology groups appearing in diagram \eqref{eqn:3_diagram} are one-dimensional
 at $q=q(\pi)=q(\pi')$:
\begin{equation}
\label{eqn:3_diagram}
\xymatrix@C=80pt{
H^{q}({\mathfrak{P}}, K;\sigma \otimes W^\vee)
 \ar[r]_{T \otimes \operatorname{pr}_{W}}^{\sim}
 \ar[d]_{A_{\ast}}
& 
H^{q}({\mathfrak{P}}', K';{\sigma}'\otimes {W'}^\vee)
 \ar[d]^{A_{\ast}'}
\\
H^{q}({\mathfrak{P}}, K;\pi \otimes F \otimes W^\vee)
 \ar[r]^{S \otimes \operatorname{pr}_{F} \otimes \operatorname{pr}_{W}} 
&
 H^{q}({\mathfrak{P}}', K';\pi' \otimes F' \otimes {W'}^\vee)
 \ar[d]^{B_{\ast}'}
\\
H^{q}({\mathfrak{P}}, K;\pi \otimes W_{\pi}^\vee)
 \ar[u]^{B_{\ast}}
 \ar@{.>}[r]
& 
H^{q}({\mathfrak{P}}', K';\pi' \otimes {W'_{\pi'}}^\vee)
}
\end{equation}

\end{conj}

We examine when each of these requirements holds
in the subsequent sections.

We provide sufficient conditions for the vertical morphisms $A_{\ast}$ and $B_{\ast}$
 in diagram~\eqref{eqn:3_diagram} to be isomorphisms,
 in Sections~\ref{subsec:vert_A} and \ref{subsec:vert_B}, respectively. 
 These results concern representations of individual groups and do not involve the restriction $G \downarrow G'$.
 The proof makes use of existing techniques related to translation functors for single groups.

 The core of Conjecture~\ref{conj:250629} is the commutativity of diagram~\eqref{eqn:ST_comm},  
for which we have developed new techniques in Section~\ref{sec:SBO_translation} concerning symmetry breaking under translation.

With these preparations, we provide supporting evidence in Theorem~\ref{thm:non-van}.

\bigskip
\subsection[Morphisms between (P, K)-Cohomologies under Translation]{Morphisms between $({\mathfrak{P}}, K)$-Cohomologies under Translation}
\label{subsec:vert_A}
~
\medskip
\newline
In this section, we provide a sufficient condition 
 for the vertical morphism $A_{\ast}$ 
 in  diagram~\eqref{eqn:3_diagram} to be an isomorphism.
The main result is stated in Proposition~\ref{prop:23082631}, which is proved
in Section~\ref{sec:comb_fd_rep},
building on certain combinatorial results on finite-dimensional representations.

\medskip

The symmetric pair $(G,K)$ is not a Gelfand pair
when $p,q \ge 2$, and the restriction $F|_K$ of an irreducible finite-dimensional representation $F$ of $G$
 to the subgroup $K$ is not necessarily multiplicity-free.
As a result, the tensor product representation
 $F|_K \otimes W$ generally exhibits even higher multiplicities for $W \in \operatorname{Irr}(K)$.  
Nevertheless, we have the following:

\begin{proposition}
\label{prop:23082631}
Let $\delta \in \operatorname{Deco}(p,q)$
 and $\lambda, \widetilde \lambda \in {\mathbb{Z}}_{\operatorname{int}}^{p+q} \cap C(\delta)$, 
 see~\eqref{eqn:Cdelta}
 for the notation.  
Let $F$ be an irreducible finite-dimensional representation of $G$
 with extremal weight $-\lambda+\widetilde \lambda$.  
Define $\pi:=\pi_{\lambda}$
 and $\sigma :=\pi_{\widetilde \lambda}$
 to be discrete series representations of $G$, with Harish-Chandra parameters $\lambda$ and $\widetilde \lambda$, respectively, and with the same decoration $\delta$.
Let $W_{\pi}$  and $W_{\sigma} \in \widehat K$ be the corresponding coherent parameters, as defined in \eqref{eqn:cohlmd}.  
\par\noindent
{\rm{(1)}}\enspace
There are isomorphisms of $G$-modules:
\[
\psi_{\lambda}^{\widetilde \lambda}(\pi)
=P_{\widetilde \lambda}(\pi \otimes F) \simeq \sigma,
\quad
\psi_{\widetilde \lambda}^{\lambda}(\sigma)
=P_{\lambda}(\sigma \otimes F^\vee) \simeq \pi.
\]
\newline{\rm(2)}\enspace
If $\tau \in \operatorname{Irr}(G)$ has $\mathfrak{Z}(\mathfrak{g})$-infinitesimal character $\lambda$
satisfying 
\[
P_{\widetilde{\lambda}}(\tau \otimes F) \simeq \sigma,
\]
then $\tau$ is isomorphic to $\pi_\lambda$.
\par\noindent
{\rm{(3)}}\enspace
Assume either of the following conditions holds:
\begin{equation}
\label{eqn:lmdij}
  \min_{1 \le i \le p} (\lambda_i -\widetilde \lambda_i)
 \ge \max_{1 \le j \le q} (\lambda_{p+j} -\widetilde \lambda_{p+j})
\end{equation}
 or 
\begin{equation}
\label{eqn:lmdji}
  \max_{1 \le i \le p} (\lambda_i -\widetilde \lambda_i)
 \le \min_{1 \le j \le q} (\lambda_{p+j} -\widetilde \lambda_{p+j})
\end{equation}
Then the natural $G$-morphisms
\[
         \sigma
  \simeq P_{\widetilde\lambda}(\pi \otimes F) 
  \hookrightarrow \pi\otimes F
\]
induce an isomorphism of $({\mathfrak{P}}, K)$-cohomology:
\begin{equation}
\label{eqn:PK_isom}
   H^j({\mathfrak{P}}, K; \sigma \otimes W_{\sigma}^\vee)
   \overset \sim \rightarrow
   H^j({\mathfrak{P}}, K; \pi\otimes F \otimes W_{\sigma}^\vee).  
\end{equation}
Furthermore, 
 the projection
\[
   \pi \otimes F \to P_{\widetilde\lambda}(\pi \otimes F)
 \simeq \sigma
\]
of $G$-modules induces the inverse map in \eqref{eqn:PK_isom}.  

These cohomology spaces are one-dimensional 
 if $j=q(\delta)$, 
 and vanish otherwise;
 see Section \ref{subsec:q}
 for the definition of $q(\delta)$.  
\end{proposition}

\begin{remark}
\label{rem:23082614}
{\rm{(1)}}\enspace
The assumption \eqref{eqn:lmdij} or \eqref{eqn:lmdji}
 fits with Theorem~\ref{thm:230816}
 which constructs an elementary pair from 
 a non-elementary coherent pair.  
\par\noindent
{\rm{(2)}}\enspace
We do not need an assumption that the extremeal weight $-\lambda+\widetilde \lambda$ is dominant.  
\end{remark}

\begin{proof}
[Proof of Proposition \ref{prop:23082631}]
(1)\enspace
Since both $\lambda$ and $\widetilde \lambda$ lie 
 in the same chamber $C(\delta)$, 
 the spectral sequence \eqref{eqn:Rbspec} in Lemma~\ref{lem:Rspecseq}) of $({\mathfrak{g}}, K)$-modules
\[
   P_{\widetilde \lambda}({\mathcal{R}}_{\mathfrak{b}}^j({\mathbb{C}}_{\lambda}
 \otimes F_i/F_{i-1}))
\Rightarrow 
   {\mathcal{R}}_{\mathfrak{b}}^j({\mathbb{C}}_{\lambda}) \otimes F
\]
collapses, 
 where $0=F_0 \subset F_1 \subset \cdots \subset F_m=F$ is a ${\mathfrak{b}}$-stable filtration
 such that the unipotent radical ${\mathfrak{u}}$
 acts trivially on $F_i/F_{i-1}$.
This yields the first isomorphism
$P_{\widetilde \lambda}(\pi \otimes F) \simeq \sigma$.  

By switching the roles of $(\pi, \lambda)$ and $(\sigma, \widetilde{\lambda})$,
we obtain the second isomorphism
$P_{\lambda}(\sigma \otimes F^\vee) \simeq \pi$,
since the contragredient representation $F^\vee$ has an extremal weight $\lambda-\widetilde{\lambda}$.  
\par\noindent
{\rm{(2)}}\enspace
If follows from the natural isomorphism
\[
\operatorname{Hom}_G(\sigma\otimes F^\vee, \tau) \simeq \operatorname{Hom}_G(\sigma, \tau\otimes F)
\]
that we have an induced isomorphism
\[
\operatorname{Hom}_G(P_\lambda(\sigma\otimes F^\vee), \tau) \simeq \operatorname{Hom}_G(\sigma, P_{\widetilde\lambda}(\tau\otimes F)).
\]
The left-hand side is equal to $\operatorname{Hom}_G(\pi,\tau)$ by the first statement (1), while the right-hand side is equal to $\operatorname{Hom}_G(\sigma, \sigma)=\C \operatorname{id}$ by the assumption. 
Therefore, the two irreducible representations $\tau$ and $\pi=\pi_\lambda^\delta$ must be isomorphic.
\par\noindent
{\rm{(3)}}\enspace
We take another (coarser) $({\mathfrak{P}}, K)$-stable filtration
 $0=F_0 \subset F_1 \subset \cdots \subset F_k=F$
 such that unipotent radical ${\mathfrak{p}}^+$ acts trivially 
 on each $K$-module $F_i/F_{i-1}$.  
Accordingly, 
 there is a spectral sequence
\begin{equation}
\label{eqn:VG437}
H^j({\mathfrak{P}}, K; \sigma \otimes F_i/F_{i-1} \otimes W_{\sigma}^\vee)
 \Rightarrow
H^j({\mathfrak{P}}, K; \sigma \otimes F \otimes W_{\sigma}^\vee).  
\end{equation}

On the other hand, the assumption \eqref{eqn:lmdij} or \eqref{eqn:lmdji} implies
the multiplicity-one result:
\[ [W_{\pi}:W_{\sigma} \otimes F|_K]=1,
\]
as will be shown in
 Proposition \ref{prop:23082610} below.  

The isomorphism \eqref{eqn:PK_isom} then follows from the vanishing results in degree $j \ne q(\delta)$;
 see Proposition~\ref{prop:dbc}.

 Finally, we derive point (3) by applying point (2) to $F^\vee \otimes \sigma$.
\end{proof}

\subsection{Inequalities on Convex Polygons}
~~~\newline
In dealing with $({\mathfrak{P}}, K)$-cohomologies in the context of translation functors,  
it becomes necessary to analyze the tensor product of two irreducible \emph{finite-dimensional} representations,  
one of $G$ and the other of $K$.

In general, the natural morphism \eqref{eqn:PK_isom} between $({\mathfrak{P}}, K)$-cohomology spaces  
is not expected to be an isomorphism.  
This is because the tensoring functor $\otimes$ and the forgetful functor (i.e., restriction $G \downarrow K$)  
may jointly increase the multiplicities of irreducible components in $F \otimes W_\sigma^\vee$.

This section presents some elementary convex analysis that will be used to establish a sufficient condition for the multiplicity-one property  
of a specific irreducible component in $F \otimes W_\sigma^\vee$  
(see Proposition~\ref{prop:23082610}).

\begin{definition}
Let $D$ be a finite group acting on a real vector space $V$.  
For a fixed vector $v \in V$,
we denote by $\operatorname{Conv}(D; v)$ the convex hull of the finite set $\{w v : w \in D\}$. 
\end{definition}

We begin with a classical example. 
Recall that $W_G$ denotes the Weyl group.
\begin{lemma}
\label{lem:wtconv}
Let $F$ be an irreducible finite-dimensional representation  
of $\mathfrak{g}$ with highest weight $\mu$.  
Then the set of weights of $F$ is contained in $\operatorname{Conv}(W_G; \mu)$.
\end{lemma}
It is well-known that any linear functional on $V$ attains its maximum or minimum on the convex set $\operatorname{Conv}(D; v)$ at one of its vertices $w v$, for some $w \in D$.

In what follows, we use coordinates of the form $t = (t^{(p)}, t^{(q)})$,  
corresponding to the direct sum decomposition
${\mathbb{R}}^{p+q} = {\mathbb{R}}^p \oplus {\mathbb{R}}^q$.

\begin{lemma}
\label{lem:convqt}
Let $\nu = (\nu_1, \dots, \nu_{p+q}) \in {\mathbb{R}}^{p+q}$  
be a weakly decreasing sequence; that is, $\nu_1 \ge \cdots \ge \nu_{p+q}$.  

Suppose that
\begin{equation}
\label{eqn:Convtnu}
 t \in \operatorname{Conv}({\mathfrak{S}}_{p+q}; \nu),
\end{equation}
and
\begin{equation}
\label{eqn:convqt}
 \nu^{(p)} \in \operatorname{Conv}({\mathfrak{S}}_{p}; t^{(p)})
 \quad \text{and} \quad
 \nu^{(q)} \in \operatorname{Conv}({\mathfrak{S}}_{q}; t^{(q)}).
\end{equation}
Then it follows that
\[
 t^{(p)} \in {\mathfrak{S}}_{p} \cdot \nu^{(p)}
 \quad \text{and} \quad
 t^{(q)} \in {\mathfrak{S}}_{q} \cdot \nu^{(q)}.
\]
\end{lemma}
\begin{proof}
Without loss of generality, we may assume that  
$t_1 \ge \cdots \ge t_p$ and $t_{p+1} \ge \cdots \ge t_{p+q}$.  
Our goal is to show that $t = \nu$.  

By considering the maxima of linear functionals, we obtain the following chain of inequalities:
\begin{alignat*}{2}
  t_1 &\ge \nu_1 &&\ge t_1 \\
  t_1 + t_2 &\ge \nu_1 + \nu_2 &&\ge t_1 + t_2 \\
           & \vdots && \\
  t_1 + \cdots + t_p &\ge \nu_1 + \cdots + \nu_p &&\ge t_1 + \cdots + t_p.
\end{alignat*}
The left-hand inequalities follow from \eqref{eqn:convqt}  
and the assumption $t_1 \ge \cdots \ge t_p$,  
while the right-hand inequalities follow from \eqref{eqn:Convtnu}  
and the assumption $\nu_1 \ge \cdots \ge \nu_{p+q}$.  
Thus, we conclude that $t^{(p)} = \nu^{(p)}$.

Similarly, by considering minima of linear functionals, we obtain:
\begin{alignat*}{2}
  t_{p+q} &\ge \nu_{p+q} &&\ge t_{p+q} \\
  t_{p+q-1} + t_{p+q} &\ge \nu_{p+q-1} + \nu_{p+q} &&\ge t_{p+q-1} + t_{p+q} \\
           & \vdots && \\
  t_{p+1} + \cdots + t_{p+q} &\ge \nu_{p+1} + \cdots + \nu_{p+q} &&\ge t_{p+1} + \cdots + t_{p+q},
\end{alignat*}
and hence $t^{(q)} = \nu^{(q)}$.  

Therefore, under the assumption that $t_1 \ge \cdots \ge t_p$ and $t_{p+1} \ge \cdots \ge t_{p+q}$,  
we have shown that $t = \nu$.
\end{proof}

\medskip

\subsection{Certain Multiplicity-One Results of Specific Components in Finite-Dimensional Representations}
\label{sec:comb_fd_rep}
~~~
\newline
This section provides multiplicity-one results necessary for proving the isomorphisms of the natural morphisms \eqref{eqn:PK_isom} between $(\mathfrak{P}, K)$-cohomologies, using the PRV conjecture (now a theorem).

The main part of the following proposition is the second assertion.
\begin{lemma}
\label{lem:230082314}
Let $\nu=(\nu_1, \dots, \nu_{p+q}) \in {\mathbb{Z}}^{p+q}$ be a weight satisfying the standard
dominance condition: $\nu_1 \ge \cdots \ge \nu_{p+q}$.
Let $w^{(p)} \in {\mathfrak{S}}_p$, $w^{(q)} \in {\mathfrak{S}}_q$.  
Suppose that both $s \in {\mathbb{Z}}^{p+q}$ and $s+(w^{(p)} \nu^{(p)}; w^{(q)} \nu^{(q)})$ are
 $\Delta^+({\mathfrak{k}})$-dominant integral weights.  
\par\noindent
{\rm{(1)}}\enspace
Let $\tau \in \widehat K$ be an irreducible subrepresentation of the restriction $F^G(\nu)|_K$.  
Then the following three conditions on $\tau$
 are equivalent:
 \begin{itemize}
  \item[(i)]\enspace
  $[F^K(s + (w^{(p)} \nu^{(p)}; w^{(q)} \nu^{(q)})) : F^K(s) \otimes \tau] \ne 0$.
  \item[(ii)]\enspace
  $[F^K(s + (w^{(p)} \nu^{(p)}; w^{(q)} \nu^{(q)})) : F^K(s) \otimes \tau] = 1$.
  \item[(iii)]\enspace
  $\tau = F^K(\nu)$.
\end{itemize}

\noindent
{\rm{(2)}}\enspace
\[
[F^K(s+ (w^{(p)} \nu^{(p)}; w^{(q)} \nu^{(q)})): F^K(s) \otimes F^G(\nu)|_K] =1.
\]
\end{lemma}
\begin{proof}
(1)\enspace
We write $K = K_1 \times K_2$, 
with $K_1 = U(p)$ and $K_2 = U(q)$.  
Let $t \in \mathbb{Z}^{p+q}$ be the highest weight
of $\tau = \tau_1 \boxtimes \tau_2 \in \widehat{K}$, 
i.e., 
$t_1 \ge \cdots \ge t_p$
and  
$t_{p+1} \ge \cdots \ge t_{p+q}$.  

The implication (ii) $\Rightarrow$ (i) is immediate.  

We next prove the implication (i) $\Rightarrow$ (iii).  
Since $\tau$ occurs in $F^G(\nu)$, 
we have $t \in \operatorname{Conv}(\mathfrak{S}_{p+q}; \nu)$
by Lemma~\ref{lem:wtconv}.

On the other hand, 
if (i) holds,
then there exist weights $v^{(p)}$ of $\tau_1 \in \widehat{K_1}$ 
and $v^{(q)}$ of $\tau_2 \in \widehat{K_2}$ 
such that 
\[
s^{(p)} + w^{(p)} \nu^{(p)} = s^{(p)} + v^{(p)}, \quad
s^{(q)} + w^{(q)} \nu^{(q)} = s^{(q)} + v^{(q)}, 
\]
i.e.,
\[
w^{(p)} \nu^{(p)} = v^{(p)}, \quad
w^{(q)} \nu^{(q)} = v^{(q)}. 
\]
Hence, $\nu^{(p)} \in \operatorname{Conv}(\mathfrak{S}_p; t^{(p)})$
and $\nu^{(q)} \in \operatorname{Conv}(\mathfrak{S}_q; t^{(q)})$.  
It then follows from Lemma~\ref{lem:convqt} that 
\[
t^{(p)} \in \mathfrak{S}_p \cdot \nu^{(p)}, \quad
t^{(q)} \in \mathfrak{S}_q \cdot \nu^{(q)}. 
\]
Since $t_1 \ge \cdots \ge t_p$ and $t_{p+1} \ge \cdots \ge t_{p+q}$, we conclude that $t = \nu$.  
Thus, (i) $\Rightarrow$ (iii) is proved.  

Third, to prove the implication (iii) $\Rightarrow$ (ii), we consider $\tau = \tau_1 \boxtimes \tau_2$, according to the direct product decomposition $K = K_1 \times K_2$. Then the implication (iii) $\Rightarrow$ (ii) reduces to verifying the following two equalities:
\[
[F^{K_1}(s^{(p)} + w^{(p)} \nu^{(p)}) : F^{K_1}(s^{(p)}) \otimes F^{K_1}(\nu^{(p)})] = 1,
\]
\[
[F^{K_2}(s^{(q)} + w^{(q)} \nu^{(q)}) : F^{K_2}(s^{(q)}) \otimes F^{K_2}(\nu^{(q)})] = 1,
\]
which follow from the affirmative solution to the PRV conjecture (see~\cite{Ku}).

\par\noindent
(2)\enspace
For $\tau \in \widehat{K}$, set
\[
m(\tau) := \dim \operatorname{Hom}_K(\tau, F^G(\nu)|_K),
\]
i.e., the multiplicity of $\tau$ in the restriction $F^G(\nu)|_K$.  
Then we have:
\begin{multline*}
  [F^K(s+(w^{(p)} \nu^{(p)}, w^{(q)} \nu^{(q)})) : F^K(s) \otimes F^G(\nu)|_K]
\\
  =
  \sum_{\tau \in \widehat{K}} 
  m(\tau)
  [F^K(s+(w^{(p)} \nu^{(p)}, w^{(q)} \nu^{(q)})) : F^K(s) \otimes \tau].  
\end{multline*}

By part (1), the right-hand side is equal to $m(\tau)$
for $\tau := F^K(\nu)$, and vanishes for all other $\tau$.
Moreover, since $\nu$ is an extremal weight 
of the irreducible $G$-module $F^G(\nu)$,
we have $m(\tau) = 1$.  
Hence, the assertion follows.
\end{proof}

We are now ready to prove the following statement.

\begin{proposition}
\label{prop:23082610}
Suppose we are in the setting of Proposition~\ref{prop:23082631}.  
Then we have
\[
  [W_{\sigma} : W_{\pi} \otimes F|_K] = 1.  
\]
\end{proposition}

\begin{proof}
[Proof of Proposition~\ref{prop:23082610}]
Let $\nu \in \mathbb{Z}^{p+q}$ be the rearrangement 
of $-\lambda + \widetilde{\lambda}$ in non-increasing order, 
{\it i.e.}, $\nu_1 \ge \cdots \ge \nu_{p+q}$.  
Then, by the definition of $F$, $\nu$ is the highest weight of 
the irreducible representation $F$ of $G$; i.e.,  $F = F^G(\nu)$.  

Moreover, there exist $w^{(p)} \in \mathfrak{S}_p$ 
and $w^{(q)} \in \mathfrak{S}_q$ such that
\[
-\lambda^{(p)} + \widetilde{\lambda}^{(p)} = w^{(p)} \nu^{(p)}
\quad \text{and} \quad
-\lambda^{(q)} + \widetilde{\lambda}^{(q)} = w^{(q)} \nu^{(q)},
\]
by assumption~\eqref{eqn:lmdij} or~\eqref{eqn:lmdji}.  

Now the proposition follows from Lemma~\ref{lem:230082314}~(2), 
with $s = \lambda - \rho_{\mathfrak{g}}$.  
\end{proof}

\subsection{Morphisms between \texorpdfstring{$({\mathfrak{P}}, K)$}{(P, K)}-Cohomologies: Change of Coefficients}
\label{subsec:vert_B}

In this section, we provide a sufficient condition 
 for the vertical morphism $B_{\ast}$ 
 in diagram~\eqref{eqn:3_diagram} to be an isomorphism.
The main result of this section is stated in Proposition~\ref{prop:PK_spect_seq}.

We consider the tensor product $F \otimes W$ 
 as a $({\mathfrak{P}}, K)$-module, 
 where the $G$-module $F$ is regarded
 as a $({\mathfrak{P}}, K)$-module via restriction.  
For future reference, 
 we work in a slightly more general setting,
 replacing $F$ with an arbitrary finite-dimensional $({\mathfrak{P}}, K)$-module $V$ satisfying a certain multiplicity-free property as described below.  
\begin{prop}
\label{prop:PK_spect_seq}
Let $W_\pi$ be a coherent parameter for $\pi \in \operatorname{Disc}(G)$.
Let $V$ be a finite-dimensional $({\mathfrak{P}}, K)$-module
 such that 
\[
   \dim \operatorname{Hom}_K(W_{\pi}^\vee, V)=1.  
\]
Then there exists a natural isomorphism:
\[
   H^{\ast}({\mathfrak{P}}, K;\pi \otimes W_\pi^\vee)
  \overset \sim \longrightarrow
  H^{\ast}({\mathfrak{P}}, K;\pi\otimes V).  
\]
\end{prop}

\begin{proof}[Proof of Proposition~\ref{prop:PK_spect_seq}]
We recall a general fact: for any short sequence of $({\mathfrak{P}}, K)$-modules
\[
   0 \to U_1 \to U_2 \to U_3 \to 0, 
\]
 there is an associated long exact sequence in cohomology:
\[
  \cdots \to 
  H^i({\mathfrak{P}}, K;U_1)
  \to 
  H^i({\mathfrak{P}}, K;U_2)
  \to 
  H^i({\mathfrak{P}}, K;U_3)
  \to 
  H^{i+1}({\mathfrak{P}}, K;U_1)
  \to \cdots.
\]

We now take a filtration of $({\mathfrak{P}}, K)$-modules
\[
   0=V_0 \subset V_1 \subset \cdots \subset V_m = V
\]
 such that each successive quotient $V_j/V_{j-1}$ is irreducible as a $K$-module.  
By our assumption, 
  there exists a unique index $k$ such that 
$V_k / V_{k-1}$ is isomorphic to $W_{\pi}^{\vee}$.

It follows from Proposition~\ref{prop:dbc} that
$H^i({\mathfrak{P}}, K;\pi \otimes V_j/V_{j-1})=\{0\}$
 for all pairs $(i,j)$ except for 
 $(i,j)= (q(\delta),k)$.
 Hence, the long exact sequence associated to the short exact sequence $0 \to V_{j-1} \to V_j \to V_j/V_{j-1} \to 0$ induces the following natural isomorphisms:
\begin{multline*}
  H^{\ast}({\mathfrak{P}}, K;\pi \otimes W_{\pi}^\vee)
  \simeq
  H^{\ast}({\mathfrak{P}}, K;\pi \otimes V_{k}/V_{k-1})
  \overset \sim \leftarrow
  H^{\ast}({\mathfrak{P}}, K;\pi \otimes V_k)
\\
  \overset \sim \rightarrow
  H^{\ast}({\mathfrak{P}}, K;\pi \otimes V_{k+1})
  \overset \sim \rightarrow
  \cdots
  \overset \sim \rightarrow
   H^{\ast}({\mathfrak{P}}, K;\pi \otimes V_{m})
  =
  H^{\ast}({\mathfrak{P}}, K;\pi \otimes V).  
\end{multline*}
Thus, 
 the proposition is shown.  
\end{proof}

Return to diagram \eqref{eqn:B_fP},
if 
\[
  [W_{\sigma} : W_{\pi} \otimes F|_K] = 
  [W'_{\sigma'} : W'_{\pi'} \otimes F'|_K] = 1.  
\] 
(see e.g., Proposition~\ref{prop:23082610} for a sufficient condition),
then it follows from Proposition~\ref{prop:PK_spect_seq}
 that there are natural isomorphisms:
\[
\begin{alignedat}{3}
B_*   &\colon\quad & 
H^{\ast}({\mathfrak{P}}, K;\, \pi \otimes W_{\pi}^\vee) 
&\xrightarrow{\;\sim\;} 
& H^{\ast}({\mathfrak{P}}, K;\, \pi \otimes F \otimes W_\sigma^\vee), \\[1ex]
B'_*  &\colon\quad &
H^{\ast}({\mathfrak{P}}', K';\, \pi' \otimes F' \otimes W_{\sigma'}'^\vee) 
&\xrightarrow{\;\sim\;}
& H^{\ast}({\mathfrak{P}}', K';\, \pi' \otimes W_{\pi'}).
\end{alignedat}
\]

\begin{remark}
\label{rem:PK_spect_seq}
We may apply Proposition~\ref{prop:PK_spect_seq} in the context of Conjecture~\ref{conj:250629}, taking $V = S^r(\mathbb{C}^{p+q}) \otimes W_\sigma^\vee$, where $S^r(\mathbb{C}^{p+q})$ denotes the $r$th symmetric tensor of
the standard representation, and $W_\sigma$ is the coherent parameter of some discrete series representation $\sigma$.

In this case, $S^r(F)$ decomposes as a direct sum of irreducible $K$-modules via restriction $G \downarrow K$:
\[
  S^r(F) \simeq \bigoplus_{a + b = r} S^a(\mathbb{C}^p) \boxtimes S^b(\mathbb{C}^q),
\]
and hence the $K$-module $S^r(F) \otimes W_\sigma^\vee$ is multiplicity-free, by Pieri's rule.
\end{remark}

\subsection{An Unconditional Non-Vanishing Result}

In this section, we discuss a non-vanishing result for the natural morphisms between $(\fP,K)$-cohomologies associated with a pair $(\pi, \pi')$ of discrete series representations,
in the setting that $(\pi, \pi')$ forms \emph{non-elementary coherent pairs}. 
Conjecture~\ref{conj:250629} is formulated in that setting. In this section, we provide supporting evidence in Theorem \ref{thm:non-van}.

Let $(\pi, \pi')$ be a coherent pair
 with associated signature $[\delta \delta']$.  
As recalled in Corollary \ref{cor:coh_to_elem}, 
 there exists a discrete series representation $\sigma$ of $G$
 such that $(\sigma, \pi')$ constitutes an elementary pair
 having the same signature $[\delta \delta']$.  
Hence, the following definition makes sense:
\begin{equation}
\label{eqn:250702}
  d(\pi, \pi'):=\operatorname{min}
\{|\pi-\sigma|: \text{$(\sigma, \pi')$ is an elementary pair 
 with signature $[\delta \delta']$}\},   
\end{equation}
 where, 
 for two discrete series representations $\pi$
 and $\sigma$ of $G$ with Harish-Chandra parameters $\lambda$ and $\mu$, respectively, 
 we define 
\[
   |\pi-\sigma|:=\underset{i=1}{\overset{p-1}\sum}|\lambda_i-\mu_i|.  
\]
Note
 that $d(\pi, \pi')>0$
 if and only if $(\pi, \pi')$ is non-elementary.  

\begin{theorem}
\label{thm:non-van}
Conjecture~\ref{conj:250629} holds for any coherent pair $(\pi, \pi')$ with \newline  $d(\pi, \pi')=1$.  
\end{theorem}

The proof of Theorem~\ref{thm:non-van} is based on a new technique \lq\lq{translation of symmetry breaking}\rq\rq\
 for the pair $(G,G')$ developed in Section~\ref{sec:SBO_translation}, 
 together with classical techniques on translations of representations of individual groups $G$ or $G'$.

The general case will be addressed in a forthcoming paper.

The rest of this section is devoted to proving Theorem~\ref{thm:non-van}.  
The assumption $d(\pi, \pi')=1$ implies
 that there exists a discrete series representation $\sigma$ of $G$
 with Harish-Chandra parameter $\mu \in C(\delta)$
 such that $(\sigma, \pi')$ is an elementary pair
 and 
\[
   \text{$\mu=\lambda+e_i$ \quad or \quad $\mu=\lambda-e_i$} 
\]
for some $i$ ($1 \le i \le p+q$).  
We treat the case:
\begin{equation}
\mu=\lambda+e_i \in C(\delta);
\end{equation}
the other case $\mu=\lambda-e_i$ is similar.  

We recall the notations used in Conjecture \ref{conj:250629}
 applied to the setting where $(F, F'):=({\mathbb{C}}^{p+q}, {\bf{1}})$ and $\sigma':=\pi'$.  
In particular, $\operatorname{pr}_{F} \colon F \to F'$ denotes the projection to the trivial $G'$-module $F'={\bf{1}}$.

Let $T \colon \sigma \to \sigma'$ be a non-trivial symmetry breaking operator between the elementary pair $(\sigma, \sigma')$, which is unique up to scalar multiplication.
Let $W := W_{\sigma} \in \widehat K$ and $W':={W'}_{\sigma'} \in \widehat{K'}$ be the coherent parameters (Definition~\ref{def:W_lambda}), 
 which satisfy $[W|_{K'}:W']=1$ as shown in Theorem~\ref{thm:elementary_pair_equiv},
 and let $\operatorname{pr}_{W} \colon W^{\vee} \to W'^{\vee}$ be a non-trivial $K'$-homomorphism.  

We take a nontrivial symmetry-breaking operator
$S \colon \pi \to \pi'$ between the coherent pair $(\sigma, \sigma')$.
Let $W_{\pi} \in \widehat K$ and $W':={W'}_{\pi'} \in \widehat{K'}$ denote the coherent parameters of $\pi$ and $\pi'$, respectively.

Since both $\lambda$ and $\mu$ 
  belong to the same regular dominant chamber $C(\delta)$, it follows from Proposition~\ref{prop:prireg} that
\[
\pi\simeq \psi_{\mu}^{\lambda}(\sigma) \simeq
  P_{\lambda}(\sigma),
  \quad
\sigma\simeq \psi_{\lambda}^{\mu}(\pi) \simeq
  P_{\mu}(\pi).
\]
In particular, $\dim\operatorname{Hom}_{G}(\sigma, \pi \otimes F)=1$.

We now apply Theorem~\ref{thm:23081404} to the symmetry breaking operator $S\colon \pi \to \pi'$ in conjunction with $\operatorname{pr}_{i, +} \in
\operatorname{End}_{G}(\pi \otimes F)$.
If the pair of Harish-Chandra parameters $(\lambda, \lambda')$ for the discrete series representations $(\pi, \pi')$  satisfies
\begin{equation}
\label{eqn:lmd_cond1}
\lambda_i + \frac 1 2 \not \in \{\lambda_j':j=1, \dots, \widehat{p}, \dots, p+q\}, 
\end{equation}
then Theorem~\ref{thm:23081404} ensures the non-vanishing of the composition map
\[
(S \otimes \operatorname{pr}_{F})\circ \operatorname{pr}_{i, +}.
\]
Since 
$\operatorname{pr}_{i, +}$
acts as a non-zero scalar multiplication on the subspace $P_\mu(\pi\otimes F)$,
the projection operator $\operatorname{pr}_{i, +}$
induces a non-zero $G$-homomorphism
\[
A\colon \sigma \hookrightarrow \pi\otimes F
\]
via the isomorphism $P_\mu(\pi \otimes F) \sim \sigma$.
In turn, the composition
 $(S \otimes \operatorname{pr}_{F})\circ A$
induces a {\emph{non-zero}} symmetry breaking operator
from $\sigma$ to $\pi'=\sigma'$,
which must take the form $c T$ for some non-zero constant $c$.
By renormalizing $S$, we may assume without loss of generality that $c=1$.

Since $\sigma'=\pi'$ and $F' ={\bf{1}}$, 
we can take $A' \colon \sigma' \to \pi' \otimes F'$ to be the identity map.
Thus, we have obtained a commutative diagram \eqref{eqn:ST_comm}, 
{\it{i.e.,}} 
the identity
\[
 (S \otimes \operatorname{pr}_F) \circ A = A' \circ T
\]
in $\operatorname{Hom}(\sigma, \pi'\otimes F')$.

Passing to $({\mathfrak{P}}, K)$-cohomologies, 
 we obtain the commutative diagram \eqref{eqn:fP_TS_diagram}.  
 The virtical homomorphism
 in diagram \eqref{eqn:fP_TS_diagram},
i.e., the homomorphism
\begin{equation*}
A_*\colon
   H^j({\mathfrak{P}}, K; \sigma \otimes W^\vee)
   \rightarrow
   H^j({\mathfrak{P}}, K; \pi\otimes F \otimes W^\vee).  
\end{equation*}
is bijective by Proposition \ref{prop:23082631} applied to $\tilde{\lambda}=\mu$,
since the assumption \eqref{eqn:lmdij}
 or \eqref{eqn:lmdji} there is clearly satisfied
 when $\mu=\lambda \pm e_i$.

We now present diagram~\ref{eqn:3_diagram}
in Conjecture~\ref{conj:250629} again:
\begin{equation*}
\xymatrix@C=80pt{
H^{\ast}({\mathfrak{P}}, K;\sigma \otimes W^\vee)
 \ar[r]_{T \otimes \operatorname{pr}_{W}}^{\sim}
 \ar[d]_{A_{\ast}}
& 
H^{\ast}({\mathfrak{P}}', K';{\sigma}'\otimes {W'}^\vee)
 \ar[d]^{A_{\ast}'}
\\
H^{\ast}({\mathfrak{P}}, K;\pi \otimes F \otimes W^\vee)
 \ar[r]^{S \otimes \operatorname{pr}_{F} \otimes \operatorname{pr}_{W}} 
&
 H^{\ast}({\mathfrak{P}}', K';\pi' \otimes F' \otimes {W'}^\vee)
 \ar[d]^{B_{\ast}'}
\\
H^{\ast}({\mathfrak{P}}, K;\pi \otimes W_{\pi}^\vee)
 \ar[u]^{B_{\ast}}
 \ar@{.>}[r]
& 
H^{\ast}({\mathfrak{P}}', K';\pi' \otimes {W'_{\pi'}}^\vee)
}
\end{equation*}

Since $(\sigma, \sigma')$ is an elementary pair, 
 we already know
 that the horizontal map $T \otimes \operatorname{pr}_{W \to W'}$ induces 
 an isomorphism:
\[
\xymatrix@C=70pt{
  H^{\ast}({\mathfrak{P}}, K;\sigma\otimes W^\vee)
 \ar[r]_{T \otimes \operatorname{pr}_{W}}^{\sim}
 & 
H^{\ast}({\mathfrak{P}}', K';
{\sigma}^{\prime} \otimes W'^\vee).  
}
\]

On the other hand, it follows from
Proposition~\ref{prop:PK_spect_seq} 
 and Remark~\ref{rem:PK_spect_seq}
 that the left vertical arrow; i.e., 
 the morphisms
\[
B_*   \colon 
H^{\ast}({\mathfrak{P}}, K;\, \pi \otimes W_{\pi}^\vee) \rightarrow
 H^{\ast}({\mathfrak{P}}, K;\, \pi \otimes F \otimes W_\sigma^\vee), 
\]
is bijective.

Since $\mathbb{F}'$ is the trivial one-dimensional representation
and  since $\sigma' = \pi'$,
the vertical maps $a'_*$ and $B'_*$ in the right-hand side corner are
clearly bijective.

Thus, 
all the natural morphisms in diagram~\ref{eqn:3_diagram}
are bijective between relative cohomologies, thereby establishing Theorem~\ref{thm:non-van}.

\newpage

\part{Global applications}\label{global}

\noindent {\bf Warning:}  Because the authors come from different mathematical traditions, we have decided that Part \ref{global} will use a different convention for $\overline{\partial}$-cohomology than the rest of the paper.  In the next three sections, we will compute coherent cohomology by a complex of relative $(\fP^-,K)$-cohomology, as in \cite{H90}.  This means we have to change various conventions.  We will be able to cite the results of the preceding sections, provided we work with the dictionary in Table~\ref{tab:250503}.  Instead of using the convention of the upper left corner of that table, we work with the convention in the lower right.  This preserves the cohomological degree $q = q(\delta)$ at the cost of dualizing the discrete series representation $\pi = \pi_\lambda^\delta$.  In particular we quote the theorem in Example \ref{ex:dual} rather than 
Theorem~\ref{thm:elementary_pair_equiv}.  We apologize for any possible confusion.

\section{Coherent Cohomology of Shimura Varieties}\label{Shimura7}

\subsection{Generalities}\label{gener}

Let 	$(\G,X)$ be a Shimura datum; $\G$ is a connected reductive group over $\QQ$ and $X$ is a homogeneous space for $\G(\RR)$, consisting in a $\G(\RR)$-conjugacy
class of homomorphisms $h$ of the real group $\IS := R_{\CC/\RR}\G_{m,\CC}$ to $\G_\RR$.   The homomorphisms $h$ satisfy a short list of axioms chosen by Deligne
(see, for example, \cite[\S 1.1]{BHR}) that guarantee that for each $h \in X$ the centralizer $K_h \subset \G(\RR)$ of $h(\IS(\RR))$ is contained in a maximal compact modulo center subgroup of $\G(\RR)$ and contains a maximal compact modulo center subgroup of the identity component $\G(\RR)^0 \subset \G(\RR)$.  The adjoint representation of $K_h$ decomposes $\fg = Lie(\G)$ as the sum of three subspaces
\begin{equation}\label{HCdecomp} \fg = \fk_h \oplus \fp^+_h \oplus \fp^-_h 
\end{equation}
where $\fp^+_h$ and $\fp^-_h$ are abelian subalgebras of $\fg$; in the case of the unitary groups of \S \ref{subsec:ugps} this is just the Harish-Chandra decomposition~\eqref{eqn:HCD}.

We fix one $h \in X$ for the moment and drop the subscript.   Choose a maximal torus $T \subset K$, with Lie algebra $\ft$ and a system of positive roots of $\ft$ containing the roots
of $\ft$ on $\fp^-$, as in Section~\ref{subsec:ugps}.  The sum $\fP^- := \fk \oplus \fp^-$ is a maximal parabolic subalgebra.
Let $\IP \subset \G$ be the corresponding parabolic subgroup, and let $\hat{X} = \G/\IP$ be the flag variety.  The group $\G$ acts on the parabolic subalgebras of $\fg$ by the adjoint representation, and the map $x \mapsto Stab_{\G}(\fk_x \oplus \fp^-_x)$ defines a $\G(\RR)$-equivariant embedding 
$$\beta:  X \hookrightarrow \hat{X}(\CC)$$
which defines the $\G(\RR)$-invariant holomorphic structure on $X$.  The Shimura variety  $S(\G,X)$ is a canonical $\G(\af)$-equivariant model, defined
 over a certain number field $E(\G,X)$, of the projective limit 
$$\varprojlim_{K_f \subset \G(\af)}  G(\QQ)\backslash X \times \G(\af)/K_f$$ 
of quasiprojective complex algebraic varieties.  

Let $Rep(\IP)$ denote the category of ind-finite algebraic representations of $\IP$; it is canonically equivalent to the category of ind-finite
$\G$-equivariant vector bundles on $\hat{X}$.  The latter category has a canonical $E(\G,X)$-rational structure.  There is a tensor functor 
$$ [\bullet]:  Rep(\IP) \ra \{\text{ $\G(\af)$-equivariant vector bundles over } S(\G,X)\}$$
defined by identifying the representation $(\rho,\CV)$ of $\IP$ with the $\G$-equivariant vector bundle $\CV_{\hat{X}}$ over $\hat{X}$, and then defining
$$[\CV] =  \varprojlim_{K_f \subset \G(\af)}  \G(\QQ)\backslash \CV_{\hat{X}} \times \G(\af)/K_f.$$

We are concerned with two specific kinds of representations of $\IP$, and their tensor products.  An irreducible representation
of $K_h$ will be denoted by $W = W_\lambda$, where $\lambda-\rho_\mathfrak{g}:  \ft \ra \CC$ is the highest weight of $W$, as in Definition \ref{def:W_lambda}, where in the applications $\lambda$ is the Harish-Chandra parameter of a discrete series. 
This inflates to a representation of $\IP$, also denoted $W$, through the projection of $\IP$ on its Levi quotient $K_h$.
As above, we identify $W$ with the corresponding $\G$-equivariant vector bundle $W_{\hat{X}}$ and let $[W]$ be the automorphic vector bundle on $S(\G,X)$.

An irreducible representation of $\G$ is denoted $F$; by restriction to $\IP$ it defines a $\G$-equivariant vector bundle $F_{\hat{X}}$, and thus an automorphic
vector bundle $[F]$ on $S(\G,X)$.  Such a bundle is endowed with a flat connection $\nabla_F:  [F] \ra [F]\otimes \Omega^1_{S(\G,X)}$.  

\subsection{Unitary Group Shimura Varieties}\label{ugsv}
The following material is adapted, with slight changes, from \cite{GHL}.   Let $\CK$ be a totally imaginary quadratic extension of the totally real field $\CKp$, with
$[\CKp:\QQ] = 2$; let $c \in Gal(\CK/\CKp)$ be complex conjugation. 
Let $V$ be an $n$-dimensional vector space over $\CK$, endowed with a hermitian form relative to $c$; let $H = R_{\CKp/\QQ}U(V)$ be the corresponding
unitary group.   Let $\Sigma^+ = Hom(\CKp,\RR)$ and let $\Sigma \subset Hom(\CK,\CC)$ be a CM type.  Thus
$$Hom(\CK,\CC) = \Sigma \coprod \Sigma\cdot c$$
and  restriction to $\CKp$ defines a bijection $\Sigma \leftrightarrow \Sigma^+$.  Thus any $v \in \Sigma^+$ extends uniquely to an element of $\Sigma$,
also denoted $v$.  Then the given hermitian form on $V$ defines on
$V_v := V \otimes_{\CK,v} \CC$  a hermitian form with signature $(r_v,s_v)$; we let $H_v = U(V_v)$.

We
define a Shimura datum $(H,Y_V)$ as in \cite[\S 2.2]{H21}. Explicitly, the base point $h_V  \in Y_V$ is given as a homomorphism 
$\IS \ra H_\RR \isoarrow \prod_{v \in \Sigma^+} H_v$, which on complex points is given by

\begin{equation}\label{yVz} h_{V,v}(z) = \begin{pmatrix}  (z/\bar{z})I_{r_v} & 0 \\ 0 & I_{s_v} \end{pmatrix} \end{equation}
The following lemma is then obvious: We record it here in order to define parameters for automorphic vector bundles in the next sections.  

\begin{lemma}\label{stab}  
Let $h \in Y_V$.   Its stabilizer $K_h=:K_{H,\infty}$ in $H_\infty = \prod_{v \in \Sigma} H_v$ is isomorphic to $\prod_{v \in \Sigma^+} U(r_v)\times U(s_v)$.
\end{lemma}
We also write $\prod_{v \in \Sigma}$ to emphasize that $U(r_v)\times U(s_v)$ and $U(s_v)\times U(r_v)$ play different 
(complex conjugate) roles in the theory of Shimura varieties. When necessary, we write $\fP = \prod_{v \in \Sigma} \fP_v$,
$K_h = \prod_{v \in \Sigma} K_v$.

The reflex field $E(H,Y_V)$ is the subfield of $F^{Gal}$ determined as the stabilizer of the cocharacter $\kappa_V$ with $v$-component  $\kappa_{V,v}(z) = \begin{pmatrix}  zI_{r_v} & 0 \\ 0 & I_{s_v} \end{pmatrix}$.  In particular, if there is $v_0 \in S_\infty$ such that $s_{v_0} > 0$ but $s_v = 0$ for $v \in S_\infty \setminus \{v_0\}$ -- the type of unitary groups that we will be mainly interested in later -- then $E(H,Y_V)$ is the subfield $\imath_{v_0}(F) \subset \CC$. 

We will also fix the following notation: Let $V'\subset V$ be a non-degenerate subspace of $V$ of codimension 1. We write $V$ as the orthogonal direct sum $V' \oplus V_1$ for some $1$-dimensional hermitian space $V_1$, and consider the unitary groups $H':=R_{\CKp/\QQ}U(V')$ and $H^\flat:=R_{\CKp/\QQ}(U(V')\times U(V_{1}))$ over $\QQ$. Obviously, there are natural inclusions $H'\subset H''\subset H$, and a homomorphism of Shimura data
\begin{equation}\label{inclusionmap} 
(H^\flat,Y_{V'}\times Y_{V_1}) \hookrightarrow (H,Y_V).
\end{equation}
The space $V_1$ is completely determined by the inclusion of $V'$ in $V$; it is ignored in the GGP conjectures but needs to be
recalled in order to define period relations consistently.

Thus we are using the letter $H$ to denote a unitary group, viewed as algebraic group over $\QQ$; $G$ is reserved for  unitary groups over $\RR$, as before.

In the applications, the base point $h_V \in Y_V$ will always be taken to be the image of a CM point, as in \cite[\S 2.2]{H21}.

\subsection{Cohomology of automorphic vector bundles on \texorpdfstring{$S(H,Y_V)$}{S(H,YV)}.}\label{avbH}
Notation for automorphic vector bundles is as in \S \ref{gener}.  Thus representations of $K_h$ are denoted $W$.  With respect to the factorization in Lemma \ref{stab}, we can write $W = \otimes_{v \in \Sigma} W_v$, where each $W_v$ is an irreducible representation of $U(r_v) \times U(s_v)$.  Similarly, representations of $H$, denoted $F$, can be written $F = \otimes_{v \in \Sigma} F_v$, where $F_v$ is an irreducible algebraic representation of the Lie group $U(r_v,s_v)$, or equivalently
of its extension to $GL(n,\CC)$.

We thus have automorphic vector bundles of the form $[W]$ and $[F]$ over $S(H,Y_V)$.  To define the coherent cohomology of Shimura varieties properly
one needs to introduce the theory of toroidal compactifications $_{K_f}S(H,Y_V) \hookrightarrow _{K_f}S(H,Y_V)^{tor}$
and their relations as the level subgroups vary.  As a substitute for cuspidal cohomology we work with the interior cohomology
 $H^*_!(S(H,Y_V),[W])$, defined as in  \cite[(3.3)]{H21} as 
 \begin{equation}\label{cohcoh}
\begin{aligned}
& H^*_!(S(H,Y_V),[W]) :=  \\ &\varinjlim_{K_f} Im[H^*(_{K_f}S(H,Y_V)^{tor},[W]^{sub}) \ra H^*(_{K_f}S(H,Y_V)^{tor},[W]^{can})],
\end{aligned}
\end{equation} 
We refer to \cite[\S 3.2]{H21} for  notation and explanations.

The local results on $\db$ cohomology in
Section~\ref{subsec:basic_Pk_cohomology}\footnote{Which we interpret by means of the dictionary in Table \ref{tab:250503}!}
have immediate consequences for the interior cohomology.  Let $\CA(H)$, $\CA(H')$ denote the spaces of automorphic forms on $H(\QQ)\backslash H(\bA)$ and $H'(\QQ)\backslash H'(\bA)$, respectively.  Let
$\CA_0(H)$ and $\CA_0(H')$ denote the subspaces of cusp forms.  Let $F$ and $W$ be as above.  The methods of \cite{H90} provide the following inclusions:

\begin{equation}\label{automorphic} H^*(\fP^-,K_h;\CA_0(H)\otimes W) \hookrightarrow H^*_!(S(H,Y_V),[W]);
\end{equation}
\begin{equation}\label{automorphicF}
H^*(\fP^-,K_h;\CA_0(H)\otimes F\otimes W) \hookrightarrow H^*_!(S(H,Y_V),[F]\otimes [W]).
\end{equation}
 Let $\lambda = (\lambda_v)_{v \in \Sigma}$ be the
Harish-Chandra parameter of a tensor product $\sigma^\vee = \otimes_{v \in \Sigma} \sigma_v^\vee$ of discrete series representations of $\prod_v U(r_v,s_v)$, 
$$\pi_\lambda = \otimes_v \pi_{\lambda_v}; ~~ W = W_\lambda = \otimes_v W_{\lambda_v}.$$   
For $? = W, F\otimes W$, we let $H^*_t(S(H,Y_V),[?])$ %$H^*(S(H,Y_V),[?]) \subset H^*_!(S(H,Y_V),[?])$ 
denote the
image of the subspace of $H^*(\fP,K_h;\CA_0(H)\otimes ?)$ generated by cuspidal automorphic representations 
with archimedean component  in the discrete series.  

For $W$ sufficiently regular it is known that \eqref{automorphic} is an isomorphism, and we will always assume this to be the case.  

\begin{remark}\label{interior}  The inclusion 
$H^*_t(S(H,Y_V),[W]) \subset H^*_!(S(H,Y_V),[W])$ is an equality for $W = W_\lambda$ with sufficiently regular $\lambda$
(see for example \cite[Corollary 5.3.3]{H90}).  In particular, for sufficiently regular $\lambda$ the complex vector space $H^*_t(S(H,Y_V),[W_\lambda])$ has a natural rational structure over a number field $E_\lambda$ which contains the reflex field for $S(H,Y_V)$ and over which the automorphic vector bundle $[W_\lambda]$ is defined.  Moreover, the Galois group $Gal(\Qbar/\QQ)$ permutes the Shimura data $S(H,Y_V)$ and $[W_\lambda]$, and thus acts on the $H^*_t(S(H,Y_V),[W_\lambda])$, as explained in \cite{H13}. 

More generally, by applying appropriate operators $A:  \CA_0(H) \ra \CA_0(H)$ in the group algebra of $H(\af)$ one can guarantee that the map
$$H^*(\fP^-,K_h;Im(A)\otimes F\otimes W) \hookrightarrow A(H^*_!(S(H,Y_V),[F]\otimes [W]))$$
derived from \eqref{automorphicF}
is an isomorphism.  For example, one can let $A$ denote a pseudocoefficient of a discrete series representation of $H(\CK_v)$ for some prime $v$ of $\CK$ split over $\CKp$; or it can be an operator that picks out an eigenspace for the unramified Hecke algebra whose base change to $GL(n+1)_\CK$ is cuspidal.  Thus we will always proceed as if the left-hand side of \eqref{automorphic} or \eqref{automorphicF} has a natural rational structure.
\end{remark}

The following Theorem is the global consequence of Proposition~\ref{prop:dbc}.  

\begin{theorem}\label{cohc}  %Let $W = W_\Lambda$, $\Lambda = (\Lambda_v)_{v \in \Sigma}$. 
For each $v \in \Sigma^+$ choose a Harish-Chandra parameter $\lambda_v$ of a discrete series representation $\pi_{\lambda_v}$ of $H_v$,
so that
$\lambda = (\lambda_v)_{v \in \Sigma}$ is a Harish-Chandra parameter for $H_\infty$.  

There is a unique irreducible representation $W_{\tlambda} = \otimes_v W_{\lambda_v}$ of $K_h$ and an integer $q = \sum_v q_v$ with the following property.  
$$H^i_t(S(H,Y_V),[W_{\alpha}]) \neq 0 \Rightarrow \alpha = \lambda \text{ and } i = q.$$
Moreover, $W_{\tlambda}$ has highest weight $\lambda - \rho_{\mathfrak{g}}$, where $\rho_{\mathfrak{g}}$ is the half sum of the positive roots for the system containing the roots in $\fp^+$. 

Moreover, suppose $\lambda$ is sufficiently regular, so that \eqref{automorphic} is an isomorphism.  Then $W_{\tlambda}$ is linked to $\lambda$ in the following sense:
 \begin{equation}\label{isotypic}
\begin{aligned}
H^\bullet_t(S(H,Y_V),[W_{\tlambda}]) &=  H^{q}_!(S(H,Y_V),[W_{\tlambda}]) \\
&= \bigoplus_{\Pi = \pi_\lambda^\vee \otimes \Pi_f} H^{q}(\fP^-,K_h; \pi_\lambda^\vee \otimes W_{\lambda})\otimes \Pi_f
\end{aligned}
\end{equation}
where $\Pi$ runs over cuspidal cohomology.
\end{theorem}

\begin{remark}  Note the dualization of $\pi_\lambda$ in the sum on the right hand side of \eqref{isotypic}, which is imposed by the dictionary in Table \ref{tab:250503}.
\end{remark}

\begin{proof}  The first part of the theorem follows immediately from Proposition~\ref{prop:dbc} and the definition of $H^*_t$.  The integer $q_v$ is determined by the property that
$$H^{q_v}(\fP_v^-,K_v;\pi_v^\vee\otimes W_{\lambda_v}) \neq 0.$$
The last statement follows from Remark \ref{interior}.
\end{proof}

\subsection{Cohomology of Reducible Automorphic Vector Bundles}

Recall that we always assume $\lambda$ is sufficiently regular.  We let $F$ be an irreducible representation of $H$.  When $\pi$ is a cuspidal automorphic representation, we write $\pi = \pi_\infty \otimes \pi_f$ with $\pi_\infty$ the archimedean and $\pi_f$ the finite part.  
For $ ? = W, F \otimes W$, the graded vector space $H^\bullet_!(S(H,Y_V),[?])$ decomposes as a direct sum of irreducible representations of $H(\af)$, with finite multiplicities, and we write
$$H^\bullet(S(H,Y_V),[?])[\pi_f]$$
for the $\pi_f$-isotypic subspace.  

Let $\delta$ be a decoration, $\tlambda = (\tlambda_v)_{v \in \Sigma}$ as in Theorem \ref{cohc}, and assume $W = W_{\tlambda}$.   Let
$H^\bullet_{\tlambda}(S(H,Y_V),[?])$ denote
the subspace
$$ \bigoplus_{\pi_{\tlambda}^\vee\otimes \pi_f \subset \CA_0(H)} H^\bullet(S(H,Y_V),[?])[\pi_f] \subset H^\bullet_!(S(H,Y_V),[?]).$$
When $? = W$ this is computed by Theorem \ref{cohc}.  As for the case $? = F\otimes W$, we have the following Proposition:

\begin{prop}\label{redu} Suppose $\pi = \pi_{\tlambda}^{\delta,\vee}$, $\sigma = \pi_\ttl^{\delta,\vee}$ and $F$ satisfy the inequalities of Proposition \ref{prop:23082631}.   
Write $W_\pi = W_{\tlambda}$, $W_\sigma = W_\ttl$.  Then
 the  isomorphism of $(\fP^-,K)$-cohomology determines an isomorphism of graded $H(\af)$-modules
$$H_{\tlambda}^\bullet(S(H,Y_V),[F]\otimes [W_\ttl]) \isoarrow H^\bullet_!(S(H,Y_V),[W_\lambda])$$
concentrated in degree $q(\tlambda)$.
This isomorphism is rational over the number field $E_{\tlambda}$ of Remark \ref{interior}.
\end{prop}
\begin{proof}  By \eqref{automorphicF} the left hand side is computed by
\begin{equation}\label{dsum}
\bigoplus_{\pi_f}\bigoplus_{\pi_\infty\otimes \pi_f \subset \CA_0(H)} H^\bullet(\fP^-,K_h;\pi_\infty\otimes F\otimes W_\ttl),
\end{equation}
where the first sum is taken over $\pi_f$ such that $\pi_0 := \pi_{\tlambda}^\vee\otimes \pi_f$ is cuspidal automorphic.  Now suppose $\pi = \pi_\infty\otimes \pi_f$ is cuspidal automorphic, with $\pi_\infty$ as yet to be determined.  By \cite{La}, $\pi$ and $\pi_0$ both admit base change to cohomological automorphic representation $BC(\pi_0), BC(\pi)$ of $GL(n)_\CK$.   Since $\pi_0$ and $\pi$ both have the same finite part, so do $BC(\pi_0)$ and $BC(\pi)$, so by Jacquet--Shalika strong multiplicity one $BC(\pi_0) = BC(\pi)$.  This implies, again by \cite{La}, that $\pi_\infty$ and $\pi_{\tlambda}^\vee$ are cohomological with the same 
infinitesimal character.  By the regularity hypothesis they belong to the same discrete series $L$-packet. 

But now we can decompose $\pi_\infty \otimes F$ as a sum according to infinitesimal character, and we have
$$H^\bullet(\fP^-,K_h;\pi_\infty\otimes F\otimes W_\ttl) \isoarrow H^\bullet(\fP^-,K_h;P_{\ttl}(\pi_\infty\otimes F)\otimes W_\ttl).$$
Now Proposition \ref{prop:23082631} (2) applies to show that $\pi_\infty \otimes \pi_f$ only contributes to \eqref{dsum} if 
$\pi_\infty \isoarrow \pi_{\lambda}^\vee$.  
Hence $\pi \isoarrow \pi_0$; and by multiplicity one for $H$, they must be equal.  Finally, by Proposition \ref{prop:dbc} again, 
$$H^\bullet(\fP^-,K_h;P_{\ttl}(\pi_\infty\otimes F)\otimes W_\ttl) = H^\bullet(\fP^-,K_h;\pi_{\ttl}^\vee\otimes W_\ttl)$$
is concentrated in degree $q(\tlambda)$ and is of dimension $1$.
\end{proof}

\section{Restriction of coherent cohomology}
\label{sec:rest_cohomology}
Let $H = R_{\CKp/\QQ}U(V)$, $H' = R_{\CKp/\QQ}U(V')$, and $H^\flat$ be as in \S \ref{ugsv}.  

\subsection{Elementary Pairs}

We now choose discrete series representations $\pi = \pi_{\tlambda}$, $\pi' = \pi'_{\tlambda'}$, of $H_\infty$ and $H'_\infty$, respectively, as in the statement of
Theorem \ref{cohc}, so that for each $v \in \Sigma$
the pair $(\pi_{\tlambda,v}, \pi'_{\tlambda',v})$ is an elementary pair as in Definition \ref{def:elementary}.  The choice of Weyl chamber $C_v$ as in Theorem \ref{cohc} for
each $\pi_{\tlambda_v}$ and  $\pi'_{\tlambda',v}$ is unique and is therefore dropped from the notation.     Let $U(V_1)$ and $H^\flat$ be as in \S \ref{ugsv}.  We define
the representation $W_{\lambda,\lambda'}^\flat = W_{\lambda'}\otimes \xi(\lambda,\lambda')$ of $H^\flat$, where $\xi(\lambda,\lambda')$ is the unique character of $U(V_1)$ with the property that the restriction of 
$W_{\lambda,\lambda'}^\flat$ to the center $Z_H$ of $H$, which is a subgroup of $H^\flat$, coincides with the character of $Z_H$ on $W_{\tlambda}$.
Let $q = \sum_v q(\lambda_v) = \sum_v q(\lambda'_v)$.

\begin{theorem}\label{elemSh}   With the above notation, the homomorphism
$$R_\tlambda:  H^q_!(S(H,Y_V),[W_{\tlambda}]) \ra H^q_!(S(H^\flat,Y_{V'}\times Y_{V_1}),[W_{\lambda,\lambda'}^{\flat}])$$
is rational over the reflex field of $[W_{\lambda,\lambda'}^{\flat,\vee}]$ and commutes with Galois automorphisms in the obvious way:
For any $\sigma \in Gal(\bar{\QQ}/\QQ)$, $\sigma(R_\tlambda)$ is canonically equal to the restriction map
$$H^q_!(\sigma(S(H,Y_V)),\sigma([W_{\tlambda}])) \ra H^q_!(\sigma(S(H^\flat,Y_{V'}\times Y_{V_1})),\sigma([W_{\lambda,\lambda'}^{\flat}]))$$
where $\sigma(S(H,Y_V))$, $\sigma(S(H^\flat,Y_{V'}\times Y_{V_1}))$, and the rest are isomorphic to the Shimura varieties determined by Langlands's conjecture on conjugation of Shimura varieties (see \cite{Milne}).
%{\color{red}  Insert the action of $Gal(\bar{\QQ}/\QQ)$ on all the data.}

Moreover, let $c \in H^q_!(S(H,Y_V),[W_{\tlambda}])$ be a non-zero cohomology class.  Then up to translating $S(H^\flat,Y_{V'}\times Y_{V_1})$ by an element $\eta \in H(\af)$, 
the image of $c$ in $H^q_!(S(H^\flat,Y_{V'}\times Y_{V_1})\cdot \eta,\eta^*[W_{\lambda,\lambda'}^{\flat}])$ does not equal zero.
\end{theorem}

\begin{proof}  The first part is a consequence of the theory of canonical models of automorphic vector bundles.  The second follows from the Burger-Sarnak method, by the arguments in  \cite[\S 2]{HL} and \cite[\S 4]{H14}.
\end{proof}

\subsection{General Coherent Pairs}  

In this section we assume Conjecture \ref{conj:250629}.   For each  $v \in \Sigma$, let $(\pi_v,\pi'_v)$ be a coherent pair with parameters $(\ttl_v,\ttl'_v)$ that are far enough from the walls for Conjecture \ref{conj:250629} to be valid, and write $\pi = \otimes_{v \in \Sigma}\pi_v, \pi' = \otimes_{v \in \Sigma} \pi'_v$.   Let $W_{\ttl} = \otimes_v W_{\ttl_v}$, $W'_{\ttl'} = \otimes_v W'_{\ttl_v}$.

For each $v$, we let $(\sigma_v, \sigma'_v)$ be an elementary pair,  We let $(F_v,F'_v)$ be a pair of irreducible (algebraic) representations of $H(\CK_v)$ and $H'(\CK_v)$, respectively, that satisfy the conditions of Conjecture \ref{conj:250629}.  Let $[F]$ denote the  (flat) automorphic vector bundle on $S(H,Y_V)$ corresponding to $\{F_v\}$; the bundle corresponding to $\{F'_v\}$ is trivial and is omitted from the notation in what follows.   We let $\tlambda_v, \tlambda'_v$ be the Harish-Chandra parameters of $\sigma_v$, $\sigma'_v$, as in Proposition 
\ref{prop:23082631}, and let $W_v = W_{\lambda_v}$ (resp. $W'_v = W'_{\tlambda_v}$) be the representations of  $U(r_v) \times U(s_v)$ (resp. $U(r_v-1) \times U(s_v)$ as in (1) of Conjecture \ref{conj:250629}.   Let $W = \otimes_v W_{\tlambda_v}$, $W' = \otimes_v W_{\tlambda_v}^{\prime}$.

We consider the set
$$S(\pi) = \{\pi_f~|~\pi\otimes \pi_f \subset \CA_0(H)\} $$
and define
$$H(\sigma,\pi) = 
\bigoplus_{\pi_f \in S(\pi)} H^q(\fP^-,K; \sigma^\vee\otimes W)\otimes \pi_f.$$
The map
$$\sigma^\vee \isoarrow P_\ttl(\pi^\vee\otimes F) \hookrightarrow \pi^\vee\otimes F$$
determines an isomorphism
\begin{equation}
\label{exo}
H(\sigma,\pi) \isoarrow H^q_\ttl(S(H,Y_V), [F]\otimes [W]),
\end{equation}
by the calculation of Proposition \ref{redu}.

\begin{theorem}\label{mainthm}  With notation and assumptions as in \S \ref{conj_PK}, 
%\ref{regularcorr}, 
let $W = W_\tlambda$,
$W^{\flat,\prime} = W^\flat_{\lambda,\tlambda'}$, $q = q(\delta)$.  We assume Conjecture \ref{conj:250629}.  Then there is a commutative diagram 

\begin{equation}\label{Fdiag:cohSh}
\xymatrix{
   H^{q}_\ttl(S(H,Y_V), [F]\otimes [W])
     \ar[r]^{T \circ \operatorname{pr}}
     &
   H^q_!(S(H^\flat,Y_{V'}\times Y_{V_1}),[W^{\flat,\prime}])
     \ar[d]^{=}
\\
    H(\sigma,\pi)
     \ar@{^{(}-_>}[u]  \ar[r]
      &H^q_!(S(H^\flat,Y_{V'}\times Y_{V_1}),[W^{\flat,\prime}])
}
\end{equation}
The bottom horizontal arrow is induced by the natural restriction, as in Theorem \ref{elemSh}, and the vertical arrows are those induced from (4) of Conjecture \ref{conj:250629}.

\end{theorem}

In particular, it follows from Theorem \ref{thm:non-van} that:

\begin{cor}\label{mainthma}

Suppose, for each $v$, $F_v$ is either the standard representation or its dual, and $F'_v \subset F_v$ is the subspace fixed by $H'(\mathcal{K}_v)$.  Then there is a commutative diagram \eqref{Fdiag:cohSh}, where the top arrow is induced by the natural restriction of coefficients. 
\end{cor}

\section{Period Invariants and Applications}\label{specialvalues}

In this section we simplify notation and assume $\CKp = \QQ$, so that $\CK$ is an imaginary quadratic field with a chosen complex embedding.   All the results of this section have the natural generalizations to CM fields.

\subsection{Periods for \texorpdfstring{$(\mathfrak{P}^-,K)$}{(P-,K)}-Cohomology}

Notation is as in \S \ref{avbH}.  Let $\pi \subset \CA_0(H)$  contribute non-trivially to
$H^*(\fP^-,K_h;\CA_0(H)\otimes F\otimes W)$, which we identify as in Remark \ref{interior} with the space $H^*_!(S(H,Y_V),[F]\otimes [W])$, which has a natural rational structure $\pi_\CK$ over a number field -- in this case over $\CK$.  In other words, $\pi_\CK$ is an admissible $(\fg,K_h)\times G(\af)$-module over $\CK$ with a natural $\CK$-rational embedding $\pi_\CK \hookrightarrow \pi$ such that the induced homomorphism
$$\pi_\CK\otimes_\CK \CC \rightarrow \pi$$
is an isomorphism.  

It is essentially proved that $\pi$ has multiplicity one in $\CA_0(H)$ \cite{AIGKMS}, and we will assume this to be the case.  Then it follows as in \cite[Proposition 3.17]{H13} the representation $\pi$ has a rational structure over a number field, compatible with the canonical rational structure on  $H^*_!(S(H,Y_V),[F]\otimes [W])$.  This rational structure depends on rational structures on the Lie algebra of $H$ and our choice of CM point $h_V \in Y_V$ as well as  on the fiber at $y$ of $[F]\otimes [W]$.  We assume in what follows that these choices have been made consistently.  The following lemma guarantees that the rational structure is intrinsic to $\pi$ and does not depend on its realization in a specific coherent cohomology space.

For simplicity, we only consider rational structures over $\Qbar$, and denote by $\pi_\Qbar \subset \pi$ the subspace of rational vectors.
Then as in \cite[(3.13)]{H13} we can define an element $Q(\pi) \in \CC^\times$, well-defined up to algebraic multiples, such that, for any $v, w \in \pi_\Qbar$ the quotient
\begin{equation}\label{Qpi}
Q(\pi)^{-1}\langle v, w \rangle \in \Qbar,
\end{equation}

where $\langle \bullet,\bullet \rangle$ is the $L_2$ norm for Tamagawa measure on $H(\bA)$.  To justify the use of this invariant in the next section, we need to show that it doesn't depend on the cohomological realization, in the following sense:

\begin{lemma}\label{Frational} Let  $[F]\otimes [W]$ be an automorphic vector bundle on $S(H,Y_V)$, attached to an irreducible representation $F$ of $G$ and an irreducible representation $W$ of $K_h$.  Let $\pi_\sigma$ be a discrete series representation, let $\{\pi_{\sigma_\delta}\}$ denote the discrete series $L$-packet containing $\pi_\sigma$, as $\delta$ runs over
the set of decorations.   Suppose the irreducible representation $W_\sigma$ of $P_{h_V}$ occurs with multiplicity one in $F \otimes W$.  Let $\pi \subset \CA_0(H)$ be a cuspidal automorphic representation with $\pi_\infty \isoarrow \pi_\sigma^\vee$.  Suppose every irreducible constituent of $\pi_\sigma^\vee \otimes F$ is in the discrete series, and let $\Delta(F\otimes W)$ denote the set of $\delta$ such that $W_{\sigma_\delta}$ occurs in the Jordan-H\"older decomposition of $F\otimes W$; let $m(\sigma,\delta)$ denote its multiplicity.  Then
\begin{enumerate}  
\item  There is a canonical isomorphism of $[\pi_f]$-isotypic components
$$H^*_!(S(H,Y_V),[F]\otimes [W])[\pi_f]  \isoarrow \bigoplus_{\delta \in \Delta(F \otimes W)} m(\sigma,\delta)H^{q(\sigma_\delta)}_!(S(H,Y_V),[W_{\sigma_\delta}])[\pi_f]$$
and both spaces are canonically isomorphic to
$$\bigoplus_{\delta \in \Delta(F \otimes W)} m(\sigma,\delta) H^{q(\sigma_\delta)}(\fP^-,K_h;\pi_\sigma^\vee \otimes W_{\sigma_\delta}) \otimes \pi_f$$
\item  The right-hand side of the last diagram in Part (1) is isomorphic to $\pi_f$, compatibly with the rational structure on $\pi$.
\item  The isomorphism in Part (1) preserves the rational structure on $\pi$.  In particular, $Q(\pi)$ depends only on $\pi$ and not on the specific realization in  coherent cohomology.
\end{enumerate}
\end{lemma}
\begin{proof}  Part (1) follows from
Proposition~\ref{prop:PK_spect_seq}.
Indeed, the stable base change of $\pi$ is an automorphic representation $\Pi$ of $GL(n)_\CK$, and by strong multiplicity one it is determined by $\pi_f$.  It follows from the properties of stable base change, as in \cite{HLa} (justified by \cite{AIGKMS}), that $\pi_f$ determines the $L$-packet of $\pi_\infty$.  Thus  if $\pi'_\infty\otimes \pi_f \subset \CA_0(H)$, $\pi'_\infty$ is necessarily one of the $\pi_{\sigma_\delta}^\vee$, and Part (1) is a special case of Proposition~\ref{prop:PK_spect_seq}.
%Corollary \ref{JH}{ii}.

Part (2) follows from Proposition~\ref{prop:dbc}; the statement about rationality is a consequence of the constructions in \cite{H13}.

We write $W_\sigma$ instead of $W_{\sigma_\delta}$.  To prove Part (3), let $\mathcal{W} \subset F \otimes W$ be the smallest $P_{h_V}$-invariant subspace containing $W_\sigma$ as a Jordan-H\"older constituent, and let
$[\mathcal{W}]$ denote the corresponding automorphic vector bundle on $S(H,Y_V)$.  Then the homomorphism
$$H^q_!(S(H,Y_V), [\mathcal{W}])[\pi_f] \rightarrow H^*_!(S(H,Y_V),[F]\otimes [W])[\pi_f]$$
is rational over $\Qbar$ (and even over the field of definition of $\pi_f$), because the rationality of maps of coefficients implies the rationality of maps of corresponding automorphic vector bundles.   Our choice of $\mathcal{W}$ implies that it maps surjectively to $W_\sigma$, and our hypothesis implies that $W_\sigma$ does not appear in the kernel of this map, nor in the cokernel of the inclusion of $\mathcal{W}$ in $F\otimes W$.  It follows that there is a surjective map
$$\begin{aligned}
H^q_!(S(H,Y_V), [\mathcal{W}])[\pi_f] \ra &H^q_!(S(H,Y_V), [W_\sigma])[\pi_f] \\ \isoarrow &H^q(\fP^-,K_h;\pi_\sigma^\vee \otimes W_\sigma)\otimes \pi_f
\end{aligned}$$
and this is again rational over $\Qbar$.

Now the construction of the rational structure on $\pi$ is determined by the rationality recalled in Part (ii) and by a certain choice of rational structure on the enveloping algebra $U(Lie(H))$, which is called the de Rham rational structure in \cite{H13}.  The claim then follows from the claim that the de Rham rational structure of $U(Lie(H))$ preserves the rational structure of the fiber at $h_V$ of the automorphic vector bundles $[F]$ and $[W]$ as well as all the automorphic subbundles of $[F]\otimes [W]$, determined by the canonical models of automorphic vector bundles.  But the de Rham rational structure on $U(Lie(H))$ is defined by viewing $U(Lie(H))$ as the fiber at $h_V$ of an automorphic vector bundle, so the claim is a consequence of the functoriality of the construction of automorphic vector bundles and its compatibility with tensor operations.
\end{proof}

\subsection{Motives, Deligne's Periods, and Critical Values of \texorpdfstring{$L$}{L}-functions}

\subsubsection{Motives as Realizations}

For Grothendieck, a (pure) motive over the field $\mathbb{Q}$ is an object of an abelian category, yet to be defined, obtained from smooth projective algebraic varieties over $\QQ$ and algebraic correspondences, to which one can assign cohomology groups of different types and comparison maps between them.  For the purposes of this paper we only care about the cohomology groups, since this is all the theory of automorphic forms is equipped to study, and specifically the cohomology of Shimura varieties with coefficients in local systems.  
More precisely, a (pure) motive $M$ over $\mathbb{Q}$, with coefficients in a field $E$, of rank $n$ and weight $w$, $n \in \mathbb{N}, w \in \mathbb{Z}$, consists of the following data:

\begin{enumerate}
\item An $E$-vector space $M_B$ (the \emph{Betti realization}) of dimension $n$, equipped with an $E$-linear involution
\[
F_\infty: M_B \to M_B;
\]

\item An $E$-vector space $M_{dR}$ (the \emph{de Rham realization}) of dimension $n$, equipped with a decreasing Hodge filtration
\[
F^iM_{dR} \subset F^{i-1}M_{dR}, \quad i \in \mathbb{Z},
\]
with
\[
F^iM_{dR} = 0 \quad \text{for } i \gg 0; \qquad 
F^iM_{dR} = M_{dR} \quad \text{for } i \ll 0;
\]

\item An $E \otimes \mathbb{C}$-linear \emph{comparison isomorphism}
\[
I_\infty: M_B \otimes \mathbb{C} \to M_{dR} \otimes \mathbb{C};
\]

\item For each non-archimedean completion $E_\lambda$ of $E$, a continuous representation
\[
\rho_{\lambda,M}: \operatorname{Gal}(\overline{\mathbb{Q}}/\mathbb{Q}) \to \operatorname{Aut}(M_\lambda),
\quad \text{where } M_\lambda := M_B \otimes_E E_\lambda,
\]
such that the $\rho_{\lambda,M}$ form a compatible family of $\lambda$-adic representations of weight $w$. This means that if $p$ is a prime unramified for $\rho_{\lambda,M}$, then the eigenvalues of $\rho_{\lambda,M}(\mathrm{Frob}_p)$ (geometric Frobenius) are algebraic numbers $\alpha$, independent of $\lambda$, such that $|\iota(\alpha)| = p^{w/2}$ for every embedding $\iota: \overline{\mathbb{Q}} \hookrightarrow \mathbb{C}$.
\end{enumerate}

\begin{remark}
\leavevmode
\begin{enumerate}
\item The anti-linear involution $F_\infty \otimes c$ on $M_B \otimes \mathbb{C}$, where $c$ is complex conjugation, induces via $I_\infty$ a Hodge structure on $M_{dR}\otimes \mathbb{C}$:
\[
M_{dR}\otimes \mathbb{C} \cong \bigoplus_{p+q = w} M^{p,q},
\]
such that
\[
(F_\infty \otimes c)(M^{p,q}) = M^{q,p}, 
\qquad 
F^iM_{dR}\otimes \mathbb{C} = \bigoplus_{p \geq i} M^{p,q}.
\]

\item For each prime $p$, the compatibility of the $\rho_{\lambda,M}$ allows us to define an Euler factor $L_\sigma(s,M)$ for each embedding $\sigma: E \hookrightarrow \mathbb{C}$.  
If $\rho_{\lambda,M}$ is unramified at $p$, then the characteristic polynomial
\[
P_p(M,T) := \det(1 - \rho_{\lambda,M}(\mathrm{Frob}_p)T)
\]
has coefficients in $E$, and we set
\[
L_{p,\sigma}(s,M) = \sigma(P_p)(M, p^{-s})^{-1}.
\]
In the ramified case, one restricts $\mathrm{Frob}_p$ to the invariants $(M_\lambda)^{I_p}$ under inertia $I_p$, and defines
\[
P_p(M,T) = \det(1 - (\mathrm{Frob}_p|_{(M_\lambda)^{I_p}})T), 
\qquad
L_{p,\sigma}(s,M) = \sigma(P_p)(M,p^{-s})^{-1}.
\]
Finally, the Euler product is defined by
\[
L_\sigma(s,M) = \prod_p L_{p,\sigma}(s,M),
\qquad
L(s,M) = \big[s \mapsto (L_\sigma(s,M), \ \sigma: E \hookrightarrow \mathbb{C}) \in E \otimes \mathbb{C}\big].
\]
If $M$ is pure of weight $w$, then the Euler product for $L(s,M)$ converges absolutely for $\Re(s) > \tfrac{w}{2} + 1$.  

\end{enumerate}
\end{remark}

Let $\pi$ be an automorphic representation $\pi$ of $H$ that contributes to the coherent cohomology, as in the previous sections, and let 
$\Pi$ denote its base change to $GL(n+1)_\CK$.  Then one defines a motive $M(\Pi)$, in the sense of the list of data above, of weight $n$, in the following way.  The $L$-function $L(s,M(\Pi))$ is just the standard automorphic $L$-function of $\Pi$, usually shifted so that the center of symmetry is the point $s = \frac{n+1}{2}$; in the discussion below, however, we will use the standard (unitary) normalization, so that the center of symmetry is the point $s = \frac{1}{2}$.  The other cohomology groups are realized on a Shimura variety $H_0$ attached to a unitary group of signature $(1,n)$.  In most cases $\Pi$ descends to a collection ($L$-packet) of automorphic representations of $H_0$, which contribute classes $\omega_i$ in the coherent cohomology of the Shimura variety attached to $H_0$, which are represented by automorphic forms $f_i$.  This Shimura variety has a canonical rational structure (canonical model) over the field $\CK$, as does its coherent cohomology.   The automorphic forms $f_i$ are eigenfunctions for Hecke operators, and the coefficient field $E$ of $M(\Pi)$ in the above definition is generated by its eigenvalues.  All of the structures above, and in the section \ref{DC} below can then be computed directly in terms of the $f_i$ and their Hecke eigenvalues.  See \cite[\S 3]{GHL} for more details.

\subsubsection{Critical values of $L$-functions}  

Let $M$ be a pure motive of weight $w$. It is conjectured that the function $L(s,M)$ admits a \emph{meromorphic continuation} to $\mathbb{C}$, with poles determined by the invariants of the action of $\Gamma_\mathbb{Q} = \operatorname{Gal}(\overline{\mathbb{Q}}/\mathbb{Q})$ on the $M_\lambda$, and that it satisfies a \emph{functional equation} of the form
\[
L_\infty(s,M)L(s,M) = \varepsilon(s,M)\, L_\infty(1-s, M^\vee) \, L(1-s, M^\vee),
\]
where $L_\infty(s,M)$ is a certain product of archimedean Euler factors (shifts of $\Gamma$-functions multiplied by exponentials), determined by the Hodge structure of $(M_B,M_{dR})$; $\varepsilon(s,M)$ is an entire function without zeros; and $M^\vee$ is the dual motive (obtained from $M$ by dualizing all realizations).

\begin{definition}[Critical values]
An integer $m \in \mathbb{Z}$ is \emph{critical} for $M$, and $L(m,M)$ is a \emph{critical value} of $L(s,M)$, if neither $L_\infty(s,M)$ nor $L_\infty(1-s,M^\vee)$ has a pole at $s = m$.
\end{definition}

For example, when $M = \mathbb{Q}(0)$ is the trivial motive, $L(s,M) = \zeta(s)$ is the Riemann zeta function. The critical integers for $\zeta(s)$ are the positive even integers and the negative odd integers. The formulas for these values of $\zeta(s)$ have been known since the 18th century. The above definition is equivalent to Deligne's \cite[Definition 1.3]{De}. At the time of Deligne's article, Shimura had already expressed the critical values of several classes of $L$-functions, attached to automorphic forms, as algebraic multiples of automorphic periods. Deligne's paper was the first to propose a general formula for critical values, up to algebraic factors.  

\subsubsection{Deligne's Conjecture}\label{DC}

Assume $w \geq 0$; if $w$ is even, always assume that the eigenvalue of $F_\infty$ on $M^{w/2,w/2}$ is $-1$. Let $M_B^\pm \subset M_B$ denote the eigenspaces of $F_\infty$ with eigenvalues $\pm 1$. Then the Hodge filtration always has a step $F^+M_{dR}$ (resp. $F^-M_{dR}$) of dimension $d^+(M) = \dim M_B^+$ (resp. $d^-(M) = \dim M_B^-$). Define
\[
M_{dR}^\pm = M_{dR}/F^{\mp}M_{dR},
\]
and the maps
\[
I^{\pm}: M_B^\pm \otimes \mathbb{C} \hookrightarrow M_B \otimes \mathbb{C} 
\xrightarrow{I_\infty} M_{dR} \otimes \mathbb{C} 
\twoheadrightarrow M_{dR}^\pm \otimes \mathbb{C}.
\]
Hodge symmetry implies that $I^\pm$ is an isomorphism of $E \otimes \mathbb{C}$-modules, and therefore its determinant $c^\pm(M)$, computed in any $E$-rational bases of $M_B^\pm$ and $M_{dR}^\pm$, is an element of $(E \otimes \mathbb{C})^\times$, well defined up to multiplication by $E^\times$.

On the other hand, we saw that $L(s,M)$ can be written as a function with values in $E \otimes \mathbb{C}$. If $\alpha, \beta \in E \otimes \mathbb{C}$, with $\beta \in (E \otimes \mathbb{C})^\times$, we write $\alpha \sim \beta$ if $\alpha/\beta \in E \subset E \otimes \mathbb{C}$.  

\begin{conj}[Deligne]\label{Deligneconj}
Let $M$ be a motive over $\mathbb{Q}$ with coefficients in $E$, and let $m \in \mathbb{Z}$ be a critical value of $L(s,M)$. Then
\[
L(m,M) \sim (2\pi i)^{d^{(-1)^m}(M)\cdot m}\, c^{(-1)^m}(M).
\]
\end{conj}

When $M = M(\Pi)$ and $M' = M(\Pi')$ with $\Pi$ (resp. $\Pi'$) the base change of an automorphic representation $\pi$ of $H$ (resp. an automorphic representation $\pi'$ of $H'$) that contributes to coherent cohomology, the left-hand side $L(s,M\otimes M')$ of Conjecture \ref{Deligneconj} is just the Rankin--Selberg $L$-function $L(s,\Pi \times \Pi')$, again up to a shift to take the weight into account.   Meanwhile, Deligne's periods $c^{\pm}(M\otimes M')$ that appear on the right hand side of the conjecture are computed in \cite{H13a, GHL} terms of the Petersson norms of holomorphic modular forms on different inner forms of $H$ and $H'$, all of which have the same base change $\Pi$ and $\Pi'$.  In other words, the calculations in \cite{H13a}, \cite{GHL}, and elsewhere express Deligne's {\it motivic period} in terms of concrete {\it automorphic periods}.  This provides the theoretical bridge between Deligne's conjecture and the automorphic theory.  The practical bridge is provided by the expression of special values of automorphic $L$-functions in terms of automorphic periods.  In the present case, the Ichino-Ikeda-N. Harris formula provides such an expression for the central value of the Rankin--Selberg $L$-function $L(s,\Pi\times \Pi')$, as we'll see in \S \ref{CentralRS} below.

\subsection{Automorphic periods}

Let $\pi$ and $\pi'$ be cuspidal automorphic representations of $H = U(V)$ and $H' = U(V')$, respectively.  Suppose there is a non-trivial pairing of $(Lie(H'),K')\times H'(\af)$-modules:
\begin{equation}\label{HomHp}
\operatorname{Hom}_{(Lie(H'),K')\times H'(\af)}(\pi \otimes \pi',\CC) \neq 0.
\end{equation}
In that case the GGP conjecture implies that the dimension of the left-hand side of \eqref{HomHp} is of dimension $1$.  We suppose both $\pi$ and $\pi'$ have discrete series components at the archimedean primes, so both contribute to coherent cohomology of their respective Shimura varieties.  Thus they have canonical rational structures, as in \cite{H13}, over some number field $E(\pi,\pi')$.  Then we define the {\it GGP period}
$$P(\pi,\pi') \in (E(\pi,\pi')\otimes \CC)^\times,$$
as in \cite[Definition 5.15]{H13} (where they were called Gross-Prasad periods), so that, for any $f \in \pi, f' \in \pi'$, both rational over $E(\pi,\pi')$, the $(Hom(E(\pi,\pi'),\CC)$-tuple
$$P_H(f,f') := (\int_{[G']}  j(f)(\iota(h))j(f')(h) dh)_{j \in Hom(E(\pi,\pi'),\CC)}$$
satisfies
\begin{equation}\label{PHff}
P(\pi,\pi')^{-1}P_H(f,f') \in E(\pi,\pi')^\times.
\end{equation}
The invariant $P(\pi,\pi')$ measures the failure of the $L_2$-pairing to respect the $E(\pi,\pi')$-rational structures of $\pi$ and $\pi'$.

Here we need to point out that we have inconspicuously changed our conventions:  the discussion in previous sections was based on non-trivial homomorphisms in $Hom_{H'(\CKp_v)}(\pi_v,\pi'_v)$, whereas \eqref{IINHformula} is based on $H'(\mathbf{A})$-invariant {\it pairings} between $\pi$ and $\pi'$.  This explains the appearance of the contragredient in the following discussion.

\subsubsection{Automorphic periods and cup products}\label{autocup}   %Let $F \simeq \RR$ be any completion of the totally real number field $\CKp$ and let
Let $\sigma_1$ and $\sigma_2$ be discrete series representations of $G = H_\infty := H(\CKp\otimes_\QQ \RR)$ and $G' = H'_\infty := H'(\CKp\otimes_\QQ \RR)$, respectively.   Say $(\sigma_1,\sigma_2^\vee)$ is a {\it GGP pair} if
\begin{equation}\label{GGPpair}
\operatorname{Hom}_{G'}(\sigma_1^\vee \otimes \sigma_2^\vee, \CC) \neq 0.
\end{equation}
In other words, $\operatorname{Hom}_{G'}(\sigma_1|_{G'},\sigma_2^\vee) \neq 0$, so that there is a symmetry-breaking operator from $\sigma_1$ to $\sigma_2^\vee$.  The dualization of the second term is more natural for applications to cohomological pairings.
We use the same terminology if $F \simeq \RR$ is any completion of $\CKp$ and we write $G = H(F)$, $G' = H'(F)$, and $\sigma_1$, $\sigma_2$ are any irreducible representations of $G$ and $G'$. 
In either case, the dimension of the space in \eqref{GGPpair} is $1$.

In what follows we let $K \subset G = H_\infty$ and $K' \subset K \cap G'$, where $G' = H'_\infty$, be maximal compact subgroups, as before.  We consider a GGP pair $(\sigma_1,\sigma_2^\vee)$,  an irreducible finite-dimensional representation $F$ of $G$, with $F' = F^{G'}$ one-dimensional,  an irreducible representation $W$ of $K$, and an irreducible $K'$-subspace $W' \subset W^\vee$.  We also choose cuspidal automorphic representations $\pi$ and $\pi'$ of $H$ and $H'$ with
archimedean components $\sigma_1^\vee$ and $\sigma_2^\vee$, respectively, and isomorphisms
$$\pi \isoarrow \sigma_1^\vee \otimes \pi_f; ~~ \pi' \isoarrow \sigma_2^\vee \otimes \pi'_f.$$  
We assume
\begin{enumerate}
\item  $\dim H^a(\fP^-,K; \sigma_1^\vee\otimes F \otimes W) = 1$ and $H^{a'}(\fP^-,K; \sigma_1^\vee\otimes F \otimes W) = 0$ for $a' \neq a$;
\item  $\dim H^b(\fP^{\prime,-},K'; \sigma_2^\vee\otimes W^{\prime}) = 1$ and $H^{b'}(\fP^{\prime,-},K'; \sigma_2^\vee \otimes W^{\prime}) = 0$ for $b' \neq b$;
\item  There is a pairing  
\begin{equation}\label{artificial}
H^a(\fP^-,K; \sigma_1^\vee\otimes F \otimes W)\otimes \pi_f \otimes H^b(\fP^{\prime,-},K'; \sigma_2^\vee\otimes W^{\prime})\otimes \pi'_f \ra \CC
\end{equation}
induced by the $L_2$ pairing of  $\pi$ and $\pi'$,  any choice of  $\Qbar$-rational $G'$-invariant pairings 
$$F \otimes F^{\prime,\vee} \ra \Qbar, ~~~ W \otimes W' \ra \Qbar;$$
and any choice of $\Qbar$-rational pairing
$$\wedge^a (\fp^-)^\vee \otimes \wedge^b (\fp^{\prime,-})^\vee \ra \Qbar.$$
\end{enumerate}
When $(\sigma_1,\sigma_2^\vee)$ is a coherent pair, the hypotheses imply that $a+b = \dim Y_{V'}$ and the pairing in the last line is the projection on top $K'$-invariant differentials.  

\begin{lemma}\label{autocup2}  Let $\phi \in H^a(S(H,Y_V),[F]\otimes [W])$ and $\phi' \in H^b(S(H',Y_{V'}), [W'])$ be $\Qbar$-rational coherent cohomology classes.  Let $\phi \cup \phi' \in \CC$ denote the image of $\phi\otimes \phi'$ obtained by applying the pairing \eqref{artificial} to the inclusion \eqref{automorphicF}.  

Then $P(\pi,\pi')^{-1}\cdot \phi \cup \phi'  \in \Qbar$.  
\end{lemma}
\begin{proof}  We can write the ``cup product pairing'' $(\phi,\phi') \mapsto \phi \cup \phi'$ as a linear combination with $\Qbar$ coefficients of $L_2$-pairings of elements of $\pi$ and $\pi'$, so the claim is clear.
\end{proof}

In the following lemma we apply Lemma \ref{autocup2} to the situation where $(\sigma_1,\sigma_2^\vee)$ is a coherent pair.

\begin{lemma}\label{Ppipi1} Suppose the $L_2$-pairing coincides with cup product in coherent cohomology, up to a factor in $\bar{\QQ}$, as in the statement of
Lemma \ref{autocup2}.  This is the case when $(\pi_\infty,\pi_\infty^{\prime,\vee})$ forms an elementary pair.  Let $\xi_\pi$ and $\xi_{\pi'}$ denote the central characters of $\pi$ and $\pi'$, respectively, and define the Hecke character
$$\xi(\pi,\pi') = \xi(\pi)^{-1}\cdot \xi(\pi').$$
Viewing $\xi(\pi,\pi')$ as a CM Hecke character of $U(V_1)$, we may define the period invariant $p(\xi(\pi,\pi')) \in \CC^\times$ so that 
$$p(\xi(\pi,\pi'))^{-1}\cdot \xi(\pi,\pi')$$
defines a $\bar{\QQ}$ rational section of the automorphic vector bundle defined by the algebraic character $\xi(\lambda,\lambda')$ of $U(V_1)$.

Then we can take $P(\pi,\pi') = p(\xi(\pi,\pi'))$.

More generally, suppose $(\pi_\infty,\pi_\infty^{\prime,\vee})$ 
is a coherent pair that satisfies Conjecture  \ref{conj:250629}.  Then again we can take 
$P(\pi,\pi') = p(\xi(\pi,\pi'))$.
\end{lemma}
\begin{proof}  The fact that the $L_2$-pairing with respect to Tamagawa measure coincides with cup product in coherent cohomology follows from \cite[Proposition 3.8]{H90}, up to an algebraic factor.  Indeed, as normalized in \cite{H90}, the cup product pairing equals $(2\pi i)^{-n}$ times the integral against a measure that is the product of an algebraic differential form at archimedean places with a measure on the finite ad\`eles that restricts to a rational number on an open compact subgroup; here $n$ is the dimension of $Sh(H',Y_{V'})$.  But this measure is a rational multiple of $(2\pi)^n$ times the Tamagawa measure.

Moreover, the cup product pairing involves contraction of pairings of finite-dimensional representations; as explained in \cite[Proposition 2.6.12]{H97} and in the proof of Lemma \ref{autocup2}, these do not affect the rationality of the pairing.  The first claim then follows as in \cite[\S 4]{H14}.  The second claim then follows from Lemma \ref{Frational}.
\end{proof}

\subsubsection*{$Q$-periods in the holomorphic case}

Fix one of the complex embeddings $\imath$ in the chosen CM type $\Sigma$.   Let $(i,n-i)$ be the signature of the complex hermitian space $V\otimes_{\CK,\imath}\CC$, and suppose the signature is $(0,n)$ at the other elements of $\Sigma$.  Let $\pi$ be a cuspidal automorphic representation of $H$, and suppose $\pi_\imath$ is in the holomorphic discrete series.  Let $\Pi$ denote the (cohomological) automorphic representation of $GL(n)_\CK$ obtained by stable quadratic base change from $\pi$.  Following the notation of \cite{GHL} we write
\begin{equation}\label{PipI}
P^{(i)}(\Pi,\imath) = Q(\pi).
\end{equation}
Implicit in the notation is the assumption that $Q(\pi_1) = Q(\pi_2)$ if $\pi_1$ and $\pi_2$ are both holomorphic discrete series representations of the chosen $U(V)$ with the same base change to $GL(n)$.  This is known at least if the infinitesimal character of $\Pi$ is sufficiently regular.  We use this notation when we rewrite the special value formula in Corollary \ref{rewrite}.

\subsection{Central Values of Rankin--Selberg \texorpdfstring{$L$}{L}-functions}\label{CentralRS}

Let $\pi$ and $\pi'$ be cuspidal automorphic representations of $H = U(V)$ and $H' = U(V')$, respectively.  
Let $\Pi = BC(\pi)$, $\Pi' = BC(\pi')$ be the 
representations of $GL(n)_\CK$, $GL(n-1)_\CK$ obtained from $\pi$ and $\pi'$ by quadratic base change.

The {\it Ichino-Ikeda-N. Harris conjecture}, now a theorem (due to W. Zhang, Chaudouard, Zydor, Y. Liu, H. Xue, X. Zhu, Beuzart-Plessis---see \cite{BPCZ} for references), expresses the central value
$$L(\frac{1}{2},\Pi \times \Pi')$$
in terms of period integrals of automorphic forms in $\pi \otimes \pi'$.
Fix a finite set $S$ of places of $\QQ$, including $\infty$ and all places at which $\pi$ and $\pi'$ are not spherical, and define
\begin{equation}\label{421} \mathcal{L}^S(\pi,\pi'): = \frac{L^S(\frac{1}{2},\Pi\otimes \Pi')}{L^S(1,\pi,Ad)L^S(1,\pi',Ad)} .\end{equation}
Here $L^S(s,\bullet)$ denotes the $L$-function of $\bullet$ with the factors in $S$ removed.

We say $f \in \pi$ is {\it factorizable} if, for one choice (and therefore any choice) of isomorphism
$$r:  \pi \isoarrow \bigotimes{}'_v \pi_v,$$
where $v$ runs over places of $\CKp$, the image of $f$ under $r$ is a pure tensor.
Then  there is an integer $\beta$, depending on the $L$-packets containing $\pi$ and $\pi'$, such that, for factorizable $f \in \pi$ and $f' \in \pi'$ 
\begin{equation}\label{IINHformula}
\mathcal{P}(f,f') = 2^{\beta}\Delta_{H} \mathcal{L}^S(\pi,\pi') \ \prod_{v \in S} I^*_v(f_v,f'_v).
\end{equation}
Here $\Delta_H$ is an elementary $L$-function, each $I^*_v(f_v,f'_v)$ is a normalized local integral of matrix coefficients, and
$$\mathcal{P}(f,f'):=\cfrac{|P_{H}(f,f')|^{2}}{\langle f,f\rangle \langle f',f'\rangle}.$$

Here
and $\langle f,f\rangle, \langle f',f'\rangle$ are Petersson norms.
The integrals are taken with respect to Tamagawa measure.

For each $v \in S$, the GGP conjecture picks out  inner forms $H^*_v, H^{\prime,*}_v$, and unique members $\pi^*_v, \pi^{\prime,*}_v$ of the
local $L$-packets defined by $\Pi_v, \Pi'_v$ such that
$$I^*_v(\bullet,\bullet) \not\equiv 0 \Rightarrow H_v = H^*_v, H'_v = H^{\prime,*}_v, \pi_v = \pi^*_v, \pi'_v = \pi^{\prime,*}_v.$$   
Moreover, Theorem \ref{nonv} (see below) guarantees that the local integrals do not vanish for appropriate choices of functions.

We assume we have picked $H' \subset H$, $\pi$, $\pi'$ so that the local pairings do not vanish identically.  In particular, we assume the root number of the $L$-function $L(s,\Pi \times \Pi')$ is $+1$, so there is no local sign obstruction.  

Now suppose $(\pi_\infty,\pi_\infty^{\prime,\vee})$ is a coherent pair of discrete series representations.  By Lemmas \ref{PHff} and \ref{autocup2}, we can rewrite the formula of \eqref{421}:  
%{\color{red}  Add details in what follows, in particular for the coefficient fields $E(\Pi)$, $E(\Pi')$.  Also we have to replace $H'$ by $H^\flat$.}

\begin{theorem}\label{IIperiods}  Suppose the pair $(\pi_\infty,\pi_\infty^{\prime,\vee})$ satisfies Conjecture \ref{conj:250629}; in particular, suppose the pair satisfies the hypotheses of Corollary \ref{mainthma}.  Then the test functions $f$ and $f'$ can be chosen to be coefficients of normalized coherent cohomology classes over $\Qbar$, respectively, as in Lemma \ref{autocup2}, so that 
$$\mathcal{P}(f,f') \sim \frac{p(\xi(\pi,\pi'))}{Q(\pi)Q(\pi')},$$
Here $Q(\pi)$ and $Q(\pi')$ are defined as in \eqref{Qpi} and we write $\alpha \sim \beta$ for $\alpha, \beta \in \CC$ if either both sides are zero or if the ratio of the two is an algebraic number.  (We can be more precise about the number field, as in \cite{H13}, but for simplicity we leave this for a future paper.)   Thus we have for following formula for the central critical value, up to scalars in $\Qbar$:
$$L(\frac{1}{2},\Pi \times \Pi') \sim \frac{p(\xi(\pi,\pi'))\cdot L^S(1,\pi,Ad)L^S(1,\pi',Ad)}{I^*_\infty\cdot\Delta_H \cdot Q(\pi)Q(\pi')}.$$
Here $I^*_\infty = \prod_{v \mid \infty} I^*_v(f_v,f'_v)$ for some local data $f_v, f'_v$ chosen so that $I^*_\infty \neq 0$.

\end{theorem} 

In the above we have used Theorem \ref{nonv}, \eqref{Ivstar}, and Corollary \ref{nonv_ar} from the next section in order to remove the non-archimedean local terms from \eqref{IINHformula} (up to non-vanishing algebraic factors) and move the archimedean local integral to the denominator.
It is straightforward to replace the relation up to $\Qbar$-factors with a result up to factors in $E(\Pi)\cdot E(\Pi')$.

\subsubsection{Non-vanishing of local integrals}

The following Theorem is part of the proof of the Ichino--Ikeda conjecture.  

\begin{theorem}\label{nonv}  Let $(\sigma_1,\sigma_2^\vee)$ be a GGP pair, and suppose $\sigma_1$ and $\sigma_2$ are both in the discrete series.  For $v_i \in \sigma_i$, $w_i \in \sigma_i^\vee$, and let 
$$f_{(v_1,w_1)}(h) = \langle \sigma_1(h)v_1,w_1 \rangle,~~
f'_{(v_2,w_2)}(h') = \langle \sigma_2^\vee(h)w_2,v_2 \rangle$$
denote the matrix coefficients of $\sigma_1$ and $\sigma_2^\vee$, respectively, with $h \in H(F), h' \in H'(F)$.

Then there exist vectors $v_i, w_j$ as above such that
$$I(f_{(v_1,w_1)},f'_{(v_2,w_2)}) :=  \int_{H'(F)} f_{(v_1,w_1)}(h')f'_{(v_2,w_2)}(h')dh' \neq 0.$$

\end{theorem}

This is a special case of Theorem 6.4.1 of \cite{SV}.  In the case of the Ichino--Ikeda conjecture for unitary groups this was proved in the non-archimedean case in 
\cite[Th\'eor\`eme 14.3.1]{BP16} and in general in \cite[Theorem 5]{BP20}. It also follows from point 2 of \cite[Theorem 2.3]{He} for discrete series of real unitary groups.

It is known that the cuspidal automorphic representations $\pi, \pi'$ that we are considering---for example, in Lemma \ref{Ppipi1} and in the following discussion---are tempered locally everywhere \cite{C}.  Thus we can apply Theorem \ref{nonv}
freely in what follows.  
For non-archimedean places this suffices for our applications.  We choose vectors $f, f'$ in \eqref{IINHformula} that are factorizable over the finite ad\`eles of H, and such that at all non-archimedean $v$ the local integrals 
\begin{equation}\label{Ivstar}
I_v^*(f_v,f'_v) \in \bar{\QQ}^*.
\end{equation}
The compatibility of algebraicity with the property of defining coherent cohomology classes rational over the field of definition is proved in \cite[Lemma 4.1.9]{H13a}.  In particular, the local factors at non-archimedean $v \in S$ in \eqref{IINHformula} can be ignored in what follows.

\subsubsection{Archimedean local integrals}\label{archloc}

We need a stronger result for real places $v$:  we need to know that the local factor $I_v^*(f_v,f'_v)$ doesn't vanish when $f_v \in \pi_v$ and $f'_v \in \pi^{\prime}$ are appropriate vectors in the support of the coherent cohomology classes attached to $\pi$ and $\pi'$ respectively.   We need to explain the relation between the pairing in \eqref{IINHformula} and the cup product in coherent cohomology, as in \cite{H14}.

Thus let $\phi \in H^q_\lambda(S(H,Y_V),[F]\otimes [W_\ttl])$ and
$\phi' \in [H^{q}_!(S(H^\flat,Y_{V'}\times Y_{V_1}), [W']]^\vee$ (note the Serre dualization).  We can realize such a $\phi'$ as the class defined by the automorphic representation $\pi^{\prime}$ such that $(\pi_\infty,\pi^{\prime,\vee}_\infty)$ is a GGP pair.  As in the proof of Lemma \ref{autocup2}, we can write these classes in terms of sums of factorizable vectors $f \in \pi$ and $f' \in \pi^{\prime}$, respectively, and replace the cup product $\phi \cup \phi'$ by the period $\mathcal{P}(f,f')$.

The following is standard:

\begin{lemma}\label{mat}  Let $F$ be a completion of $\CKp$ and let $(\sigma_1,\sigma_2^\vee)$ be a tempered GGP pair of representations of $H(F)$ and $H'(F)$; let
$\Lambda^\vee$ and $\Lambda$ denote generators of the $1$-dimensional spaces
$$\mathcal{V} := \operatorname{Hom}_{H'(F)}(\sigma_{1}^\vee \otimes \sigma^\vee_2,\CC)$$
and 
$$\mathcal{V}' :=\operatorname{Hom}_{H'(F)}(\sigma_1 \otimes \sigma_2,\CC).$$
Let  $f_i \in \sigma_i$, $f_i^\vee \in \sigma_i^\vee$, $i = 1, 2$, and assume that
$$\Lambda(f_1\otimes f_2) \neq 0; ~~~~~~ \Lambda^\vee(f_1^\vee \otimes f_2^\vee) \neq 0.$$
For $\psi_i$ a matrix coefficient of $\sigma_i$, $i = 1, 2$, define
$$I(\psi_1,\psi_2) = \int_{H'(F)}\psi_1(\iota(h'))\psi_2(h')dh'.$$
For $i = 1, 2$, let $\varphi_i$ be the matrix coefficient of $\sigma_i$
attached to the pair $(f_i,f_i^\vee)$.  Then
$I(\varphi_1,\varphi_2) \neq 0.
$
\end{lemma}

\begin{proof}  We know by Theorem \ref{nonv} that the integral $I$ doesn't vanish identically.  It thus defines a non-zero multiple of the generator $\Lambda\otimes \Lambda^\vee$ of the one-dimensional space $\mathcal{V}\otimes \mathcal{V}'$.  The Lemma then follows from the choice of $f_i$ and $f_i^\vee$.
\end{proof}

Let $v \in S$ be an archimedean place of $\CKp$.  Comparing Lemma \ref{mat} to the formula \eqref{IINHformula}, we have $f_1 = f_v \in \pi_v = \pi_{\lambda_v}^\vee$, $f_2 = f'_v \in \pi_v^{\prime} = \pi_{\lambda'_v}^\vee$, with $f_i^\vee$ the complex conjugate of $f_i$.  
In order to prove that the term $I_\infty^*$ in the denominator in Theorem \ref{IIperiods} doesn't vanish, we need to choose $f'_v$ in the minimal $K'_v$ type $\mu_{\lambda_v'}^\vee$ of $\pi_v^{\prime}$ and $f_v$ in a $K_v$-type $\mu$ of $\pi_{\lambda_v}^\vee$ that pairs non-trivially with $\mu_{\lambda'}^\vee$ but that also contributes to the cohomology class  
$$c_v(\pi) \in H^{q_v(\lambda_v)}(\fP^-,K_v;\pi_{\lambda_v}^\vee\otimes F \otimes W_\sigma)$$
that is the local component of the contribution of $\pi$ to
$H^q(S(H,Y_V), [F]\otimes [W])$.  

We thus write the cocycle representing $c_v(\pi)$ as a sum
$\sum_j \omega_j \otimes \mathbf{f}_j$ where $\mathbf{f}_j$ is the basis of the space $F$ as in Setting \ref{set:Fpi} -- now written in bold face to avoid confusion with the test functions $f, f'$ -- and $\omega_j$ is an element of the tensor product of $\pi_v^\vee \otimes (\wedge^{q_v}\fp_v^-)^\vee \otimes W_v$.  

Let $\Lambda$ be a generator of the one-dimensional space
$$\mathcal{V} = \operatorname{Hom}_{H'(F)}(\pi_v\otimes \pi_v^{\prime},\CC) $$
Let $\mathbf{f}_0$ be the generator of the $H'(F)$-fixed subspace $F'$.  Then the commutativity of \eqref{eqn:3_diagram}, which is valid under our hypotheses by Theorem \ref{thm:non-van}, implies that there is a $K$-type $\mu$ that contributes to the coefficient $\omega_0$ of $c_v(\pi)$ that pairs non-trivially, under $\Lambda$, with the minimal $K'$-type $\mu_{\lambda'}^\vee$ of $\pi_v^{\prime}$.  Choose $f_v \in \mu$ and $f'_v \in \mu_{\lambda'}^\vee$ such that
$\Lambda(f_v\otimes f'_v) \neq 0$.  Then it follows from Lemma \ref{mat} that

\begin{cor}\label{nonv_ar}  With the above choice of $f_v$ and $f'_v$ for archimedean places $v$, the term $I^*_\infty$ in the denominator of Theorem \ref{IIperiods} doesn't vanish.
\end{cor}

It follows from Lemma \ref{autocup2} that, with this choice of archimedean $f_v$, the relation
$$\mathcal{P}(f,f') \sim \frac{p(\xi(\pi,\pi'))}{Q(\pi)Q(\pi')}$$
of Theorem \ref{IIperiods} is valid.
This completes the justification of Theorem \ref{IIperiods}.

\subsubsection{Expressing the central value in terms of periods} 
In \cite[Theorem 5.15]{GHL}, the denominator $L^S(1,\pi,Ad)L^S(1,\pi',Ad)$ has been expressed, up to algebraic factors, as a product of Petersson norms of holomorphic automorphic forms on unitary groups of different signatures. Theorem 6.22 of \cite{GH} gives a different expression in terms of the Betti--Whittaker coefficients of the base change of $\pi$ and $\pi'$ to automorphic representations of $GL(n)$ and $GL(n-1)$, respectively.  Both results are only valid under certain conditions.  The expression in \cite{GHL} is more complete but assumes $\pi$ and $\pi'$ are sufficiently regular; the result of \cite{GH} does not require a regularity condition but includes an  archimedean factor that had not been determined at the time, and also assumed some hypotheses that guaranteed the validity of the Ichino--Ikeda--N. Harris conjecture of \cite{NH}.  These conditions are no longer necessary, since the conjecture of \cite{NH} has been proved in full, but \cite{GH} only considers the case where $\CK$ is imaginary quadratic.

We can therefore rewrite the right hand side of the expression in Theorem \ref{IIperiods} purely in terms of such periods.  We begin with the statement of the result of \cite{GHL}.\footnote{The preprint \cite{GHL} is in the process of being replaced by a new article, including the same theorem, that should appear shortly.}

\begin{theorem}\label{Asai thm}
Let $\pi$ be as in Theorem \ref{IIperiods}. We assume that the base change $\Pi$ of $\pi$ to $GL(n)_\CK$ is $5$-regular and admits descent to unitary groups of all possible signatures. 
Then, one has
\begin{equation}\label{Asai L-value 1}
L^{S}(1,\pi,Ad)\sim (2\pi i)^{dn(n+1)/2}\prod\limits_{\imath \in \Sigma}\prod\limits_{0\leq i\leq n}P^{(i)}(\Pi,\imath).
\end{equation}
Here $P^{(i)}(\Pi,\imath)$ is defined as in \eqref{PipI}.  

Interpreted as families, this relation is equivariant under the action of ${\rm Aut}(\C/F^{Gal})$.
\end{theorem}

\begin{remark}  The condition ``$5$-regular'' measures how far the infinitesimal character of $\Pi$ (or $\pi$) is from the walls.  Such a strong condition is unnecessary if we use the expression
$Q(\pi^{(i)})$, where $\pi^{(i)}$ is the descent of $\Pi$ to a unitary group of signature $(i,n-i)$ at all archimedean places $\imath$, instead of its factorization $\prod\limits_{\imath \in \Sigma}P^{(i)}(\Pi,\imath)$, which is only known under a regularity hypothesis.  

The existence of $\pi^{(i)}$ is guaranteed by the descent hypothesis.

In \cite{GHL} the formula is stated for $L(1,\Pi,As^{(-1)^n})$, but this Asai $L$-function for $GL(n)_\CK$ is the same as the automorphic $L$-function $L(1,\pi,Ad)$ for $H$.
\end{remark}

Next, we copy the expression of $L^S(1,\pi,Ad)$ from \cite[Theorem 6.22]{GH}.  Here we assume $\CK$ is imaginary quadratic:
\begin{equation}\label{pBW}
L(1,\pi,Ad) \sim a\left(\pi_{\infty}\right) \mathcal{G}\left(\varepsilon_{\mathcal{K},f}\right)^{g(n)}p_{BW}(\Pi). 
\end{equation}
Here $p_{BW}$ is the Betti--Whittaker period, denoted $p(\Pi)$ in \cite{GH}, $\mathcal{G}\left(\varepsilon_{\mathcal{K},f}\right)$ is the Gauss sum attached to $\CK$, $g(n)$ is a certain integer, and $a\left(\pi_{\infty}\right)$ is a complex constant coming from our normalization of the archimedean Whittaker functional for $\Pi$. 

The article \cite{JLS} computes the archimedean factor explicitly for a specific normalization.

Now we can rewrite the expression for $L(\frac{1}{2},\Pi \times \Pi')$ in terms of these period invariants
\begin{cor}\label{rewrite} 
In the situation of Theorem \ref{IIperiods}, we have the following expression for the central value:
$$L(\frac{1}{2},\Pi \times \Pi') 
\sim (2\pi i)^{dn^2}\frac{p(\xi(\pi,\pi'))\cdot\prod\limits_{\imath \in \Sigma}[\prod\limits_{0\leq i\leq n}P^{(i)}(\Pi,\imath)\cdot
\prod\limits_{0\leq j\leq n-1}P^{(j)}(\Pi',\imath)]}{I^*_\infty\cdot\Delta_H \cdot Q(\pi)Q(\pi')}.$$

At least when $\CK$ is imaginary quadratic, we also have

$$L(\frac{1}{2},\Pi \times \Pi') 
\sim \frac{p(\xi(\pi,\pi'))\cdot a\left(\pi_{\infty}\right)a(\pi'_\infty) \mathcal{G}\left(\varepsilon_{\mathcal{K},f}\right)^{g(n)+g(n-1}}{p_{BW}(\Pi)p_{BW}(\Pi')I^*_\infty\cdot\Delta_H \cdot Q(\pi)Q(\pi')}.$$

\end{cor}

\newpage

\part{Some  Examples and a Review of Notation}
\label{part:4}

This chapter offers a small selection of concrete examples and computations, intended to help clarify some of the definitions and notations that appear in Part~\ref{part:1} and \ref{part:translation} of the paper. A list of (some of) the first occurrences of these notations is provided in section 11.
We also explore the 8 different possibilities of defining 
coherent cohomology and give a guide to translate them into our conventions.

  \section{Interleaving Patterns and Symmetry Breaking}
  \label{sec:10}

The coherent signature discussed in this paper is a special case of GGP interlacing patterns (Definition~\ref{def:interlace}). Several equivalent definitions of the coherent signature are given in Proposition~\ref{prop:23022216}. In this section, we enumerate GGP interleaving patterns and specify when each of the various properties related to them is satisfied in the low-dimensional case.
In Section~\ref{subsec:Part4_Branching}, we illustrate with examples which parts of GGP interlacing patterns give rise to \emph{coherent pairs} $(\pi, \pi')$, or more strongly, \emph{elementary pairs} $(\pi, \pi')$.

\subsection{Counting Interleaving Patterns for 
 \texorpdfstring{$GL_{n+1} \downarrow GL_n$}{GL(n+1) -> GL(n)}}
~~\par

Let $(G,G')=(U(p,q), U(p-1,q))$.
The total number of GGP interleaving patterns is given by
\[
\binom{p+q}{q}\binom{p+q-1}{q}.
\]
Among them,
there are $\binom{p+q-1}{q}$ coherent signatures (Definition~\ref{def:coh_deco}).

Thus, if we consider all the pairs  $(G,G')=(U(p,q), U(p-1,q))$ with 
\[
p+q=n+1,
\]
the total numbers of coherent pairs and GGP interleaving patterns are given, respectively, by
\begin{equation}
\label{eqn:23020409_250812}
2^{n}=
\sum_{q=0}^{n} \binom{n}{q} 
\
\text{ and }
\
\binom{2n+1}{n}
=
\sum_{q=0}^{n} \binom{n+1}{q} \binom{n}{q}.  
\end{equation}
 
\subsection{List of GGP Interlacing Patterns for Small \texorpdfstring{$p$}{p} and \texorpdfstring{$q$}{q}}
\hspace*{2em}\par
Let $[\delta \delta']$ be a GGP interlacing pattern,
where $\delta \in \operatorname{Deco}(p,q)$
and $\delta' \in \operatorname{Deco}(p-1,q)$.
One of the defining conditions for $[\delta \delta']$ to form a \emph{coherent pair} is that Condition~($W$) (Definition~\ref{def:23022211}) is satisfied for the signature $[\delta \delta']$; equivalently, this holds if the reduction procedure (Definition~\ref{def:reduction}) terminates at the singleton $+$ (see Proposition~\ref{prop:23022216}). 

On the other hand, Condition~(q) (regarding the equality of degrees, $q(\delta) = q(\delta')$) and Condition~($\mu$) (concerning minimal $K$-types), see Definition~\ref{def:23022211}, are strictly weaker than the condition of being a coherent signature and do not necessarily imply it.

We illustrate these distinctions through explicit examples and examine how the reduction procedure described in Definition~\ref{def:reduction} applies in practice.

$\bullet$\enspace
$U(2,1) \downarrow U(1,1)$

Let us begin by enumerating the GGP interlacing patterns for $(p,q)=(1,1)$, and provide a table indicating which of the given conditions are satisfied for each pattern.

There are six GGP interlacing patterns $[\delta\delta']$ in total.
Among them, three satisfy $q(\delta)=q(\delta')$, 
and among these
 $\binom{2}{1}=2$ cases correspond to coherent pairs---that is, those satisfying Condition~$(W)$, or equivalently, for which the reduction procedure terminates at the singleton $+$.

\begin{table}[H]%[htb]
\begin{center}
\begin{tabular}{c|c|c|c|c|c|c|c|c|c}
%\hline
\text{Case}
&\text{interleaving}
&$\delta$
&$q(\delta)$
&$\delta'$
&$q(\delta')$
&$q(\delta)=q(\delta')$
&($\mu$)
&(W)
&\text{reduction}
\\
\hline
\text{I}
&$+ \op + - \oi$
&$++-$
&$2$
&$\op\oi$
&$1$
&
&\text{Yes}
&\text{}
&$+ - \oi$
\\[.3em]
\text{II}
&$+ \op \oi - +$
&$+-+$
&$1$
&$\op\oi$
&$1$
&\text{Yes}
&\text{Yes}
&\text{Yes}
&$+$
\\[.3em]
\text{III}
&$+-+ \op \oi$
&$+-+$
&$1$
&$\op\oi$
&$1$
&\text{Yes}
&\text{}
&\text{}
&$+-+\op\oi$
\\[.3em]
\text{IV}
&$+-\oi \op +$
&$+-+$
&$1$
&$\oi\op$
&$0$
&
&\text{Yes}
&\text{}
&$+-\oi$
\\[.3em]
\text{V}
&$\oi\op +-+$
&$+-+$
&$1$
&$\oi\op$
&$0$
&
&\text{}
&\text{}
&$\oi \op +-+$
\\[.3em]
\text{VI}
&$\oi -+ \op +$
&$-++$
&$0$
&$\oi\op$
&$0$
&\text{Yes}
&\text{Yes}
&\text{Yes}
&$+$
\\
\end{tabular}
\end{center}
\caption{Interleaving patterns for $(U(2,1), U(1,1))$}
\label{tab:U21}
\end{table}

The subtle interrelations among these conditions, as shown in the above table, are discussed in
Proposition~\ref{prop:23022216} and Remark~\ref{rem:equiv_coherent_sign}

$\bullet$\enspace 
$U(2,2) \downarrow U(1,2)$

In the case $(p,q) = (1,2)$, there are $18 = \binom{4}{2} \binom{3}{2}$ GGP interleaving patterns $[\delta \delta']$. 
Among these, we present in the table only the five cases where the degrees satisfy $q(\delta) = q(\delta')$.

Among these five, there are $\binom{3}{2} = 3$ coherent signatures---that is, those satisfying Condition~$(W)$.

\begin{table}[H]
\begin{center}
\begin{tabular}{c|c|c|c|c|c|c|c|c}
&\text{interleaving}
&$\delta$
&$q(\delta)$
&$\delta'$
&$q(\delta')$
&($\mu$)
&($W$)
&\text{reduction}
\\
\hline
&$+ \op \oi - \oi - +$
&$+--+$
&$2$
&$\op\oi\oi$
&$2$
&\text{Yes}
&\text{Yes}
&$+$
\\[.3em]
&$ \oi - + \op \oi - +$
&$-+-+$
&$1$
&$\oi\op\oi$
&$1$
&\text{Yes}
&\text{Yes}
&$+$
\\[.3em]
&$\oi\op\oi-+-+$
&$-+-+$
&$1$
&$\oi\op\oi$
&$1$
&\text{}
&\text{}
&$\oi\op\oi-+-+$
\\[.3em]
&$\oi-+-+\op \oi$
&$-+-+$
&$1$
&$\oi\op\oi$
&$1$
&\text{}
&\text{}
&$+-+\op \oi$
\\[.3em]
&$\oi-\oi-+\op +$
&$--++$
&$0$
&$\oi\oi\op$
&$0$
&\text{Yes}
&\text{Yes}
&$+$
\\
\end{tabular}
\end{center}
\end{table}

$\bullet$\enspace$U(3,2) \downarrow U(2,2)$

In the case $(p,q) = (2,2)$, there are $60 = \binom{5}{2} \binom{5}{2}$ GGP interleaving patterns $[\delta \delta']$. 
Among these, we present in the table only the seventeen cases where the degrees satisfy $q(\delta) = q(\delta')$.

Among these seventeen, there are $\binom{4}{2} = 6$ coherent signatures---that is, those satisfying Condition~$(W)$.

\begin{table}[H]
\begin{center}
\begin{tabular}{c|c|c|c|c|c|c|c}
%\text{Case}
&\text{interleaving}
&$G$
&$q$
&$G'$
&$q'$
&\text{reduction}
\\
\hline
&$\oi - \oi - + \op + \op +$
&$--+++$
&$0$
&$\oi\oi\op\op$
&$0$
&$+$
&$(\mu)$
\\[.3em]
&$ \oi \op \oi - + - + \op +$
&$-+-++$
&$1$
&$\oi\op\oi\op$
&$1$
&$\oi \op \oi - + - +$
&
\\[.3em]
&$\oi-+\op\oi-+\op+$
&$-+-++$
&$1$
&$\oi\op\oi\op$
&$1$
&$+$
&$(\mu)$
\\[.3em]
&$\oi-+-+\op \oi\op +$
&$-+-++$
&$1$
&$\oi\op\oi\op$
&$1$
&$+-+\op \oi$
&
\\[.3em]
&$\oi-+\op+\op\oi- +$
&$-++-+$
&$2$
&$\oi\op\op\oi$
&$2$
&$+$
&$(\mu)$
\\[.3em]
&$\oi-+\op+-+\op\oi$
&$-++-+$
&$2$
&$\oi\op\op\oi$
&$2$
&$+-+\op\oi$
&
\\[.3em]
&$+\op\oi-\oi-+\op+$
&$+--++$
&$2$
&$\op\oi\oi\op$
&$2$
&$+$
&$(\mu)$
\\[.3em]
&$+\op\oi\op\oi-+-+$
&$+-+-+$
&$3$
&$\op\oi\op\oi$
&$3$
&$\oi\op\oi-+-+$
&
\\[.3em]
&$+-\oi\op+\op\oi-+$
&$+-+-+$
&$3$
&$\oi\op\op\oi$
&$2$
&$+-\oi\op+\op\oi-+$
&$(\mu)$
\\[.3em]
&$+\op\oi-+\op\oi- +$
&$+-+-+$
&$3$
&$\op\oi\op\oi$
&$3$
&$+$
&$(\mu)$
\\[.3em]
%I
&$+-\oi-+\op+\op\oi$
&$+--++$
&$2$
&$\oi\op\op\oi$
&$2$
&$+-\oi-+\op+\op\oi$
&
\\[.3em]
&$+\op\oi-+-+\op\oi$
&$+-+-+$
&$3$
&$\op\oi\op\oi$
&$3$
&$+-+\op\oi$
&
\\[.3em]
&$+-+\op\oi\op\oi-+$
&$+-+-+$
&$3$
&$\op\oi\op\oi$
&$3$
&$+-+\op\oi\op\oi-+$
&
\\[.3em]
&$+\op+\op\oi-\oi-+$
&$++--+$
&$4$
&$\op\op\oi\oi$
&$4$
&$+$
&$(\mu)$
\\[.3em]
&$+-+\op\oi-+\op\oi$
&$+-+-+$
&$3$
&$\op\oi\op\oi$
&$3$
&$+-+\op\oi-+\op\oi$
&\\[.3em]
&$+-+-+\op\oi\op\oi$
&$+-+-+$
&$3$
&$\op\oi\op\oi$
&$3$
&$+-+-+\op\oi\op\oi$
&
\\[.3em]
%II
&$+-+\op+\op\oi-\oi$
&$+-++-$
&$4$
&$\op\op\oi\oi$
&$4$
&$+-+\op+\op\oi-\oi$
&%$-$
\\
\end{tabular}
\end{center}
\end{table}

\subsection{Branching Laws: Coherent Pairs and Elementary Pairs}
\label{subsec:Part4_Branching}
~~~\par

In the following, we present simple examples of branching laws and examine, for pairs $(\pi, \pi')$ such that $\operatorname{Hom}_{G'}(\pi|_{G'}, \pi') \neq 0$, which parts correspond to \emph{coherent pairs} (Definition~\ref{def:coherent_pair})
and which to \emph{elementary pairs} (Definition~\ref{def:elementary}).

For this purpose, we introduce an additional notation: for a discrete series representation $\pi$ of $G = U(p,q)$ with Harish-Chandra parameter $\lambda$, we recall that $\mu_\lambda$ denotes the minimal $K$-type parameter (\emph{Blattner parameter}), as defined in \eqref{eqn:min_K_type_delta}.

Our previous notation was $\pi = \pi_\lambda$, or $\pi_\lambda^\delta$ when emphasizing the decoration $\delta \equiv \delta(\lambda) \in \operatorname{Deco}(p,q)$. In what follows, we also use the notation $\pi_{[\mu_\lambda]}^\delta$, which is convenient for detecting elementary pairs.

\subsubsection{Branching Laws of Holomorphic Discrete Series Reprentations}

As noted in Convention~\ref{holo}, the decorations of holomorphic discrete series representations of $G = U(p,q)$ and $G' = U(p-1,q)$ are given respectively by
\[
\underbrace{-\cdots-}_{q}\underbrace{+\cdots+}_{p}
\quad \text{and} \quad
\underbrace{-\cdots-}_{q}\underbrace{+\cdots+}_{p-1}.
\]
We saw in Section~\ref{holorep} that the only GGP interleaving pattern $[\delta\delta']$ compatible with this decoration is
\begin{equation}
\oi ~ - ~ \oi ~\dotsb ~ - ~ \oi ~ - ~ + ~ \op ~ + ~ \op ~ +~ \dotsb ~\op ~ +
\end{equation}

$\bullet$ \enspace $U(2,1) \downarrow U(1,1)$

Holomorphic discrete series representations of $G=U(2,1)$ are parametrized by $\lambda \in \mathbb{Z}^3$ subject to
 $\lambda_3 > \lambda_1 > \lambda_2$, and the corresponding Blattener parameter is given by $\mu_\lambda=(\lambda_1-1, \lambda_2, \lambda_3+1)$.
Accordingly, we denote this discrete series representation $\pi$ by 
\[
   \pi_{\lambda_1, \lambda_2, \lambda_3}^{-++} 
  =
   \pi_{[\lambda_1-1, \lambda_2, \lambda_3+1]}^{-++}
\]

The irreducible decomposition of the restriction $U(2,1) \downarrow U(1) \times U(1,1)$
is discretely decomposable, and can be described in terms of Harish-Chandra parameters as
\[
\pi_{\lambda_1, \lambda_2, \lambda_3}^{-++}|_{G'}
= \bigoplus_{n=0}^{\infty} \bigoplus_{\ell=1}^{\lambda_1 - \lambda_2}
\mathbb{C}_{\lambda_1 - n - \ell} \boxtimes \pi_{\lambda_2 - \tfrac{1}{2} + \ell,\, \lambda_3 + n + \tfrac{1}{2}}^{\oi \op}.
\]
The same branching law is described in terms of Blattner parameters as follows:
\[
\pi_{[\lambda_1-1, \lambda_2, \lambda_3+1]}^{-++}|_{G'}
=\bigoplus_{n=0}^{\infty} \bigoplus_{\ell=1}^{\lambda_1-\lambda_2}
{\mathbb{C}}_{\lambda_1-n-\ell} \boxtimes 
\pi_{[\lambda_2-1 + \ell, \lambda_3+n+1]}^{\oi \op}, 
\]
Hence, for any $\pi'$ such that $\operatorname{Hom}_{G'}(\pi|_{G'}, \pi') \neq 0$, the pair $(\pi, \pi')$ is a coherent pair in this case. Moreover, it is an elementary pair if and only if $n = 0$---that is,
\[
 \pi'= \pi_{\lambda_2 - \tfrac{1}{2} + \ell,\, \lambda_3 + \tfrac{1}{2}}^{\oi \op}
\quad \text{for some } 1 \le \ell \le \lambda_1 - \lambda_2.
\]

\vskip 1pc
$\bullet$ $U(p,q) \downarrow U(p-1,q)$

Likewise, for general $p,q$,
 if $\pi$ is a holomorphic discrete series representation, then
any pair $(\pi, \pi')$ such that
$\operatorname{Hom}_{G'}(\pi|_{G'}, \pi') \neq 0$, forms a coherent pair.
We determine which of these pairs are elementary
in the branching law. 

Holomorphic discrete series representations of $G = U(p,q)$ are parametrized by the
Harish-Chandra parameter
\[
\lambda = (\lambda_1, \dots, \lambda_{p+q}) \in \left(\mathbb{Z} + \tfrac{p+q+1}{2} \right)^{p+q}
\]
subject to the strict inequalities
\[
\lambda_{p+1} > \dots > \lambda_{p+q} > \lambda_1 > \dots > \lambda_p,
\]
which correspond to the decoration $\delta=-\dotsb-+\dots+$.
The associated Blattner parameter is given by
\[
\mu_\lambda = \left(\lambda_i - \tfrac{1}{2}(p+q+1) + i \,;\, \lambda_{p+j} + \tfrac{1}{2}(p - q - 1) + j\right),
\]
where the indices $i$ and $j$ range over $1 \le i \le p$ and $1 \le j \le q$, respectively.
Here, we adopt the notation
$\mu = (\mu_i \,;\, \mu_j)$
to represent $\mu = (\mu_1, \dots, \mu_p, \mu_{p+1}, \dots, \mu_{p+q})$.

Then, the representation $\pi_\lambda$ decomposes discretely as a direct sum of holomorphic discrete series representations
of $G'=U(p-1,q)$:
\[
\pi_\lambda|_{G'} \simeq
\bigoplus_{\lambda'}
\pi'_{\lambda'}
\]
where each $\pi'_{\lambda'}$ has Harish-Chandra parameter
\[
\lambda'=(\lambda_{1}', \cdots, \lambda_{p-1}', \widehat{\lambda'_p};\lambda_{p+1}', \cdots, \lambda_{p+q}') 
\in \left(\mathbb{Z} + \tfrac{p+q}{2} \right)^{p+q}
\]
ranges over the set of parameters satisfying the interleaving condition \eqref{eqn:holo_branch}.

The holomorphic discrete series representation $\pi'_{\lambda'}$ of $G'=U(p-1,q)$
has Blattner parameter
\[
\mu'\equiv \mu_{\lambda'}'=(\mu_1', \cdots, \mu_{p-1}', \widehat{\mu'_p}; \mu_{p+1}', \cdots, \mu_{p+q}'),
\]
where
\[
\mu'_i=\lambda_i - \tfrac{1}{2}(p+q) + i 
\, (1 \le i \le p-1)
;\, \mu'_{p+j}=\lambda'_{p+j} + \tfrac{1}{2}(p - q - 2) + j \, (1 \le j \le q).
\]

The interleaving condition \eqref{eqn:holo_branch} for the Harish-Chandra parameters can equivalently be expressed in terms of the following inequalities involving the Blattner parameters $(\mu, \mu')
\in {\mathbb Z}^{p+q} \times \mathbb{Z}^{p+q-1}$:
\begin{equation}
\label{eqn:23012313}
\begin{cases}
&\mu_{p+1}'\ge \mu_{p+1} \ge \mu_{p+2}'\ge \mu_{p+2} \ge  \cdots 
 \ge \mu_{p+q}'\ge \mu_{p+q}, 
\\ 
&\mu_{p+q} \ge \mu_{1} + p+q-1, 
\\
&\mu_{1} \ge \mu_{1}'\ge \cdots \ge \mu_{p-1}'\ge \mu_{p}
\end{cases}
\end{equation}
In particular, $(\pi_\lambda, \pi'_{\lambda'})$ satisfies the condition for being an \emph{elementary pair} if and only if
$\mu_{p+j}'=\mu_{p+j}$ $(1 \le j \le q)$, or equivalently,
\[
 \lambda_{p+j}=\lambda'_{p+j}-\tfrac{1}{2}
 \textrm{ for all }
 1 \le j\le q
\]
in addition to the interleaving condition \eqref{eqn:holo_branch}.

\subsubsection{Symmetry Breaking of Non-holomorphic Discrete Series : $U(2,1) \downarrow U(1,1)$}
\label{subsubsec:SB_u21}

\medskip
The symmetry breaking of non-holomorphic discrete series representations first appears in the case 
$(G,G') =(U(2,1),U(1,1))$.
We examine this case.

A discrete series representation of $G=U(2,1)$
that is neither holomorphic nor anti-holomorphic 
is associated with the decoration $\delta=+-+$.
Let $\lambda=(\lambda_1, \lambda_2, \lambda_3) \in {\mathbb{Z}}^3$ satisfy the corresponding inequality
 $\lambda_1 > \lambda_3 > \lambda_2$. 
The associated non-holomorphic discrete
series representation has the following parameters: 
\begin{alignat*}{2}
&\text{The Harish-Chandra parameter} \qquad 
&& \tlambda=(\lambda_1, \lambda_2, \lambda_3), 
\\
&\text{The Blattner parameter}
&&
\mulmd=(\mu_1, \mu_2, \mu_3)=(\lambda_1, \lambda_2; \lambda_3).  
\end{alignat*}
We denote this representation by 
\[
\pi_{\lambda_1, \lambda_2, \lambda_3}^{+-+}
 =\pi_{[\mu_1, \mu_2;\mu_3]}^{+-+}.  
\]

\vskip 1pc

In what follows, we fix $\lambda = (\lambda_1, \lambda_2, \lambda_3) \in \mathbb{Z}^3$ satisfying $\lambda_1 > \lambda_3 > \lambda_2$, and consider the support of the restriction of the non-holomorphic discrete series $\pi_{\lambda_1, \lambda_2, \lambda_3}^{+-+}$ to $G' = U(1,1)$, where we vary the Harish-Chandra parameter $\lambda'=(\lambda'_1, \widehat{\lambda'_2}; \lambda'_3)$ of $G'$ with $\lambda_1', \lambda_3' \in \mathbb{Z} + \tfrac{1}{2}$.

\begin{alignat*}{2}
\operatorname{Supp}\big(\pi_{\lambda_1, \lambda_2, \lambda_3}^{+-+} \big|_{G'}\big)
&=
\{\pi_{\lambda_1', \lambda_3'}^{\op \oi} 
: \lambda_1 > \lambda_1' > \lambda_2' > \lambda_3 \; (>\lambda_2)\}
\qquad
&&\text{(Case II)}
\\
&\phantom{=}\cup
\{\pi_{\lambda_1', \lambda_3'}^{\op \oi} 
: (\lambda_1 > \lambda_3 >)\; \lambda_2 > \lambda_1' > \lambda_3'\}
\qquad
&&\text{(Case III)}
\\
&\phantom{=}\cup
\{\pi_{\lambda_1', \lambda_3'}^{\oi \op} 
: (\lambda_1 >)\; \lambda_3 > \lambda_3' > \lambda_1' > \lambda_2\}
\qquad
&&\text{(Case IV)}
\\
&\phantom{=}\cup
\{\pi_{\lambda_1', \lambda_3'}^{\oi \op} 
: \lambda_3' > \lambda_1' > \lambda_1 \; (>\lambda_3 > \lambda_2)\}
\qquad
&&\text{(Case V)}.
\end{alignat*}

The same formula can be described in terms of Blattner parameters. In this case,
we fix $\mu = (\mu_1, \mu_2, \mu_3) \in \mathbb{Z}^3$ satisfying $\mu_1 > \mu_3 > \mu_2$, and describe the support of the restriction of the non-holomorphic discrete series $\pi_{[\mu_1, \mu_2, \mu_3]}^{+-+}$ to $G' = U(1,1)$, where we vary Blattner parameters $(\mu_1', \mu_3') \in \mathbb{Z}^2$.

\begin{alignat*}{2}
\operatorname{Supp}(\pi_{[\mu_1, \mu_2; \mu_3]}^{+-+}|_{G'})
&=
\{\pi_{[\mu_1';\mu_3']}^{\op \oi}: 
    \mu_1 - 1 \ge \mu_1' - 1 > \mu_3' \ge \mu_3 \; (\ge \mu_2) \}
\qquad
&&\text{(Case II)}
\\
&\phantom{=}\cup\;\{\pi_{[\mu_1';\mu_3']}^{\op \oi}:
    (\mu_1 > \mu_3 >)\; \mu_2 \ge \mu_1' > \mu_3' + 1 \}
\qquad
&&\text{(Case III)}
\\
&\phantom{=}\cup\;\{\pi_{[\mu_1';\mu_3']}^{\oi \op}:
    (\mu_1 >)\; \mu_3 \ge \mu_3' > \mu_1' + 1 \ge \mu_2 + 1 \}
\qquad
&&\text{(Case IV)}
\\
&\phantom{=}\cup\;\{\pi_{[\mu_1',\mu_3']}^{\oi \op}:
    \mu_3' - 1 > \mu_1' \ge \mu_1 \; (>\mu_3 > \mu_2) \}
\qquad
&&\text{(Case V)}.
\end{alignat*}

In particular, Case~II corresponds to the situation where $(\pi, \pi')$ forms a coherent pair. Among these, the pair is elementary if and only if $\lambda_3 = \lambda_3' - \tfrac{1}{2}$.

The degrees satisfy $q(\delta) = q(\delta') = 1$ in both Case~II and Case~III; that is, when $(\delta, \delta') = (- + -,\, \op \oi)$.

\section{Guide to Notation}
\label{sec:11}
{ This section is a guide to the notation used  in the article.  The last part discusses different conventions for coherent cohomology of discrete series representations}

\subsection{Combinatorial Data for Discrete Series Representations}
We summarize some notation introduced in Section~\ref{sec:convent}
 of the main paper.  

Let $(G,K)=(U(p,q), U(p) \times U(q))$.
In our convention of \eqref{eqn:GK_embed},
 the holomorphic tangent space of the hermitian symmetric space $G/K$ at the origin
 is identified with the subalgebra 
 $\mathfrak{p}^+ \subset \mathfrak{g}$,
 and we set
 \[
   {\mathfrak{P}}\equiv {\mathfrak{P}}^+:={\mathfrak{p}}^+ + {\mathfrak{k}}, 
\quad
{\mathfrak{P}}^-:={\mathfrak{p}}^- + {\mathfrak{k}}.
\]
Using a standard basis $e_1, \dots, e_{p+q}$ of the dual of the fundamental Cartan subalgebra $\mathfrak t \subset \mathfrak{g}$, we write
\begin{align*}
\Phi^{c,+} \equiv\, \Delta^+({\mathfrak{k}})=&\{e_i-e_j:\text{$1 \le i<j \le p$ or $p+1 \le i < j \le p+q$}\},
\\
\Phi^{nc,+} \equiv\, \Delta({\mathfrak{p}}^+)
=&\{e_i-e_j:\text{$1 \le i \le p$ and $p+1 \le j \le p+q$}\},
\\
\Phi^+ \equiv\, \Delta^+({\mathfrak{g}})=&\Phi^{c,+} \cup \Phi^{nc,+}. 
\end{align*}
Accordingly, the half-sum of positive roots \( \rho_{\mathfrak{g}} \), as given in \eqref{eqn:rho_G}, is
\[
\rho_{\mathfrak g}
  =  \rho_c+\rho_n
   =\frac{1}{2}(p+q-1, \dots, 1-p-q) \in C_+.
\]

We denote by $F^K(\mu)$ the irreducible finite-dimensional representation of $K$ with highest weight $\mu$ with respect to $\Phi^{c,+}$.

The discrete series representations of $G=U(p,q)$ fall into $\binom{p+q}{q}$ distinct series. These series can be indexed combinatorially, either by the set of decorations (denoted by $\operatorname{Deco}(p,q)$) or by the set of Kostant representatives $W^\mathfrak{k} \subset W \simeq \mathcal S_{p+q}$
defined in \eqref{eqn:Kostant_W}. 
There is a natural bijection between these two indexing sets:
\[
 \operatorname{Deco}(p,q) \simeq  W^{\mathfrak{k}}, \quad \delta \leftrightarrow w. 
\]
The choice of positive root system described above serves as the reference point for this bijection
(Lemma~\ref{lem:decobij}).

A discrete series representation of $ G = U(p,q) $, corresponding to the Harish-Chandra parameter $ \lambda $, is denoted by $\pi_\lambda$.
The following invariants are associated with this representation:
\begin{alignat*}{5}
\delta&=\delta(\lambda) \in \operatorname{Deco}(p,q) &: 
&\quad \text{decoration determined by } \lambda 
&&\quad (\text{see \eqref{eqn:decobij}}),
\\
q(\delta) &= \# \Phi^{nc,+}(\lambda) 
&: 
&\quad \text{degree of } \pi 
&&\quad (\text{Section~\ref{subsec:q}}),
\\
\mu_{\lambda} &= \lambda + w \rho_{\mathfrak{g}} - 2\rho_c 
&: 
&\quad \text{minimal $K$-type parameter} 
&&\quad (\text{Definition~\ref{eqn:min_K_type_delta}}),
\\
W_{\lambda} &= F^K(\lambda - \rho_{\mathfrak{g}}) 
&: 
&\quad \text{coherent parameter} 
&&\quad (\text{Definition~\ref{def:W_lambda}}).
\end{alignat*}

Here, 
$
\Phi^{nc,+}(\lambda) \equiv \Phi^{nc,+}(\delta)
:= \{ \alpha \in \Delta(\mathfrak{p}^+) : \langle \alpha, \lambda \rangle > 0 \}
$ is determined solely by the decoration $\delta$, or equivalently, by the Kostant representative $w\in W^{\mathfrak{k}}$ (Lemma~\ref{lem:q}):
\[\Phi^{nc,+}(\lambda)=
\Delta(\mathfrak{p}^+) \cap w \Delta^+(\mathfrak{g})
= \{ e_i - e_j : \text{$+_i$ appears to the left of $-_j$ in } \delta \}.
\]
The representation is also denoted by $ \pi_{\lambda^+}^\delta $, where $\lambda^+ $ is the reference parameter (Definition~\ref{def:ref_deco}) defined by
$
\lambda^+ = w^{-1} \lambda \in C_+.
$
Although the decoration $\delta$ is determined by $\lambda$, we also write $\pi_\lambda^\delta$, when we wish to emphasize the decoration $\delta$ explicitly.
The parameter set of Harish-Chandra parameters for each decoration $\delta \in \operatorname{Deco}(p,q)$ is
$\mathbb{Z}_{\operatorname{int}}^{p+q}(\delta) := \left(\mathbb{Z} + \tfrac{p+q+1}{2}\right)^{p+q} \cap w C_+$ (see~\eqref{eqn:Z_int_delta}).

\begin{remark}
\label{rem:q_W_choice}
The degree $q(\delta)$ and the coherent parameter $W_\lambda$ depend only on the representation $\pi$, and hence we sometimes denote them by $q(\pi)$ and $W_\pi$, respectively.
However, both $q(\pi)$ and $W_\pi$ are sensitive to the choice of the set of noncompact roots $\Delta(\mathfrak{p}^+)$, or equivalently, on the 
complex structure of the hermitian symmetric space $G/K$.
This dependence is detailed in Table~\ref{tab:250503}.
\end{remark}

\bigskip 
We summarize below some useful formulas in connection with our convention.

Here is a relation between the Blattner parameter $\mu_\lambda$ of the minimal $K$-type and the coherent parameter $\lambda -\rho_{\mathfrak{g}}$ for $W_\lambda$, as stated in Lemma~\ref{lem:24042109}:
\[
\mu_\lambda - (\lambda - \rho_{\mathfrak{g}}) =
w \rho_{\mathfrak g} + \rho_{\mathfrak{g}} - 2\rho_c
=\sum_{\alpha \in \Phi^{nc,+}(\delta)} \alpha.
\]
The irreducible decomposition of the $K$-module
$\Lambda({\mathfrak{p}}^+)$ is given in \eqref{eqn:25050326}:
\[
\Lambda({\mathfrak{p}}^+)
=\bigoplus_{j=0}^{pq}
\Lambda^j({\mathfrak{p}}^+)
   \simeq
  \bigoplus_{j=0}^{pq} \ \bigoplus_{\substack{\delta \in \operatorname{Deco}(p,q)\\ q(\delta)=j}} 
  \  F^K \left(\sum_{\alpha \in  \Phi^{nc, +}(\delta)} \alpha\right).  
\]

\subsection{ Convention for \texorpdfstring{$({\mathfrak{P}}, K)$}{(P, K)}-Cohomologies}
~
\newline
There are several conventions for $(\fP, K)$-cohomology found in the literature. These conventions vary depending on the choice of complex structure on the hermitian symmetric space $G/K$ (in our case, the holomorphic tangent bundle is identified with $G\times_K {\fp^+}$),  the definition of a maximal parabolic subalgebra $\fP$ (whether one takes $\fP^+=\mathfrak{k}+\fp^+$ or $\fP^-=\mathfrak{k}+\fp^-$), and whether the coefficients are given by a representation $\pi$ or its contragredient $\pi^\vee$ as coefficients.
Thus, in theory, there are 
$2\times 2 \times 2=8$ conceivable conventions.

Needless to say, the results proved in Part~\ref{part:translation} hold under all eight conceivable conventions.
Although the convention for $(\mathfrak{P}, K)$-cohomology used in Part 2 is the same as that in \cite{H14} and differs from those in some other parts of the literature, we explain how the theorems in Part 2 can be translated precisely into each of the eight conceivable conventions.

For the reader’s convenience, we provide a correspondence among the different conventions.
First, we fix a complex structure on $G/K$, and define the nilpotent Lie algebra $\mathfrak{p}^+$ so that the holomorphic tangent bundle is identified with $G \times_K \mathfrak{p}^+$.
\begin{quote}
(In our paper, $\mathfrak{p}^+$ is defined via the root system, but since the final results are independent of coordinates, the same choice can equivalently be made from the complex structure.)
\end{quote}
With this choice, four conceivable conventions remain.

With the notation introduced in Section~\ref{sec:convent},  
we summarize below the conditions under which the $(\fP, K)$-cohomologies are non-vanishing.
In the table, 
the top-left entry reflects the convention adopted in Sections~\ref{sec:coherent} and~\ref{sec:Pk} of Part~\ref{part:translation} where $\mathfrak{P}=\mathfrak{P}^+$, whereas the bottom-right entry corresponds to the convention used in Part~\ref{global}.
\renewcommand{\arraystretch}{1.2} 
\begin{table}[hbtp]
\centering
\begin{tabular}{c|c|c}
& $\pi_{\lambda}^{\delta} \otimes W^{\vee}$
& $(\pi_{\lambda}^{\delta})^{\vee} \otimes W$
\\[.5em]
\hline
 $H^{j}({\mathfrak{P}}^+, K;\hphantom{M})$\quad
& $\begin{matrix} W=W_{\lambda}\\  j=q(\delta) \end{matrix}$
& $\begin{matrix} W=W_{\lambda}^-\\  j=pq- q(\delta) \end{matrix}$
\\[.5em]
\hline
 $H^{j}({\mathfrak{P}}^-, K;\hphantom{M})$\quad
& $\begin{matrix} W=W_{\lambda}^-\\  j=pq-q(\delta) \end{matrix}$\qquad
& $\begin{matrix} W=W_{\lambda} \\  j= q(\delta) \end{matrix}$
\\
\hline
  \end{tabular}
\caption{Non-vanishing $(\fP, K)$-cohomologies}%Table \ref{tab:250503}}
\label{tab:250503}
\end{table}
\renewcommand{\arraystretch}{1.0}  

\begin{quote}
Since we use the duality between the geometric realization of the discrete series representation $\pi$ via Dolbeault cohomology and the $(\mathfrak{P}, K)$-cohomology, we transfer the results of Part 2 for
\[
H^{j}(\mathfrak{P}^+, K; \pi \otimes W^\vee)
\]
to
\[
H^{j}(\mathfrak{P}^-, K; \pi^\vee \otimes W).
\]
The table above provides the exact parameters.
In particular, the degree of non-vanishing cohomology remains the same.
\end{quote}

In the left-lower block and the right-upper block of the above table, we use the following definition:
\[
   W_{\lambda}^-:=F^K(\lambda-\rho_c+\rho_n).
\]
instead of our previous definition 
$
   W_{\lambda}=F^K(\lambda-\rho_{\mathfrak g}) = F^K(\lambda-\rho_c-\rho_n).  
$

Let $\sigma_0$ denote the longest element in $W_K$. 
Let $\pi$ be a discrete series representation with Harish-Chandra parameter $\lambda$ and decoration $\delta$.
In passing from $\pi$ to its contragredient $\pi^\vee$, we observe that $\pi^\vee$ is associated with the Harish-Chandra parameter $\lambda^\vee$, and the decoration $\delta^\vee$, and the Blattner parameter $\mu_{\lambda^\vee}$ for the minimal $K$-type, which are defined as follows:
\begin{alignat*}{5}
\lambda^\vee&:= -\sigma_0 \lambda &:
&\quad \text{Harish-Chandra parameter of $\pi^\vee$} &&%
\\
\delta^\vee& \in \operatorname{Deco}(p,q) &: 
&\quad \text{the reverse of decoration of } \delta 
&& 
\\
\mu_{\lambda^\vee} &= -\sigma_0 \mu_\lambda 
&: 
&\quad \text{Blattner parameter for the minimal $K$-type.} 
&&\quad 
\end{alignat*}

On the other hand, as illustrated in the lower-left and lower-right blocks of Table~\ref{tab:250503},
the values of the degree $q(\pi^\vee)$ and the coherent parameter $W_{\pi^\vee}$ depend on the choice of parabolic subalgebras involved in the computation of the relative cohomology.

For example, the equivalent conditions for elementary pairs $(\pi, \pi')$ stated in Theorem~\ref{thm:elementary_pair_equiv} are transferred from the formulation in terms of $(\mathfrak{P}, K)$-cohomology in the upper-left block to $(\mathfrak{P}^-, K)$-cohomology in the lower-right block, as formulated below.
\begin{example}
[Dual of Theorem~\ref{thm:elementary_pair_equiv}]
\label{ex:dual}
Let $\pi = \pi_\lambda$ and $\pi' = \pilmdp$
be discrete series representations of $G$ and $G'$, respectively, with Harish-Chandra parameters $\lambda$
and $\lambda'$.
Let $\delta$ and $\delta'$ denote the decorations of $\pi$ and $\pi'$, respectively.
Let $W_\lambda$ and $W'_{\lambda'}$ denote the coherent parameters of $\pi$ and $\pi'$, respectively
(Definition~\ref{def:W_lambda}).
Assume that the signature $[\delta \delta']$ associated with $\lambda$ and $\lambda'$ is coherent
(see Proposition~\ref{prop:23022216} for six equivalent definitions).
Then, the following four conditions are equivalent:
\begin{itemize}
    \item[{\rm{(i)}}]
    $[\mu_{\lambda^\vee} : \mu'_{{\lambda'}^\vee}] \neq 0$.
    
    \item[{\rm{(ii)}}]
    In addition to (i), $\operatorname{Hom}_{G'}(\pi^\vee|_{G'}, {\pi'}^\vee) \neq \{0\}$.
    
    \item[{\rm{(iii)}}]
    In addition to (ii), $[W_\lambda : \Wlmdp] \neq 0$.
    
    \item[{\rm{(iv)}}]
    In addition to (iii), the restriction map
    $T\otimes \operatorname{pr}_{W} $
    induces an isomorphism
    $$
 H^\ell(\mathfrak{P}^-, K; \pi^\vee \otimes W_{\lambda}) \overset{\sim}\longrightarrow
    H^\ell(\mathfrak{P}^{\prime, -}, K'; {\pi'}^\vee \otimes {W'_{{\lambda}'}}),
    $$

where $T \in \operatorname{Hom}_{G'}(\pi^\vee|_{G'}, {\pi'}^\vee)$ and $\operatorname{pr}_{W} \in \operatorname{Hom}(W_\lambda|_{K'}, {W'_{\lambda'}})$ are non-zero homomorphisms as guaranteed by (ii) and (iii), respectively, each unique up to scalar multiplication.
The cohomology groups are one-dimensional when
$\ell=q(\delta) = q(\delta')$, and vanish in all other degrees..
\end{itemize}
\end{example}

To verify the table above, we also use the following identities:
\begin{align*}
    %\label{eqn:24042110}
  \Phi^{nc, +}(\delta^{\vee})
  &=\Delta({\mathfrak{p}}^+) \setminus \sigma_0\Phi^{nc, +}(\delta), 
\\
   (W_{\lambda}^-)^{\vee} &\simeq  
   F^K(-\sigma_0\lambda-\rho_{\mathfrak g})
  \quad  \text{($K$-isomorphism)}
\end{align*}

The following example is a key checkpoint for understanding the conventions used throughout the paper.
\begin{example}
In order to clarify our convention, we consider the simplest case: a holomorphic discrete series representation, whose decoration is given by
$\delta=(-\cdots -+\cdots+)$.  
Then $\Phi^{nc, +}(\delta)=\emptyset$, 
 $q(\delta)=0$, 
 and $\mu_{\lambda}=\lambda-\rho_{\mathfrak g}$
 by Lemma \ref{lem:24042109}.  
The highest weight of any $K$-type occuring in $\pi_{\lambda}^{\delta}$
 is of the form
\[
   \mu=\mu_{\lambda} + \sum_{\alpha \in \Delta({\mathfrak{p}}^-)} m_{\alpha} \alpha.
\]
Hence, we have a $K$-isomorphism
 $(\pi_{\lambda}^{\delta})^{{\mathfrak{p}}^+} \simeq \mu_{\lambda}$,
 and natural isomorphisms of $(\fP,K)$-cohomologies:  
\[
   H^0({\mathfrak{P}}, K;\pi_{\lambda}^{\delta} \otimes W_{\lambda}^{\vee})
   =
   ((\pi_{\lambda}^{\delta})^{\mathfrak{p}_+} \otimes W_{\lambda}^\vee)^K
   =
  \operatorname{Hom}_K(W_{\lambda}, \mu_{\lambda}) = {\mathbb{C}}.  
\]
\end{example}

\end{document}